\documentclass{elsarticle}
  
\usepackage[margin=1in]{geometry}

\usepackage{bbm}
\usepackage{graphicx}
\usepackage{pstool}
\usepackage{wrapfig}
\usepackage{caption}
\usepackage{subcaption}
\usepackage[section]{placeins} 

\usepackage{fancyhdr} 
\usepackage{lastpage} 
\usepackage{extramarks} 

\usepackage{xcolor} 
\usepackage{enumerate} 
\usepackage{paralist} 
\usepackage{amsmath, amsthm, amssymb, mathtools}
\usepackage{mathabx, pifont, stmaryrd} 
\usepackage[explicit]{titlesec} 
\usepackage{etoolbox} 
\usepackage{bibentry} 
\makeatletter\let\saved@bibitem\@bibitem\makeatother 
\usepackage[colorlinks, bookmarksopen, bookmarksnumbered,
            citecolor=red,urlcolor=red]{hyperref} 
\makeatletter\let\@bibitem\saved@bibitem\makeatother 
\usepackage{rotating}

\usepackage{algorithm}
\usepackage{algorithmic}

\newtheorem{remark}{Remark}

\theoremstyle{definition}

\usepackage{bm} 
\usepackage{amsfonts} 

\DeclareMathOperator*{\argmin}{arg\,min}

\newcommand{\norm}[1]{\ensuremath{\left\| #1 \right\|}}

\newcommand{\suchthat}{\mathrel{}\middle|\mathrel{}}

\newcommand{\pder}[2]{\ensuremath{\frac{\partial #1}{\partial #2}}}

\newcommand{\Ecal}{\ensuremath{\mathcal{E}}}
\newcommand{\Fcal}{\ensuremath{\mathcal{F}}}
\newcommand{\Gcal}{\ensuremath{\mathcal{G}}}
\newcommand{\Hcal}{\ensuremath{\mathcal{H}}}
\newcommand{\Ical}{\ensuremath{\mathcal{I}}}

\newcommand{\Kcal}{\ensuremath{\mathcal{K}}}

\newcommand{\Pcal}{\ensuremath{\mathcal{P}}}

\newcommand{\Scal}{\ensuremath{\mathcal{S}}}
\newcommand{\Tcal}{\ensuremath{\mathcal{T}}}

\newcommand{\Vcal}{\ensuremath{\mathcal{V}}}
\newcommand{\Wcal}{\ensuremath{\mathcal{W}}}

\newcommand{\Ycal}{\ensuremath{\mathcal{Y}}}

\newcommand{\Gbb}{\ensuremath{\mathbb{G}}}

\newcommand{\Rbb}{\ensuremath{\mathbb{R} }}
\newcommand{\Sbb}{\ensuremath{\mathbb{S} }}

\newcommand\Rbm{{\ensuremath{\bm{R}}}}

\newcommand\rbm{{\ensuremath{\bm{r}}}}

\newcommand\ubm{{\ensuremath{\bm{u}}}}
\newcommand\vbm{{\ensuremath{\bm{v}}}}

\newcommand\xbm{{\ensuremath{\bm{x}}}}
\newcommand\ybm{{\ensuremath{\bm{y}}}}

\newcommand\phibold{{\ensuremath{\boldsymbol{\phi}}}}

\newcommand\zerobold{\ensuremath{\mathbf{0}}}

\usepackage{tikz}
\usepackage{pgfplots}
\usepackage{pgfplotstable, filecontents, booktabs}
\pgfplotsset{compat=1.9}

\usetikzlibrary{pgfplots.groupplots}
\usepgfplotslibrary{fillbetween}
\usetikzlibrary{calc,fit,matrix,arrows,automata,positioning,shapes}
\usetikzlibrary{arrows.meta}

\pgfplotsset{select coords between index/.style 2 args={
    x filter/.code={
        \ifnum\coordindex<#1\fi
        \ifnum\coordindex>#2\fi
    }
}}

\tikzset{
 invisible/.style={opacity=0},
 visible on/.style={alt={#1{}{invisible}}},
 alt/.code args={<#1>#2#3}{%
   \alt<#1>{\pgfkeysalso{#2}}{\pgfkeysalso{#3}}
 },
}

\newcommand{\colorbarMatlabParula}[5]{
\begin{tikzpicture}
\begin{axis}[
   hide axis, scale only axis,
   height=0pt, width=0pt,
   colormap={parula}{rgb255=(62,38,168) rgb255=(62,39,172) rgb255=(63,40,175) rgb255=(63,41,178) rgb255=(64,42,180) rgb255=(64,43,183) rgb255=(65,44,186) rgb255=(65,45,189) rgb255=(66,46,191) rgb255=(66,47,194) rgb255=(67,48,197) rgb255=(67,49,200) rgb255=(67,50,202) rgb255=(68,51,205) rgb255=(68,52,208) rgb255=(69,53,210) rgb255=(69,55,213) rgb255=(69,56,215) rgb255=(70,57,217) rgb255=(70,58,220) rgb255=(70,59,222) rgb255=(70,61,224) rgb255=(71,62,225) rgb255=(71,63,227) rgb255=(71,65,229) rgb255=(71,66,230) rgb255=(71,68,232) rgb255=(71,69,233) rgb255=(71,70,235) rgb255=(72,72,236) rgb255=(72,73,237) rgb255=(72,75,238) rgb255=(72,76,240) rgb255=(72,78,241) rgb255=(72,79,242) rgb255=(72,80,243) rgb255=(72,82,244) rgb255=(72,83,245) rgb255=(72,84,246) rgb255=(71,86,247) rgb255=(71,87,247) rgb255=(71,89,248) rgb255=(71,90,249) rgb255=(71,91,250) rgb255=(71,93,250) rgb255=(70,94,251) rgb255=(70,96,251) rgb255=(70,97,252) rgb255=(69,98,252) rgb255=(69,100,253) rgb255=(68,101,253) rgb255=(67,103,253) rgb255=(67,104,254) rgb255=(66,106,254) rgb255=(65,107,254) rgb255=(64,109,254) rgb255=(63,110,255) rgb255=(62,112,255) rgb255=(60,113,255) rgb255=(59,115,255) rgb255=(57,116,255) rgb255=(56,118,254) rgb255=(54,119,254) rgb255=(53,121,253) rgb255=(51,122,253) rgb255=(50,124,252) rgb255=(49,125,252) rgb255=(48,127,251) rgb255=(47,128,250) rgb255=(47,130,250) rgb255=(46,131,249) rgb255=(46,132,248) rgb255=(46,134,248) rgb255=(46,135,247) rgb255=(45,136,246) rgb255=(45,138,245) rgb255=(45,139,244) rgb255=(45,140,243) rgb255=(45,142,242) rgb255=(44,143,241) rgb255=(44,144,240) rgb255=(43,145,239) rgb255=(42,147,238) rgb255=(41,148,237) rgb255=(40,149,236) rgb255=(39,151,235) rgb255=(39,152,234) rgb255=(38,153,233) rgb255=(38,154,232) rgb255=(37,155,232) rgb255=(37,156,231) rgb255=(36,158,230) rgb255=(36,159,229) rgb255=(35,160,229) rgb255=(35,161,228) rgb255=(34,162,228) rgb255=(33,163,227) rgb255=(32,165,227) rgb255=(31,166,226) rgb255=(30,167,225) rgb255=(29,168,225) rgb255=(29,169,224) rgb255=(28,170,223) rgb255=(27,171,222) rgb255=(26,172,221) rgb255=(25,173,220) rgb255=(23,174,218) rgb255=(22,175,217) rgb255=(20,176,216) rgb255=(18,177,214) rgb255=(16,178,213) rgb255=(14,179,212) rgb255=(11,179,210) rgb255=(8,180,209) rgb255=(6,181,207) rgb255=(4,182,206) rgb255=(2,183,204) rgb255=(1,183,202) rgb255=(0,184,201) rgb255=(0,185,199) rgb255=(0,186,198) rgb255=(1,186,196) rgb255=(2,187,194) rgb255=(4,187,193) rgb255=(6,188,191) rgb255=(9,189,189) rgb255=(13,189,188) rgb255=(16,190,186) rgb255=(20,190,184) rgb255=(23,191,182) rgb255=(26,192,181) rgb255=(29,192,179) rgb255=(32,193,177) rgb255=(35,193,175) rgb255=(37,194,174) rgb255=(39,194,172) rgb255=(41,195,170) rgb255=(43,195,168) rgb255=(44,196,166) rgb255=(46,196,165) rgb255=(47,197,163) rgb255=(49,197,161) rgb255=(50,198,159) rgb255=(51,199,157) rgb255=(53,199,155) rgb255=(54,200,153) rgb255=(56,200,150) rgb255=(57,201,148) rgb255=(59,201,146) rgb255=(61,202,144) rgb255=(64,202,141) rgb255=(66,202,139) rgb255=(69,203,137) rgb255=(72,203,134) rgb255=(75,203,132) rgb255=(78,204,129) rgb255=(81,204,127) rgb255=(84,204,124) rgb255=(87,204,122) rgb255=(90,204,119) rgb255=(94,205,116) rgb255=(97,205,114) rgb255=(100,205,111) rgb255=(103,205,108) rgb255=(107,205,105) rgb255=(110,205,102) rgb255=(114,205,100) rgb255=(118,204,97) rgb255=(121,204,94) rgb255=(125,204,91) rgb255=(129,204,89) rgb255=(132,204,86) rgb255=(136,203,83) rgb255=(139,203,81) rgb255=(143,203,78) rgb255=(147,202,75) rgb255=(150,202,72) rgb255=(154,201,70) rgb255=(157,201,67) rgb255=(161,200,64) rgb255=(164,200,62) rgb255=(167,199,59) rgb255=(171,199,57) rgb255=(174,198,55) rgb255=(178,198,53) rgb255=(181,197,51) rgb255=(184,196,49) rgb255=(187,196,47) rgb255=(190,195,45) rgb255=(194,195,44) rgb255=(197,194,42) rgb255=(200,193,41) rgb255=(203,193,40) rgb255=(206,192,39) rgb255=(208,191,39) rgb255=(211,191,39) rgb255=(214,190,39) rgb255=(217,190,40) rgb255=(219,189,40) rgb255=(222,188,41) rgb255=(225,188,42) rgb255=(227,188,43) rgb255=(230,187,45) rgb255=(232,187,46) rgb255=(234,186,48) rgb255=(236,186,50) rgb255=(239,186,53) rgb255=(241,186,55) rgb255=(243,186,57) rgb255=(245,186,59) rgb255=(247,186,61) rgb255=(249,186,62) rgb255=(251,187,62) rgb255=(252,188,62) rgb255=(254,189,61) rgb255=(254,190,60) rgb255=(254,192,59) rgb255=(254,193,58) rgb255=(254,194,57) rgb255=(254,196,56) rgb255=(254,197,55) rgb255=(254,199,53) rgb255=(254,200,52) rgb255=(254,202,51) rgb255=(253,203,50) rgb255=(253,205,49) rgb255=(253,206,49) rgb255=(252,208,48) rgb255=(251,210,47) rgb255=(251,211,46) rgb255=(250,213,46) rgb255=(249,214,45) rgb255=(249,216,44) rgb255=(248,217,43) rgb255=(247,219,42) rgb255=(247,221,42) rgb255=(246,222,41) rgb255=(246,224,40) rgb255=(245,225,40) rgb255=(245,227,39) rgb255=(245,229,38) rgb255=(245,230,38) rgb255=(245,232,37) rgb255=(245,233,36) rgb255=(245,235,35) rgb255=(245,236,34) rgb255=(245,238,33) rgb255=(246,239,32) rgb255=(246,241,31) rgb255=(246,242,30) rgb255=(247,244,28) rgb255=(247,245,27) rgb255=(248,247,26) rgb255=(248,248,24) rgb255=(249,249,22) rgb255=(249,251,21) },
   colorbar horizontal,
   point meta min=#1, point meta max=#5,
   colorbar style={width=10cm, xtick={#1,#2,#3,#4,#5}}
]
\addplot [draw=none] coordinates {(0,0)};
\end{axis}
\end{tikzpicture}
}

\usepackage{setspace}

\newbool{fastcompile}
\setbool{fastcompile}{false}

\newcommand{\wslab}{s}
\newcommand{\physF}{\mathcal{F}}
\newcommand{\physS}{\mathcal{S}}
\newcommand{\physU}{U}

\newcommand{\physUnitN}{\eta}

\newcommand{\transF}{\bar{\physF}} 
\newcommand{\transS}{\bar{\physS}}
\newcommand{\transU}{\bar{\physU}}

\newcommand{\transUnitN}{\bar{\physUnitN}}

\begin{document}
\title{A space-time high-order implicit shock tracking method for shock-dominated unsteady flows}

\author[rvt1]{Charles Naudet\fnref{fn1}}
\ead{cnaudet@nd.edu}

\author[rvt1]{Matthew J. Zahr\fnref{fn2}\corref{cor1}}
\ead{mzahr@nd.edu}

\address[rvt1]{Department of Aerospace and Mechanical Engineering, University
               of Notre Dame, Notre Dame, IN 46556, United States}
\cortext[cor1]{Corresponding author}

\fntext[fn1]{Graduate Student, Department of Aerospace and Mechanical
             Engineering, University of Notre Dame}
\fntext[fn2]{Assistant Professor, Department of Aerospace and Mechanical
             Engineering, University of Notre Dame}


\begin{abstract}
High-order implicit shock tracking (fitting) is a class of high-order,
optimization-based numerical methods to approximate solutions of conservation
laws with non-smooth features by aligning elements of the computational mesh
with non-smooth features. This ensures the non-smooth features are perfectly
represented by inter-element jumps and high-order basis functions approximate
smooth regions of the solution without nonlinear stabilization, which leads
to accurate approximations on traditionally coarse meshes. In this work,
we extend implicit shock tracking to time-dependent problems using
a slab-based space-time approach. This is achieved by reformulating a
time-dependent conservation law as a steady conservation law in one
higher dimension and applying existing implicit shock tracking techniques.
To avoid computations over the entire time domain and unstructured mesh generation
in higher dimensions, we introduce a general procedure to generate conforming,
simplex-only meshes of space-time slabs in such a way that preserves features
(e.g., curved elements, refinement regions) from previous time slabs. The use
of space-time slabs also simplifies the shock tracking problem by reducing
temporal complexity. Several practical adaptations of the implicit shock tracking
solvers are developed for the space-time setting including 1) a self-adjusting temporal
boundary, 2) nondimensionalization of a space-time slab, 3) adaptive
mesh refinement, and 4) shock boundary conditions, which lead to accurate
solutions on coarse space-time grids, even for problem with complex flow
features such as curved shocks, shock formation, shock-shock and shock-boundary
interaction, and triple points.
\end{abstract}
\maketitle

\section{Introduction}
Shock-dominated flows commonly arise in science and engineering disciplines including
supersonic and hypersonic aerodynamics, astrophysics, oceanic wave modeling,
radiation hydrodynamics, and many others. Despite their prevalence, accurate and
robust simulation of shock-dominated flows remains a significant challenge for
modern computational fluid dynamics methods, which becomes more significant as
the strength of shocks increase. High-order methods, such as discontinuous
Galerkin (DG) methods \cite{2001_cockburn_RKconv, 2007_hesthaven_dgmethod}, have
received considerable attention because they are highly accurate per degree of
freedom and introduce minimal dissipation, and have been shown to provide higher fidelity
solutions for problems with propagating waves, turbulent field flow, nonlinear
interactions, and multiple scales \cite{2013_wang_highordercfd}. DG offers additional
advantages such as geometric flexibility and a high degree of parallel scalability.
Despite these advantages, high-order methods are known to lack robustness
for shock-dominated flows because high-order approximation of shocks and contact
discontinuities leads to spurious oscillations that cause a breakdown of numerical solvers. 

The most popular class of methods to suppress these oscillations is known as
\textit{shock capturing}, which includes
limiting \cite{van1979towards,1999_baumann_discontfem},
reconstruction \cite{1987_harten_eno,1994_liu_WENO}, and artificial viscosity
\cite{2006_persson_shkcapt,2010_barter_shkcapt} approaches.
These methods reduce the order of accuracy of the method near
discontinuities, and therefore globally, to first-order, which reduces
the accuracy per degree of freedom of the underlying method. This loss
of accuracy is usually remedied by aggressive $h$-refinement, particularly
in the vicinity of discontinuities. While these methods have been shown to
solve complex flows, their robustness is inherently tied to an interplay
between the amount of dissipation added to smear discontinuities and the
grid fineness, and their accuracy will always be limited to first-order due
to the smooth approximation of discontinuous features. Furthermore, these
methods are susceptible to disturbances such as
surface anomalies \cite{gnoffo_computational_2004,gnoffo_simulation_2007,candler_unstructured_2007},
spurious oscillations \cite{lee_spurious_1999,carpenter_accuracy_1999},
and carbuncles  \cite{robinet_shock_2000,pandolfi_numerical_2001,roe_carbuncles_2005}.

An alternative approach, known a \textit{shock tracking} or \textit{shock fitting},
adjusts the computational mesh to align faces of mesh elements with discontinuities
in the solution. This allows the discontinuities to be represented perfectly with the
inter-element jump in the solution basis without the need of additional stabilization,
which recovers the high accuracy per degree of freedom of high-order methods and avoids
disturbances. However, this is a difficult meshing problem that requires generating a
fitted mesh to the (unknown) shock surface. Many previous approaches employ specialized
formulations and solvers which are dimension dependent and do not easily
generalize \cite{harten1983self, glimm2003conservative, bell1982fully}
and/or are limited to relatively simple
problems \cite{shubin1981steady, shubin1982steady, rosendale1994floating}.
One noteworthy class of methods is explicit shock tracking, surveyed in
\cite{moretti2002thirty, salas2009shock}, which mostly consist of explicitly
identifying the shock and using the Rankine-Hugoniot conditions to compute its
motion and upstream/downstream states. These methods are not easily applicable
to discontinuities whose topologies are not known \textit{a priori} or undergo
substantial shock motion, making them ill-suited for many important phenomena including
curved shocks, shock intersections, shock formation, and shock-boundary interactions.

A new class of high-order methods, known as \textit{implicit shock tracking}, has
recently emerged, which includes the Moving Discontinuous Galerkin Method with Interface
Condition Enforcement (MDG-ICE)
\cite{2019_corrigan_MDGICE,2021_kercher_MDGICE_LS,2021_kercher_MDGICE}
and High-Order Implicit Shock Tracking (HOIST)
\cite{2018_zahr_hoist,2020_zahr_HOIST,huang2022robust} methods.
Like traditional shock tracking, these methods aim to align element faces with
discontinuities in the flow to recover high accuracy per degree of freedom
without stabilization and avoid nonlinear instabilities. However, instead
of attempting to explicitly mesh the unknown discontinuity surface, they
discretize the conservation law on a mesh without knowledge of discontinuities
and solve an optimization problem over the discrete solution variables and nodal
coordinates of the mesh. The solution of this optimization problem is a
discontinuity-aligned mesh and the solution of the discretized conservation law.
This critical innovation has led to a general approach that has been shown to
handle complex shock phenomena including curved and reflected shocks, shock formation,
and shock-shock interaction for viscous and inviscid, inert and reacting flows of
varying degree of complexity \cite{2019_corrigan_MDGICE, 2021_kercher_MDGICE, 2020_zahr_HOIST,2021_zahr_istreact,huang2022robust}.

Implicit shock tracking methods are formulated in terms of \textit{steady}
conservation laws and have been predominately applied to steady problems
or unsteady problems using a space-time formulation that converts an
unsteady conservation law to a steady conservation law in one
higher dimension (usually coupling the entire temporal domain).
An exception is \cite{2021_shi_methlines} where the HOIST method
is embedded in a method of lines temporal discretization that solves
a steady shock tracking problem at each time step.
As noted in \cite{2021_shi_methlines,shi2022high}, this
approach suffers from similar limitations as explicit shock fitting
methods; in particular, it is most suitable for problems with limited
shock motion and cannot handle shock-shock or shock-boundary interactions
without complex mesh modifications. Space-time formulations
\cite{lowrie1998space,sudirham2006space,klaij2006space}
do not suffer from these limitations because propagating shocks manifest
as stationary discontinuity surfaces in space-time and shock intersection
are simply triple points in space-time.

In this work, we develop a slab-based space-time shock tracking approach based
on the HOIST method to track complex shock dynamics with a computational grid
without coupling the entire temporal domain. Our work is similar to
\cite{2019_corrigan_MGDICEunsteady} that demonstrated the MDG-ICE method in
a slab-based space-time setting. The contributions of our work are:
\begin{itemize}
 \setlength\itemsep{-1pt} 
 \item systematic formulation of a traditional conservation law where the
   spatial and temporal independent variables are treated separately as
   a steady conservation law in one higher dimension with special care
   given to the relationship between the space-time terms (flux function,
   source term, numerical flux function) and the corresponding term from
   the spatial formulation;
 \item a general procedure to generate and refine a conforming, simplex-only mesh
   of a space-time slab in such a way that preserves features (e.g., curved elements,
   refinement regions, tracked discontinuities) from the previous time slab;
 \item formulation of the HOIST method over a space-time slab; and
 \item practical adaptations of the HOIST for the space-time setting including
   a self-adjusting temporal boundary (first introduced in
   \cite{2019_corrigan_MGDICEunsteady}),
   nondimensionalization of space-time slab, adaptive mesh refinement,
   and shock boundary conditions whereby boundary conditions are directly
   applied to shocks and upstream elements are removed from the space-time
   mesh (if the upstream state is known) \cite{huang2023high}.
\end{itemize}

The remainder of this paper is organized as follows. Section~\ref{sec:govern}
systematically builds up a space-time formulation of a time-dependent conservation
law from a traditional formulation that separates the temporal and spatial independent
variables, transforms the space-time conservation law to a reference domain such that
domain deformations appear explicitly, and details a high-order discontinuous Galerkin
discretization of the transformed space-time conservation law. Section~\ref{sec:slab}
introduces our approach to generate conforming, simplex-only meshes of space-time slabs
that preserve features (refinement regions, high-order elements, tracked discontinuities)
in the previous slab. Section~\ref{sec:ist} details the extension of the HOIST method
to the space-time slab setting. Section~\ref{sec:numexp} shows the proposed
method successfully tracks space-time shocks and produces accurate solutions on
a series of increasingly complex one- and two-dimensional time-dependent conservation
laws, some of which include complicated shock features such as curved shocks,
shock-shock and shock-boundary intersections, and shock formation. Finally,
Section~\ref{sec:conclude} offers conclusions.

\section{Governing equations and high-order discretization}
\label{sec:govern}
In this section, we introduce a time-dependent system of conservation laws that will be
used to model unsteady, shock-dominated flows (Section~\ref{sec:govern:spatial}) as
well as its formulation as a space-time partial differential equation that treats the
temporal dimension as another spatial dimension (Section~\ref{sec:govern:sptm}).
Then we recast the space-time conservation law on a fixed reference space-time domain
such that deformations to the physical space-time domain appear explicitly in the
flux and source term (Section~\ref{sec:govern:transf}). Finally, we introduce a
DG discretization of the transformed space-time conservation law
(Section~\ref{sec:govern:disc}). Throughout this section, we emphasize systematic
construction of space-time quantities in terms of spatial quantities.

\subsection{Time-dependent system of conservation laws}
\label{sec:govern:spatial}
Consider a general system of $m$ hyperbolic partial differential equations
posed in the spatial domain $\Omega_x \subset \Rbb^{d'}$ over the time interval
$\mathcal{T} \coloneqq (t_0, t_1) \subset \Rbb{\ge 0}$
\begin{equation} \label{eqn:gen_cons_law}
 \pder{\physU_x}{t}+ \nabla_x \cdot \physF_x (\physU_x) = \physS_x(\physU_x),
\end{equation}
where $t \in \mathcal{T}$ is the temporal coordinate,
$x = (x_1, ..., x_{d'}) \in \Omega_x$ is the spatial coordinate,
$\physU_x(\cdot,t) : \Omega_x \rightarrow \Rbb^m$ is the conservative state
implicitly defined as the solution to (\ref{eqn:gen_cons_law}),
$\physF_x : \Rbb^m \rightarrow \Rbb^{m \times d'}$
with $\physF_x : W_x \mapsto \physF_x(W_x)$ is the physical flux function where $W_x$ is used to denote the first argument of the physical flux function,
$\physS_x : \Rbb^m \rightarrow \Rbb^m$ is the physical source term,
$(\nabla_x \cdot)$ is the divergence operator on the domain $\Omega_x$ defined as
$(\nabla_x\cdot\psi)_i := \partial_{x_j} \psi_{ij}$ (summation implied on repeated
index), and $\partial \Omega_x$ is the boundary of the spatial domain (with appropriate boundary conditions prescribed). In general, the solution $U(x,t)$ may contain
discontinuities, in which case, the conservation law (\ref{eqn:gen_cons_law}) holds
away from the discontinuities and the Rankine-Hugoniot conditions hold at discontinuities.

The Jacobian of the projected inviscid flux is 
\begin{equation} \label{eqn:space_inv_proj_jacob}
  B_x : \Rbb^m \times \Sbb_{d'} \rightarrow \Rbb^{m\times m}, \qquad
	B_x : (W_x, \eta_x) \mapsto \frac{\partial [\physF_x(W_x) \eta_x]}{\partial W_x},
\end{equation}
where $\Sbb_{d'} \coloneqq \{ \eta\in\Rbb^{d'} \mid \norm{\eta} = 1\}$ and $\eta_x$ is the unit normal in the spatial setting. The eigenvalue
decomposition of the Jacobian is
\begin{equation} \label{eqn:spat_eigval_decomp}
  B_x(W_x, \eta_x) = V_x(W_x, \eta_x) \Lambda_x(W_x, \eta_x) V_x(W_x, \eta_x)^{-1},
\end{equation}
where $\Lambda_x : \Rbb^m \times \Sbb_{d'} \rightarrow \Rbb^{m\times m}$
is a diagonal matrix containing the real eigenvalues of $B_x$ and the columns of
$V_x : \Rbb^m \times \Sbb_{d'} \rightarrow \Rbb^{m\times m}$
contain the corresponding right eigenvectors.
The projected Jacobian is commonly used to define linearized Riemann solvers
\cite{leveque2002finite,toro2013riemann} for the conservation law
(\ref{eqn:gen_cons_law}), which in turn is used to define monotonic,
upwind numerical flux functions for numerical methods,
$\Hcal_x : \Rbb^m\times\Rbb^m\times\Sbb_{d'} \rightarrow \Rbb^m$, as
\begin{equation} \label{eqn:space_numflux}
 \Hcal_x : (W_x^+,W_x^-,\eta_x) \mapsto \frac{1}{2}\left(\physF_x(W_x^+)\eta_x + \physF_x(W_x^-)\eta_x\right) + \frac{1}{2}\tilde{B}_x(W_x^+,W_x^-,\eta_x)(W_x^+-W_x^-),
\end{equation}
where $\tilde{B}_x : \Rbb^m\times\Rbb^m\times\Sbb_{d'} \rightarrow \Rbb^{m\times m}$
is the Jacobian of the linearized Riemann problem \cite{toro2013riemann}. In this work,
we use Roe's numerical flux \cite{roe1981approximate} that takes
\begin{equation} \label{eqn:space_linriem}
  \tilde{B}_x(W_x^+,W_x^-,\eta_x) =  \left|B_x(\hat{U}(W_x^+,W_x^-),\eta_x)\right|,
\end{equation}
where $|B_x| = V_x |\Lambda_x| V_x^{-1}$ is the matrix absolute value of the
projected Jacobian and $\hat{U} : \Rbb^m \times \Rbb^m \rightarrow \Rbb^m$ is
the equation-specific Roe average.

\begin{remark}
The Jacobian of the projected inviscid flux is well-defined for any projection
direction using (\ref{eqn:space_inv_proj_jacob}); however, in this work, we restrict
its definition to \textit{unit} projection directions ($\Sbb_{d'}$) because the eigenvalue
decompositions can be substantially simpler under this assumption. Expressions for the
projected Jacobian and its eigenvalue decomposition for the conservation laws considered
in this work are provided in \ref{app:projjac}.
\end{remark}

\begin{remark}
Other standard numerical flux functions fit the form in (\ref{eqn:space_numflux}).
The Vijayasundaram flux \cite{vijayasundaram1986transonic} uses a linear
Riemann problem of the form
(\ref{eqn:space_numflux}) with the projected Jacobian evaluated at the arithmetic
average of the two states $\hat{U}(W_x^+,W_x^-) = (W_x^++W_x^-)/2$ and the flux
average replaced with the linearized flux average
$B_x(\hat{U}(W_x^+,W_x^-),\eta_x)\hat{U}(W_x^+,W_x^-)$ \cite{wallraff2013numerical}.
Furthermore, the local Lax-Friedrichs or Rusanov \cite{toro2013riemann}
flux takes the linearized Riemann problem matrix to be
\begin{equation}
 \tilde{B}_x : (W_x^+,W_x^-,\eta_x) \mapsto
 \max\left\{\norm{\Lambda_x(W_x^+,\eta_x)}_\infty, \norm{\Lambda_x(W_x^-,\eta_x)}_\infty\right\} I_m,
\end{equation}
where $I_m \in \Rbb^{m\times m}$ is the identity matrix. 
\end{remark}

\begin{remark}
Smoothed \cite{2020_zahr_HOIST} and entropy-fixed \cite{harten1983self} variants of
the numerical flux in (\ref{eqn:space_numflux}) are obtained by replacing the
component-wise absolute value function with a smoothed absolute value or entropy fix.
\end{remark}

\subsection{Space-time formulation} 
\label{sec:govern:sptm} 
The conservation law in (\ref{eqn:gen_cons_law}) describes a general time-dependent
system of conservation laws in a $d'$-dimensional spatial domain. Because the proposed
method tracks discontinuities over space-time slabs, we reformulate
(\ref{eqn:gen_cons_law}) as a steady conservation law in a space-time domain
\cite{lowrie1998space,sudirham2006space, klaij2006space}.
To this end, we define the space-time domain as
$\Omega \coloneqq \Omega_x \times \Tcal \subset \Rbb^d$ ($d = d'+1$) with boundary
$\partial\Omega$, and let $z=(x,t)\in\Omega$ denote the space-time coordinate.
Because the space-time domain is a Cartesian product of the spatial domain with
a time interval, we will refer to it as a \textit{space-time slab}. The boundary
of a space-time slab $\partial\Omega$ consists of three pieces: 1) the spatial
boundary, $\partial\Omega_x \times \Tcal$, 2) the bottom of the slab,
$\Omega_x \times \{t_0\}$, and 3) the top of the slab, $\Omega_x \times \{t_1\}$.
Without loss of generality, we formulate the space-time conservation law, as well
as its transformation and discretization
(Section~\ref{sec:govern:disc}) over a single arbitrary slab corresponding
to the time interval $\Tcal$; in practice, we use a sequence of space-time slabs
to cover the temporal domain of interest.

The conservation law in (\ref{eqn:gen_cons_law}) can be written as a
steady conservation law over the space-time slab as
\begin{equation} \label{eqn:sptm_claw}
\nabla \cdot \physF (\physU) = \physS(\physU),
\end{equation}
where $U : \Omega \rightarrow \Rbb^m$ is the space-time conservative vector implicitly
defined as the solution of (\ref{eqn:sptm_claw}) and related to the solution of the
spatial conservation law as $\physU : z \mapsto  \physU_x(x,t)$,
$\physF : \Rbb^m \rightarrow \Rbb^{m \times d}$ and 
$\physS : \Rbb^m \rightarrow \Rbb^{m}$ are the
space-time flux function and source term, respectively, and related
to the spatial conservation law terms as
\begin{equation} \label{eqn:sptm_flux_src}
 \physF : W \mapsto \begin{bmatrix} \physF_x(W) & W \end{bmatrix}, \qquad
 \physS : W \mapsto \physS_x(W),
\end{equation}
and $(\nabla\cdot)$ is the space-time divergence operator defined as
$\nabla\cdot\begin{bmatrix}\psi&\phi\end{bmatrix} = \nabla_x\cdot\psi + \partial_t \phi$.

The Jacobian of the space-time projected inviscid flux
$B : \Rbb^m \times \Sbb_d \rightarrow \Rbb^{m\times m}$
is
\begin{equation} \label{eqn:space_time_proj_jac}
  B : (W, \eta) \mapsto \frac{\partial [\physF (W) \eta]}{\partial W}.
\end{equation}
Any $\eta\in\Sbb_d$ can be written as $\eta = (n_x, n_t)$ where
$n_x \in \Rbb^{d'}$ and $n_t \in \Rbb$, and $n_x$ can be expanded
$n_x = \eta_x \norm{n_x}$ with $\eta_x \in \Sbb_{d'}$ ($\eta_x$
is uniquely defined as $n_x/\norm{n_x}$ in the case where $n_x \neq 0$,
otherwise it is arbitrary). From this expansion and the form of the inviscid
flux in (\ref{eqn:sptm_flux_src}), the space-time Jacobian and be related to
the original Jacobian
\begin{equation} \label{eqn:jac_sp_sptm}
 B(W, \eta) = B_x(W,\eta_x) \norm{n_x} + n_t I_m.
\end{equation}
The eigenvalue decomposition of the projected Jacobian is denoted
\begin{equation}
 B(W, \eta) = V(W,\eta) \Lambda(W,\eta) V(W,\eta)^{-1},
\end{equation}
where $\Lambda : \Rbb^m \times \Sbb_d$ is a diagonal matrix containing the
real eigenvalues of $B$ and the columns of
$V : \Rbb^m \times \Sbb_d \rightarrow \Rbb^{m\times m}$
contain the corresponding right eigenvectors. Owing to the relationship
between the projected Jacobians of the spatial and space-time inviscid
fluxes (\ref{eqn:jac_sp_sptm}), their eigenvalue decompositions are related as
\begin{equation} \label{eqn:sptm_eigval_decomp}
 \Lambda(W, \eta) = \Lambda_x(W,\eta_x)\norm{n_x} + n_t I_m, \qquad
 V(W, \eta) = V_x(W,\eta_x)
\end{equation}
Linearized Riemann solvers for the space-time conservation law
(\ref{eqn:sptm_claw}) and the associated numerical fluxes,
$\Hcal : \Rbb^m \times \Rbb^m \times \Sbb_d \rightarrow \Rbb^m$,
are obtained from the formulas in Section~\ref{sec:govern:spatial} with
the space-time quantities in place of the corresponding spatial terms
\begin{equation} \label{eqn:sptm_numflux}
 \Hcal : (W^+,W^-,\eta) \mapsto \frac{1}{2}\left(\physF(W^+)\eta + \physF(W^-)\eta\right) + \frac{1}{2}\tilde{B}(W^+,W^-,\eta)(W^+-W^-).
\end{equation}
The Jacobian of the linearized Riemann problem 
$\tilde{B} : \Rbb^m\times\Rbb^m\times\Sbb_d \rightarrow \Rbb^{m\times m}$ is
\begin{equation} \label{eqn:sptm_linriem}
  \tilde{B}(W^+,W^-,\eta) =  \left|B(\hat{U}(W^+,W^-),\eta)\right|,
\end{equation}
which can be directly written in terms of spatial quantities using
(\ref{eqn:sptm_eigval_decomp})
\begin{equation} \label{eqn:sptm_linriem2}
 \tilde{B}(W^+,W^-,\eta) =
  V_x(\hat{W},\eta_x)
  \left|\Lambda_x(\hat{W},\eta_x)\norm{n_x} + n_t I_m\right|
  V_x(\hat{W},\eta_x)^{-1},
\end{equation}
where $\hat{W} = \hat{U}(W^+,W^-)$ was introduced for brevity.

\begin{remark}
The space-time version of the Vijayasundaram flux uses a linear
Rimean problem of the form (\ref{eqn:sptm_numflux}) with the arithmetic average
$\hat{U}(W^+,W^-) = (W^++W^-)/2$ and the flux average replaced with the linearized
flux average $B(\hat{U}(W^+,W^-),\eta)\hat{U}(W^+,W^-)$ \cite{wallraff2013numerical}.
The space-time Rusanov flux takes the linearized Riemann problem matrix to be
\begin{equation}
 \tilde{B} : (W^+,W^-,\eta) \mapsto
 \max\left\{\norm{\Lambda(W^+,\eta)}_\infty, \norm{\Lambda(W^-,\eta)}_\infty\right\} I_m,
\end{equation}
which can be directly written in terms of spatial quantities using
(\ref{eqn:sptm_eigval_decomp})
\begin{equation}
 \tilde{B} : (W^+,W^-,\eta) \mapsto
 \max\left\{\norm{\Lambda_x(W^+,\eta_x)\norm{n_x}+n_t I_m}_\infty,
            \norm{\Lambda_x(W^-,\eta_x)\norm{n_x}+n_t I_m}_\infty\right\} I_m.
\end{equation}
\end{remark}

\begin{remark} \label{rem:sptm_impl}
The space-time formulation of the conservation law (\ref{eqn:sptm_claw})
can be readily obtained from spatial formulation (\ref{eqn:gen_cons_law})
with the physical flux and source term in (\ref{eqn:sptm_flux_src})
and popular numerical flux functions based on linearized Riemann solvers
in (\ref{eqn:sptm_numflux})-(\ref{eqn:sptm_linriem2}). This allows
for a general, systematic implementation of space-time conservation
laws from the standard spatial formulation, which eliminates one of
the implementation burdens associated with space-time methods.
\end{remark}

\subsection{Transformed space-time conservation law on a fixed reference domain}
\label{sec:govern:transf}
Because the proposed numerical method is based on deforming the space-time domain
to track discontinuities with the computational grid, it is convenient to recast
the space-time conservation law such that domain deformations appear explicitly.
To this end, we define $\bar\Omega \subset \Rbb^d$ as a fixed space-time reference
domain, which we require to take the form of a space-time slab, i.e.,
$\bar\Omega \coloneqq \Omega_x \times \Tcal$. Let $\Gbb$ be the collection of
diffeomorphisms from the reference domain to the physical domain, i.e., for any
$\Gcal \in \Gbb$, $\Gcal : \bar\Omega \rightarrow \Omega$ with
$\Gcal : \bar{z} \mapsto \Gcal(\bar{z})$, that preserve the space-time slab
structure of the physical domain (slab-preserving mappings will be constructed
in Section~\ref{sec:slab}). The space-time conservation law in
(\ref{eqn:sptm_claw}) can be reformulated as a conservation law over the
reference domain as \cite{2018_zahr_hoist}
\begin{equation} \label{eqn:trans1}
  \bar{\nabla} \cdot \transF(\transU; G) = \transS (\transU; g),
\end{equation}
where $\bar{U} : \bar\Omega \rightarrow \Rbb^m$ is the solution of the transformed
conservation law, $\bar\Fcal : \Rbb^m \times \Rbb^{d\times d} \rightarrow \Rbb^{m\times d}$
and $\bar\Scal : \Rbb^m \times \Rbb \rightarrow \Rbb^m$ are the transformed flux function
and source term, respectively, $\bar\nabla\cdot$ is the divergence operator in the
reference domain $\bar\Omega$ defined as
$(\nabla \cdot \psi)_i = \partial_{\bar{z}_j}\psi_{ij}$, 
and the deformation gradient $G : \bar\Omega \rightarrow \Rbb^{d\times d}$
and mapping Jacobian $g : \bar\Omega \rightarrow \Rbb$ are
\begin{equation}
 G : \bar{z} \mapsto \partial_{\bar{z}}\Gcal(\bar{z}), \qquad
 g: \bar{z} \mapsto \det G(\bar{z}).
\end{equation}
The reference domain quantities are defined in terms of the corresponding
physical domain quantities as \cite{2018_zahr_hoist}
\begin{equation} \label{eqn:trans2}
 \transU(\bar{z}) = U(\Gcal(\bar{z})),
 \qquad
 \transF : (\bar{W}; \Theta) \mapsto (\det\Theta) \Fcal(\bar{W}) \Theta^{-T},
 \qquad
 \transS : (\bar{W}; q) \mapsto q S(\bar{W}).
\end{equation}
Transformed numerical flux functions,
$\bar\Hcal : \Rbb^m\times\Rbb^m\times\Sbb_d\times\Rbb^{d\times d}\mapsto\Rbb^m$,
by virtue of being a surrogate for the inviscid flux projected in a particular
direction are related to the physical numerical flux \cite{2020_zahr_HOIST} as
\begin{equation} \label{eqn:trans3}
 \bar\Hcal : (\bar{W}^+, \bar{W}^-, \bar\eta; \Theta) \mapsto
    \sigma \Hcal(\bar{W}^+, \bar{W}^-, \sigma^{-1}(\det\Theta)\Theta^{-T}\bar\eta )
\end{equation}
where $\sigma = \norm{(\det\Theta) \Theta^{-T}\bar\eta}$ is defined for brevity.

\begin{remark} \label{rem:trans_impl}
Analogous to Remark~\ref{rem:sptm_impl}, the transformed conservation
law (\ref{eqn:trans1}) can readily be obtained from the space-time
conservation (\ref{eqn:sptm_claw}) with the transformed flux and source
term defined in (\ref{eqn:trans2}) and the transformation of any
numerical flux in (\ref{eqn:trans3}). This allows for a general,
systematic implementation of transformed space-time conservation
laws from the standard spatial formulation.
\end{remark}

\subsection{Discontinuous Galerkin discretization of the transformed space-time conservation law}
\label{sec:govern:disc}
Because we are interested in high-order shock fitting solutions, we discretize the
transformed conservation law (\ref{eqn:trans1}) with a discontinuous Galerkin method
\cite{2007_hesthaven_dgmethod}
with high-order piecewise polynomials spaces used to approximate the state
$\transU$ and domain mapping $\Gcal$.  To this end, let $\bar\Ecal_h$ represent
a discretization of the reference domain $\Omega_0$ into non-overlapping,
potentially curved, computational elements. To establish the finite-dimensional
DG formulation, we introduce the DG approximation (trial) space of discontinuous
piecewise polynomials associated with the mesh $\bar\Ecal_h$
\begin{equation}
 \Vcal_h^p \coloneqq
 \left\{
  v \in [L^2(\bar\Omega)]^m \suchthat \left. v\right|_{\bar{K}} \in [\Pcal_p(\bar{K})]^m,~\forall \bar{K}\in\bar\Ecal_h
 \right\},
\end{equation}
where $\Pcal_p(\bar{K})$ is the space of polynomial functions of degree at most
$p \geq 0$ over the domain $\bar{K}$. Furthermore, we define the space of globally
continuous piecewise polynomials of degree $q$ associated with the mesh $\bar\Ecal_h$ as
\begin{equation}
 \Wcal_h \coloneqq
 \left\{
  v \in C^0(\bar\Omega) \suchthat \left. v\right|_{\bar{K}} \in \Pcal_q(\bar{K}),~\forall \bar{K}\in\bar\Ecal_h
 \right\}
\end{equation}
and discretize the domain mapping with the corresponding vector-valued space
$\left[\Wcal_h\right]^d$. With these definitions, the DG variational problem
is: given $\Gcal_h \in \left[\Wcal_h\right]^d$, find $\bar{U}_h\in\Vcal_h^p$
such that for all $\bar\psi_h \in \Vcal_h^{p'}$, we have
\begin{equation} \label{eqn:weak1}
 r_h^{p',p}(\bar\psi_h, \transU_h; \bar\nabla\Gcal_h) = 0
\end{equation}
where $p' \geq p$ and the global residual function
$r_h^{p',p} : \Vcal_h^{p'}\times\Vcal_h^p\times[\Wcal_h]^d \rightarrow \Rbb$
is defined in terms of elemental residuals
$r_{\bar{K}}^{p',p} : \Vcal_h^{p'}\times\Vcal_h^p\times[\Wcal_h]^d \rightarrow \Rbb$
as
\begin{equation} \label{eqn:weak2}
r_h^{p', p} : (\bar\psi_h, \bar{W}_h; \Theta_h) \mapsto
 \sum_{\bar{K}\in\bar\Ecal_h} r_{\bar{K}}^{p', p} (\bar\psi_h, \bar{W}_h; \Theta_h).
\end{equation}
The elemental residuals come directly from a standard DG formulation applied to
the transformed space-time conservation law in (\ref{eqn:trans1})
\begin{equation} \label{eqn:weak3}
\begin{aligned}
 r_{\bar{K}}^{p',p} : (\bar\psi_h, \bar{W}_h; \Theta_h) \mapsto
 &\int_{\partial \bar{K}}\bar\psi_h \cdot \bar{\Hcal} ( W^+_h,  W^-_h, \transUnitN_h; \Theta_h  ) \, dS  \\
- &\int_{\bar{K}} \transF (W_h; \Theta_h) : \bar{\nabla} \bar\psi_h \, dV \\
- &\int_{\bar{K}} \bar\psi_h \cdot \transS (W_h; \det \Theta_h)) \, dV,
\end{aligned}
\end{equation}
where $\bar\eta_h : \partial\bar{K} \mapsto \Rbb^d$ is the outward unit normal to
an element $\bar{K} \in \bar\Ecal_h$, and $\bar{W}_h^+$ ($\bar{W}_h^-$) denotes the interior
(exterior) trace of $\bar{W}_h \in \Vcal_h^p$ to the element.

Finally, we introduce a basis for the test space ($\Vcal_h^{p'}$), trial space
($\Vcal_h^p$), and domain mapping space ($[\Wcal_h]^d$) to reduce the weak
formulation in (\ref{eqn:weak1})-(\ref{eqn:weak3}) to a system of nonlinear
algebraic equations. In the case where $p = p'$, we have
\begin{equation}
 \rbm : \Rbb^{N_\ubm} \times \Rbb^{N_\xbm} \rightarrow \Rbb^{N_\ubm}, \qquad
 \rbm : (\ubm, \xbm) \mapsto \rbm(\ubm,\xbm),
\end{equation}
where $N_\ubm = \mathrm{dim}(\Vcal_h^p)$ and $N_\xbm = \mathrm{dim}([\Wcal_h]^d)$,
which is the residual of a standard space-time DG method. Furthermore, we define
the algebraic enriched residual associated with a test space of degree $p' > p$
($p' = p+1$ in this work) as
\begin{equation}
 \Rbm : \Rbb^{N_\ubm} \times \Rbb^{N_\xbm} \rightarrow \Rbb^{N_\ubm'}, \qquad
 \Rbm : (\ubm, \xbm) \mapsto \Rbm(\ubm,\xbm),
\end{equation}
where $N_\ubm' = \mathrm{dim}(\Vcal_h^{p'})$, which will be used to construct
the implicit shock tracking objective function. Finally, to maintain a connection
between the algebraic and functional representation of the DG solution, we define
the operator $\Xi : \Rbb^{N_\ubm} \rightarrow \Vcal_h^p$ that maps
$\vbm\in\Rbb^{N_\ubm}$ to its representation as a function over the reference
domain $\bar\Omega$, $V_h = \Xi(\vbm) \in \Vcal_h^p$.


\section{Mesh generation of space-time slab}
\label{sec:slab}
In this section, we address construction of a conforming,
simplex-only mesh $\bar\Ecal_h$ of the space-time slab. We insist
on a simplex-only mesh because element removal via edge collapse
is well-defined in any dimension and key to robust shock tracking
\cite{huang2022robust}. Furthermore, we construct a conforming mesh
because it makes construction of the continuous space $[\Wcal_h]^d$
simple and convenient using standard finite elements where adjacent
elements share nodes. Non-conforming meshes can be used provided constraints
are added to ensure the mapped mesh remains valid (no gaps or overlaps).

One option to construct a conforming, simplex-only mesh of the space-time
slab $\bar\Omega$ is to directly use an unstructured mesh generator in
the $(d'+1)$-dimensional domain. While substantial progress has been made
in higher dimensional mesh generation
\cite{behr2008simplex,yano2012optimization,von2019simplex,anderson2023surface},
direct and automatic generation of high-order, high-quality (sliver-free) meshes
in dimensions beyond $d=d'+1=2$ remains a challenge. Automation is particularly
important to avoid user intervention as a simulation advances in time (between slabs).
Furthermore, this approach may neglect important information from the previous time slab,
e.g., the location of tracked discontinuities in the present setting. Instead, we utilize
an extrusion approach
\cite{klaij2006space,2015_wang_sptm,diosady2017tensor}
whereby the (simplicial) elements of a mesh of the spatial domain
are extruded into the temporal dimension to define space-time prismatic elements, which
are subsequently split into space-time simplices. With this approach, only a complete mesh
of the spatial domain is required along with two simple mesh operations: extrusion of a
$d'$-dimensional simplex and splitting of a $(d'+1)$-dimensional prism into simplices.
We detail these steps individually in the remainder of this section, including an approach
to split high-order elements.

\subsection{Construction of $d'$-dimensional (spatial) simplicial mesh}
\label{sec:slab:spatial}
\begin{figure}
\centering
\includegraphics[width=0.3\textwidth]{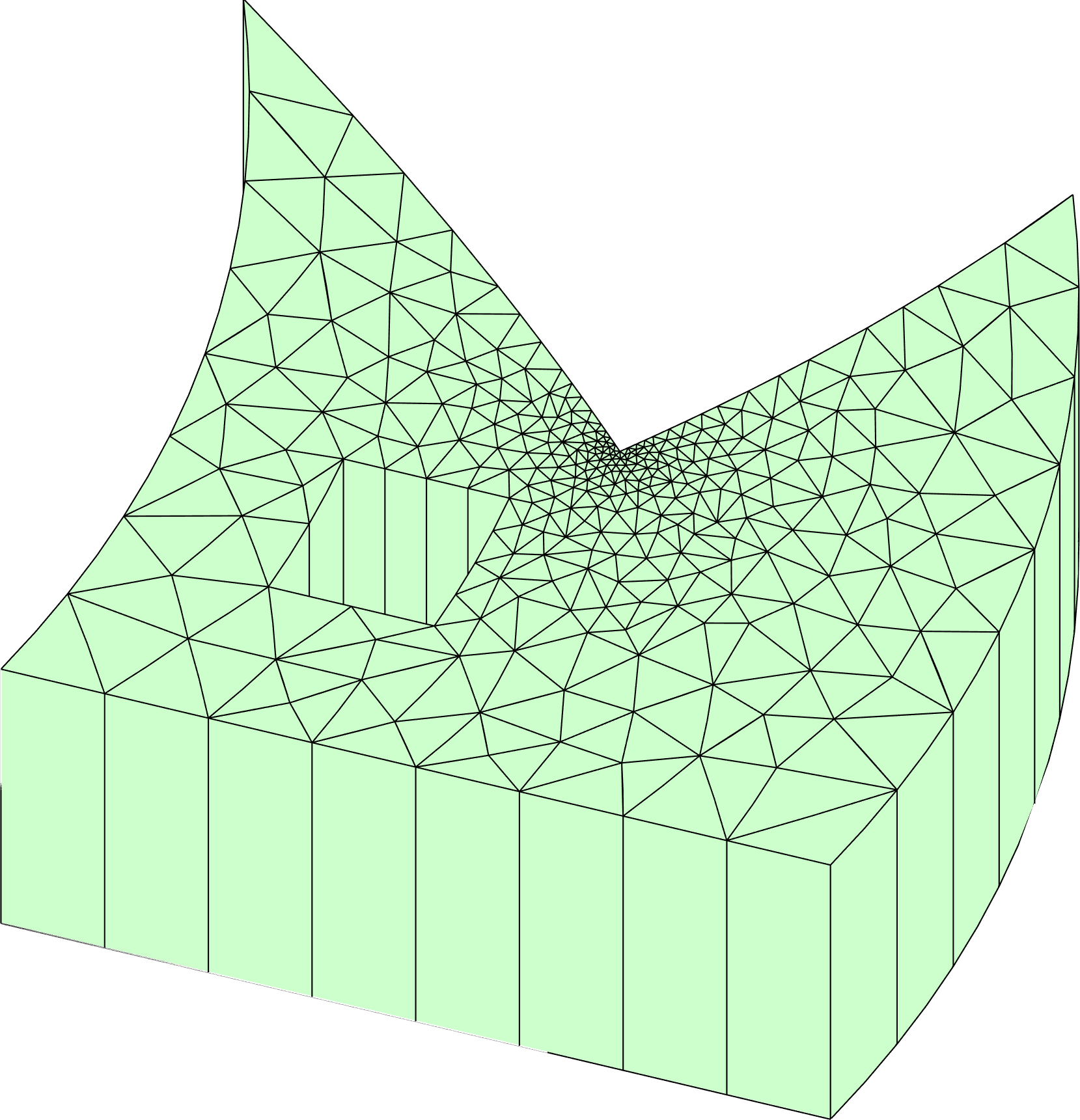} \qquad \qquad
\includegraphics[width=0.3\textwidth]{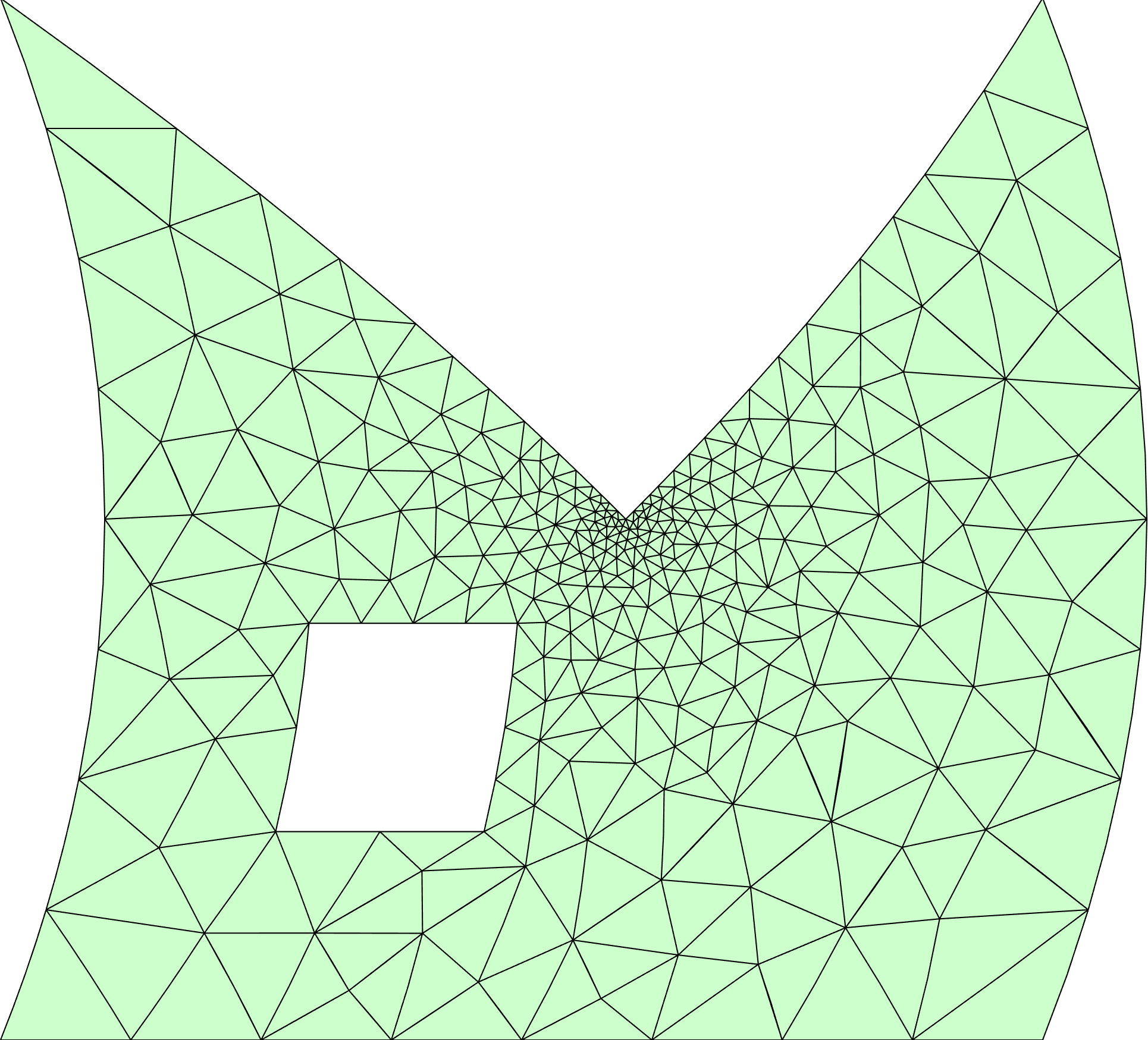}
\caption{Space-time mesh in $d'+1=3$ (\textit{left}) and $d'=2$ spatial mesh
 extracted from top of time slab (\textit{right}).}
\label{fig:extract}
\end{figure}

In the general case where multiple slabs are used to cover the temporal domain of
interest, we distinguish two cases. For the first time slab, a (potentially) curved
simplex-only mesh, $\bar\Ecal_{x,h}$, of the spatial domain $\Omega_x$ is constructed
using standard mesh generation approaches. While this can suffer from the issues
discussed in the previous section, it only needs to be performed once for an entire
time-dependent simulation (as opposed to being required for each slab) and is a
standard part of computational physics workflows. For all other time slabs, the
spatial mesh is extracted from the top boundary $\partial\Omega_x \times \{t_1\}$
of the space-time mesh of the previous time slab (Figure~\ref{fig:extract}).
This approach ensures the spatial mesh $\bar\Ecal_{x,h}$ contains all of the
features (e.g., refinement regions, tracked discontinuities, element curvature)
present in the mesh of the previous space-time slab.

\subsection{Extrusion of $d'$-dimensional (spatial) simplicial mesh}
\label{sec:slab:extrude}
\begin{figure}
\centering
\includegraphics[width=0.3\textwidth]{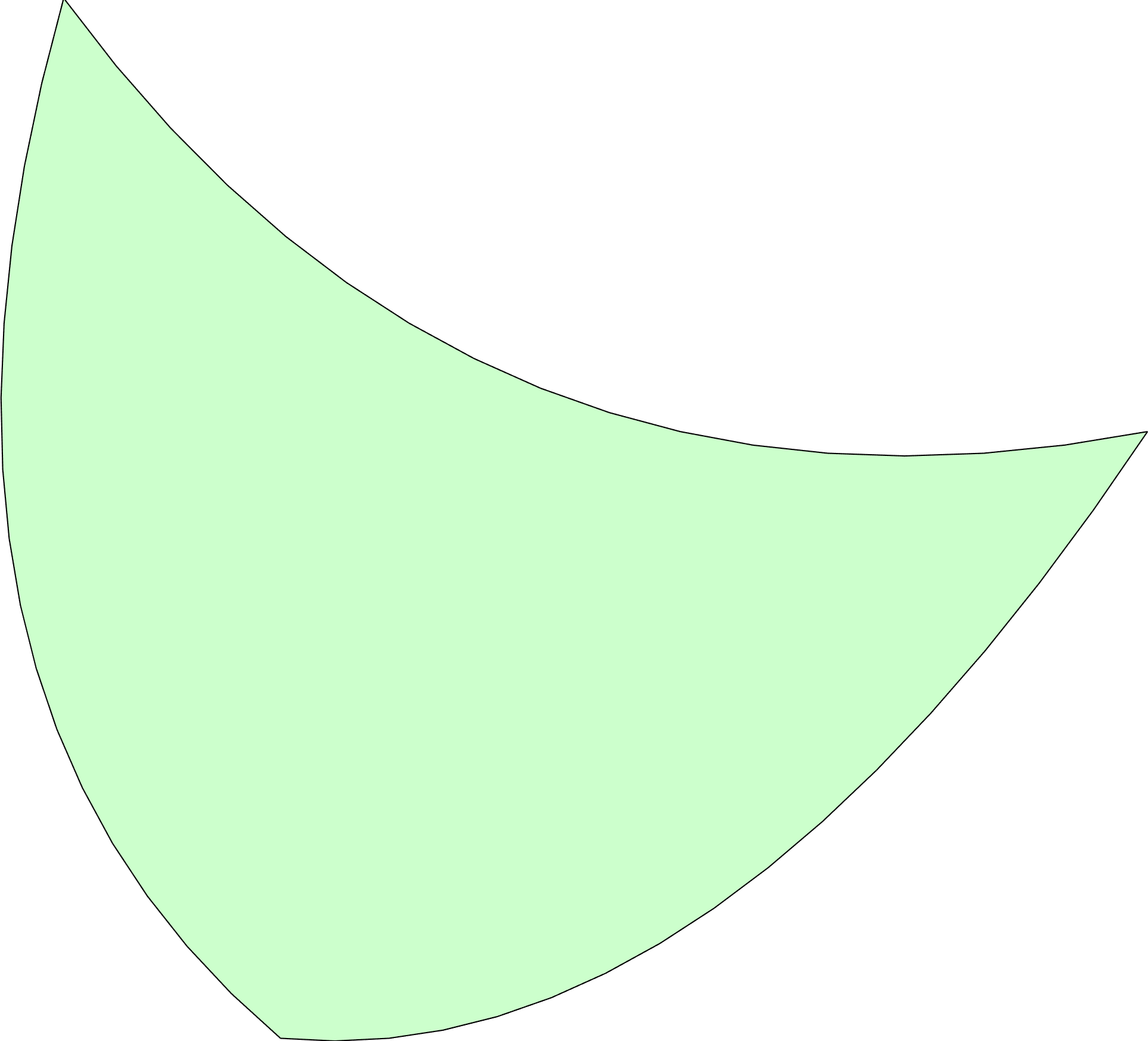} \qquad \qquad
\includegraphics[width=0.27\textwidth]{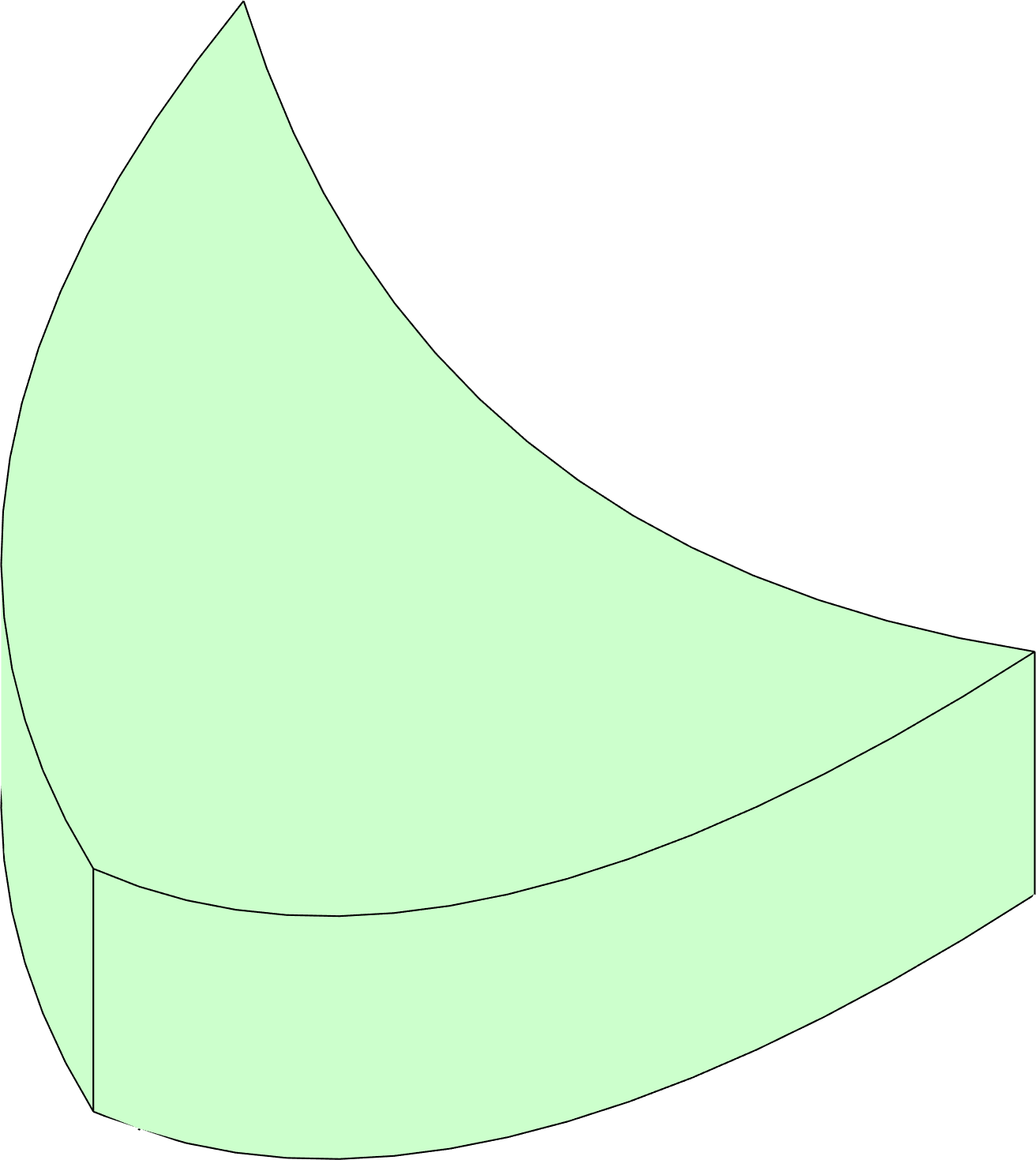}
\caption{Curved triangle in $d'=2$ spatial domain (\textit{left}) extruded
 to triangular prism in $d'+1=3$ space-time domain (\textit{right}).}
\label{fig:extrusion}
\end{figure}

Once a simplex-only mesh $\bar\Ecal_{x,h}$ of the spatial domain $\Omega_x$
is available, each spatial element is extruded in the temporal dimension
to form a space-time simplicial prism. That is, the mesh $\bar\Ecal_h$ of the
space-time slab $\bar\Omega = \Omega_x \times \Tcal$ is defined as
\begin{equation}
 \bar\Ecal_h = \bigcup_{K_x \in \bar\Ecal_{x,h}} K_x \times \Tcal,
\end{equation}
where $K_x \times \Tcal \subset \bar\Omega$ is a (potentially) curved
simplicial prism in $d$-dimensions (Figure~\ref{fig:extrusion}).
This extrusion approach preserves characteristics of the spatial mesh
throughout the space-time slab such as refinement regions and tracked
discontinuities.

\subsection{High-order space-time mesh refinement via transformations}
\label{sec:slab:split}
The final step required to create a mesh of a space-time slab is to
split the potentially curved simplicial prism into simplices. We
define a systematic procedure that splits a straight-sided version
of the elements in the parent domain and uses the isoparametric
mapping to define a high-order splitting in the reference domain.
Then the domain mapping $\Gcal$ defines a high-order splitting
in the physical domain. This is simpler than refining
potentially curved elements in the reference domain.
This approach will be used to both split prism
elements into simplices, as well as subdivide simplex elements into
smaller simplicies (for mesh refinement, Section~\ref{sec:ist:amr}) while still maintaining potentially curved high-order elements.

\begin{figure}
\centering
\input{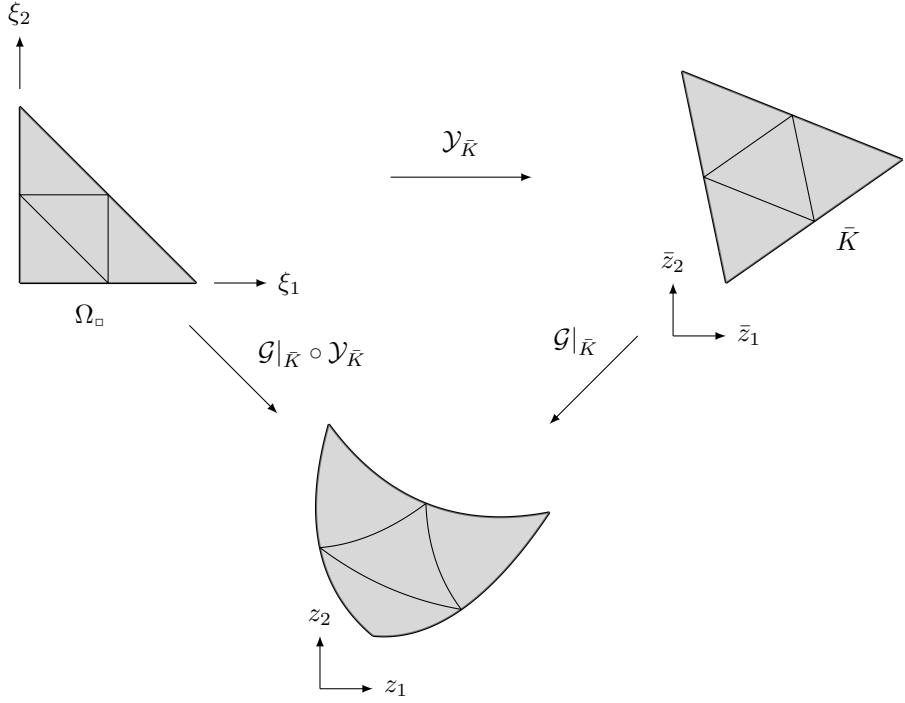}
\caption{Schematic of high-order element refinement whereby the refinement operator
 $\Pi$ is applied to the parent element, which is mapped via the isoparametric mapping
 of the unrefinement element $\Ycal_{\bar{K}}$ to define a refinement in the reference
 domain $\bar\Omega$, which is in turn mapped via the transformation $\Gcal$ to define a
 refinement in the physical domain $\Omega$.}
\label{fig:meshop}
\end{figure}

To this end, we define a parent domain $\Omega_\square\subset\Rbb^d$
along with a transformation $\Ycal_{\bar{K}} : \Omega_\square \rightarrow \bar{K}$,
$\Ycal_{\bar{K}} : \xi \mapsto \Ycal_{\bar{K}}(\xi)$ that maps
$\Omega_\square$ to the element $\bar{K} \in \bar\Ecal_h$. In this work, the parent
domain $\Omega_\square$ will either be the unit simplex $\Kcal_0^d\subset\Rbb^d$
or simplicial prism $\mathcal{SP}_0^d$, defined as
\begin{equation}
 \begin{aligned}
 \Kcal_0^d &\coloneqq \left\{ \xi \in \Rbb^d \suchthat \sum_{i=1}^d \xi_i \leq 1,~\xi_j\geq 0,~j=1,\dots,d\right\} \\
 \mathcal{SP}_0^d &\coloneqq \left\{ \xi \in \Rbb^d \suchthat \sum_{i=1}^{d-1} \xi_i \leq 1,~\xi_j\geq 0,~j=1,\dots,d-1,~-1\leq\xi_d\leq 1\right\},
 \end{aligned}
\end{equation}
and $\Ycal_{\bar{K}} \in [\Pcal_q(\Omega_\square)]^d$ is defined using a standard
isoparametric approach
\begin{equation}
 \Ycal_{\bar{K}} : \xi \mapsto \sum_{I=1}^n \bar{z}_I^{\bar{K}} \varphi_I(\xi),
\end{equation}
where $\{\bar{z}_I^{\bar{K}}\}_{I=1}^n\subset\Rbb^d$ are the nodal coordinates of 
the high-order element $\bar{K}\in\bar\Ecal_h$ in the reference domain,
$\{\varphi_I\}_{I=1}^n$ is a Lagrangian basis of $\Pcal_q(\Omega_\square)$,
and $n = \dim\Pcal_q(\Omega_\square)$. Composition of the isoparametric
mapping $\Ycal_{\bar{K}}$ with the domain mapping $\Gcal$ maps the parent element
$\Omega_\square$ to the physical domain. This setup is illustrated in
Figure~\ref{fig:meshop}.
\begin{remark}
For brevity, we use \textit{isoparametric mapping} to refer to
the mapping from parent domain ($\Omega_\square$) to reference domain element
($\bar{K}$), regardless of the relationship between $p$ (solution polynomial degree)
and $q$ (domain polynomial degree).
\end{remark}

Next, let $\Pi \in 2^{\Omega_\square}$ define a decomposition
of the parent element into subdomains such that
$\bigcup_{Z \in \Pi} Z = \Omega_\square$ and
$Z \cap Z' = \emptyset$ for $Z,Z'\in\Pi$ with $Z \neq Z'$,
and let
$\bar\Ecal_h^\mathrm{r}\subset\bar\Ecal_h$ denote the subset of
elements to be partitioned. Then a new mesh $\bar\Ecal_h'$ of
$\bar\Omega$ is obtained from the original mesh $\bar\Ecal_h$ by
composing the isoparametric mapping for each element in
$\bar\Ecal_h^\mathrm{r}$ with the refinement operator
\begin{equation} \label{eqn:refine_refmsh}
\bar\Ecal_h' =
 \left\{
  \Ycal_{\bar{K}}(Z) \suchthat \bar{K} \in \bar\Ecal_h^\mathrm{r},~Z \in \Pi
 \right\} \cup (\bar\Ecal_h\setminus\bar\Ecal_h^\mathrm{r}).
\end{equation}
Furthermore, by mapping each element of this new mesh $\bar\Ecal_h'$, a
mesh of the physical domain $\Omega$ is defined as
\begin{equation}
\Ecal_h' =
 \left\{
  \Gcal(\bar{K}') \suchthat \bar{K}' \in \bar\Ecal_h'
 \right\}.
\end{equation}
This process is illustrated in Figure~\ref{fig:meshop}.

If the original mesh $\bar\Ecal_h$ is the simplicial prism mesh produced
in Section~\ref{sec:slab:extrude}, $\bar\Ecal_h^\mathrm{r}=\bar\Ecal_h$, and $\Pi$ is
an operator that decomposes $\mathcal{SP}_0^d$ into simplices
(Section~\ref{sec:slab:split:pr2smplx}),
the resulting mesh $\bar\Ecal_h$ will be a simplex-only mesh of $\bar\Omega$.
On the other hand, if $\bar\Ecal_h$ is a simplex-only mesh of $\bar\Omega$ and
$\Pi$ is an operator that decomposes $\mathcal{K}_0^d$ into smaller simplices,
the resuting mesh $\bar\Ecal_h$ will be a refined simplex-only mesh of $\bar\Omega$.
However, in both of these cases, $\bar\Ecal_h'$ is not guaranteed to be conforming
even if $\bar\Ecal_h$ is a conforming mesh. In the remainder of this section,
we discuss strategies to ensure the new mesh $\bar\Ecal_h'$ is conforming,
which will in turn guarantee the physical mesh $\Ecal_h'$ is conforming
because the mappings $\Gcal\in[\Wcal_h]^d$ are continuous.

\begin{remark}
The expression for the refined mesh in (\ref{eqn:refine_refmsh}) is a simplified
version of the algorithm used in practice because it only uses a single parent element
decomposition $\Pi$. We usually introduce several $\Omega_\square$ decompositions,
which are used selectively to ensure a conforming mesh is produced. This can be
written as
\begin{equation} \label{eqn:refine_refmsh2}
\bar\Ecal_h' =
 \left\{
  \Ycal_{\bar{K}}(Z) \suchthat \bar{K} \in \bar\Ecal_h^\mathrm{r},~Z \in \Pi_{\bar{K}}
 \right\} \cup (\bar\Ecal_h\setminus\bar\Ecal_h^\mathrm{r}),
\end{equation}
where $\Pi_{\bar{K}}\in 2^{\Omega_\square}$ is the decomposition used for element
$\bar{K}\in\bar\Ecal_h$. Details on construction of conforming refinements are
included in Sections~\ref{sec:slab:split:pr2smplx}-\ref{sec:slab:split:smplx2smplx}.
\end{remark}

\subsubsection{Prism-to-simplex refinement}
\label{sec:slab:split:pr2smplx}
\begin{figure}
\centering
\includegraphics[width=0.2\textwidth]{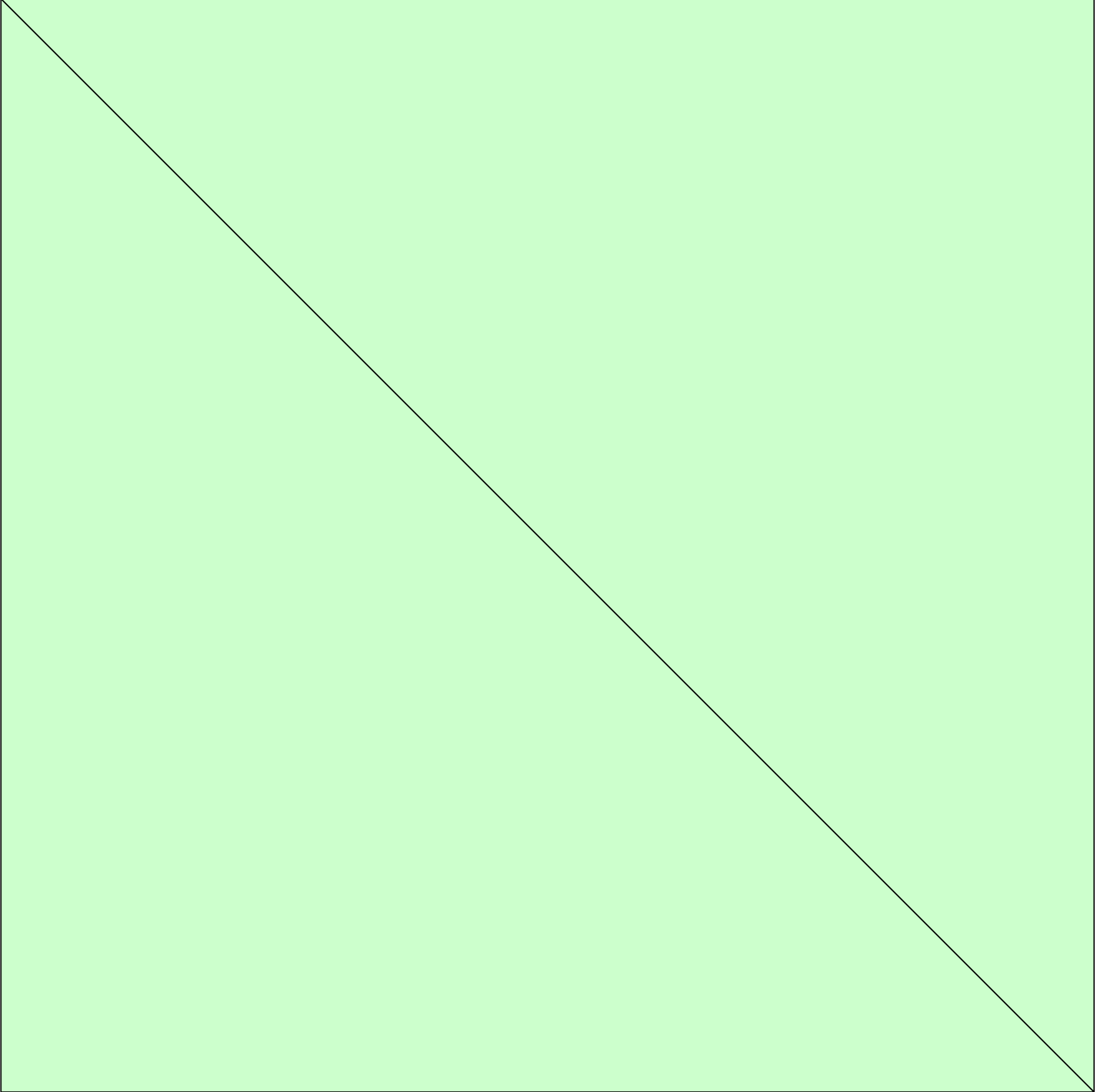} \qquad 
\includegraphics[width=0.2\textwidth]{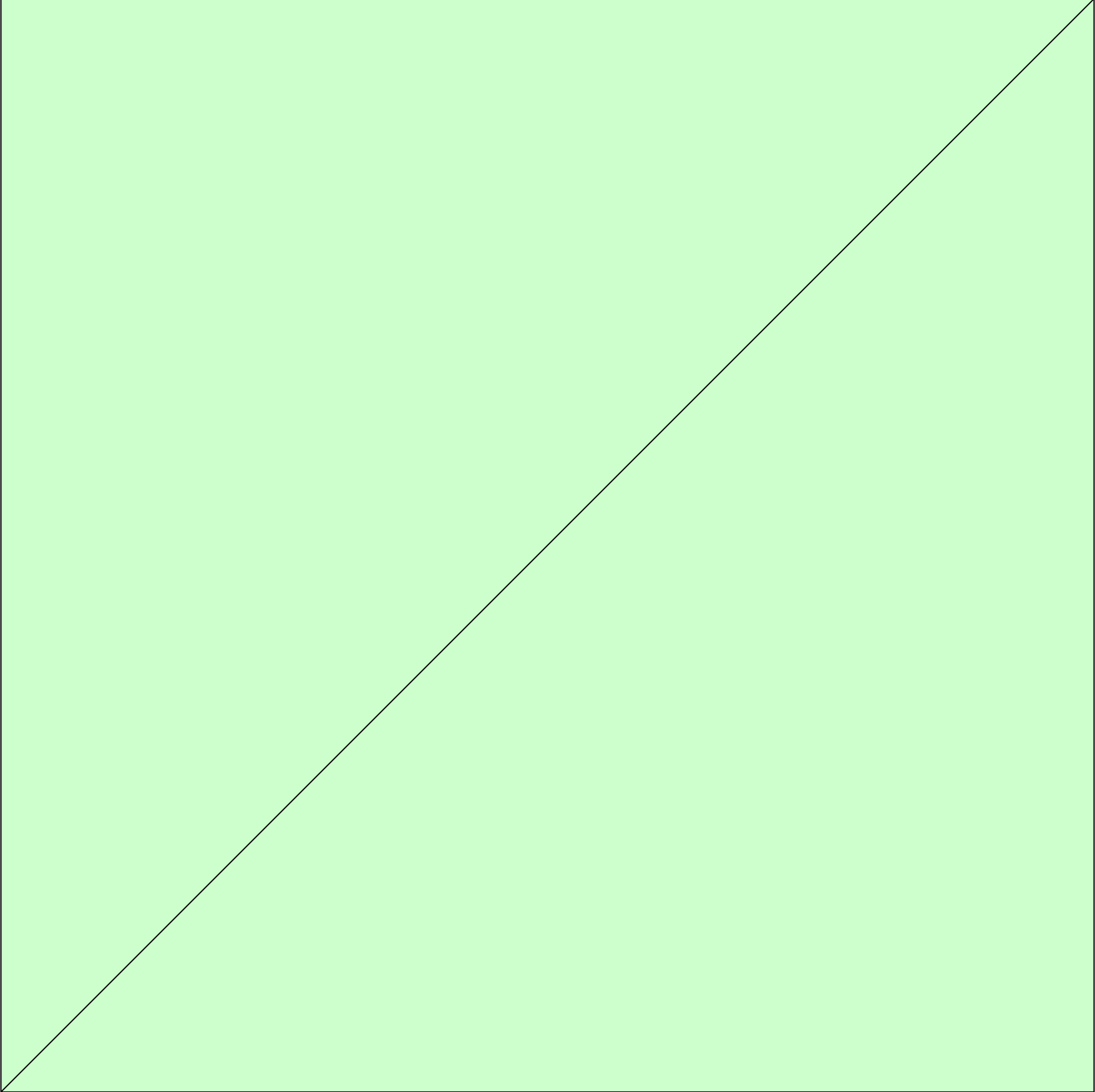} \qquad 
\includegraphics[width=0.2\textwidth]{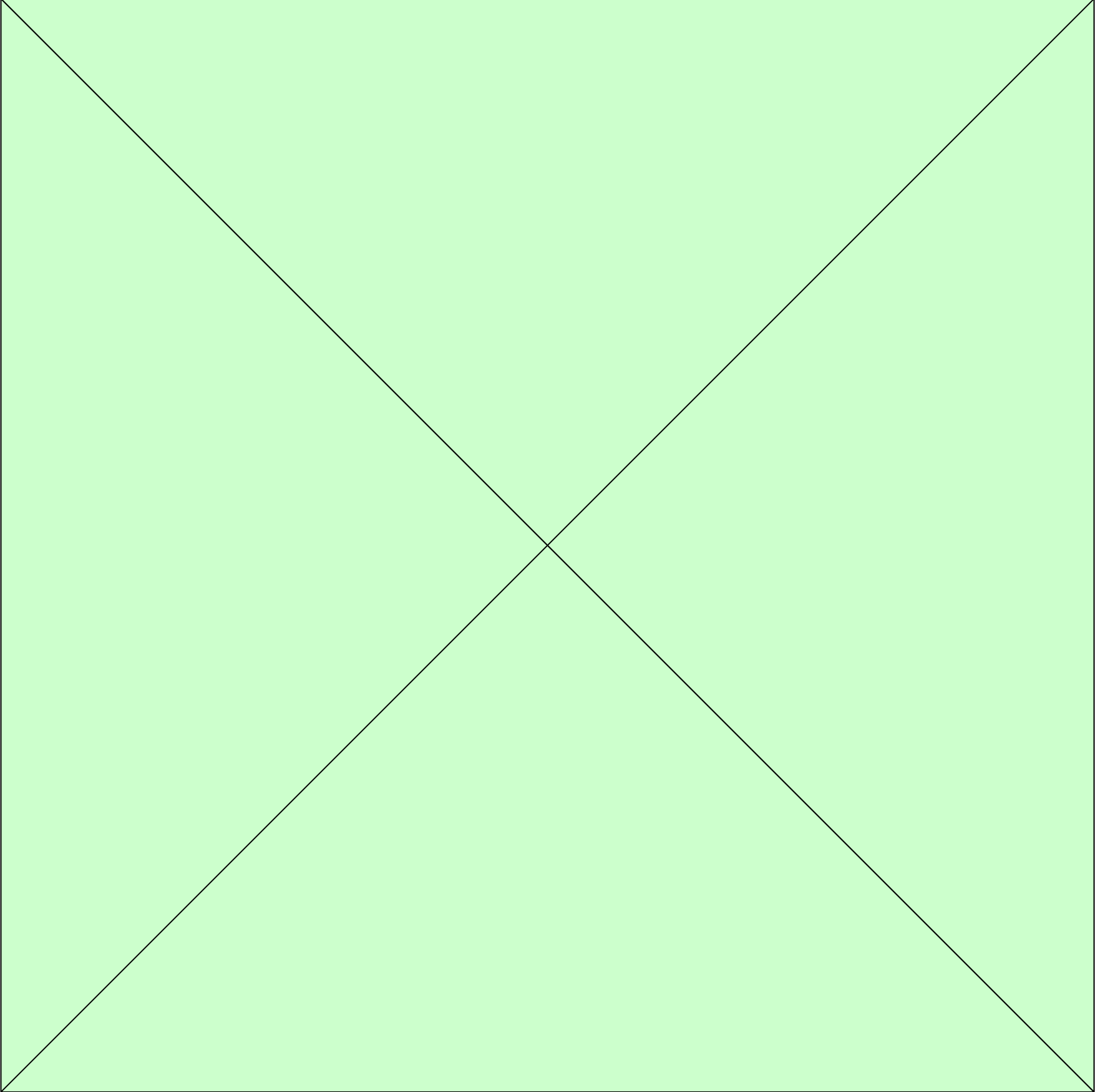}
\caption{Conforming decompositions ($d=2$) of a quadrilateral into triangles using
 a single diagonal (\textit{left}, \textit{center}), and both diagonals (\textit{right}).}
\label{fig:split_prism2d}
\end{figure}

In the $d = 2$ case, a simplicial prism is simply a
quadrilateral, which can be split with one or both
diagonal lines (Figure~\ref{fig:split_prism2d}).
Any of these splittings can be applied to
any subset of elements of the simplicial mesh,
regardless of element orientation, and the
resulting simplicial mesh is guaranteed to
be conforming.

\begin{figure}
\centering
\includegraphics[width=0.16\textwidth]{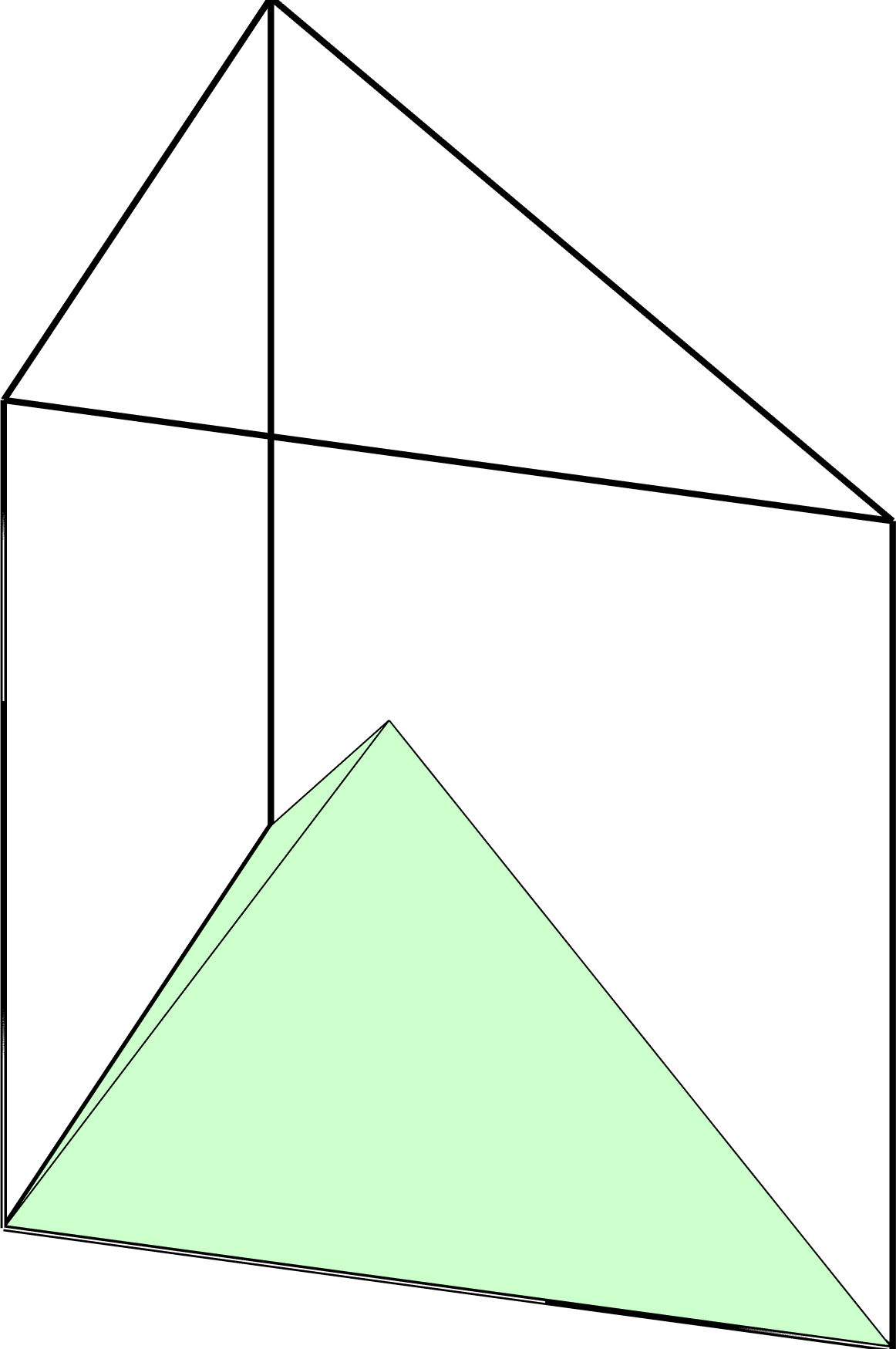} \,
\includegraphics[width=0.16\textwidth]{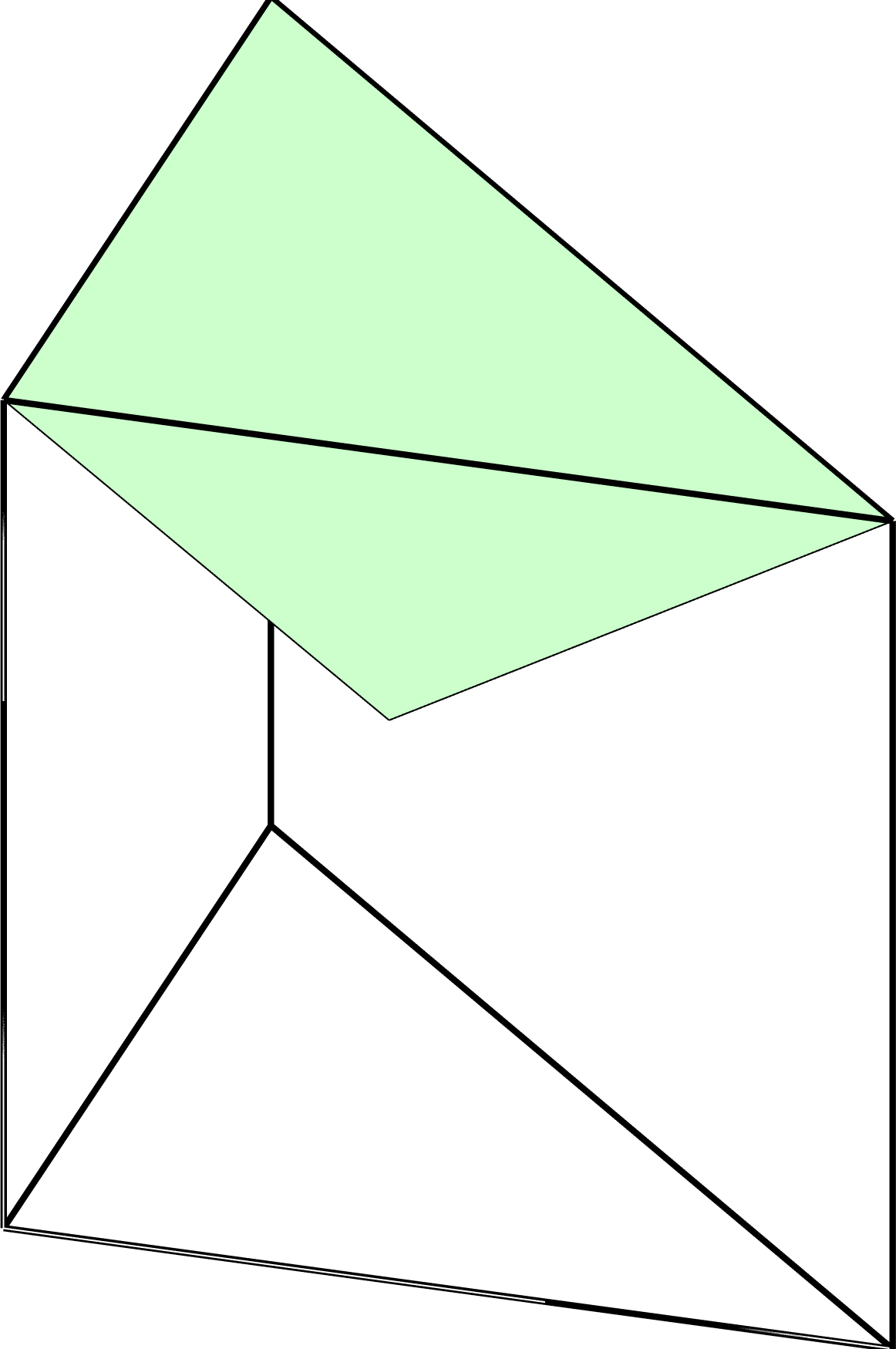} \,
\includegraphics[width=0.16\textwidth]{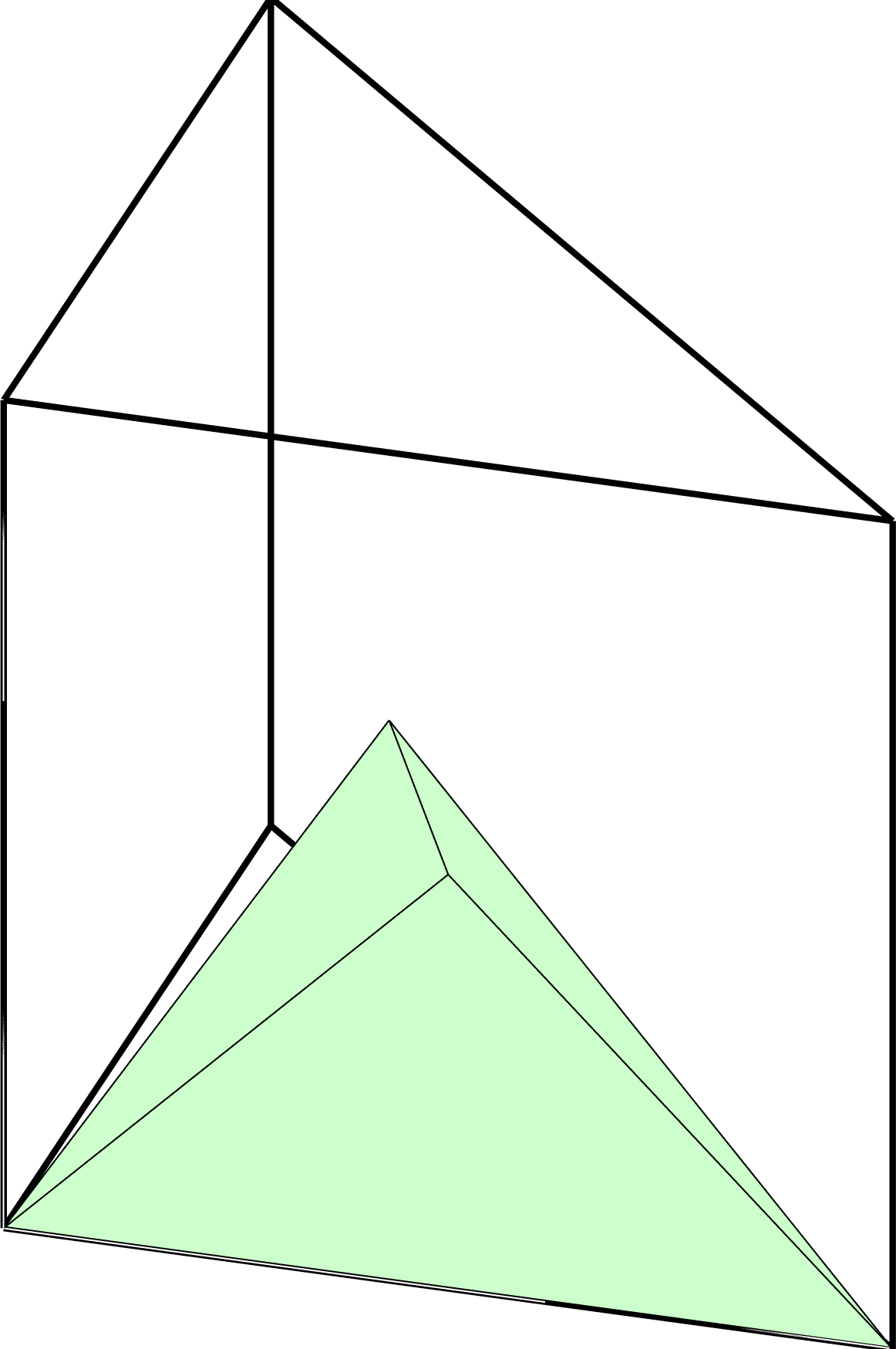} \,
\includegraphics[width=0.16\textwidth]{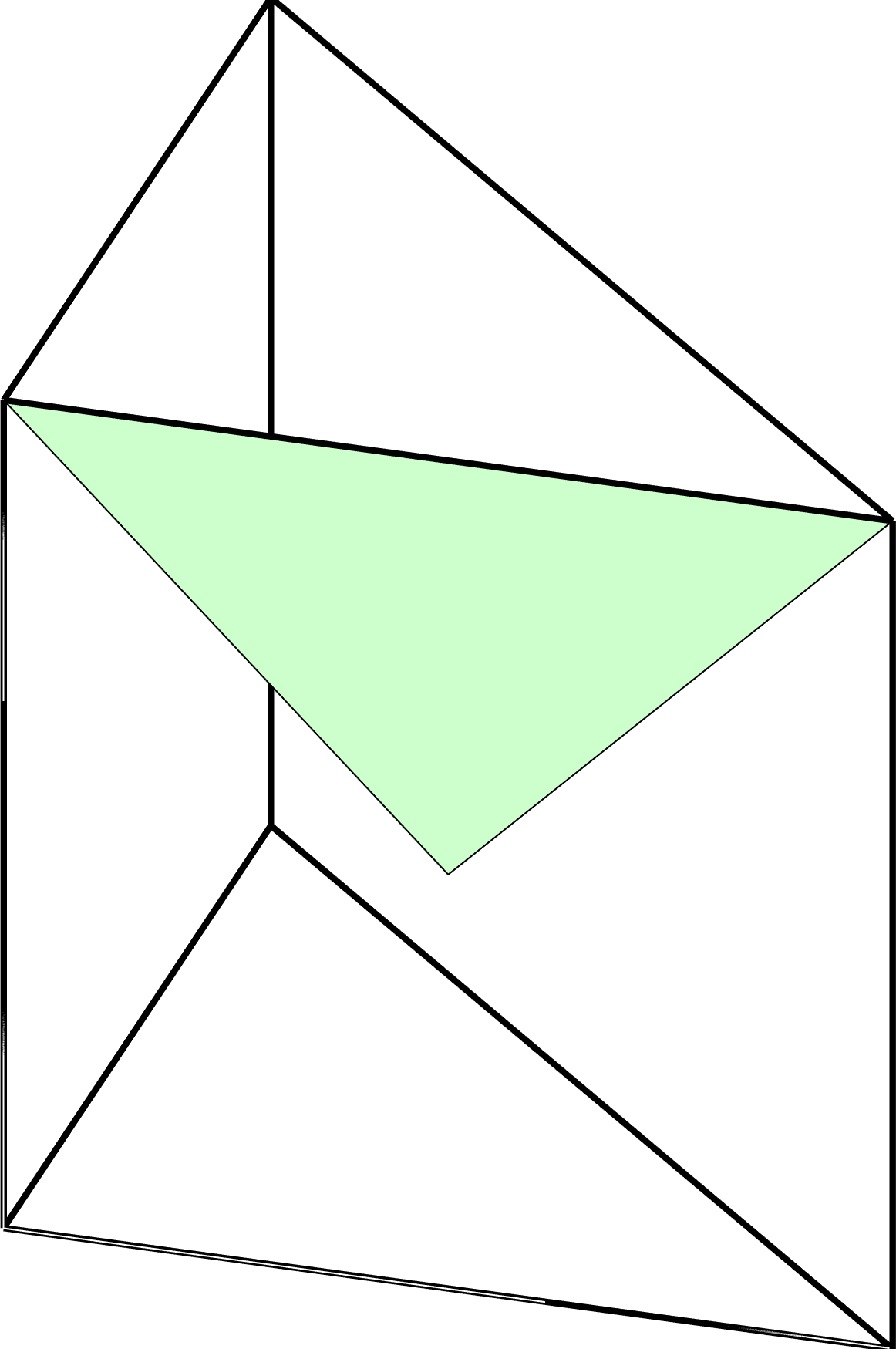} \,
\includegraphics[width=0.16\textwidth]{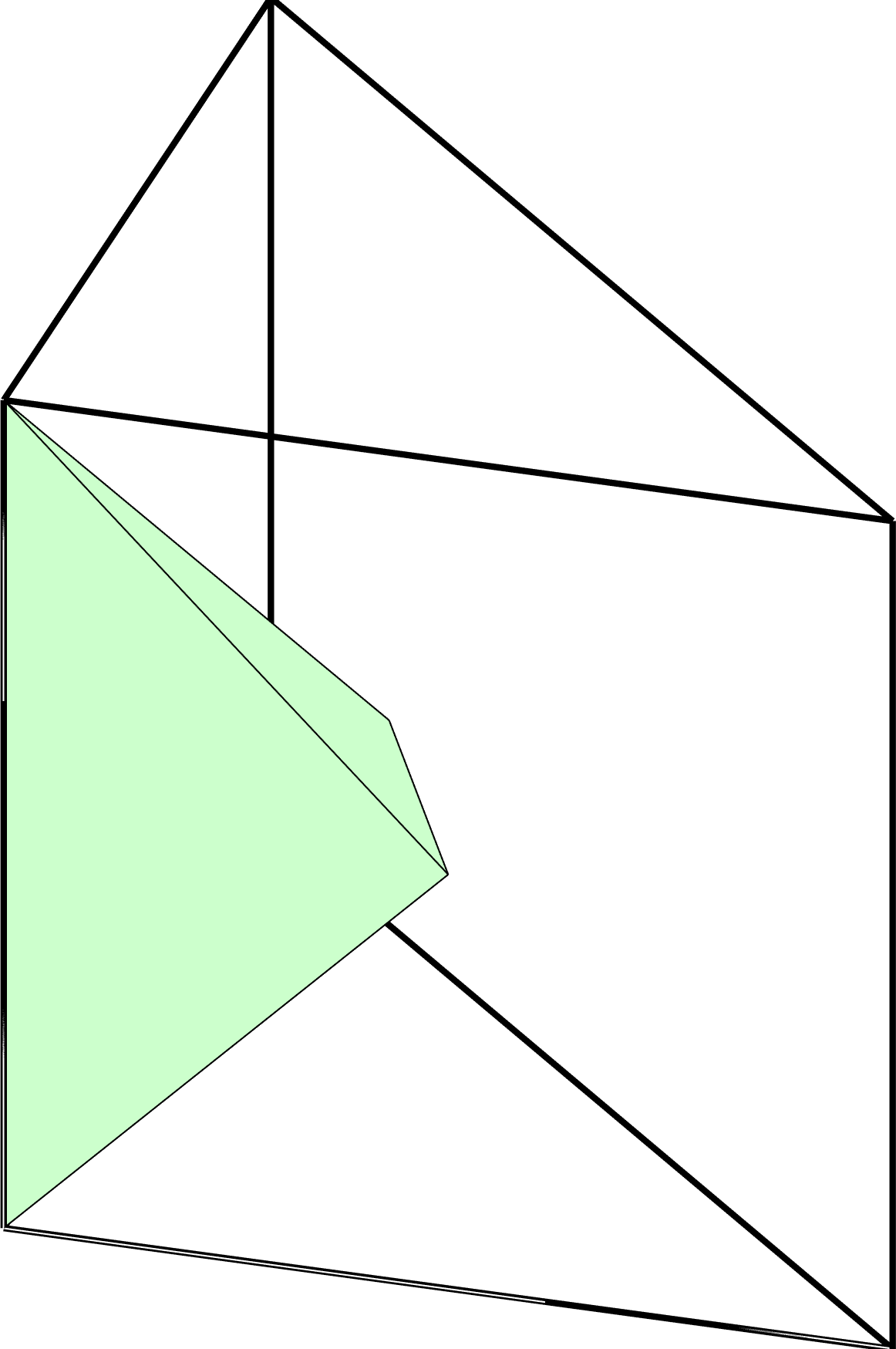} \\
\includegraphics[width=0.16\textwidth]{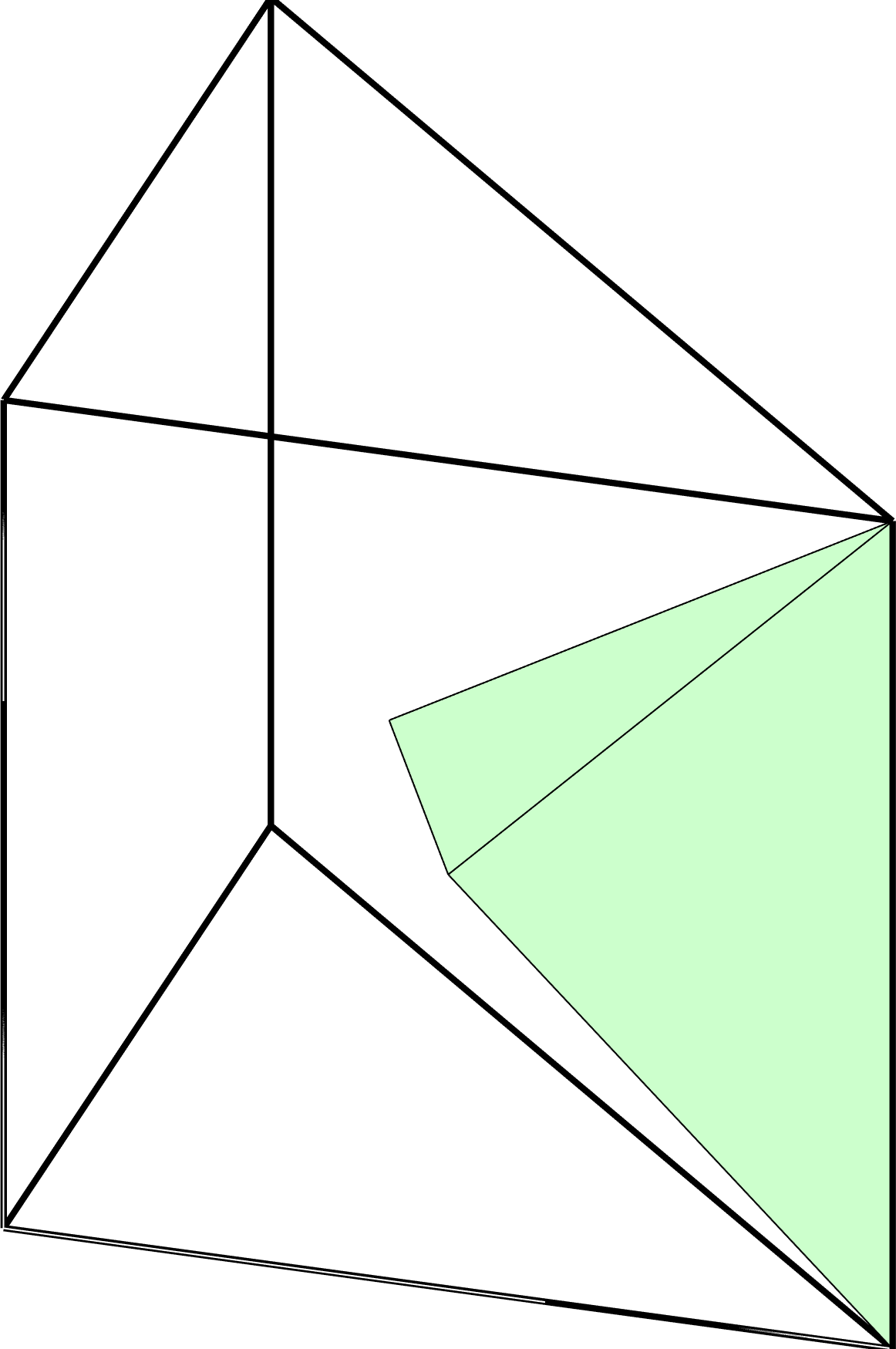} \,
\includegraphics[width=0.16\textwidth]{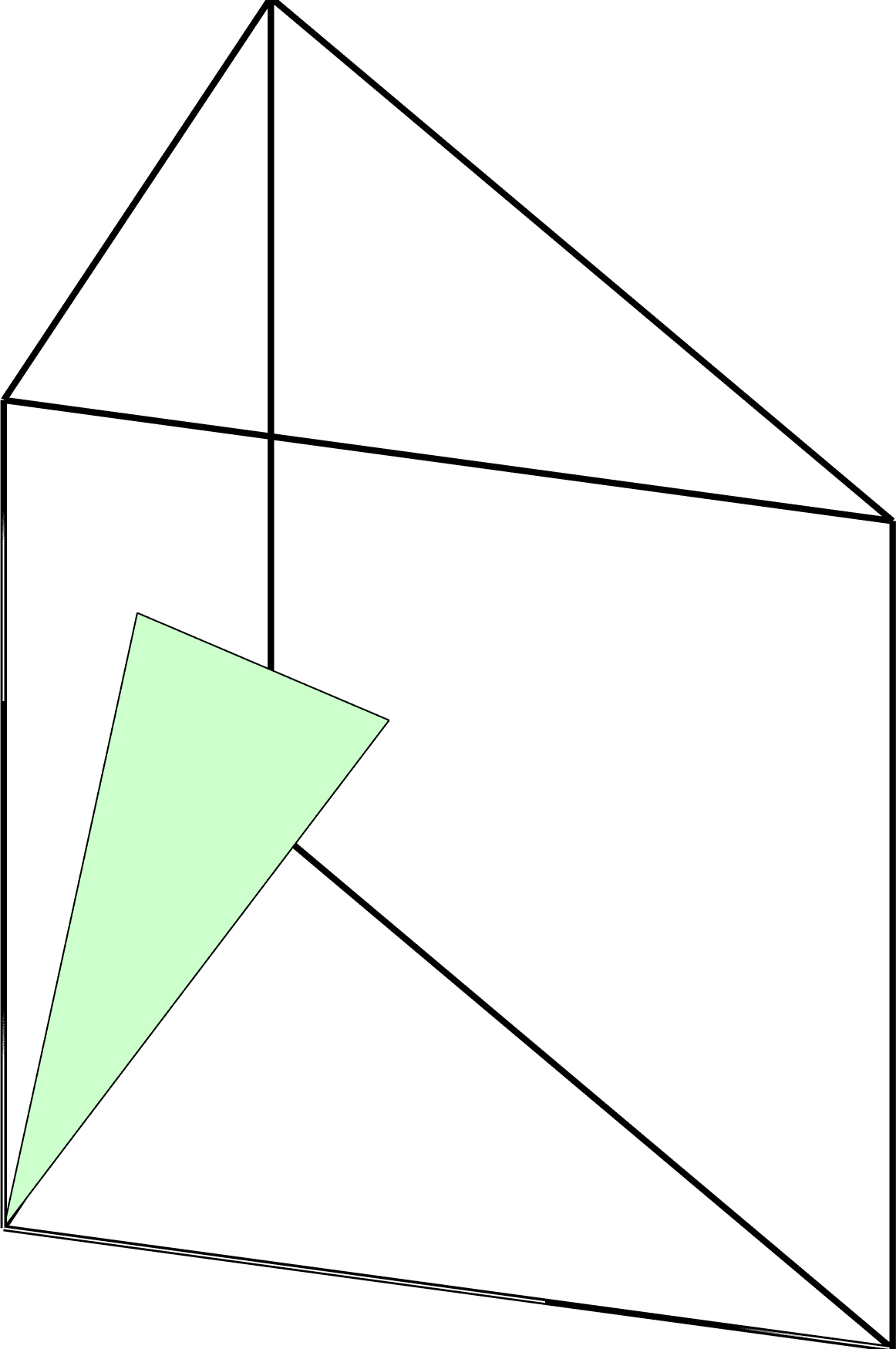} \,
\includegraphics[width=0.16\textwidth]{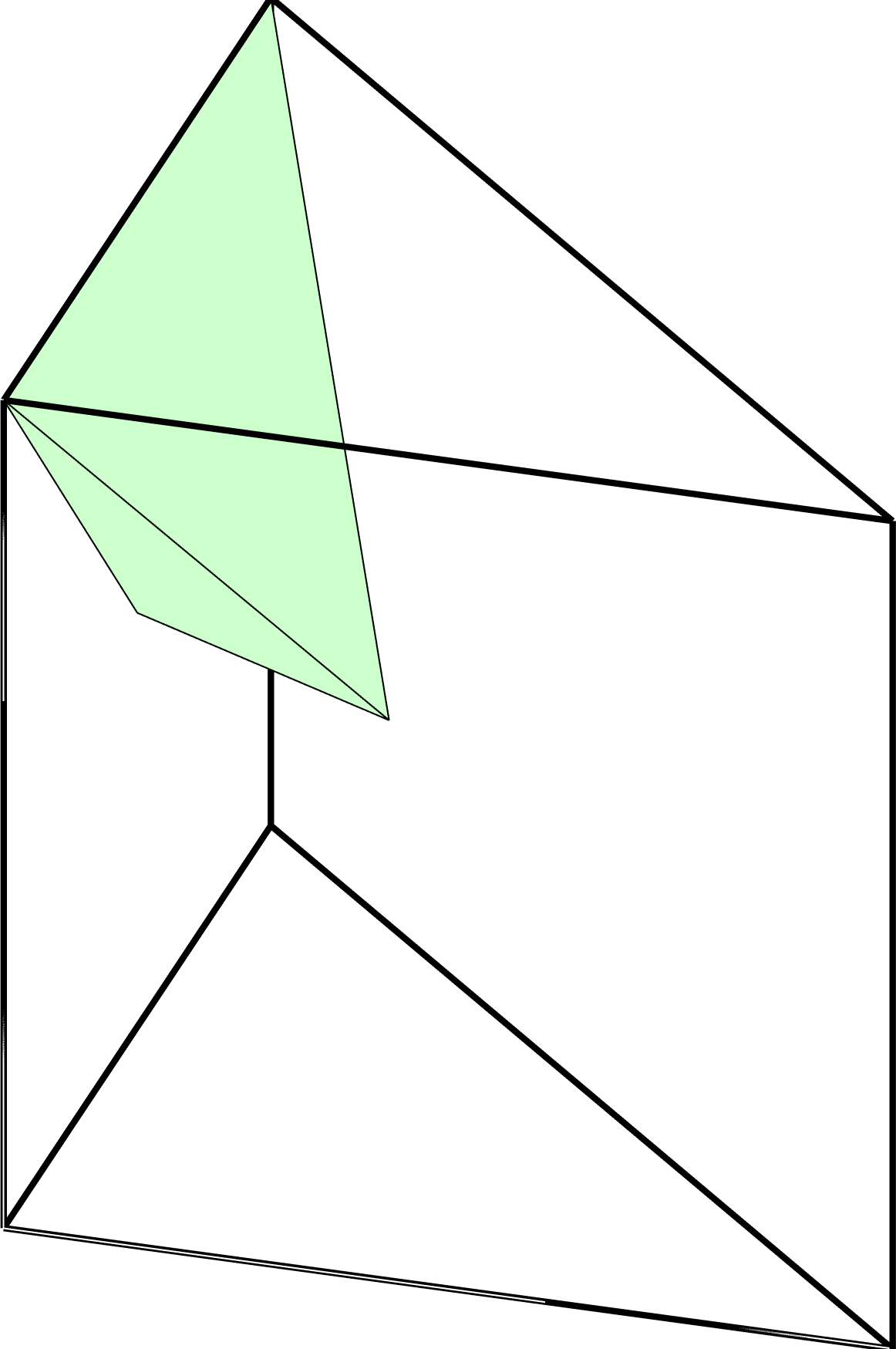} \,
\includegraphics[width=0.16\textwidth]{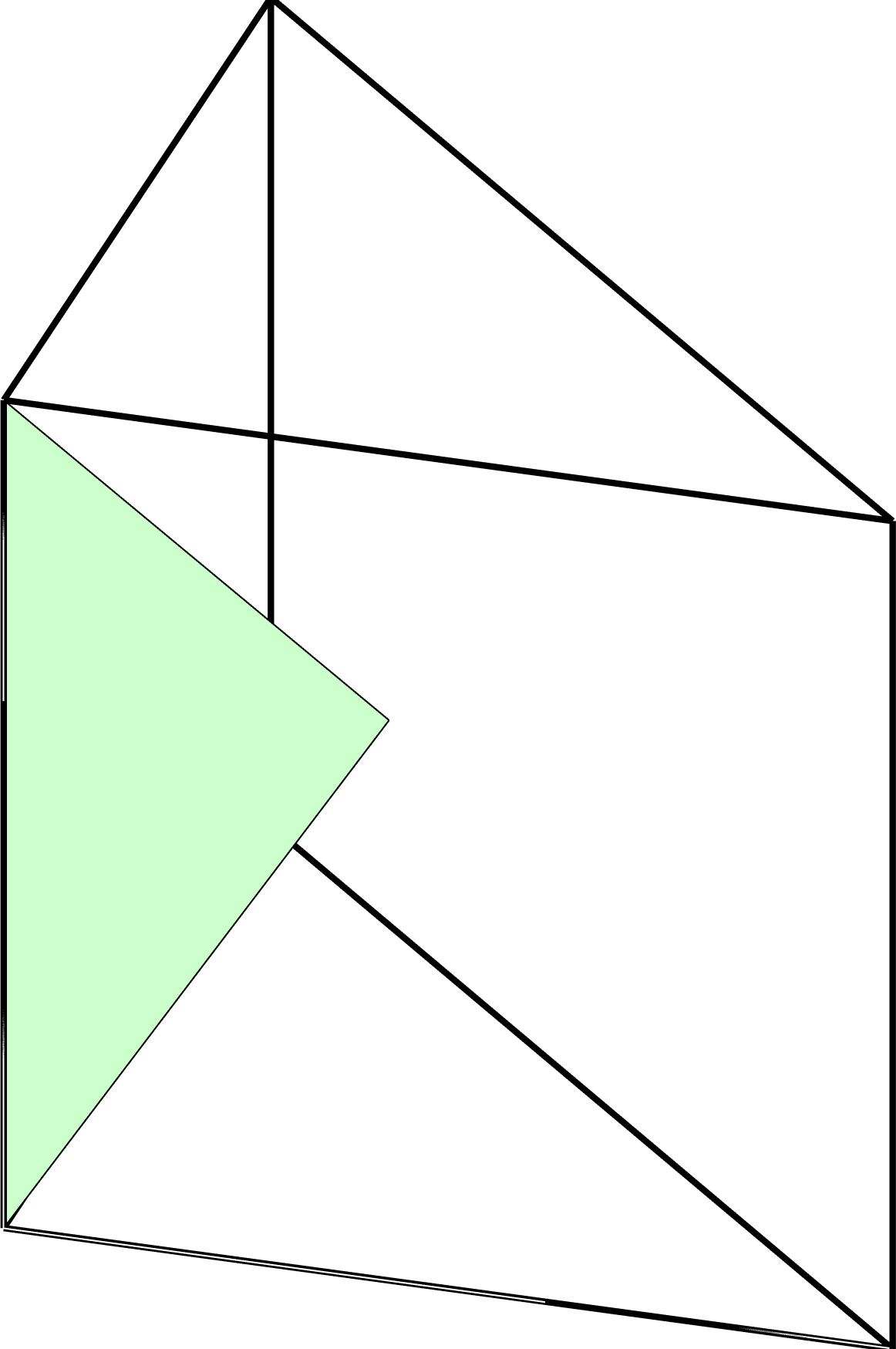} \,
\includegraphics[width=0.16\textwidth]{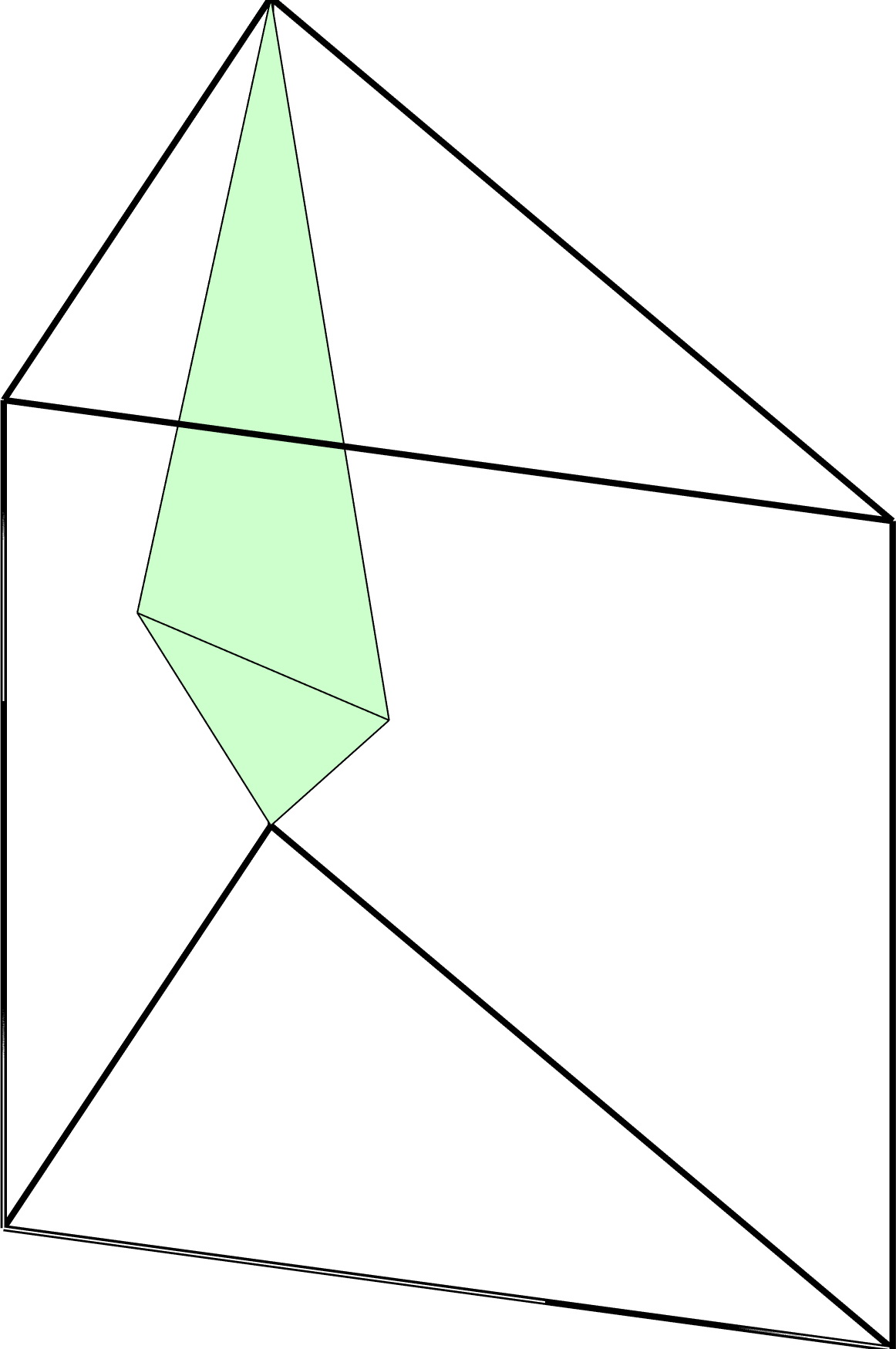} \\
\includegraphics[width=0.16\textwidth]{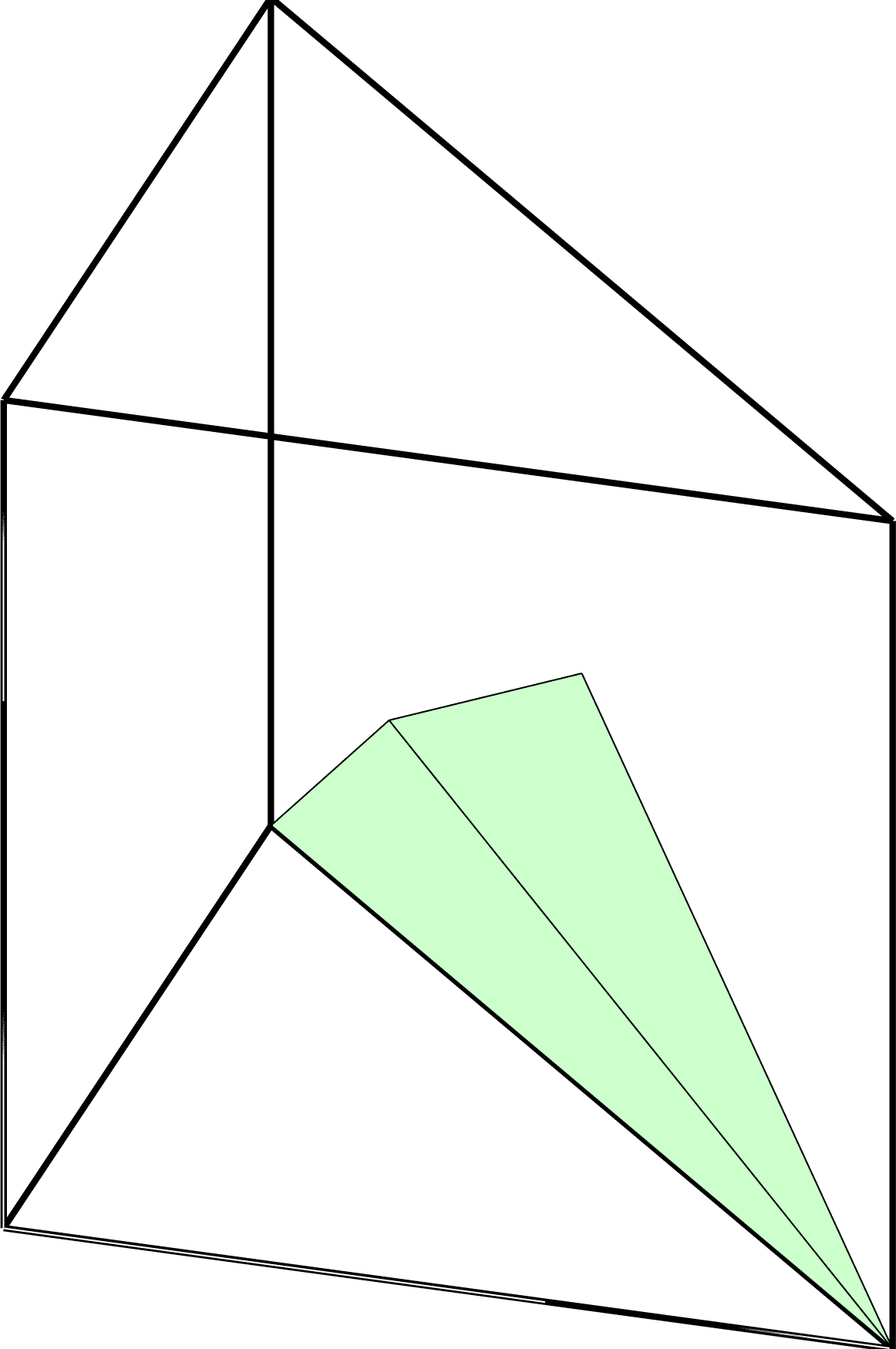} \,
\includegraphics[width=0.16\textwidth]{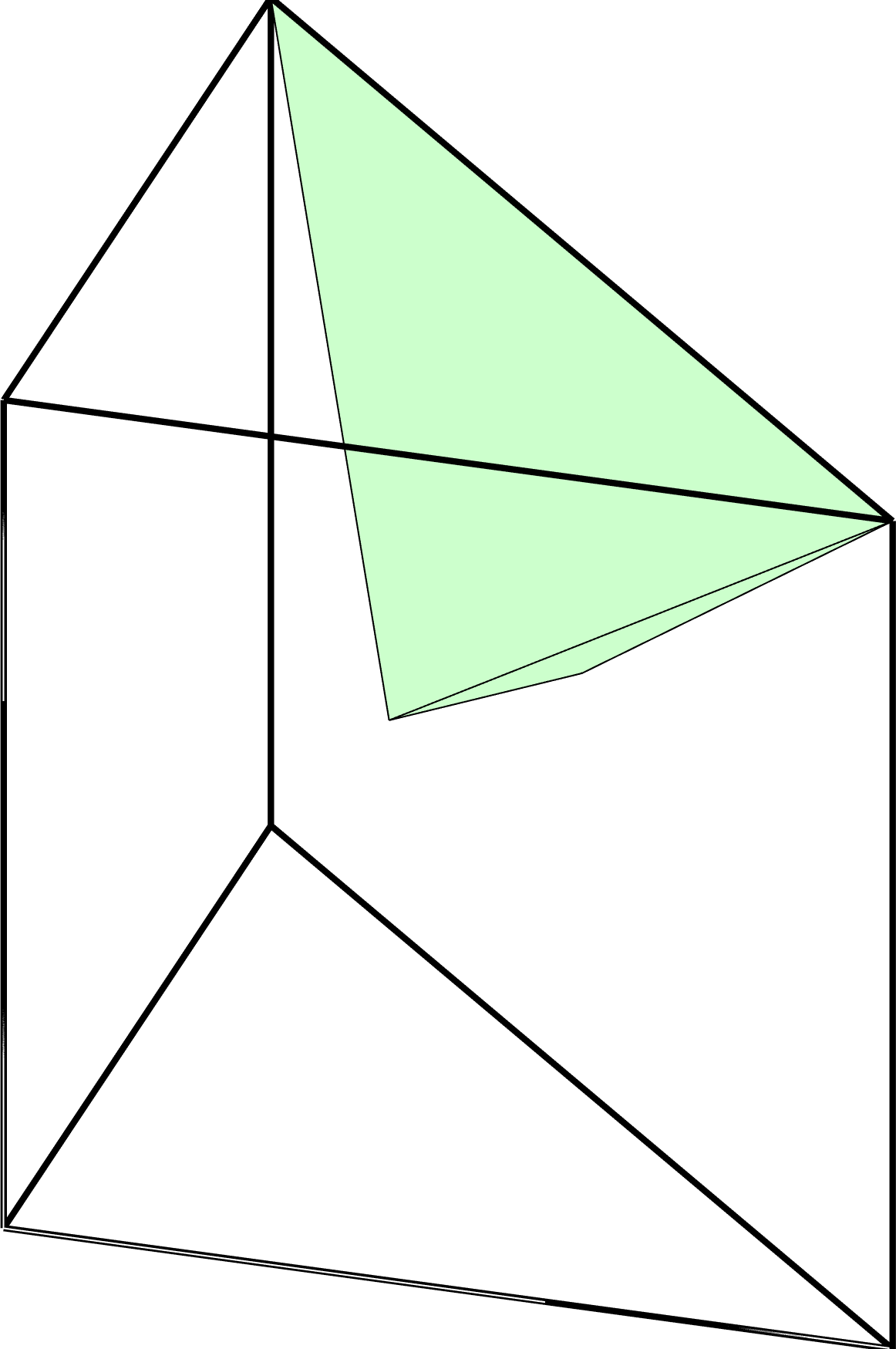} \,
\includegraphics[width=0.16\textwidth]{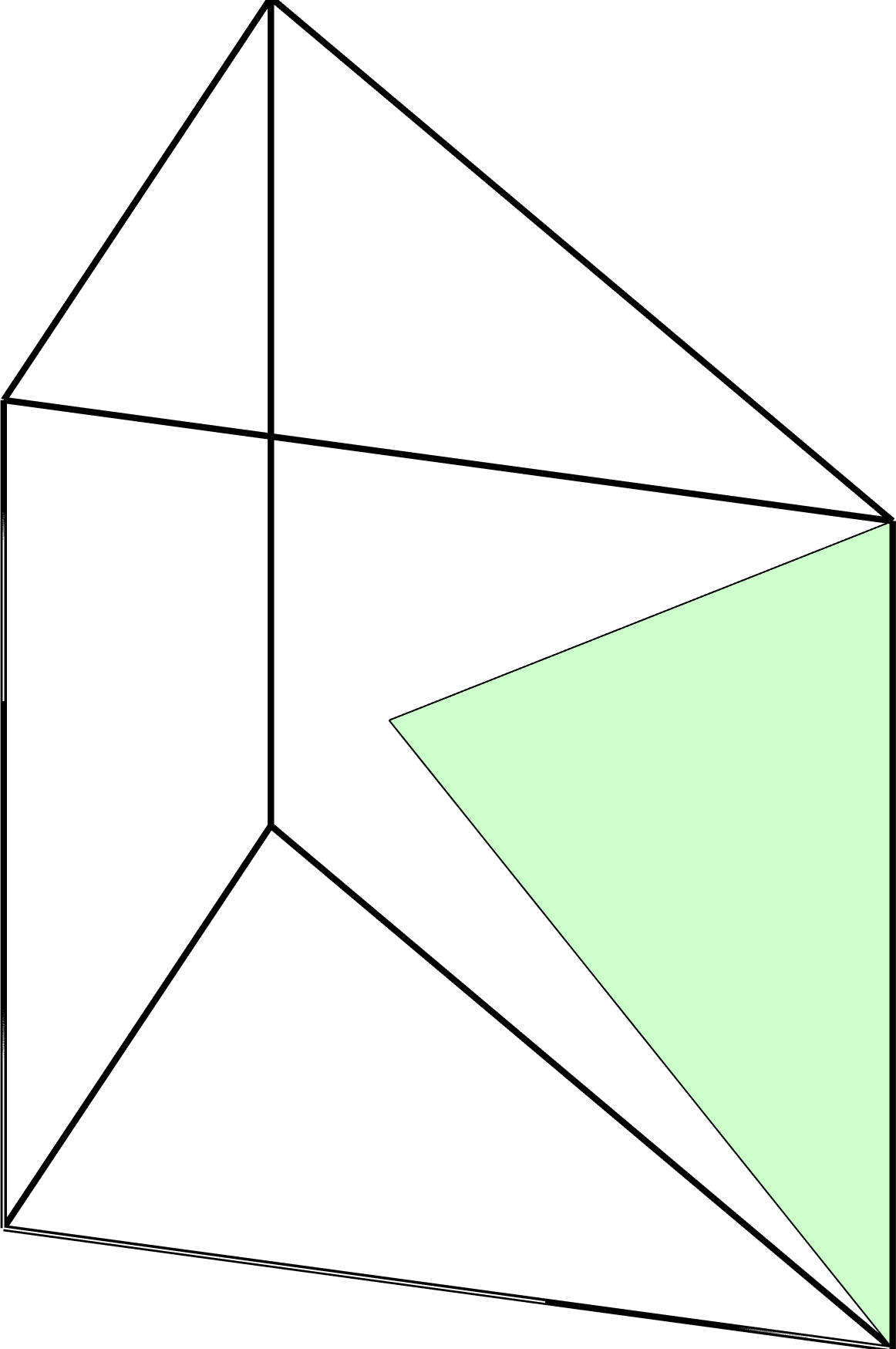} \,
\includegraphics[width=0.16\textwidth]{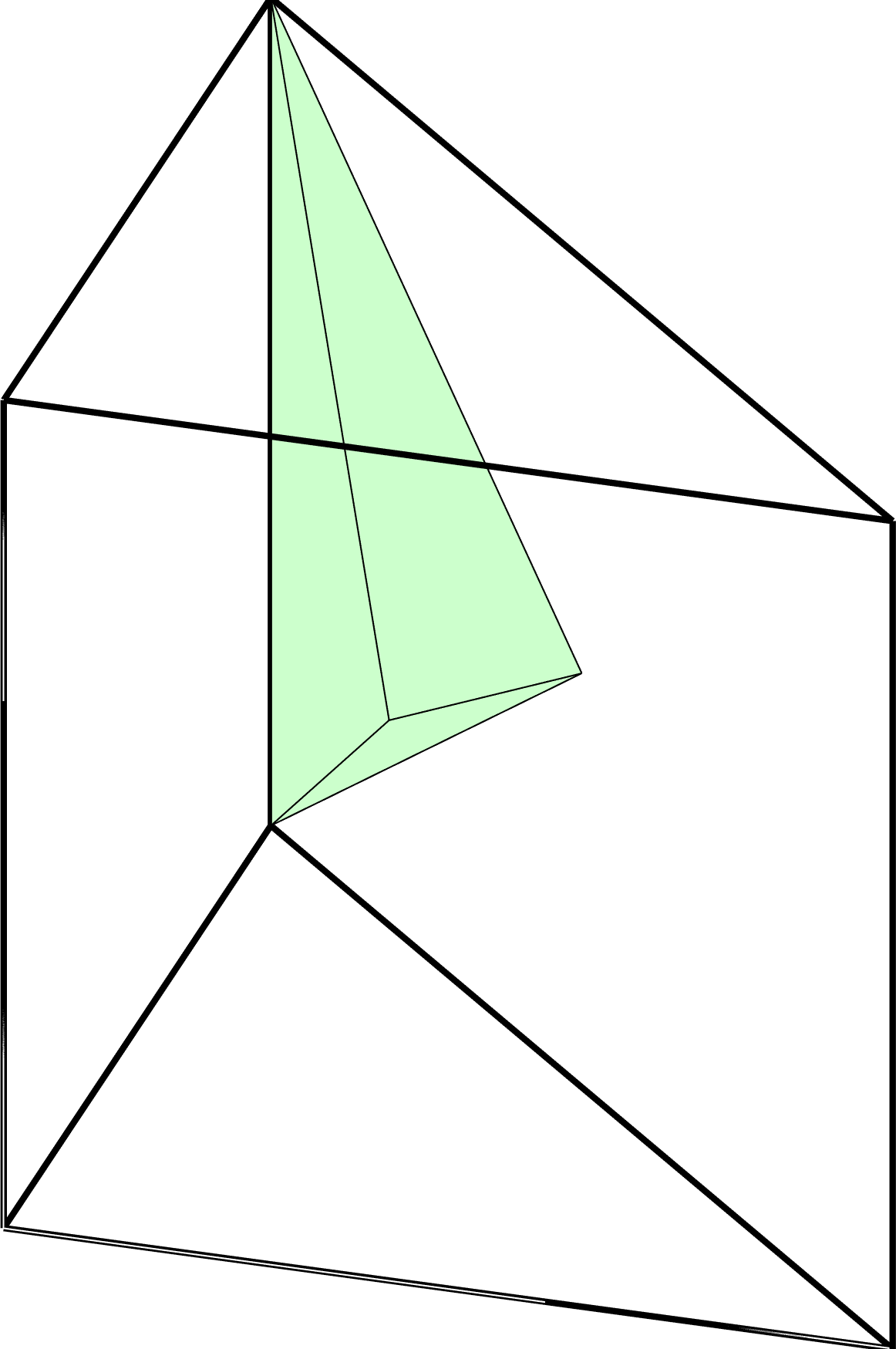} \hspace{2.5cm}
\caption{Conforming decomposition ($d=3$) of a triangular prism into $N_3 = 14$
 tetrahedra using a $d=2$ conforming splitting of the quadrilateral faces and connecting
 all face/edge nodes to a central node.}
\label{fig:split_prism3d}
\end{figure}

In the $d > 3$ case, $\mathcal{SP}_0^d$ contains $d+2$
faces: two simplex ``caps`` and $d$ simplicial prism
``sides''. We aim to define a decomposition ($\Pi$) of
$\mathcal{SP}_0^d$ such that if it is applied to all
elements of the prismatic mesh via (\ref{eqn:refine_refmsh}),
the resulting simplex-only mesh is guaranteed to
be conforming without requiring global mesh operations.
This is accomplished by decomposing $\mathcal{SP}_0^d$
in such a way that the decomposition of its $d+2$
faces into $(d-1)$-dimensional simplices is independent
of the orientation of the element. Any decomposition
of $\mathcal{SP}_0^d$ that guarantees such a
decomposition of its faces will trivially
lead to a conforming simplex-only mesh.
A dimension-independent approach based on
this concept first applies an orientation-independent
decomposition to partition each face into $(d-1)$-dimensional
simplices. Then the vertices of each $(d-1)$-dimensional
simplex are connected to a central node in the $d$-dimensional
prism to form $d$-dimensional simplices. Assuming the simplicial
caps are not refined, this procedure splits $\mathcal{SP}_0^d$
into $N_d$ simplices, where
\begin{equation}
 N_d = dN_{d-1} + 2
\end{equation}
for $d > 2$ and $N_2 = 4$, which leads to $N_3 = 14$ and $N_4 = 58$
($d=3$ case in Figure~\ref{fig:split_prism3d}).
Alternatively, a constrained Delaunay triangulation can be used to
form $N_3 = 11$ $d$-dimensional simplices without introducing a central
node. While these approaches are general and allow each element to be split
independently of its neighbors, it produces large meshes in $d > 2$.

\begin{remark}
In this work, we use the approach that splits the simplical prism into
$N_2 = 4$ and $N_3 = 14$ simplices instead of the more efficient approaches,
i.e., split with a single diagonal in $d = 2$ ($N_2 = 2$) and split without a
central node in $d = 3$ ($N_3 = 11$). Despite the creation of more simplices,
we have observed the chosen approach improves the robustness of shock fitting.
This is attributed to the fact that the chosen splitting is essentially uniform
refinement and, as such, regardless of the trajectory of discontinuities,
there will be a face that provides a reasonable initial guess to its space-time
trajectory. This could also be achieved with a shock-dependent splitting, which
was not pursued in this work.
\end{remark}

\begin{remark}
In the $d=3$ case, there are several ways to split $\mathcal{SP}_0^d$ into $N_3=4$
simplices although each element cannot be split independently. Rather, each prism
must be split into four simplices based on the partition of the neighbors to ensure
the diagonals of each face conform, which requires global mesh operations
\cite{2015_wang_sptm}.
\end{remark}


\subsubsection{Simplex-to-simplex refinement}
\label{sec:slab:split:smplx2smplx}
\begin{figure}
\centering
\includegraphics[width=0.2\textwidth]{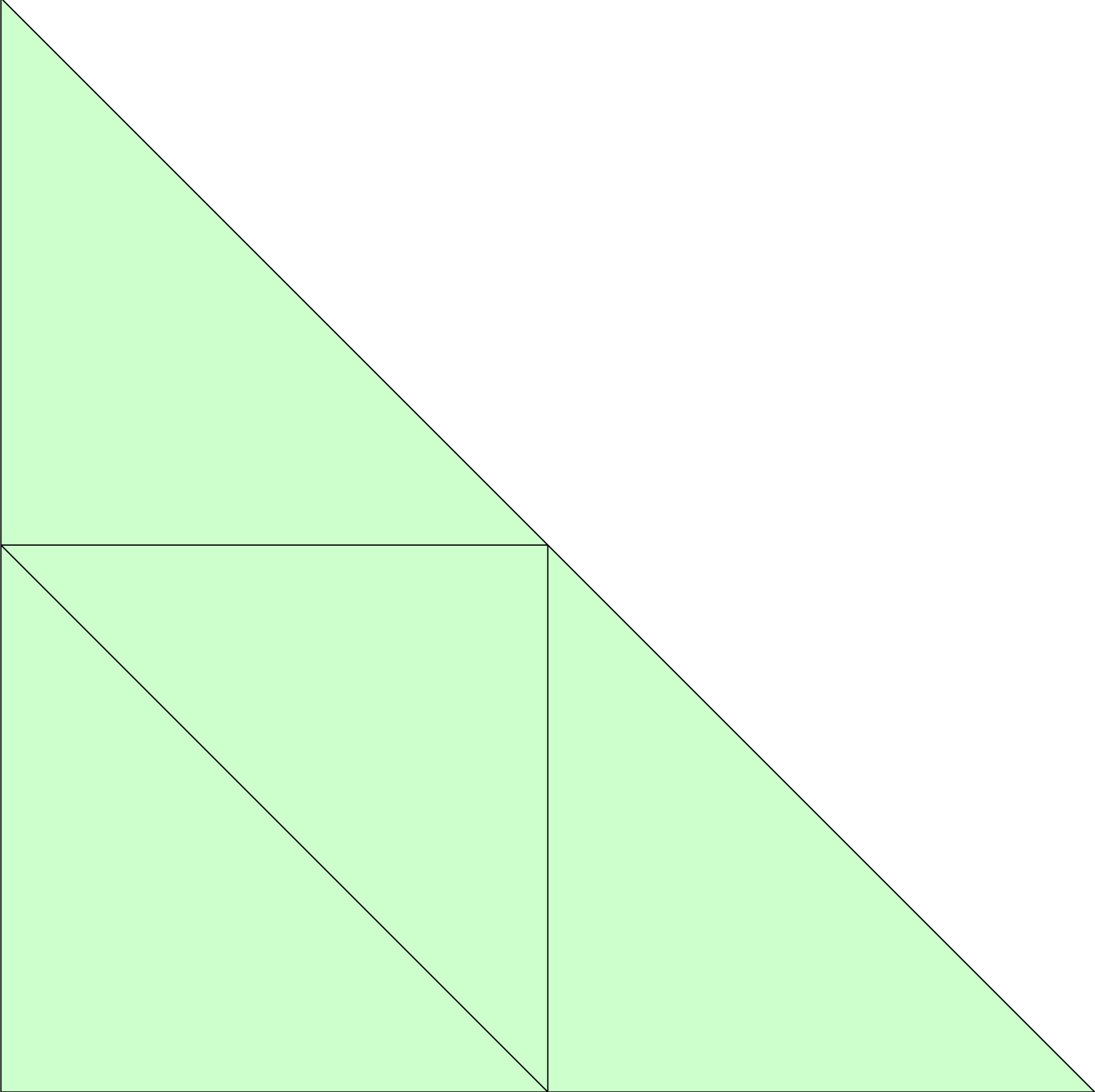} \qquad 
\includegraphics[width=0.2\textwidth]{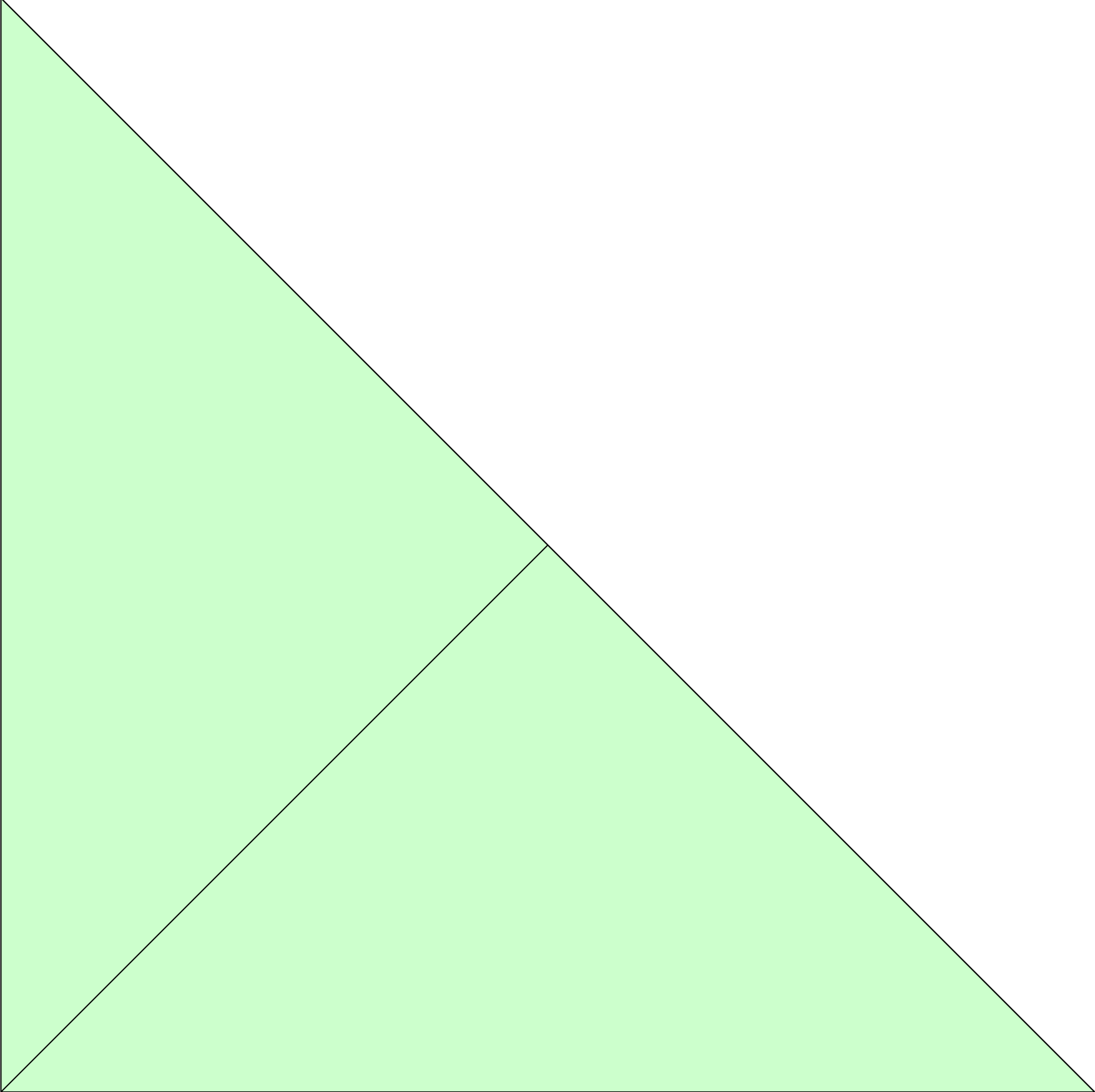} \qquad
\includegraphics[width=0.2\textwidth]{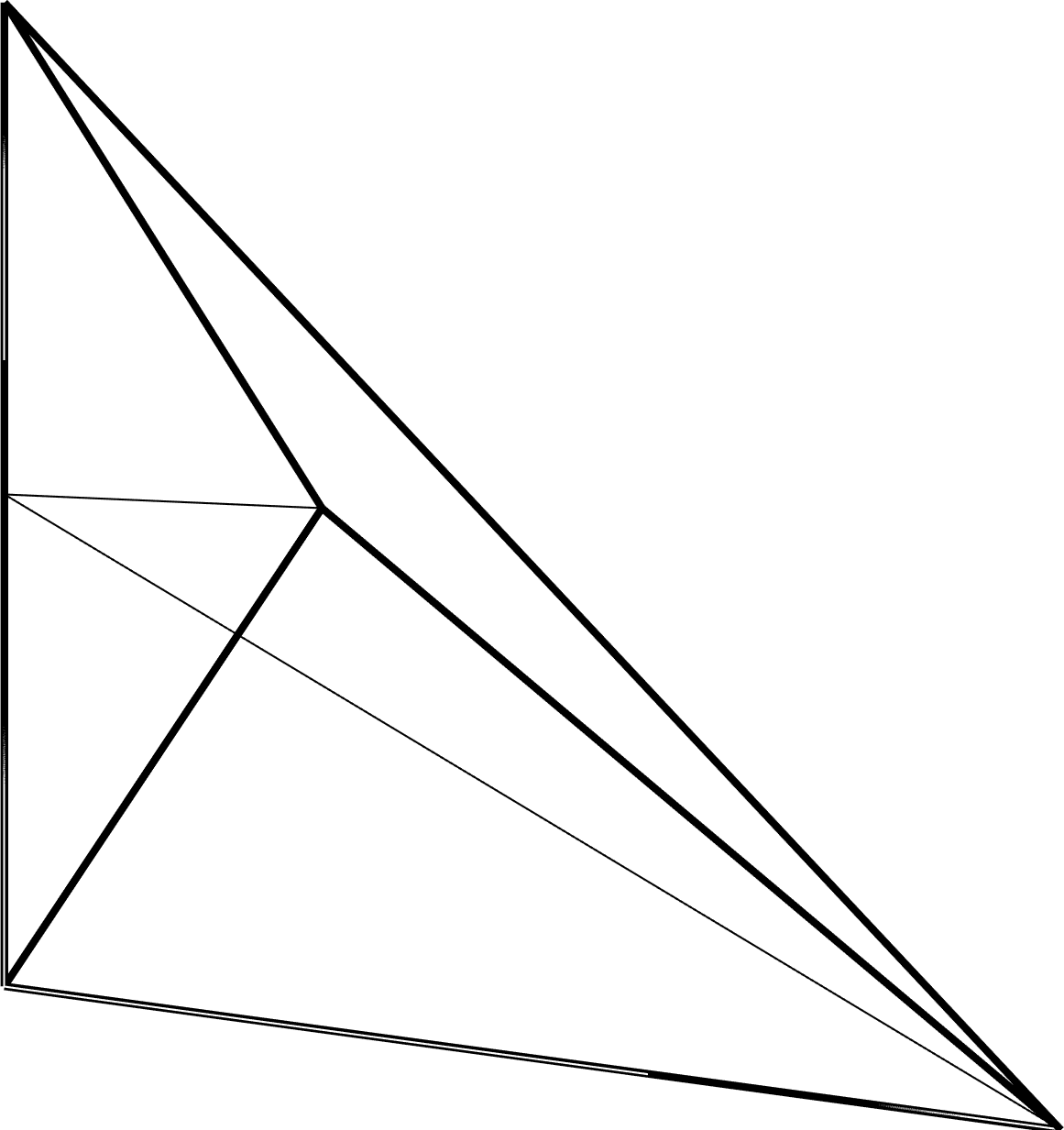}
\caption{Conforming refinement of a simplices using
 uniform refinement in $d = 2$ (\textit{left}),
 edge splitting in $d = 2$ (\textit{center}),
 and edge splitting in $d = 3$ (\textit{right}).}
\label{fig:split_smplx}
\end{figure}

We consider two approaches to refine a simplex into smaller simplices ($\Pi$),
which will be used for adaptive mesh refinement in Section~\ref{sec:ist:amr}. A general,
dimension-independent approach is known as an edge split
\cite{persson2005mesh,lohner2008applied},
whereby a new node is introduced at the midpoint of one edge of the original
simplex and the unique plane that passes through this point and the $d-1$
vertices of the simplex not associated with the original edge is used to
split the simplex into two smaller simplices (Figure~\ref{fig:split_smplx}).
To ensure a conforming mesh, all simplices that share the edge being split
must also be split using the general form of the mesh refinement in
(\ref{eqn:refine_refmsh2}). Another approach that is only practical in
$d = 2$ is to split a triangle uniformly into four smaller triangles
(Figure~\ref{fig:split_smplx}). To ensure the resulting mesh is conforming,
all neighbors of an element that is uniformly refined must also be uniformly
refined or the edge must be split \cite{bank1981adaptive, mitchell201630}.

\section{High-order implicit shock tracking over space-time slabs}
\label{sec:ist}
In this section, we extend the HOIST method \cite{2020_zahr_HOIST,huang2022robust} to
the slab-based space-time setting with special attention
given to the challenges and opportunities inherent in
this setting. In particular, we construct admissible domain mappings
over space-time slabs (Section~\ref{sec:ist:dommap}), formulate the
space-time HOIST method (Section~\ref{sec:ist:optform}), and discuss
solution transfer between space-time slabs (Section~\ref{sec:ist:transfer}).
Additionally, we describe an approach to scale time slabs in the temporal
direction (Section~\ref{sec:ist:scaling}) such that well-conditioned space-time
grids can be obtained, even when the solution contains waves with extreme speeds.
Finally, we describe a simple approach to incorporate adaptive mesh refinement
into the shock tracking solver (Section~\ref{sec:ist:amr}), describe an opportunity
for efficiency that applies boundary conditions directly to shocks
(Section~\ref{sec:ist:shkbc}) and summarize the complete algorithm
(Section~\ref{sec:ist:summ}).

\subsection{Admissible domain mappings}
\label{sec:ist:dommap}
Following the approach in \cite{huang2022robust}, we introduce a parametrization
of the physical nodes
\begin{equation} \label{eqn:y0_trans}
  \phibold : \Rbb^{N_\ybm} \rightarrow \Rbb^{N_\xbm}, \qquad
  \phibold : \ybm \mapsto \phibold(\ybm),
\end{equation}
where $\ybm \in \Rbb^{N_\ybm}$ are unconstrained degrees of freedom (DoFs),
which guarantees the boundaries of the space-time slab are preserved.
Using this parametrization, the unconstrained DoFs can
be moved freely without modifying the computational domain, which
simplifies the central optimization problem in the next section because
boundary-preserving constraints are unnecessary. Approaches to construct
the parametrization $\phibold$ for straight-sided \cite{huang2022robust}
and curved \cite{zahr2020r} boundaries are detailed in previous work.

\begin{figure}
\centering
\begin{tikzpicture}
\begin{groupplot} [
group style={group size = 2 by 1, horizontal sep = 0.05cm, vertical sep = 0.8cm},
title style={at={(current bounding box.north west)}, anchor=west}]
\nextgroupplot[axis equal image, width=0.59\textwidth, xtick={}, ytick={}, xticklabels={}, yticklabels={}, xmin=-0.2, xmax=1, ymin=0, ymax=1]
\addplot []
graphics [xmin=-0.2,xmax=1,ymin=0,ymax=0.7346] { 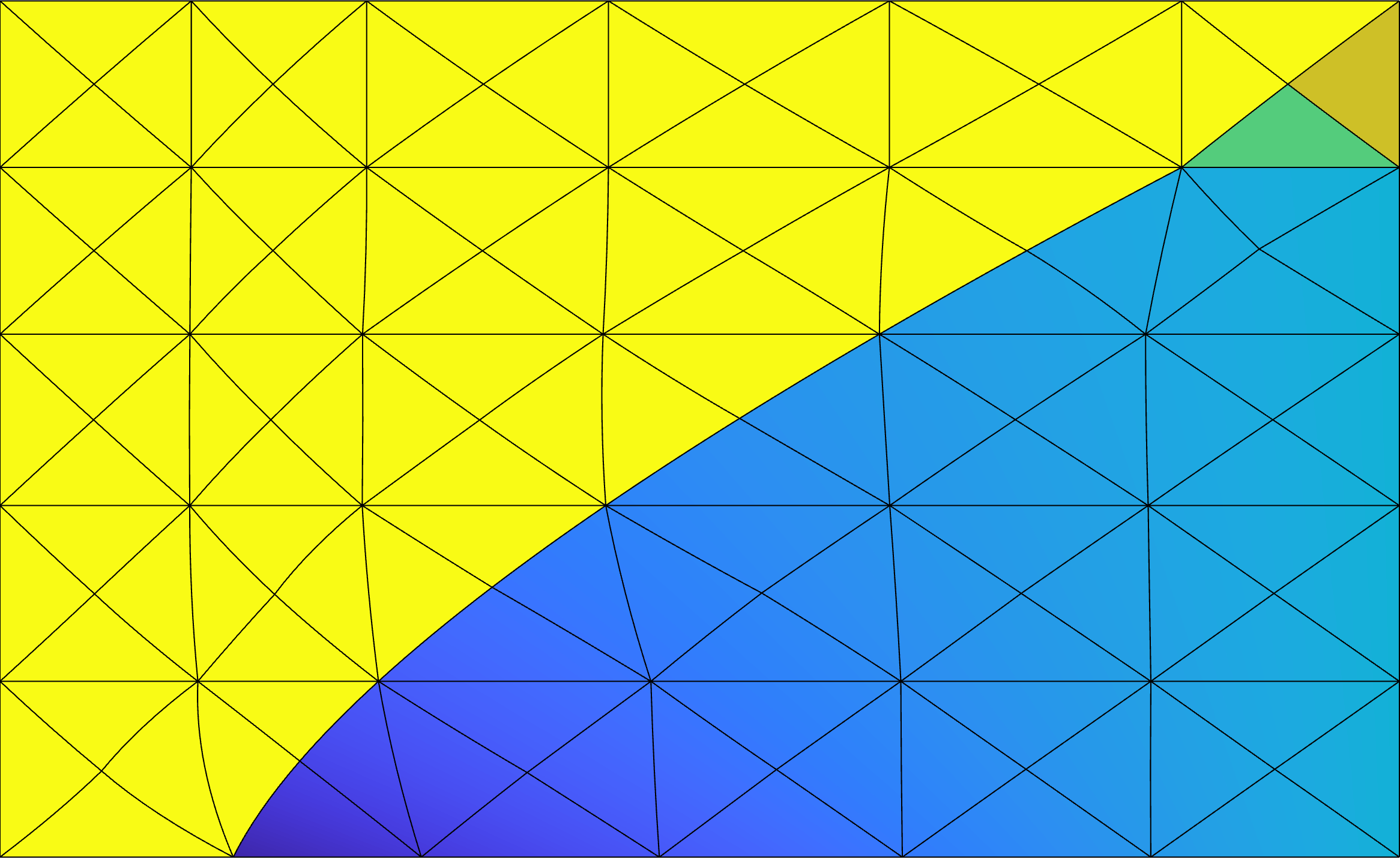};

\addplot [red, thick, dashed]
coordinates {
(-2.00000000e-01,  7.34700000e-01)
( 1.00000000e+00,  7.34600000e-01)};\label{line:pre_tracked_tbnd}

\addplot [blue, thick, dashed]
coordinates {
(-2.00000000e-01,  6.87500000e-01)
( 1.00000000e+00,  6.87500000e-01)};\label{line:tracked_tbnd}

\nextgroupplot[axis equal image, width=0.59\textwidth, xtick={}, ytick={}, xticklabels={}, yticklabels={}, xmin=-0.2, xmax=1, ymin=0, ymax=1]
\addplot []
graphics [xmin=-0.2,xmax=1,ymin=0,ymax=0.6875] { 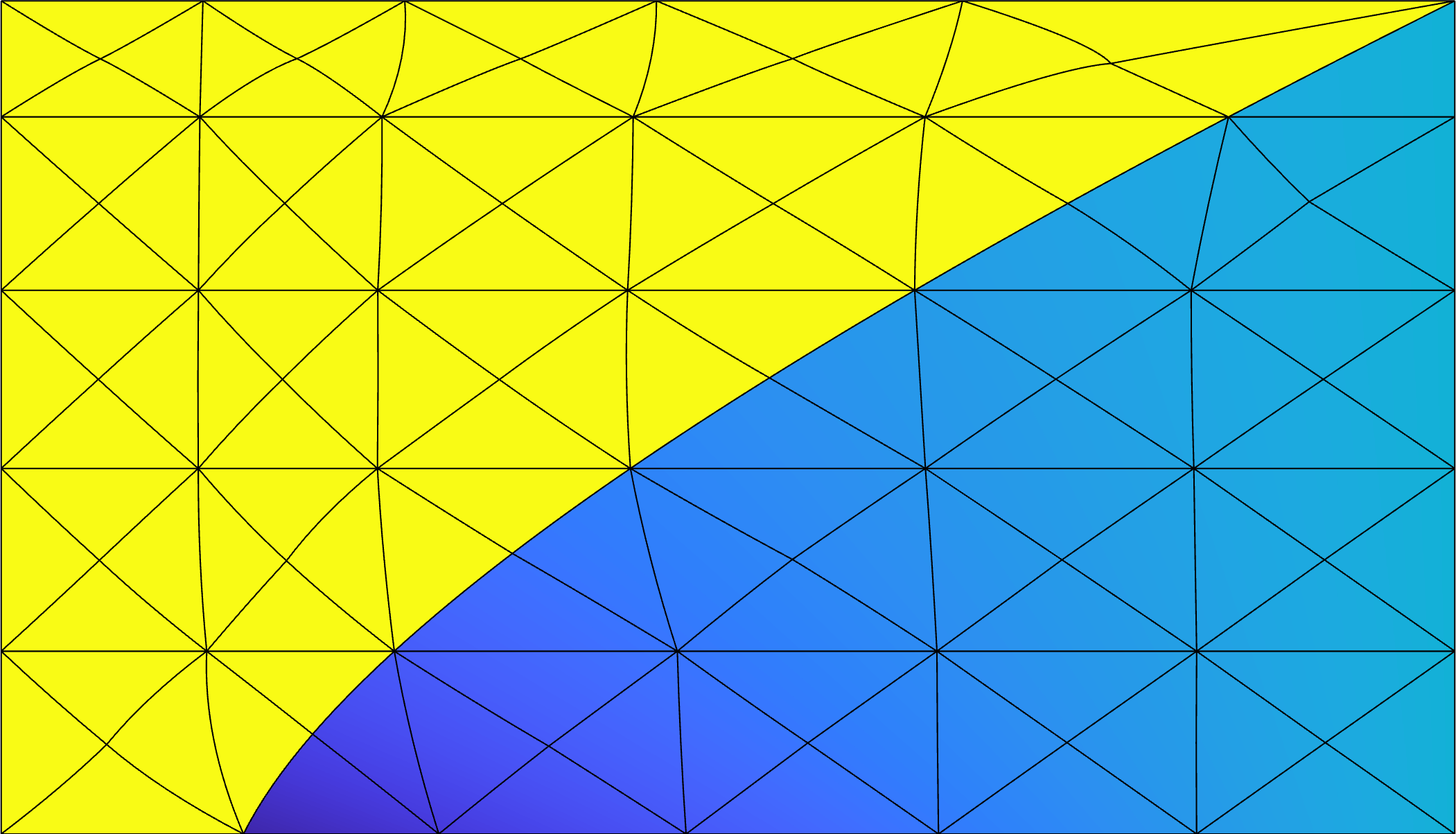};

\addplot [red, thick, dashed]
coordinates {
(-2.00000000e-01,  7.34700000e-01)
( 1.00000000e+00,  7.34600000e-01)};\label{line:pre_tracked_tbnd}

\addplot [blue, thick, dashed]
coordinates {
(-2.00000000e-01,  6.87500000e-01)
( 1.00000000e+00,  6.87500000e-01)};\label{line:tracked_tbnd}

\end{groupplot}\end{tikzpicture}
\colorbarMatlabParula{-3}{-1}{0}{2}{4}
\caption{Demonstration of translating upper boundary of space-time slab.
 With \textit{a priori} slab heights, there may be insufficient resolution
 in the mesh to track shock-boundary interactions (\textit{left}), whereas
 a translating upper boundary allows the slab height to be determined by the
 shock tracking framework such that the shock-boundary interaction is tracked
 (\textit{right}). Legend: \textit{a priori} choice of time slab height
 (\ref{line:pre_tracked_tbnd}), location of shock-boundary interaction
 (\ref{line:tracked_tbnd}).}
\label{fig:temporal_trans}
\end{figure}
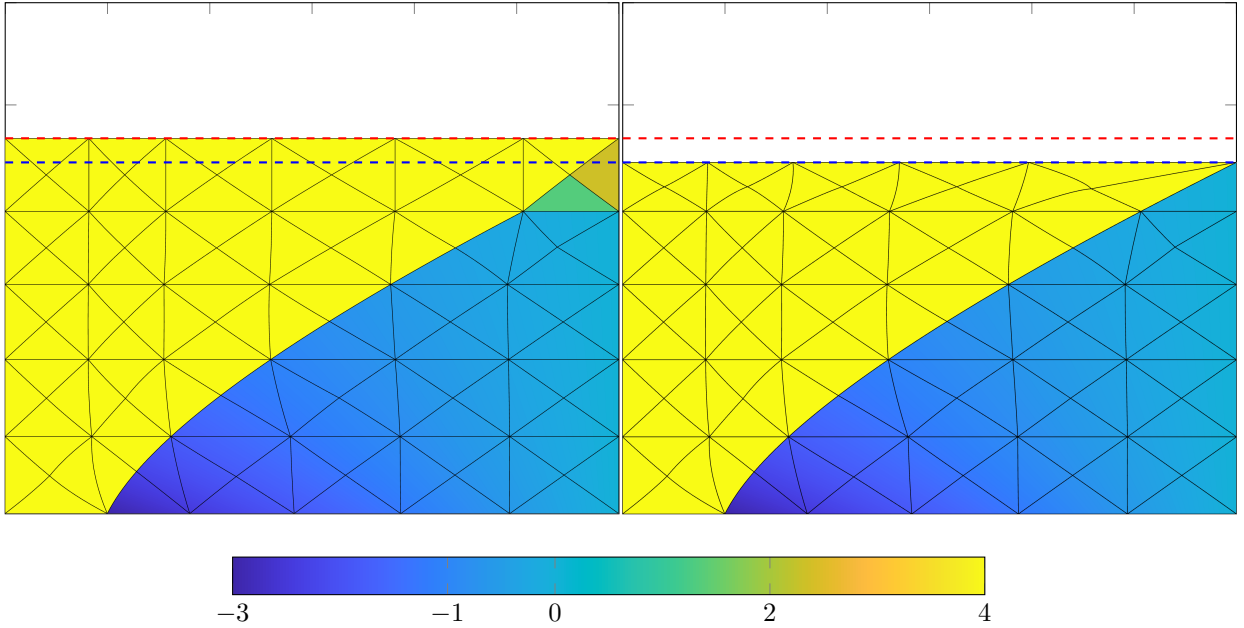

In time-dependent problems, discontinuities can interact with the spatial
boundaries, e.g., leave the domain or reflect off a wall, at unknown times.
If such an event occurs near, but not at, the top of a time slab, there will
likely be insufficient resolution in the mesh to track it. To avoid this issue,
we follow the approach proposed in \cite{2019_corrigan_MGDICEunsteady} to allow
the top of the time slab to translate so the shock-boundary interaction can be tracked
(Figure~\ref{fig:temporal_trans}). Mathematically,
we let $\Ical_\mathtt{t}\subset\{1,\dots,N_\xbm\}$ be the temporal mesh
DoFs associated with the nodes on the top of the time slab. Then
the unconstrained DoFs are written as $\ybm = (\ybm', \zeta)$
and the parametrization satisfies
\begin{equation}
 \phi_i(\ybm) = \zeta, \qquad i \in \Ical_\mathtt{t}
\end{equation}
where $\ybm' \in \Rbb^{N_\ybm-1}$ are the unconstrained DoFs
associated with the mesh degrees of freedom away from $\Ical_\mathtt{t}$
and $\zeta \in \Rbb$ is the translation of the top of the time slab.
Unlike the parametrization in \cite{huang2022robust}, this parametrization
allows the space-time slab to change; however, it guarantees the slab structure
is preserved (because the top boundary translates in the temporal direction).
By incorporating the translating boundary into the admissible domain mapping,
we are relying on the shock tracking method to position it appropriately.
More general parametrizations of the upper temporal boundary of a space-time
slab can be considered to break the slab structure for problems with events
(e.g., shock-shock or shock-boundary interactions) that occur at similar,
but not the same, time (Section~\ref{sec:numexp:euler:sod}).

\subsection{Optimization formulation}
\label{sec:ist:optform}
The HOIST method introduced in \cite{2018_zahr_hoist, 2020_zahr_HOIST} simultaneously
computes the discrete solution of the conservation law and the nodal coordinates of the
mesh that cause element faces to align with discontinuities. This is achieved through a
fully discrete, full space PDE-constrained optimization formulation with the optimization
variables taken to be the discrete flow solution and unconstrained mesh DoFs over the
time slab
\begin{equation} \label{eqn:hoist}
 (\ubm^\star,\ybm^\star) \coloneqq \argmin_{\ubm\in\Rbb^{N_\ubm},\ybm\in\Rbb^{N_\ybm}}~f(\ubm,\ybm) \quad \text{such that} \quad \rbm(\ubm,\phibold(\ybm)) = \zerobold,
\end{equation}
where $f : \Rbb^{N_\ubm} \times \Rbb^{N_\ybm} \rightarrow \Rbb$ is the objective function
and $\xbm^\star = \phibold(\ybm^\star)$ are the physical nodes of the shock-aligned
space-time slab mesh with associated flow solution $\ubm^\star$. The objective function
is composed of two terms as 
\begin{equation}
 f : (\ubm,\ybm) \mapsto
 f_\mathrm{err}(\ubm, \ybm) + \kappa ^2 f_\mathrm{msh}(\ybm_\wslab)
\end{equation}
which balances alignment of the mesh with non-smooth features and the quality of the
elements, $f_\mathrm{err} : \Rbb^{N_\ubm} \times \Rbb^{N_\ybm} \rightarrow \Rbb$
is the alignment term taken as the norm of the enriched residual
\begin{equation}
 f_\mathrm{err} : (\ubm,\ybm) \mapsto \frac{1}{2}\norm{\Rbm(\ubm,\phibold(\ybm))}_2^2,
\end{equation}
$f_\mathrm{msh}$ is a mesh distortion term defined in \cite{2020_zahr_HOIST},
and $\kappa\in\Rbb_{>0}$ is a mesh quality parameter. The mesh distortion term
is designed to blow up if the Jacobian of the mapping goes to zero, which ensures
the domain mapping with never become singular. The HOIST optimization
formulation in (\ref{eqn:hoist}) is identical to the HOIST formulation in
\cite{2020_zahr_HOIST,huang2022robust} over a space-time slab. As such,
we use the sequential quadratic programming (SQP) solver described in
\cite{2020_zahr_HOIST,huang2022robust} to solve (\ref{eqn:hoist}).

\begin{remark}
In addition to avoiding global temporal coupling,
space-time slabs also reduce the difficulty of the shock tracking problem
(\ref{eqn:hoist}), especially for problems with complex interactions that
evolve over time. This is a result of shock tracking only being performed
over a (potentially) small portion of the domain, which limits the temporal complexity.
\end{remark}

\subsection{Boundary conditions and inter-slab solution transfer}
\label{sec:ist:transfer}
Boundary conditions at the bottom of the first slab $\Omega_x\times\{t_0\}$
as well as on the spatial boundaries $\partial\Omega_x \times \Tcal$ are
specified as part of the problem description. Due to causality, boundary
conditions at the top of any slab $\Omega_x\times\{t_1\}$ are extrapolated
from the interior of the space-time domain. Finally, the boundary condition
at the bottom of all slabs except the first must be interpolated from the
top of the previous time slab. In practice, this involved interpolating
the solution in the previous time slab at the quadrature nodes along
the bottom boundary of the current slab. To ensure this interpolation
does not need to be updated at every shock tracking iteration (or
linearized for use in the SQP solver), the mapping in Section~\ref{sec:ist:dommap}
is modified such that the domain deformation is the identity map
when restricted to the bottom boundary. In practice, this is constructed
by fixing the domain parametrization to the reference mesh node along the
bottom boundary. This does not interfere with shock tracking at the
bottom of the slab because any discontinuities that intersect
the bottom of the slab would have been tracked by the previous
slab, which would be preserved in the current slab through the
extraction-based mesh generation (Section~\ref{sec:slab}).

\subsection{Nondimensionalization and scaling of space-time slab}
\label{sec:ist:scaling}
\begin{figure}
\centering
\begin{tikzpicture}
\begin{groupplot} [
  group style={group size = 1 by 3, vertical sep = 0.7cm},
  title style={at={(current bounding box.north west)}, anchor=west}]
\nextgroupplot[axis equal image, width=0.5\textwidth, xmin=0.0, xmax=10, ymin=0, ymax=0.5, xtick={0, 5, 10}, xticklabels={,,}, ytick={0.5}, ylabel={$t$}]
\addplot []
graphics [xmin=0,xmax=10,ymin=0,ymax=0.25] {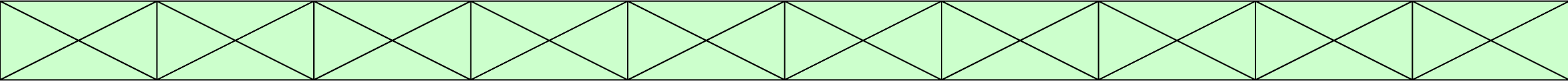};

\nextgroupplot[axis equal image, width=0.5\textwidth, xmin=0.0, xmax=10.0, ymin=0, ymax=5, xtick={0, 5, 10}, xticklabels={,,}, ytick={0, 2.5, 5}, ylabel={$t$}]
\addplot []
graphics [xmin=0,xmax=10,ymin=0,ymax=2.5] { 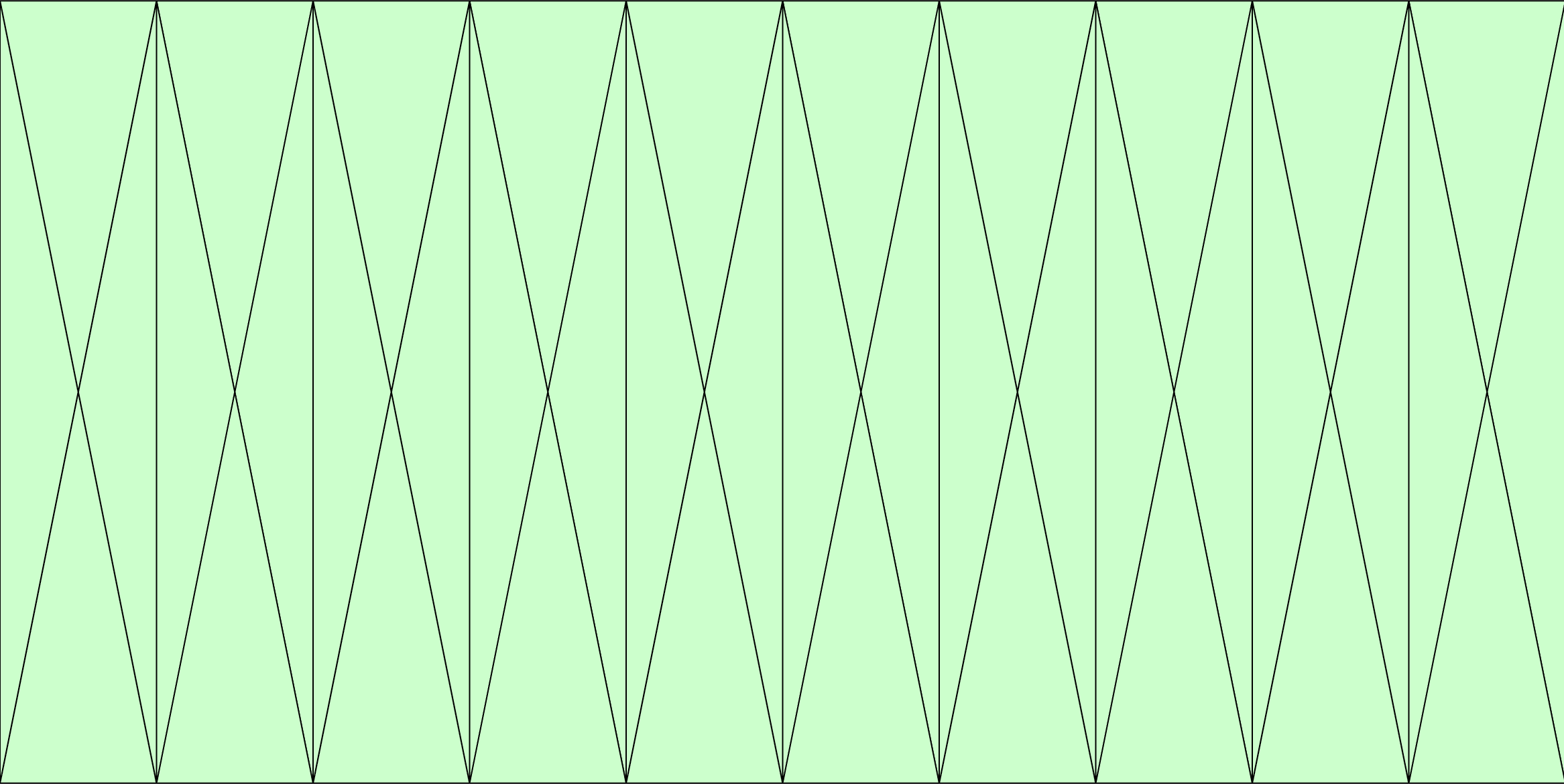};

\nextgroupplot[axis equal image, width=0.5\textwidth, xmin=0.0, xmax=10.0, ymin=0, ymax=2, xtick={0, 5, 10}, ytick={0, 1, 2}, xlabel={$x$}, ylabel={$\hat{t}$}]
\addplot []
graphics [xmin=0,xmax=10,ymin=0,ymax=1] { 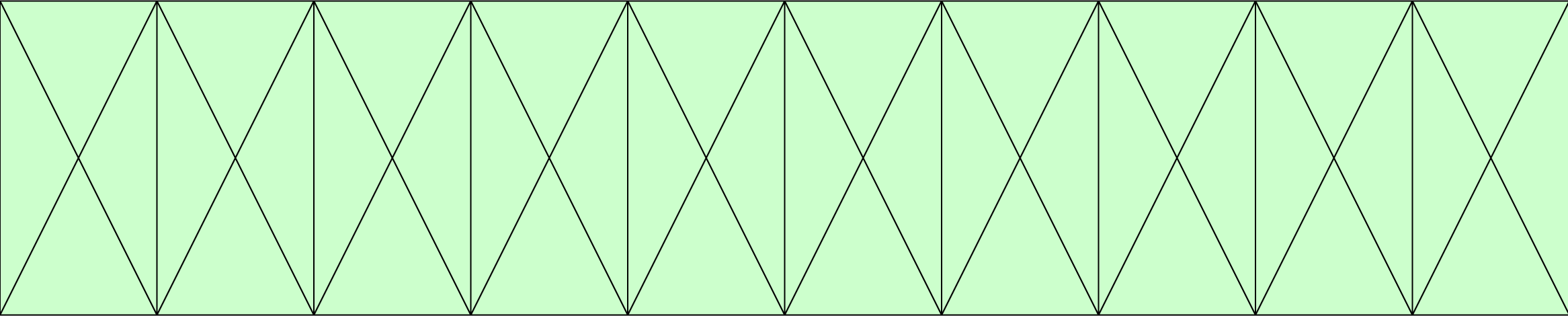};
\end{groupplot}
\end{tikzpicture}
\caption{Schematic of a spatial domain $\Omega_x = (0, 10)$ and temporal domain
 corresponding to fast waves $\Tcal=(0,0.5)$ (\textit{top}) and slower waves
 (or long-time dynamics) $\Tcal=(0,5)$ (\textit{middle}) discretized using
 10 elements in the spatial direction and 2 temporal slabs. The small temporal
 dimension compared to the spatial dimension leads to unnecessary refinement
 and poor quality elements (\textit{top}), whereas the larger temporal
 domain leads to skewed/stretched elements (\textit{midde}).
 Both can be fixed via nondimensionalization with $t^\star < 1$
 (transforms \textit{top} to \textit{bottom}) and
 $t^\star > 1$ (transforms \textit{middle} to \textit{bottom}).}
\label{fig:time_scale}
\end{figure}

There can be orders of magnitude difference between the spatial and
temporal scales, depending on the problem and unit system used.
This can lead to space-time meshes with excessive resolution in
either the spatial or temporal dimension or highly skewed meshes
(Figure~\ref{fig:time_scale}). To circumvent this issue, we
nondimensionalize the time-dependent conservation law with two goals:
1) ensure the components of the solution vector have similar magnitudes
and 2) ensure the spatial dimensions and the temporal interval have
similar magnitudes. The first is standard practice in CFD to improve
the conditioning of the Jacobian of the DG residual ($\rbm$) and the
latter is used to ensure well-conditioned space-time meshes are produced.

Let $D_x^\star\in\Rbb^{d'\times d'}$ be a diagonal matrix defining a
spatial length scale in each dimension that maps a nondimensional domain
$\hat\Omega_x$ to the dimensional domain $\Omega_x$ as
$\hat{x} \mapsto x = D_x^\star\hat{x}$. Furthermore, let
$t^\star\in\Rbb$ be a temporal length scale that maps a
nondimensional time $\hat{t}$ to the dimensional time
as $\hat{t} \mapsto t = t^\star\hat{t}$ and $C\in\Rbb^{m\times m}$
be a diagonal matrix containing the dimensionalization factor
for each solution component (derived in part from the spatial
and temporal length scales), which maps a nondimensional
state vector to the dimensional one as $\hat{W}_x \mapsto W_x = C\hat{W}_x$.
With these definitions, the conservation law in (\ref{eqn:gen_cons_law}) can be
written as a nondimensional conservation law over the nondimensional domain
$\hat\Omega_x$ as
\begin{equation} \label{eqn:claw_nondim}
 \pder{\hat{U}_x}{\hat{t}} + \hat\nabla_x \cdot \hat\Fcal_x(\hat{U}_x)
  = \Scal_x(\hat{U}_x),
\end{equation}
where $\hat{U}_x : \hat\Omega_x \rightarrow \Rbb^m$ is the nondimensional conservative
vector implicitly defined as the solution of (\ref{eqn:claw_nondim}) and related to
the solution of the original conservation law as
$\hat{U}_x(\hat{x},\hat{t}) = C^{-1} U_x(D_x^\star \hat{x}, t^\star \hat{t})$,
$\hat\Fcal_x : \Rbb^m \rightarrow \Rbb^{m \times d}$ and
$\hat\Scal_x : \Rbb^m \rightarrow \Rbb^m$ are the nondimensional
flux function and source term, respectively, and related to the original terms
as
\begin{equation}
 \hat\Fcal_x : \hat{W}_x \mapsto t^\star C^{-1} \Fcal_x(C\hat{W}_x) (D_x^\star)^{-1},
 \qquad
 \hat\Scal_x : \hat{W}_x \mapsto t^\star C^{-1} \Scal_x(C\hat{W}_x),
\end{equation}
and ($\hat\nabla_x \cdot$) is the nondimensional divergence operator. The procedure
to formulate transformed space-time conservation laws in
Sections~\ref{sec:govern:sptm}-\ref{sec:govern:transf}
proceeds without modfication with the spatial conservation
law in (\ref{eqn:claw_nondim}) in place of the conservation
law in (\ref{eqn:gen_cons_law}).

In this work, we take the spatial scaling to be the identity matrix (with
appropriate units for nondimensionalization) and focus on the choice of
temporal scaling that leads to well-conditioned space-time simplices.
In the case where the temporal dimension of the slab is much smaller than the
spatial dimensions---usually occurring with fast moving waves evolving
over small time scales---we take $t^\star < 1$ to expand the temporal
domain, effectively slowing the features. Otherwise, we take $t^\star > 1$,
which compresses the temporal domain, effectively speeding up the features
and increasing the temporal resolution of the mesh. In both situations, the
conditioning of the space-time simplices improves substantially
(Figure~\ref{fig:time_scale}). Once the characteristic length ($D_x^\star$)
and time ($t^\star$) scales are determined, the solution scaling $C$ usually
follows directly from a characteristic mass scale.

\subsection{Adaptive mesh refinement}
\label{sec:ist:amr}
In time-dependent problems, the solutions can develop complex features in localized
regions of space-time, e.g., shock formation or features with small length scales,
which calls for increased resolution using either $h$-, $p$-, or $hp$-refinement.
In this work, we use the approach described in Sections~\ref{sec:slab:split} and
\ref{sec:slab:split:smplx2smplx} to $h$-refine elements in the refinement set
$\bar\Ecal_h^\mathrm{r}$. The refinement set is the collection of elements
where the magnitude of the enriched DG residual $\Rbm$ and the Persson-Peraire
shock sensor \cite{persson2013shock, persson2006sub},
$s_{\bar{K}} : \Vcal_h^p \rightarrow \Rbb$, $\bar{K}\in\bar\Ecal_h$,
is largest. The Persson-Peraire shock sensor, which identifies elements
with oscillatory solutions, is defined as
\begin{equation}
 s_{\bar{K}} : \bar{W}_h \mapsto \log_{10} \left( \sqrt{ \frac{\int_{\bar{K}} ( \chi ( \bar{W}_h) - \mathcal{X} (\Pi_{p-1} \bar{W}_h))^2 \, dV } {\int_{\bar{K}} \chi (\bar{W}_h) ^2 \, dV} } \right),
\end{equation}
where $\Pi_{p-1} : \mathcal{V}_h^p \rightarrow \mathcal{V}_h^{p-1}$ is a projection onto
the space of piecewise polynomials of total degree $p-1$ and
$\chi : \Rbb^m \rightarrow \Rbb$ is a suitable scalar of the state. In this work
we choose $\chi : (W_1, ..., W_m) \mapsto W_1$ and define the refinement set as
\begin{equation}
 \bar\Ecal_h^\mathrm{r} : (\ubm,\xbm) \mapsto
 \left\{ \bar{K} \in \bar\Ecal_h \suchthat \norm{\Rbm_{\bar{K}}(\ubm,\xbm)} \geq  c_1  R_\mathrm{max}(\ubm,\xbm)  \text{ or } s_{\bar{K}} ( \Xi(\ubm)) \geq c_2   s_\mathrm{max}(\ubm)   \right\}
\end{equation}
where
\begin{equation}
 R_\mathrm{max} : (\ubm, \xbm) \mapsto \max_{\bar{K}\in\Ecal_h}~\norm{\Rbm_{\bar{K}}(\ubm,\xbm)}, \qquad
 s_\mathrm{max} : \ubm \mapsto \max_{\bar{K}\in\Ecal_h}~s_{\bar{K}}(\Xi(\ubm))
\end{equation}
are the largest elemental residual norm and sensor, respectively, and $c_1, c_2 \in \Rbb_{\geq 0}$ are indicator fractions whose default values are set to be $c_1 = 0.01$ and $c_2 = 0.01$ unless otherwise specified. We interrupt the HOIST solver at iteration $k$ to refine the mesh, provided the number of elements in the mesh has not exceeded a user-defined limit, where $k$ is any HOIST iteration with the DG residual below a user-defined threshold.

\subsection{Shock boundary conditions}
\label{sec:ist:shkbc}
As detailed in \cite{huang2023high}, shock tracking provides an opportunity to
reduce elements in the mesh that is not available to shock capturing methods.
In situations where a discontinuity completely separates a region in which
the flow is \textit{known} from the downstream region, all elements in the
region where the flow is known can be removed and the known flow used as
an essential boundary condition applied directly to the shock surface.
Such a strategy is available even if the shock position is not known by
defining the mesh parametrization $\phibold$ to allow arbitrary deformations
to the shock-boundary and using the implicit shock tracking framework to
determine its precise location. In the space-time setting, this approach
can only be leveraged if a space-time discontinuity separates a region
of space-time where the flow is known from the remainder of the space-time
domain. This is applied by creating a spatial mesh of only the region where
the flow is unknown at the initial condition, extruding it to form a space-time
slab, applying a boundary condition that sets the flow solution to the known value
on the elements that separate the regions where the flow is unknown vs. known,
and using implicit shock tracking to position the shock-boundary appropriately.
This technique is demonstrated for the accelerating shock problem
(Section~\ref{sec:numexp:burg:acc}) and Shu-Osher problem
(Section~\ref{sec:numexp:euler:shuosher}).

\subsection{Summary}
\label{sec:ist:summ}
The complete slab-based space-time HOIST method for shock-dominated unsteady
conservation laws is summarized in Algorithm~\ref{alg:hoist}.
\begin{algorithm}
 \caption{Slab-based space-time HOIST method}
 \label{alg:hoist}
 \begin{algorithmic}[1]
  \REQUIRE Reference spatial mesh $\bar\Ecal_{x,h}$, desired time duration $\tau$, nominal temporal slab height $\Delta t$
  \ENSURE Shock-aligned meshes and corresponding DG solution
  \STATE $t = 0$
  \WHILE{$t \leq \tau$}
   \STATE \textbf{Construct mesh of space-time slab}: Extrude $\bar\Ecal_{x,h}$ over
      time interval $(t,t+\Delta t)$ (Section~\ref{sec:slab:extrude}) and split into
      simplices to form mesh of space-time slab $\bar\Ecal_h'$ (Section~\ref{sec:slab:split})
   \STATE \textbf{Construct parametrization of space-time slab}: Define
      $\phibold$ according to \cite{huang2022robust} with self-adjusting
      upper boundary (Section~\ref{sec:ist:dommap}) and fixed lower
      boundary (Section~\ref{sec:ist:transfer}) if $t > 0$
   \STATE \textbf{Solution initialization}: Initialize SQP iteration with
      the unconstrained DoFs of the reference mesh and the DG($p=0$) solution
      on the reference mesh \cite{huang2022robust}
   \STATE \textbf{Implicit shock tracking on space-time slab}: Solve optimization
      problem (\ref{eqn:hoist}) for shock-aligned mesh and DG solution on current slab
   \STATE \textbf{Extract spatial mesh from top of time slab}: Update
      $\bar\Ecal_{x,h}$ by extracting the spatial mesh from the top of
      the (physical) space-time slab mesh $\Ecal_h$ (Section~\ref{sec:slab:spatial})
   \STATE \textbf{Update time}: Update $t$ from the top of the (physical)
      space-time slab mesh
  \ENDWHILE
\end{algorithmic}
\end{algorithm}


\section{Numerical experiments}
\label{sec:numexp}
In this section, we formulate six problems (three conservation laws) to
demonstrate the convergence, accuracy, and overall merit of the proposed
slab-based implicit shock tracking approach via direct comparisons to
exact solutions (if available) and popular shock capturing methods.


\subsection{Inviscid Burgers' equation}
We consider the time-dependent, inviscid Burgers' equation that governs nonlinear
advection of a scalar quantity $w : \Omega_x \rightarrow \Rbb$ through an
$d'$ -dimensional domain $\Omega_x \subset \Rbb^{d'}$
\begin{equation} \label{eqn:burg}
\pder{}{t}  w(x,t) + \nabla_x\cdot\left(( w(x,t)^2/2)\beta^T\right) = 0,
 \qquad  w(x, 0) = w_0(x)
\end{equation}
for $x \in \Omega_x$, $t\in\Tcal$, where $w_0 : \Omega_x \mapsto \Rbb$ is
the initial condition, $\beta \in \Rbb^{d'}$ is the flow direction, and the boundary conditions are consistent with the initial
condition at inflow boundaries. Burgers' equation is cast as a conservation law
of the form (\ref{eqn:gen_cons_law}) and the projected inviscid flux Jacobian
and its eigenvalue decomposition are given in \ref{app:projjac}. From
this information, the transformed space-time version of the conservation law
(including Riemann solver-based numerical flux) follows systematically from
Section~\ref{sec:govern:sptm}-\ref{sec:govern:transf}.

\subsubsection{One-dimensional accelerating shock leaving domain}
\label{sec:numexp:burg:acc}
First, we consider Burgers' equation (\ref{eqn:burg}) ($\beta=1$) with an accelerating
shock in the one-dimensional spatial domain $\Omega_x \coloneqq (-0.2, 1)$ and temporal
domain $\mathcal{T} \coloneqq (0, 1)$ produced by the initial condition
\begin{equation}
  w_0 : x \mapsto \mu_1 H(-x) - \mu_2 H(x),
\end{equation}
where $H : \Rbb \rightarrow \Rbb$ is the Heaviside function, and
$\mu_1 = 4$ and $\mu_2 = 3$ in this work. The analytical solution is
\begin{equation}
w : (x,t) \mapsto \mu_1 H(x_s(t) - x) + \frac{\mu_2 (x-1)}{1+\mu_2 t} H(x - x_s(t))
\end{equation}
 where $x_s : \mathcal{T} \rightarrow \Rbb$ is the shock speed
\begin{equation}
 x_s : t \mapsto (\mu_1 / \mu_2 + 1)(1 - \sqrt{1 + \mu_2 t}) + \mu_1 t.
\end{equation}
For the case considered, the shock leaves the domain before the end of the
time interval of interest. Optimal convergence of the HOIST method for this
problem (single slab) has been demonstrated in previous work \cite{huang2022robust} and
will not be repeated here. The purpose of this numerical experiment is to
demonstrate the space-time slab-based HOIST method and shock boundary condition
approach for a simple problem, examine convergence of the optimization solver, and
investigate performance with a different number of slabs.

\begin{figure}
	\centering
 	\raisebox{-0.5\height}{\input{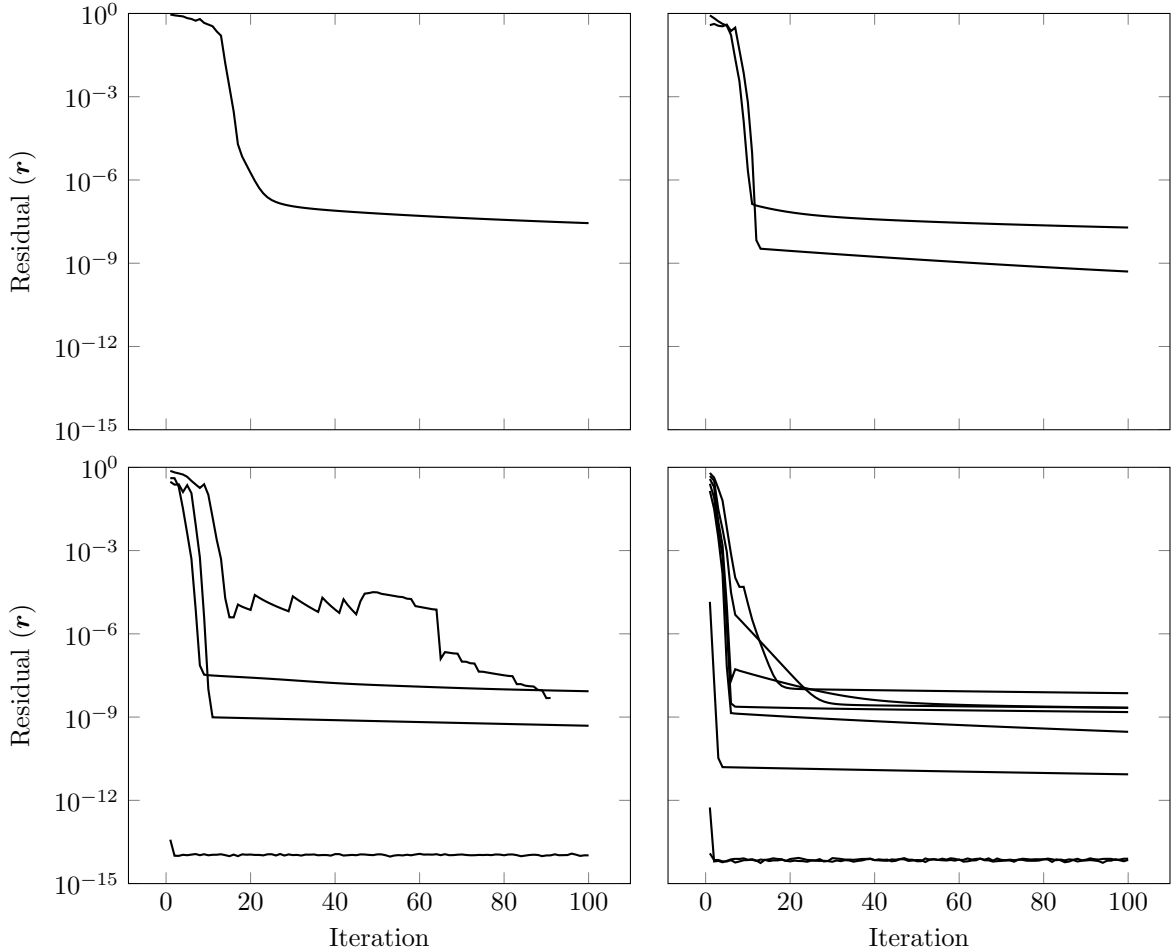}}
  \caption{Convergence of the HOIST method for the accelerating shock problem with different slab configurations with $S = 1$ (\textit{top left}), $S = 2$ (\textit{top right}),
 $S = 4$ (\textit{bottom left}), and $S = 8$ (\textit{bottom right}) slabs. Different
 lines on a single plot correspond to the convergence history of different slabs.}
 	\label{fig:iburg_nrl_conv}
\end{figure}

\begin{figure}
  \centering
  \raisebox{-0.5\height}{\begin{tikzpicture}
\begin{groupplot} [
group style={group size = 4 by 2, horizontal sep = 0.5cm, vertical sep = 0.5cm},
title style={at={(current bounding box.north west)}, anchor=west}]
\nextgroupplot[axis equal image, width=0.3\textwidth, xtick={-0.2, 1}, ytick={0, 0.5/2, 1/2}, xticklabels={}, yticklabels={0, 0.5, 1}, ylabel={$t$}, xmin=-0.1, xmax=0.5, ymin=0, ymax=0.5]
\addplot []
graphics [xmin=-0.1,xmax=0.5,ymin=0,ymax=0.06643] { 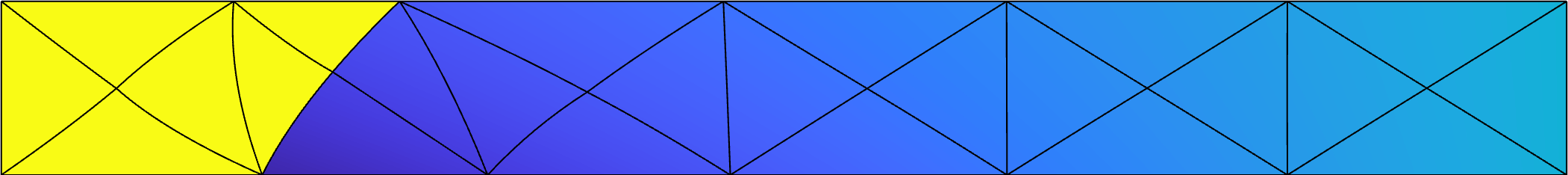};

\nextgroupplot[axis equal image, width=0.3\textwidth, xtick={-0.2/2, 1/2}, ytick={0, 1/2}, xticklabels={}, yticklabels={}, xmin=-0.1, xmax=0.5, ymin=0, ymax=0.5]
\addplot []
graphics [xmin=-0.1,xmax=0.5,ymin=0,ymax=0.13762] { 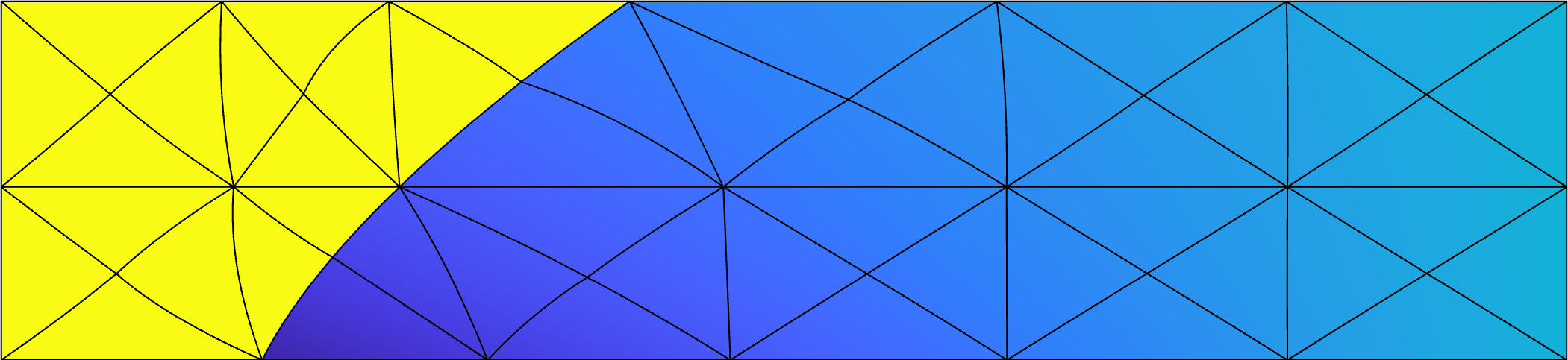};

\nextgroupplot[axis equal image, width=0.3\textwidth, xtick={-0.2/2, 1/2}, ytick={0, 1/2}, xticklabels={}, yticklabels={}, xmin=-0.1, xmax=0.5, ymin=0, ymax=0.5]
\addplot []
graphics [xmin=-0.1,xmax=0.5,ymin=0,ymax=0.206145] { 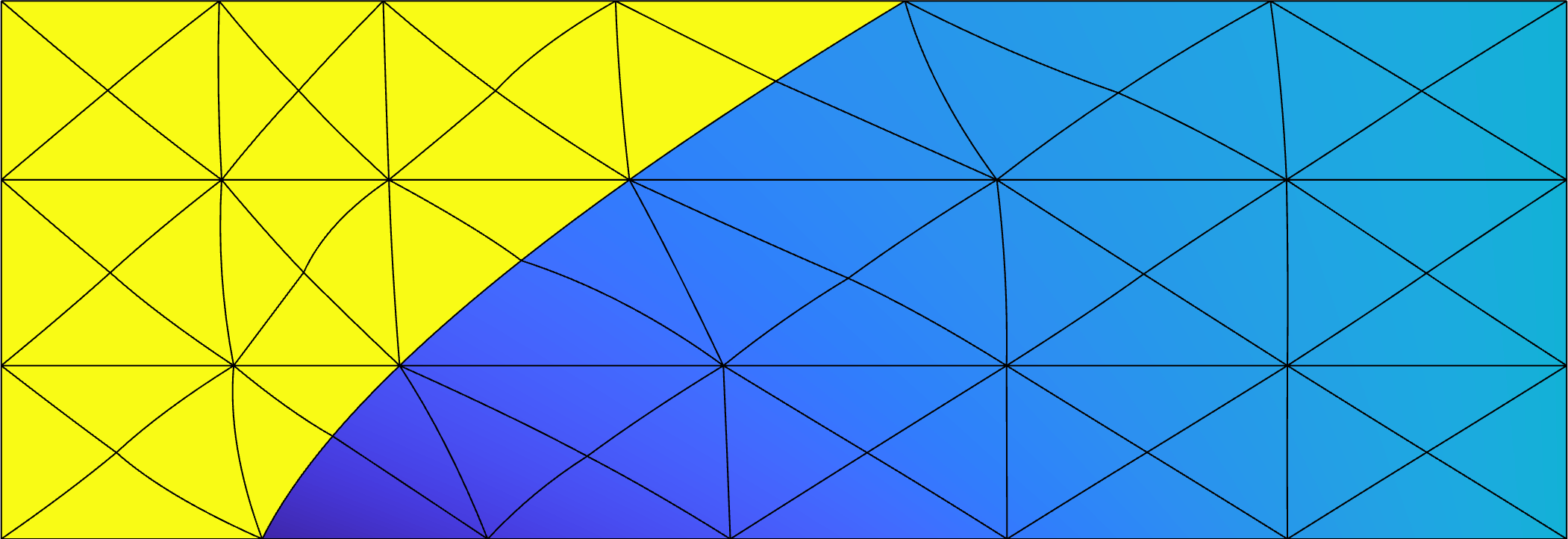};

\nextgroupplot[axis equal image, width=0.3\textwidth, xtick={-0.2/2, 1/2}, ytick={0, 1/2}, xticklabels={}, yticklabels={}, xmin=-0.1, xmax=0.5, ymin=0, ymax=0.5]
\addplot []
graphics [xmin=-0.1,xmax=0.5,ymin=0,ymax=0.13762] { 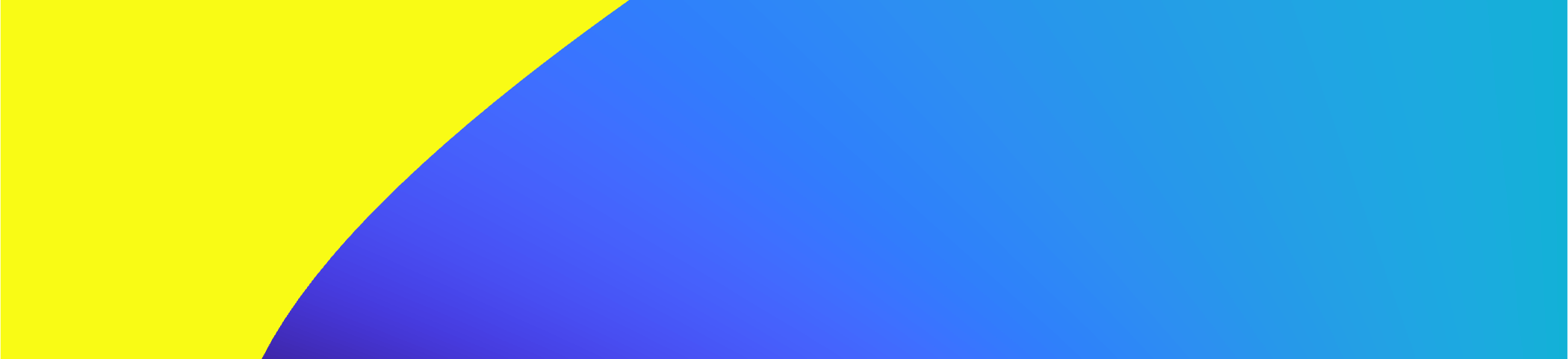};

\addplot []
graphics [xmin=-0.1,xmax=0.5,ymin=0,ymax=0.27305] { 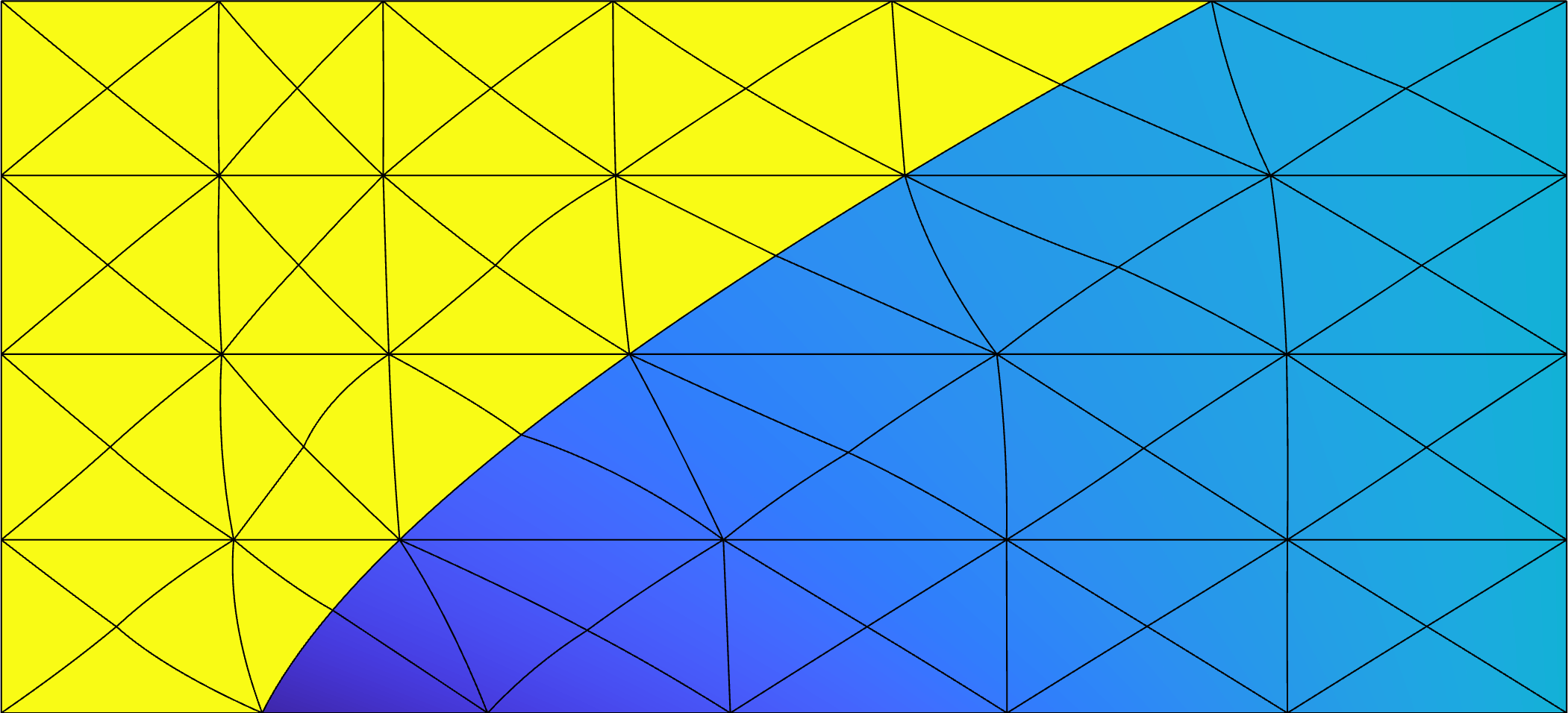};

\nextgroupplot[axis equal image, width=0.3\textwidth, xlabel={$x$}, ylabel={$t$}, xtick={-0.2/2, 0, 0.5/2, 1/2}, ytick={0, 0.5/2, 1/2}, xticklabels={-0.2, 0, 0.5, 1}, yticklabels={0, 0.5, 1}, xmin=-0.1, xmax=0.5, ymin=0, ymax=0.5]
\addplot []
graphics [xmin=-0.1,xmax=0.5,ymin=0,ymax=0.34375] { 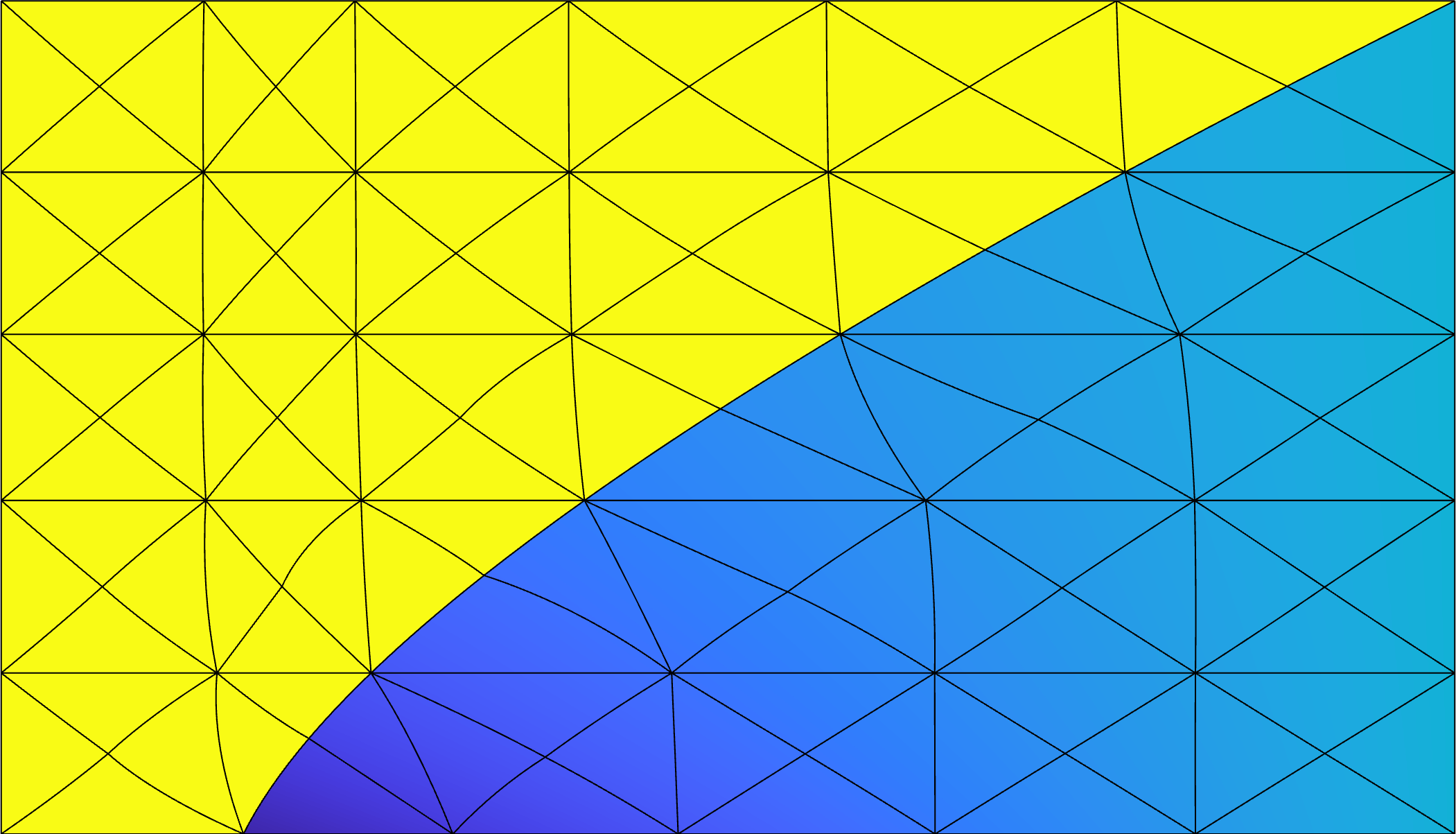};

\nextgroupplot[axis equal image, width=0.3\textwidth, xlabel={$x$}, xtick={-0.2/2, 0, 0.5/2, 1/2}, ytick={0, 1/2}, xticklabels={-0.2, 0, 0.5, 1}, yticklabels={}, xmin=-0.1, xmax=0.5, ymin=0, ymax=0.5]
\addplot []
graphics [xmin=-0.1,xmax=0.5,ymin=0,ymax=0.40863] { 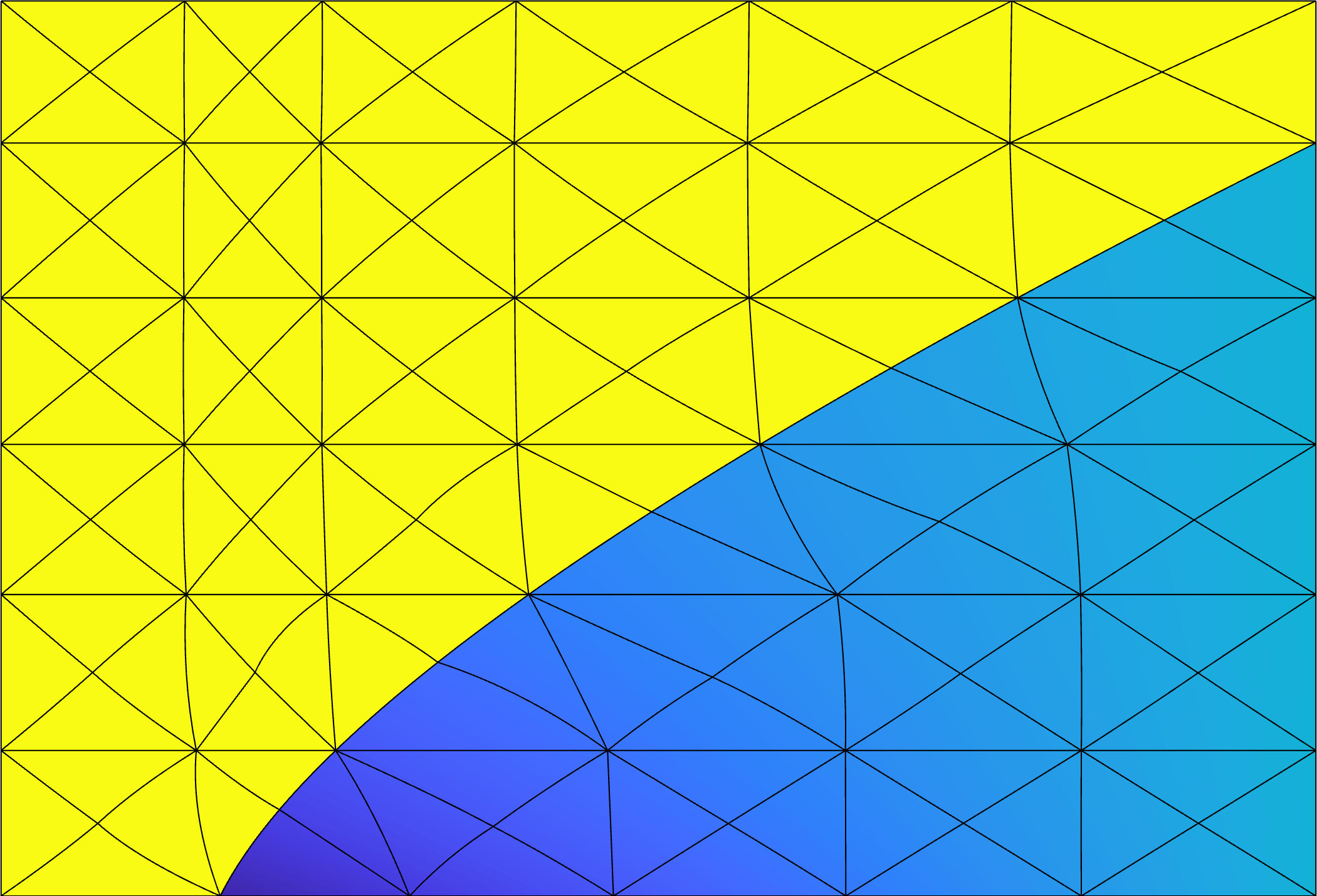};

\nextgroupplot[axis equal image, width=0.3\textwidth, xlabel={$x$}, xtick={-0.2/2, 0, 0.5/2, 1/2}, ytick={0, 1/2}, xticklabels={-0.2, 0, 0.5, 1}, yticklabels={}, xmin=-0.1, xmax=0.5, ymin=0, ymax=0.5]
\addplot []
graphics [xmin=-0.1,xmax=0.5,ymin=0,ymax=0.4735] { 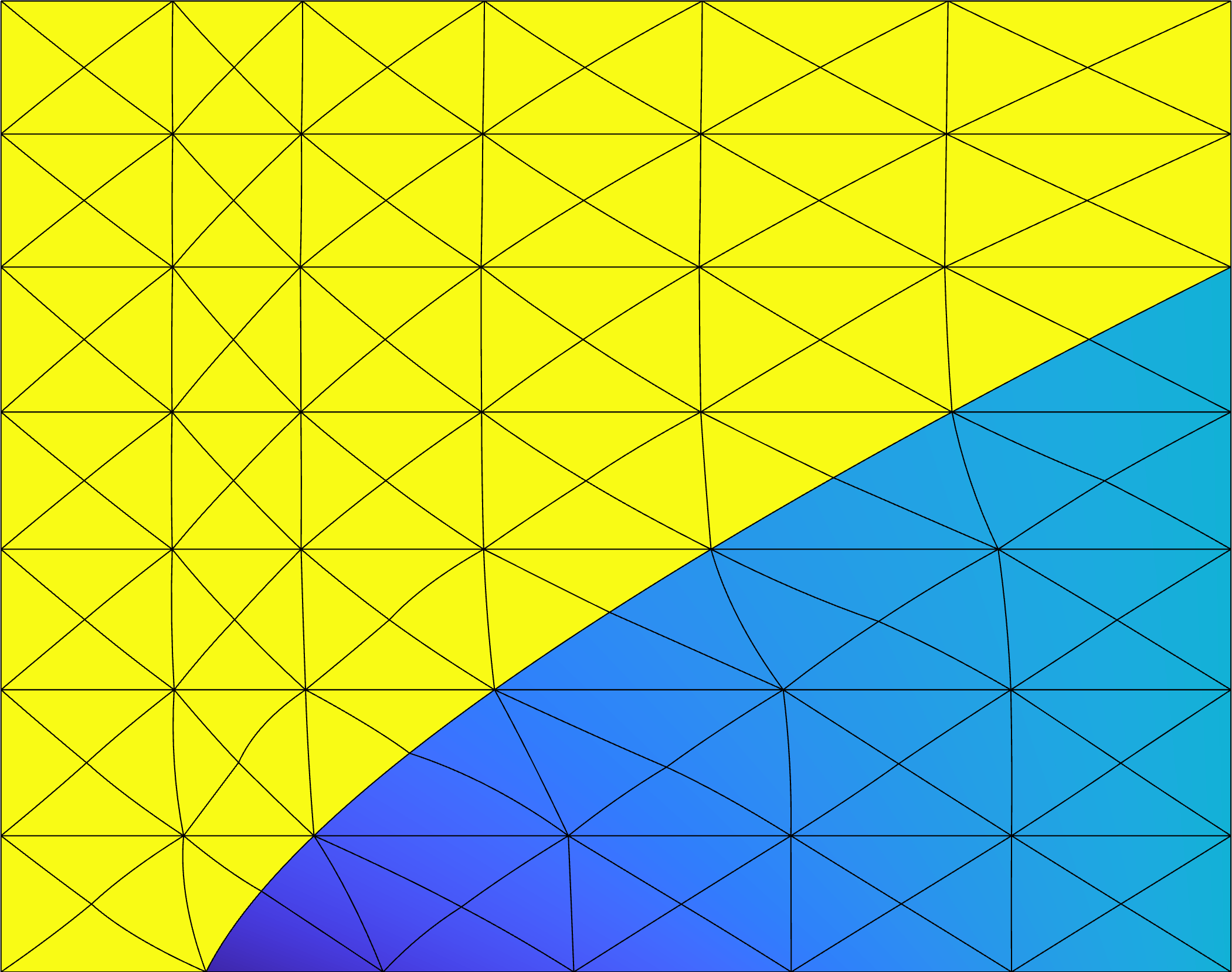};

\nextgroupplot[axis equal image, width=0.3\textwidth, xlabel={$x$}, xtick={-0.2/2, 0, 0.5/2, 1/2}, ytick={0, 1/2}, xticklabels={-0.2, 0, 0.5, 1}, yticklabels={}, xmin=-0.1, xmax=0.5, ymin=0, ymax=0.5]
\addplot []
graphics [xmin=-0.1,xmax=0.5,ymin=0,ymax=0.5] { 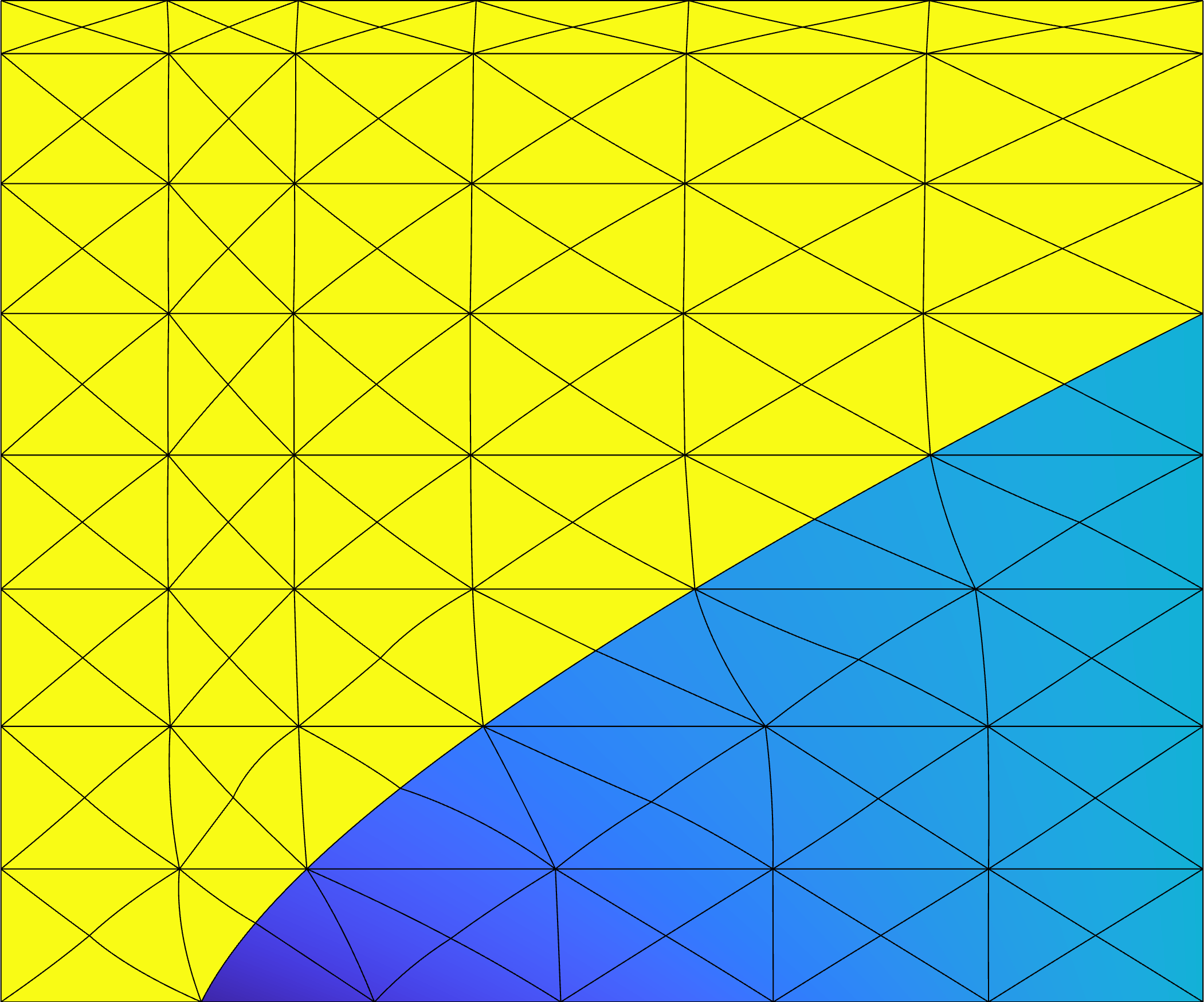};

\nextgroupplot[axis equal image, width=0.59\textwidth, xtick={-0.2/2, 0, 0.5/2, 1/2}, ytick={0, 1/2}, xticklabels={-0.2, 0, 0.5, 1}, yticklabels={}, xmin=-0.1, xmax=0.5, ymin=0, ymax=0.5]
\nextgroupplot[axis equal image, width=0.59\textwidth, xtick={-0.2/2, 0, 1/2}, ytick={0, 0.5/2, 1/2}, xticklabels={}, yticklabels={}, xmin=-0.1, xmax=0.5, ymin=0, ymax=0.5]
\nextgroupplot[axis equal image, width=0.59\textwidth, xtick={-0.2/2, 0, 1/2}, ytick={0, 0.5/2, 1/2}, xticklabels={}, yticklabels={}, xmin=-0.1, xmax=0.5, ymin=0, ymax=0.5]
\nextgroupplot[axis equal image, width=0.59\textwidth, xtick={-0.2/2, 0, 1/2}, ytick={0, 0.5/2, 1/2}, xticklabels={}, yticklabels={}, xmin=-0.1, xmax=0.5, ymin=0, ymax=0.5]
\nextgroupplot[axis equal image, width=0.59\textwidth, xtick={-0.2/2, 0, 1/2}, ytick={0, 0.5/2, 1/2}, xticklabels={}, yticklabels={}, xmin=-0.1, xmax=0.5, ymin=0, ymax=0.5]
\nextgroupplot[axis equal image, width=0.59\textwidth, xtick={-0.2/2, 0, 1/2}, ytick={0, 0.5/2, 1/2}, xticklabels={}, yticklabels={}, xmin=-0.1, xmax=0.5, ymin=0, ymax=0.5]
\nextgroupplot[axis equal image, width=0.59\textwidth, xtick={-0.2/2, 0, 1/2}, ytick={0, 0.5/2, 1/2}, xticklabels={}, yticklabels={}, xmin=-0.1, xmax=0.5, ymin=0, ymax=0.5]
\nextgroupplot[axis equal image, width=0.59\textwidth, xtick={-0.2/2, 0, 1/2}, ytick={0, 0.5/2, 1/2}, xticklabels={}, yticklabels={}, xmin=-0.1, xmax=0.5, ymin=0, ymax=0.5]
\end{groupplot}\end{tikzpicture}}
  \colorbarMatlabParula{-3}{-1}{0}{2}{4}
  \caption{HOIST solution of accelerating shock inviscid Burgers' problem over a sequence of space-time slabs ($S = 8$).}
  \label{fig:iburg_nrl_slabs_1_thru_2}
\end{figure}

We consider the convergence of the HOIST solver as a function of slab height.
To this end, we consider slab heights of $\Delta t \in \{1.2, 0.6, 0.3, 0.15\}$,
which corresponds to $S$ slabs, where $S \in \{1, 2, 4, 8\}$. Because the slab
height is optimized (Section~\ref{sec:ist:dommap}), $\Delta t$ is the initial height
as the actual slab height changes during the optimization. We use a grid with
$192$ quadratic ($p = q = 2$) space-time simplices, which corresponds to $192 / S$
elements per slab. All slabs are run for $100$ iterations and the convergence of the
DG residual is shown in Figure~\ref{fig:iburg_nrl_conv}. As more slabs are added,
the residuals converge more rapidly and to deeper tolerances, indicating the shock
tracking problem is simpler for smaller slabs. The last slab for $S = 4$ and last
three slabs for $S = 8$ correspond to times after the shock has left the domain, which
is the reason the corresponding residuals are immediately machine zero. The solution
field and space-time mesh for the $S = 8$ is shown in
Figure~\ref{fig:iburg_nrl_slabs_1_thru_2}. The curved shock is tracked by the
quadratic mesh elements until it exits the domain in the fifth slab, where the
translating boundary (Section~\ref{sec:ist:dommap}) adjusts the upper temporal
boundary such that the point where the shock exists the right boundary is
exactly tracked.

\begin{figure}
  \centering
  \raisebox{-0.5\height}{\begin{tikzpicture}
\begin{groupplot} [
group style={group size = 5 by 2, horizontal sep = 0.05cm, vertical sep = 0.4cm},
title style={at={(current bounding box.north west)}, anchor=west}]
\nextgroupplot[axis equal image, width=0.27\textwidth, xtick={-0.2, 1}, ytick={0, 0.738}, xticklabels={}, yticklabels={}, ylabel={$t$}, xmin=-0.2, xmax=1, ymin=0, ymax=0.738]
\addplot []
graphics [xmin=-0.2,xmax=1,ymin=0,ymax=0.143] { 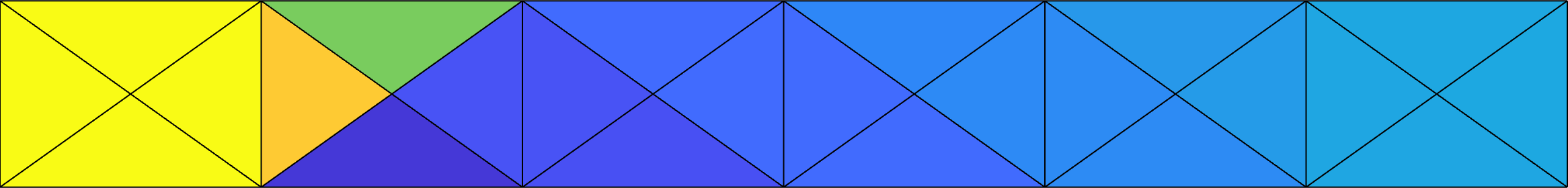};

\nextgroupplot[axis equal image, width=0.27\textwidth, xtick={-0.2, 1}, ytick={0, 1}, xticklabels={}, yticklabels={}, xmin=-0.2, xmax=1, ymin=0, ymax=0.738]
\addplot []
graphics [xmin=0,xmax=1,ymin=0,ymax=0.294] { 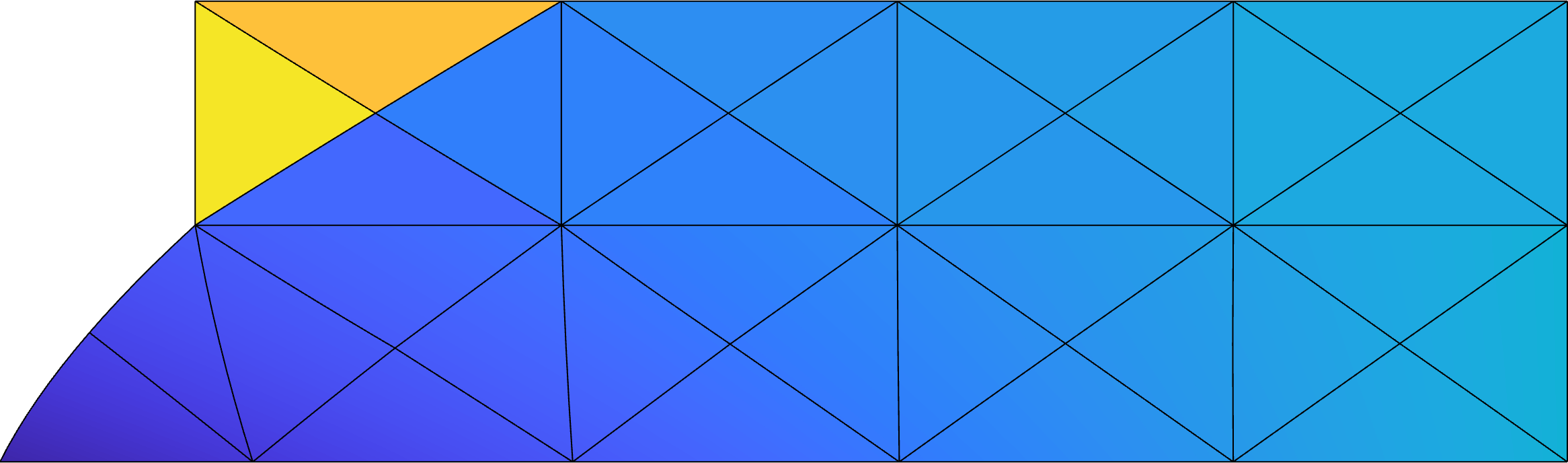};

\nextgroupplot[axis equal image, width=0.27\textwidth, xtick={-0.2, 1}, ytick={0, 1}, xticklabels={}, yticklabels={}, xmin=-0.2, xmax=1, ymin=0, ymax=0.738]
\addplot []
graphics [xmin=0,xmax=1,ymin=0,ymax=0.446] { 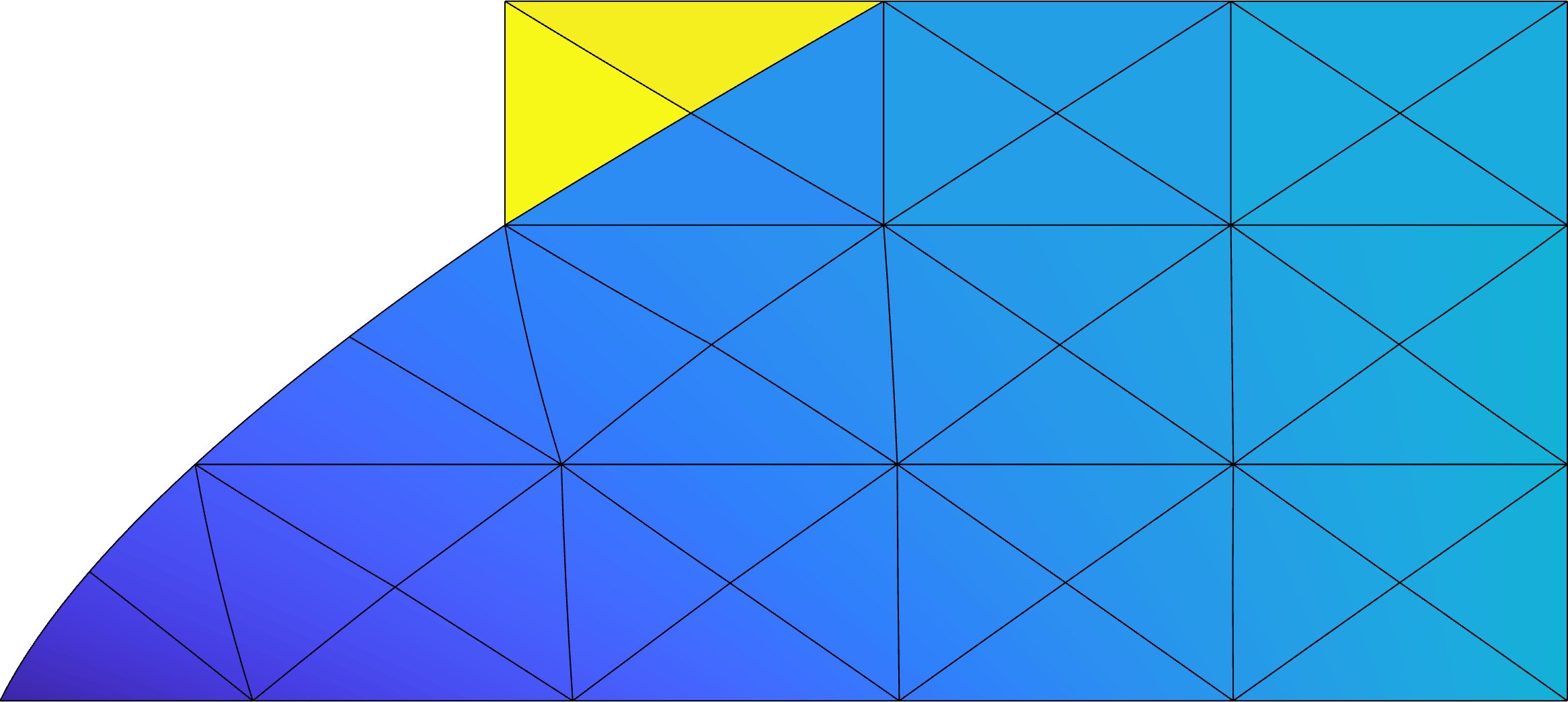};

\nextgroupplot[axis equal image, width=0.27\textwidth, xtick={-0.2, 1}, ytick={0, 1}, xticklabels={}, yticklabels={}, xmin=-0.2, xmax=1, ymin=0, ymax=0.738]
\addplot []
graphics [xmin=0,xmax=1,ymin=0,ymax=0.595] { 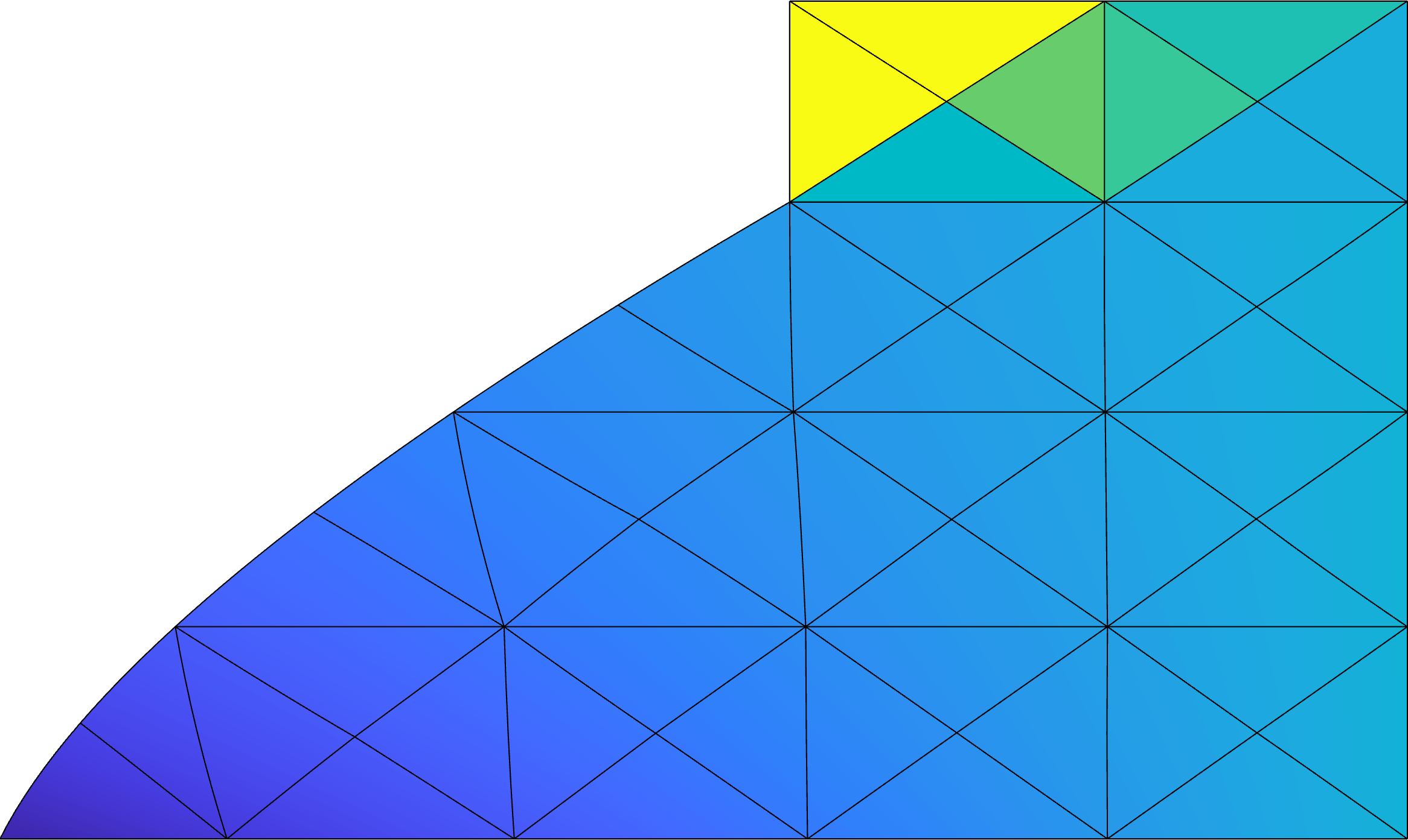};

\nextgroupplot[axis equal image, width=0.27\textwidth, xtick={-0.2, 1}, ytick={0, 1}, xticklabels={}, yticklabels={}, xmin=-0.2, xmax=1, ymin=0, ymax=0.738]
\addplot []
graphics [xmin=0,xmax=1,ymin=0,ymax=0.738] { 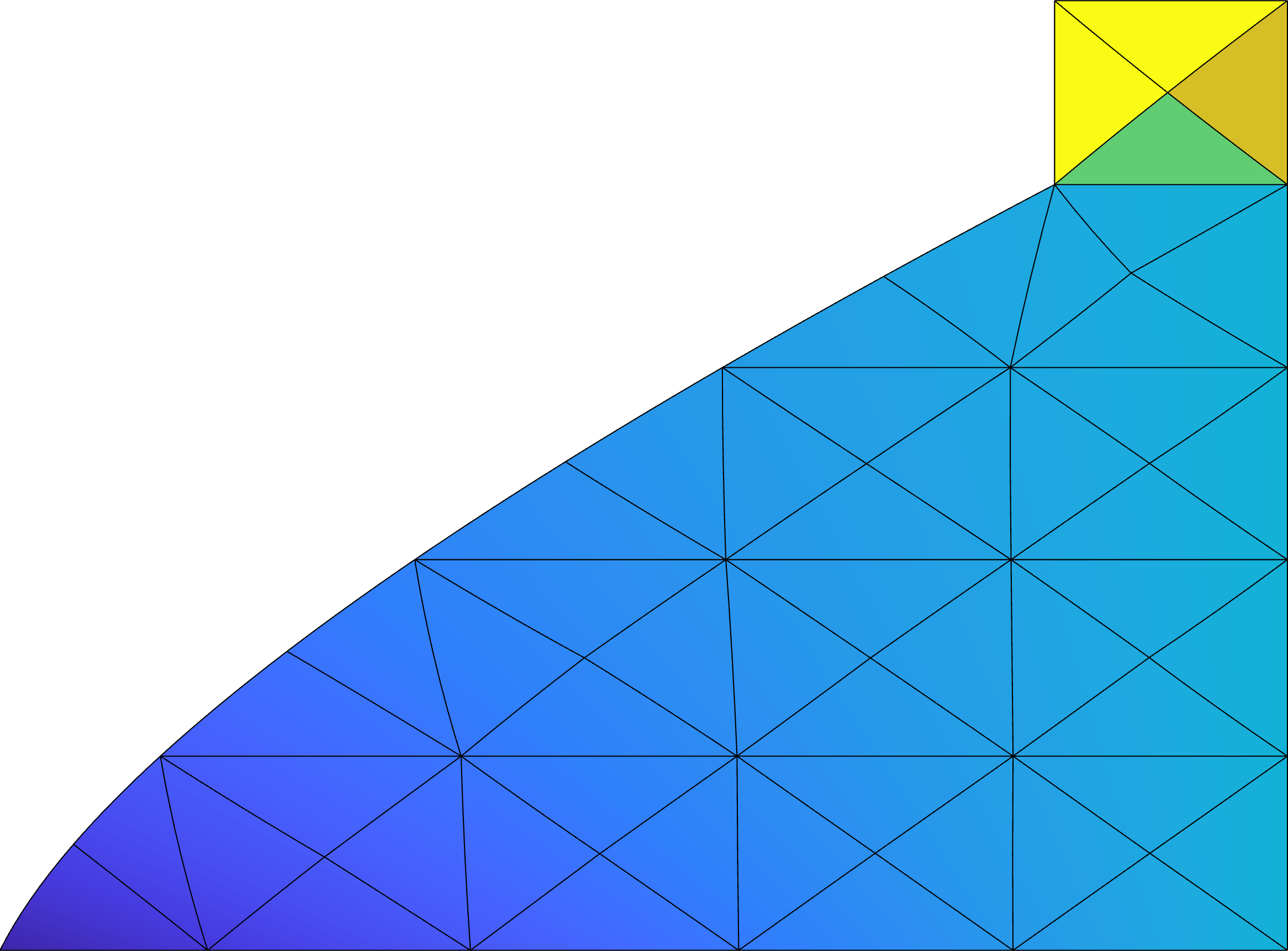};

\nextgroupplot[axis equal image, width=0.27\textwidth, xlabel={$x$}, ylabel={$t$}, xtick={-0.2, 1}, ytick={0, 0.738}, xticklabels={}, yticklabels={}, xmin=-0.2, xmax=1, ymin=0, ymax=0.738]
\addplot []
graphics [xmin=0,xmax=1,ymin=0,ymax=0.15] { 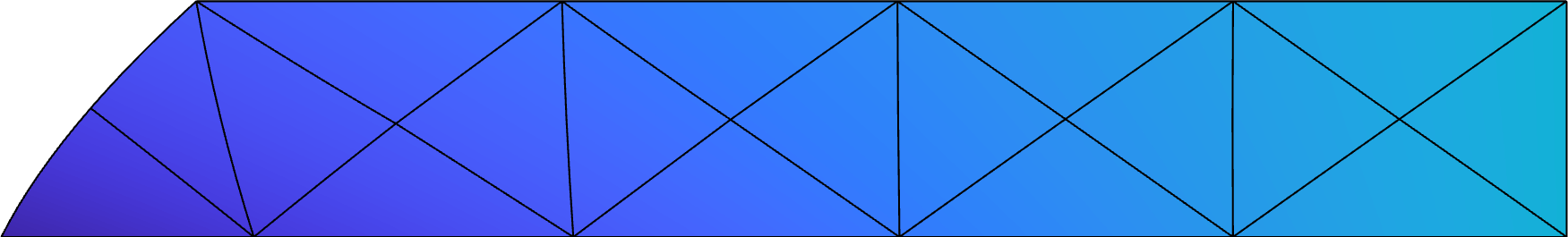};

\nextgroupplot[axis equal image, width=0.27\textwidth, xlabel={$x$}, xtick={-0.2, 1}, ytick={0, 1}, xticklabels={}, yticklabels={}, xmin=-0.2, xmax=1, ymin=0, ymax=0.738]
\addplot []
graphics [xmin=0,xmax=1,ymin=0,ymax=0.303] { 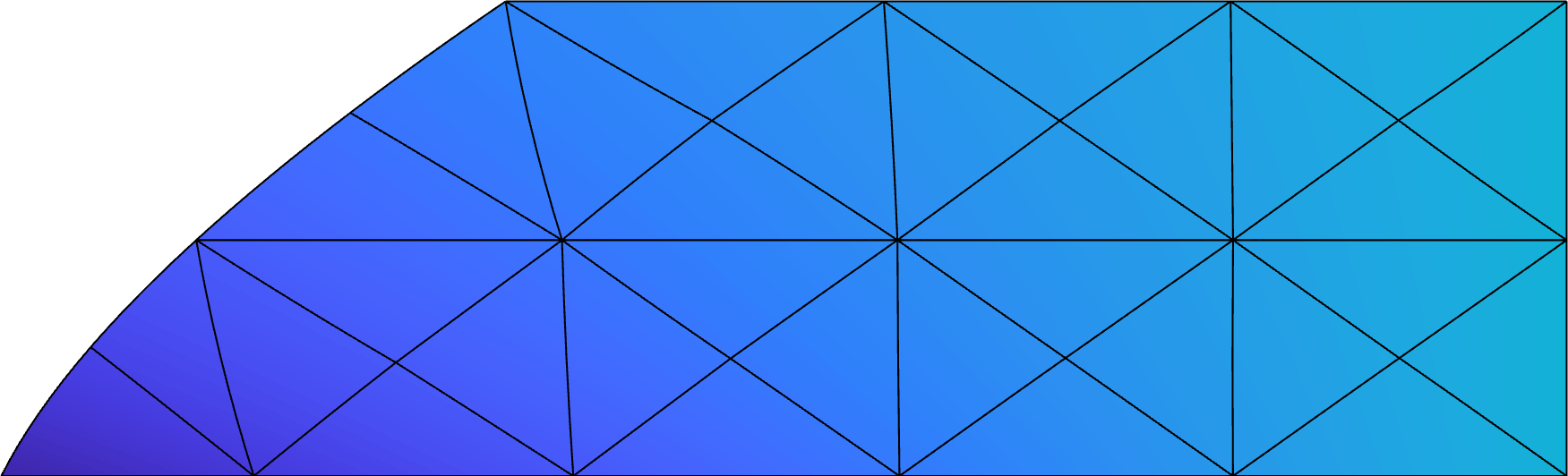};

\nextgroupplot[axis equal image, width=0.27\textwidth, xlabel={$x$}, xtick={-0.2, 1}, ytick={0, 1}, xticklabels={}, yticklabels={}, xmin=-0.2, xmax=1, ymin=0, ymax=0.738]
\addplot []
graphics [xmin=0,xmax=1,ymin=0,ymax=0.453] { 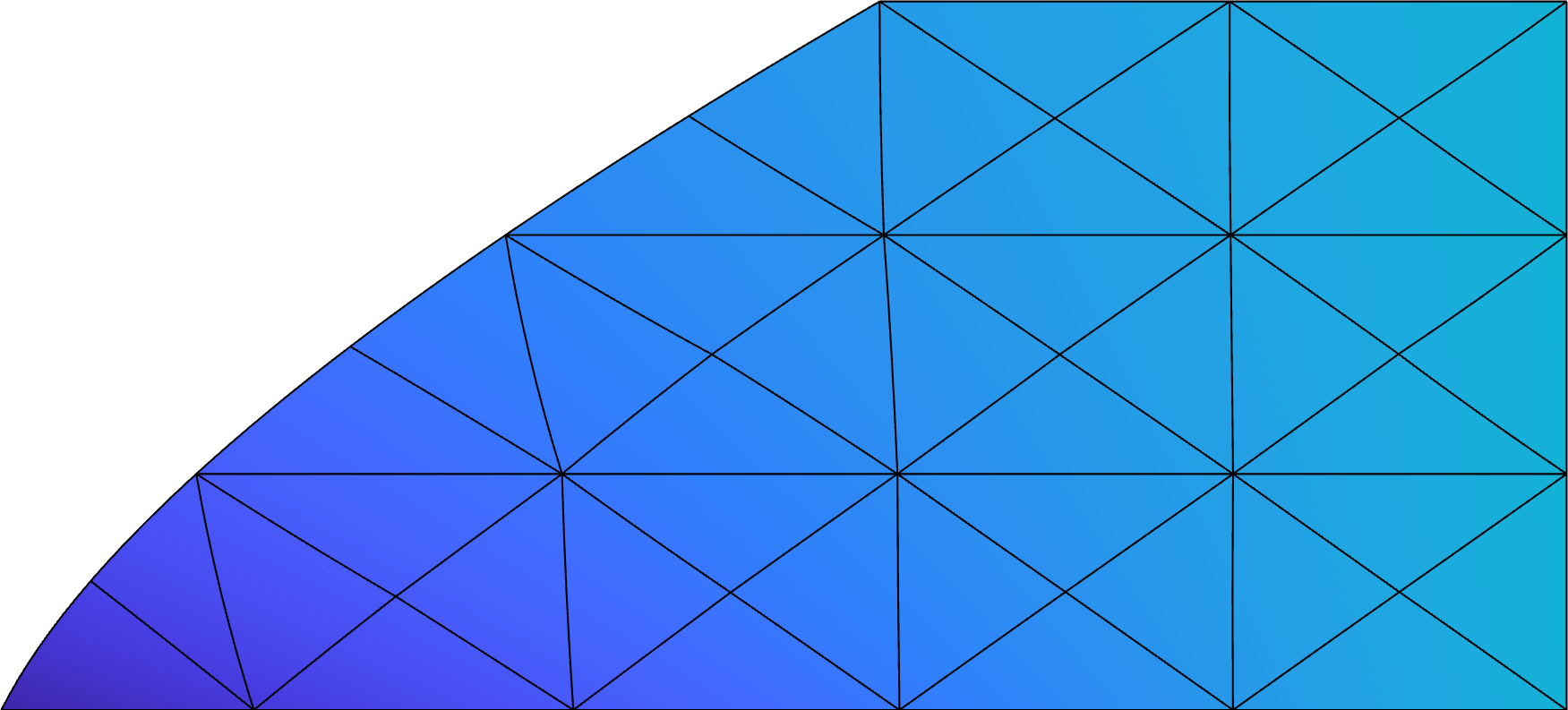};

\nextgroupplot[axis equal image, width=0.27\textwidth, xlabel={$x$}, xtick={-0.2, 1}, ytick={0, 0.5, 1}, xticklabels={}, yticklabels={}, xmin=-0.2, xmax=1, ymin=0, ymax=0.738]
\addplot []
graphics [xmin=0,xmax=1,ymin=0,ymax=0.595] { 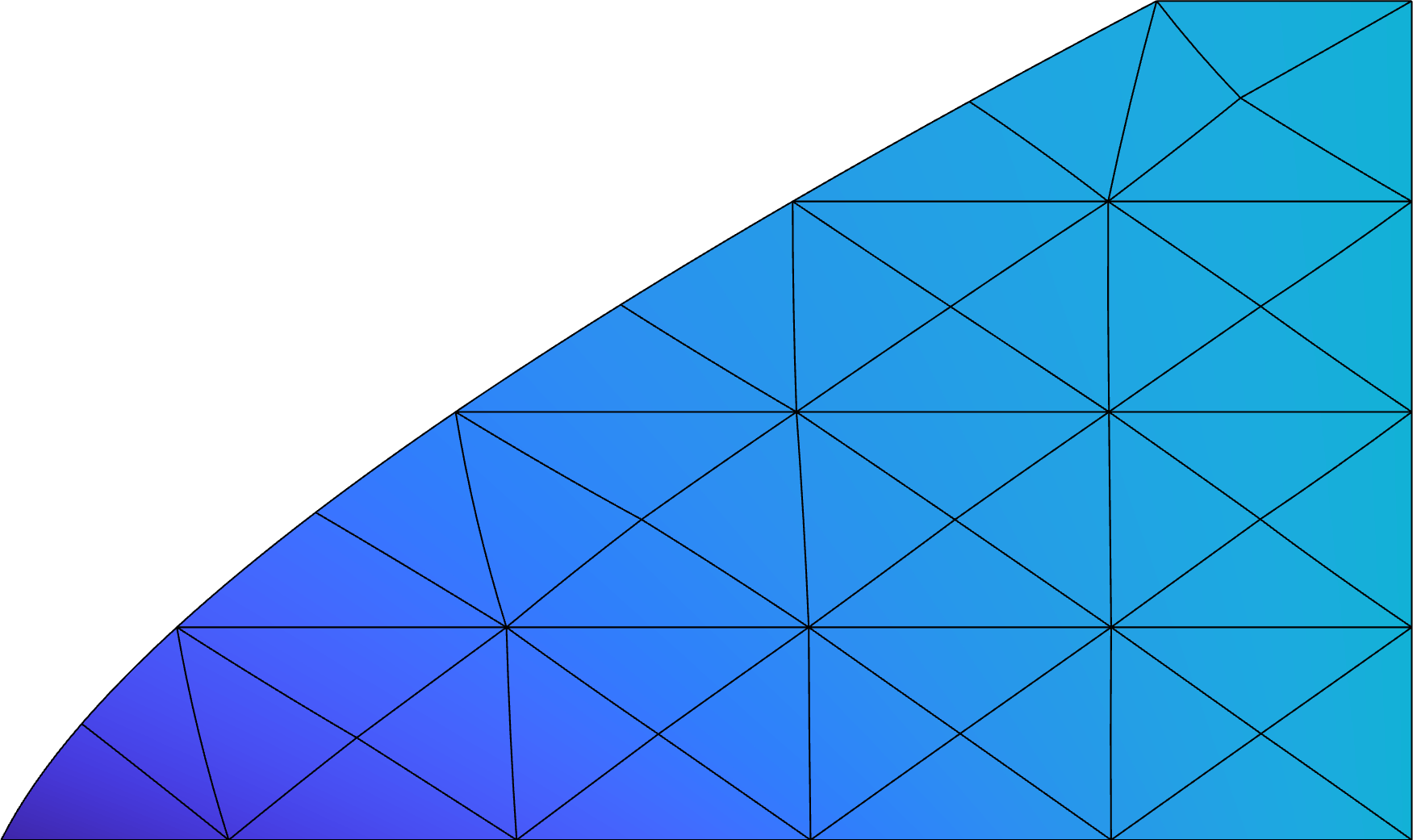};

\nextgroupplot[axis equal image, width=0.27\textwidth, xlabel={$x$}, xtick={-0.2, 1}, ytick={0, 0.5, 1}, xticklabels={}, yticklabels={}, xmin=-0.2, xmax=1, ymin=0, ymax=0.738]
\addplot []
graphics [xmin=0,xmax=1,ymin=0,ymax=0.69] { 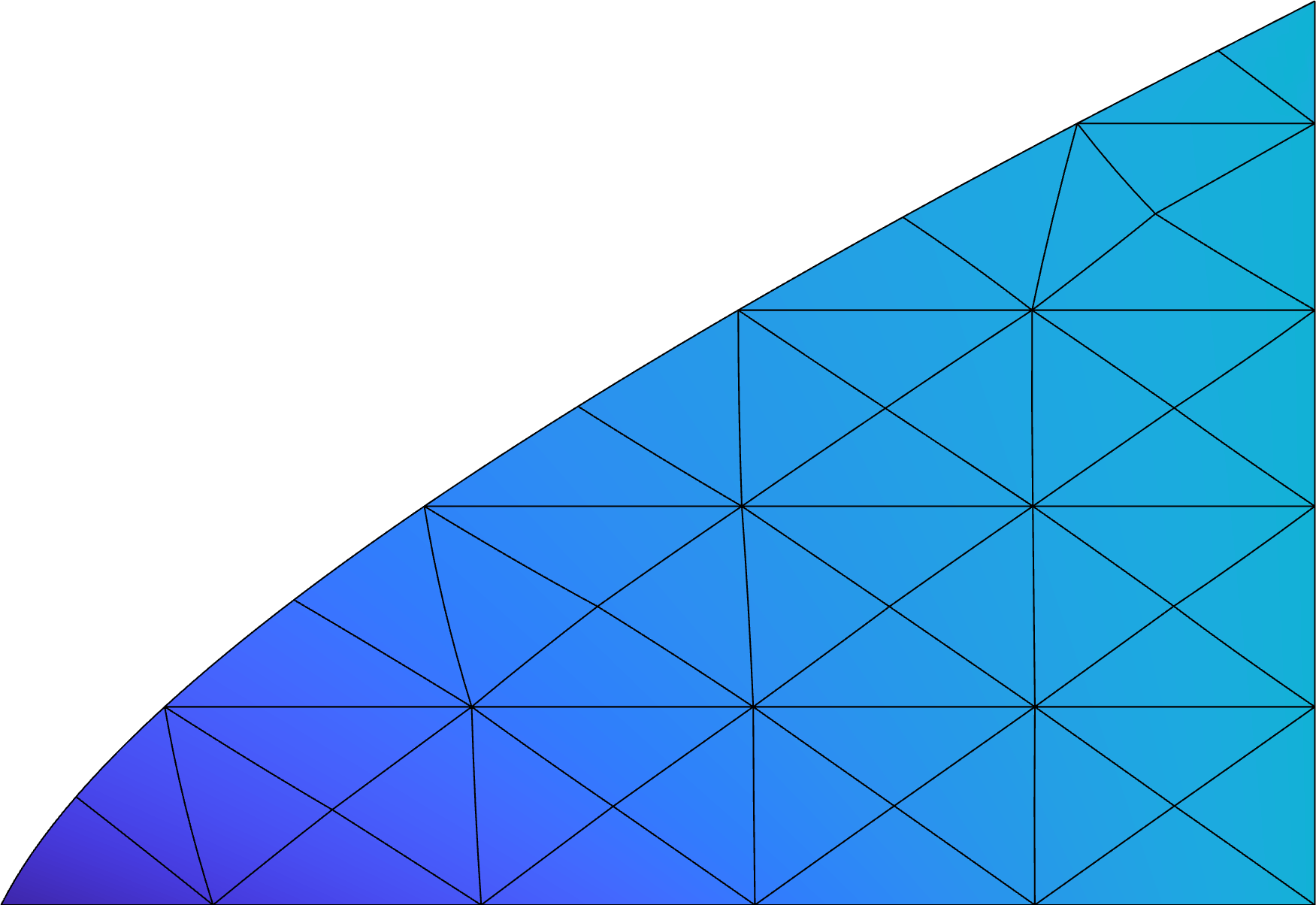};

\end{groupplot}\end{tikzpicture}}
  \caption{HOIST solution of accelerating shock inviscid Burgers' problem with the shock boundary condition approach (Section~\ref{sec:ist:shkbc}) including the initial (\textit{top row}) and converged (\textit{bottom row}) mesh and flow solution. Colorbar in Figure~\ref{fig:iburg_nrl_slabs_1_thru_2}.}
  \label{fig:iburg_nrl_slabs_bnd_rmv}
\end{figure}

Next, we solve the same problem with $S = 5$ slabs using the shock boundary condition approach \cite{huang2023high} (Section~\ref{sec:ist:shkbc}). The first space-time slab is generated by extruding spatial mesh with $6$ elements to produce a space-time mesh with $24$ elements and implicit shock tracking is used to align it with the accelerating shock. The elements to the left of the shock are removed and the remaining domain is extruded to form the next space-time slab. Because the state to the left of the shock is known, it is applied as a boundary condition to the leftmost elements. All deformation degrees of freedom associated with the shock boundary are free so it can move as necessary to track the discontinuity as it propagates in space-time. This process is repeated until the shock leaves the domain (Figure~\ref{fig:iburg_nrl_slabs_bnd_rmv}). This approach eliminates over half of the space-time elements and simplifies the shock tracking problem, which improves efficiency of the simulation.

\subsubsection{One-dimensional shock formation and merge}
\label{sec:numexp:burg:shkfrm}
Next, we consider a Burgers' equation (\ref{eqn:burg}) ($\beta = 1$) problem with
a smooth initial condition that forms two distinct shock waves that eventually merge.
The spatial domain is $\Omega_x \coloneqq (-1, 1)$ and
temporal domain is $\Tcal \coloneqq (0, 1)$ with initial condition
\begin{equation}
  w_0 : x \mapsto
   1.2 \exp\left(-\frac{(x + 0.5)^2}{0.025}\right)-\exp\left(-\frac{(x - 0.5)^2}{0.025}\right).
\end{equation}
We use this problem to demonstrate the ability of the slab-based HOIST method to
resolve shock formation and intersection. We also show mesh adaptation
(Section~\ref{sec:ist:amr}) in the HOIST method helps
adequately resolve the shock formation, i.e., the point at which
the continuous features abruptly become discontinuities (wave steepening).

\begin{figure}
	\centering
 	\raisebox{-0.5\height}{\begin{tikzpicture}
\begin{groupplot} [
group style={group size = 2 by 5, horizontal sep = 0.5cm, vertical sep = 0.5cm},
title style={at={(current bounding box.north west)}, anchor=west}]
\nextgroupplot[axis equal image, width=0.5\textwidth, xtick={-1, 0, 1}, ytick={0, 0.45, 0.9}, xticklabels={}, yticklabels={0, 0.45, 0.9}, ylabel={$t$}, xmin=-1, xmax=1, ymin=0, ymax=0.9]
\addplot []
graphics [xmin=-1,xmax=1,ymin=0,ymax=0.9] { 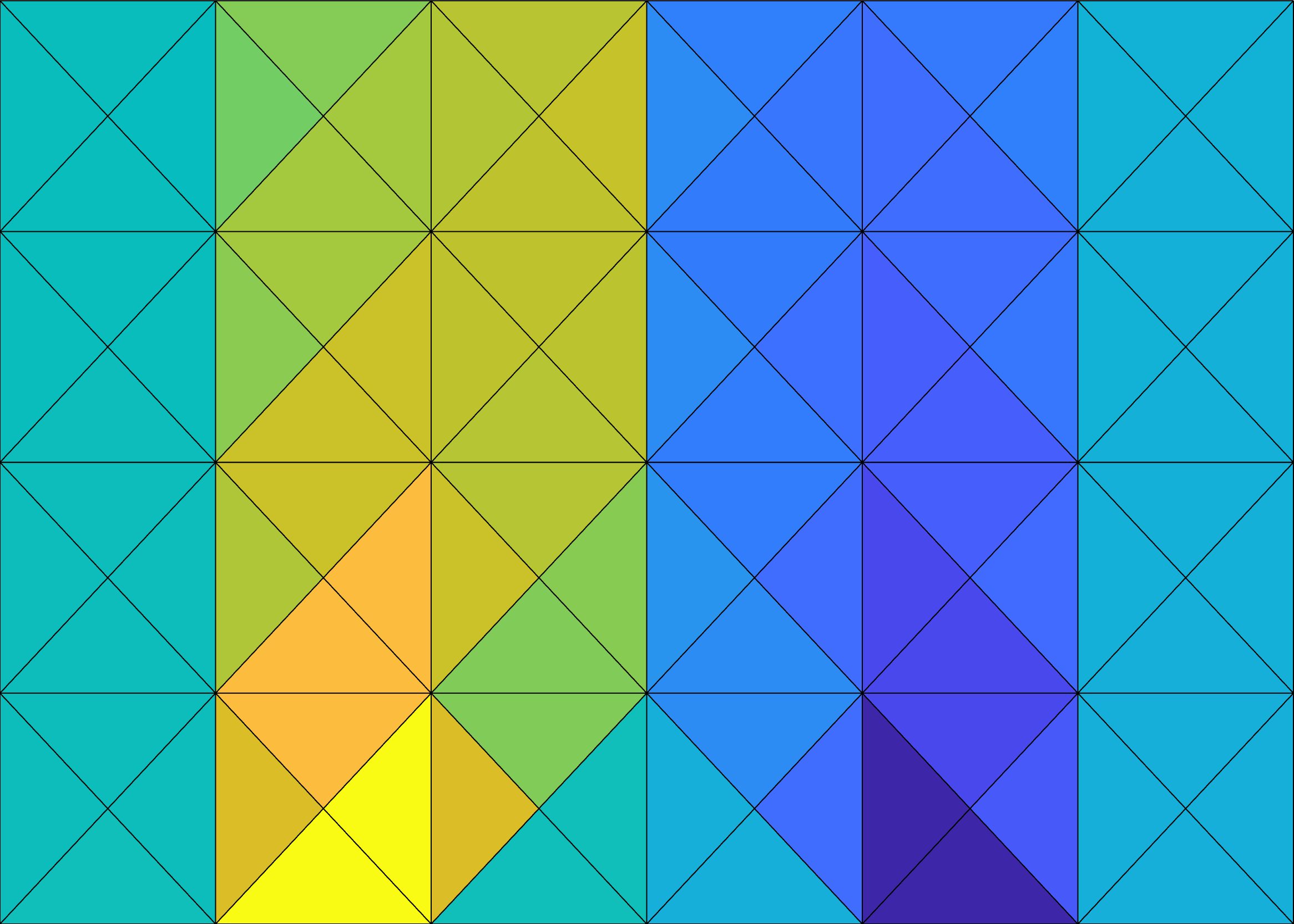};

\nextgroupplot[axis equal image, width=0.5\textwidth, xtick={-1, 0, 1}, ytick={0, 0.45, 0.9}, xticklabels={}, yticklabels={}, xmin=-1, xmax=1, ymin=0, ymax=0.9]
\addplot []
graphics [xmin=-1,xmax=1,ymin=0,ymax=0.9] { 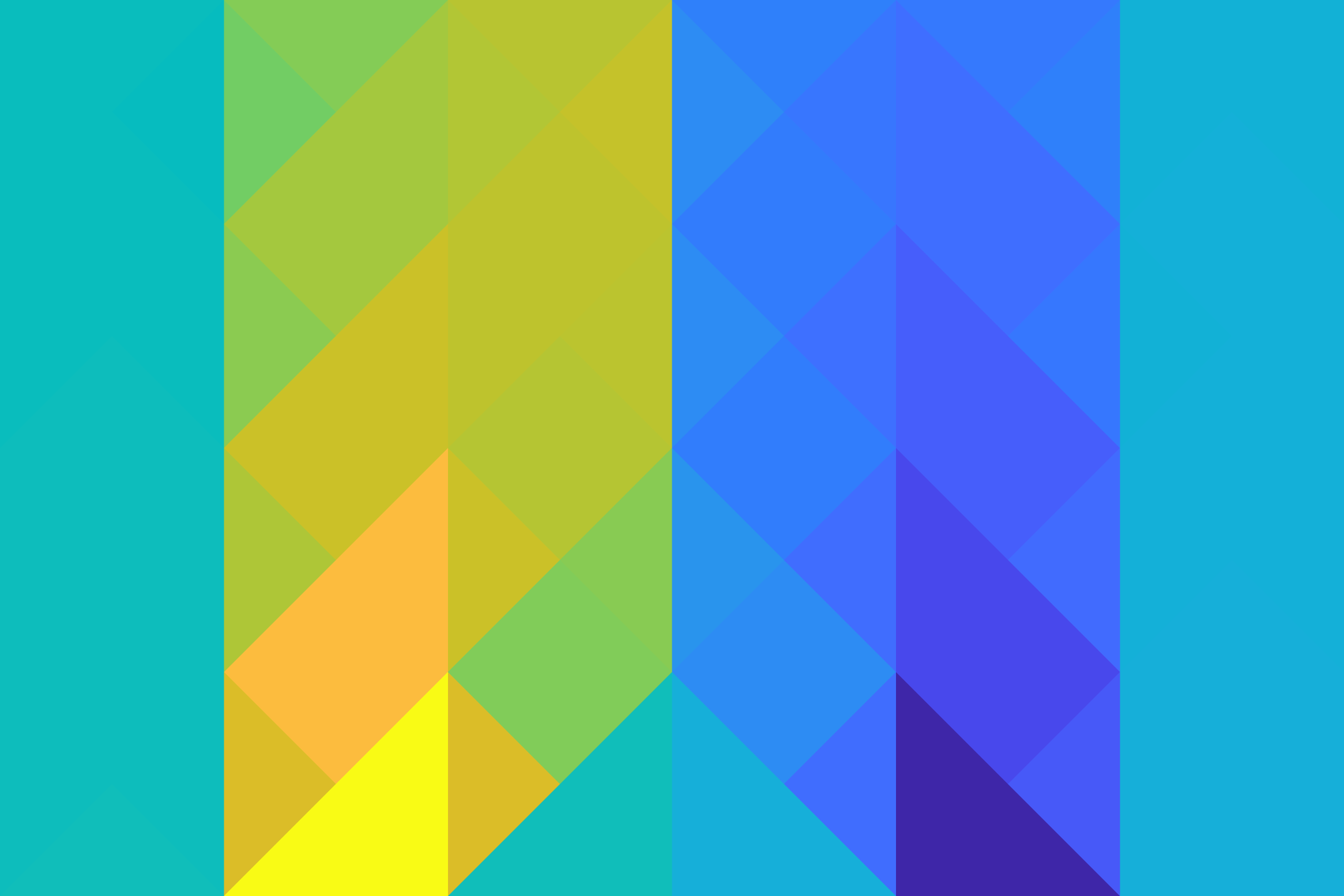};

\nextgroupplot[axis equal image, width=0.5\textwidth, xtick={-1, 0, 1}, ytick={0, 0.45, 0.9}, xticklabels={}, yticklabels={0, 0.45, 0.9}, ylabel={$t$}, xmin=-1, xmax=1, ymin=0, ymax=0.9]
\addplot []
graphics [xmin=-1,xmax=1,ymin=0,ymax=0.9] { 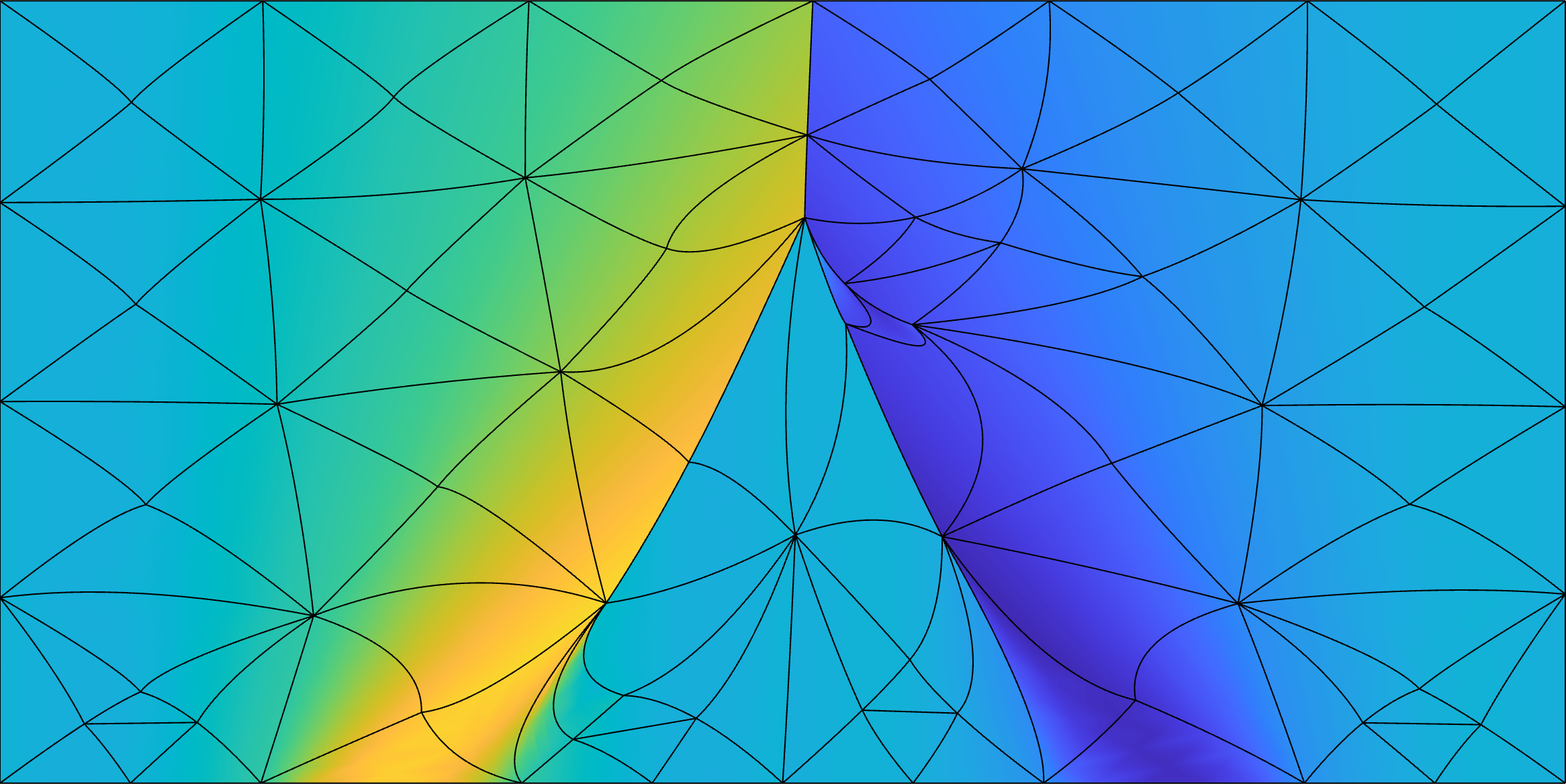};

\nextgroupplot[axis equal image, width=0.5\textwidth, xtick={-1, 0, 1}, ytick={0, 0.45, 0.9}, xticklabels={}, yticklabels={}, xmin=-1, xmax=1, ymin=0, ymax=0.9]
\addplot []
graphics [xmin=-1,xmax=1,ymin=0,ymax=0.9] { 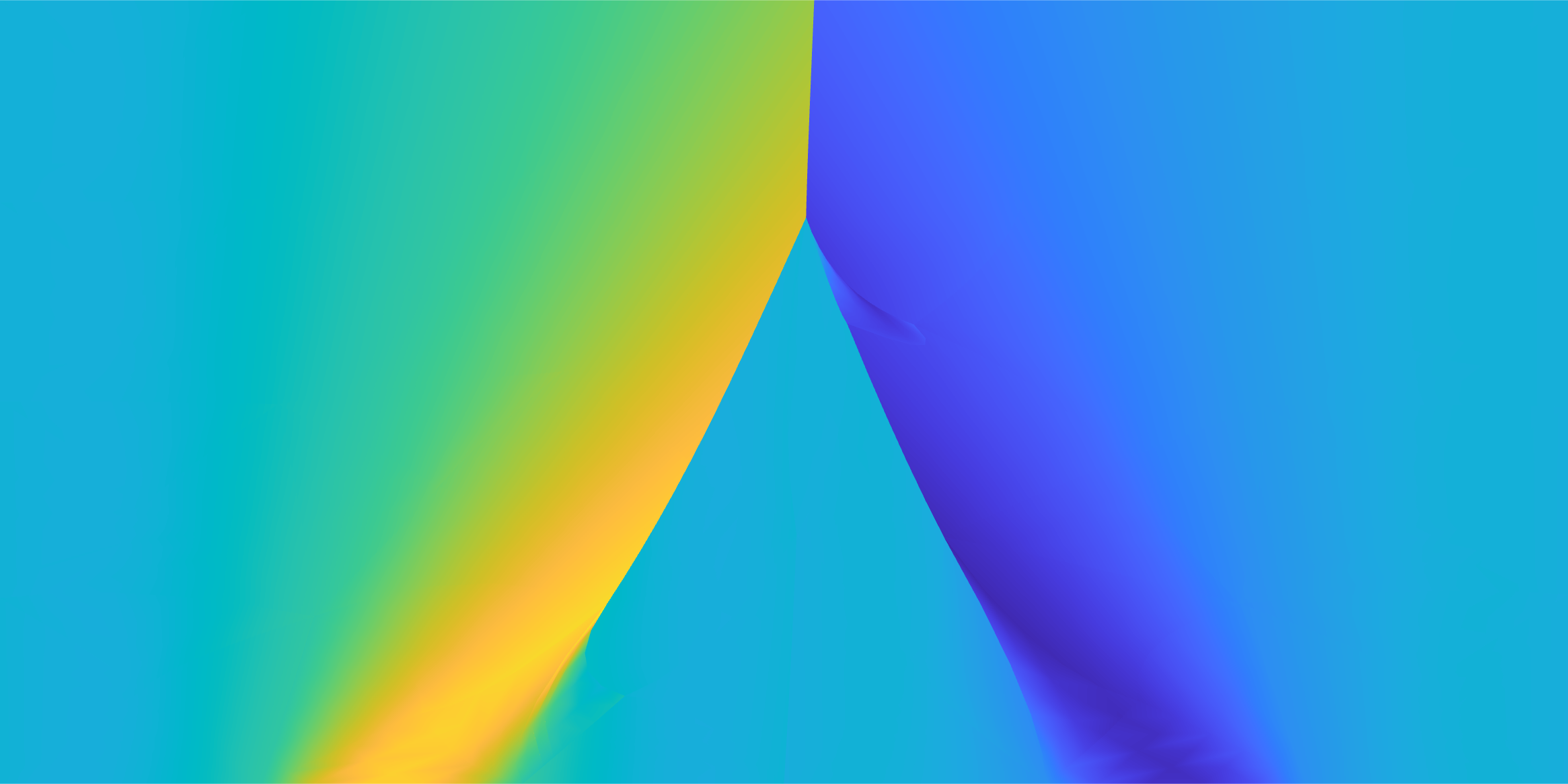};

\nextgroupplot[axis equal image, width=0.5\textwidth, xtick={-1, 0, 1}, ytick={0, 0.45, 0.9}, xticklabels={}, yticklabels={0, 0.45, 0.9}, ylabel={$t$}, xmin=-1, xmax=1, ymin=0, ymax=0.9]
\addplot []
graphics [xmin=-1,xmax=1,ymin=0,ymax=0.9] { 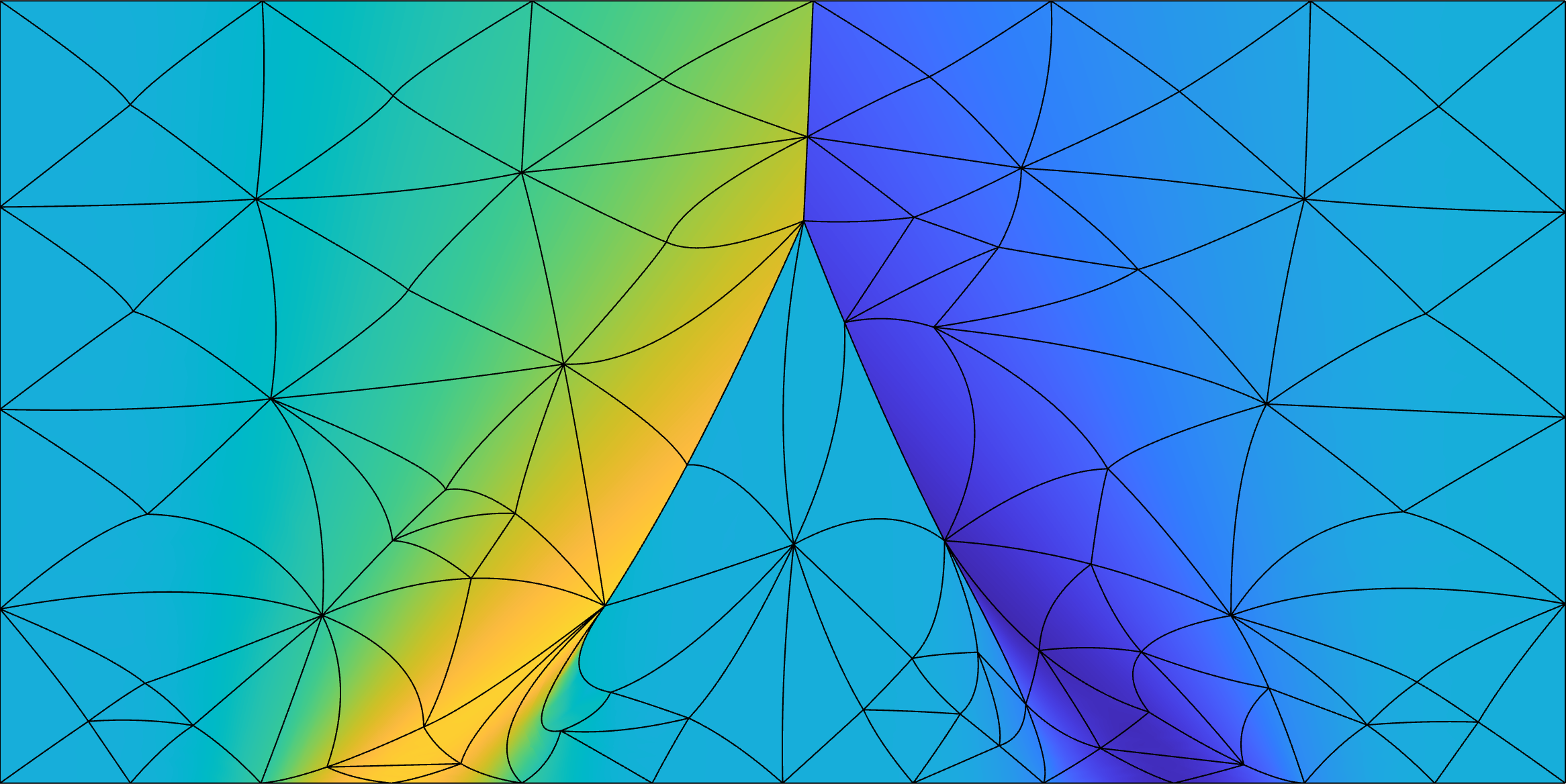};

\nextgroupplot[axis equal image, width=0.5\textwidth, xtick={-1, 0, 1}, ytick={0, 0.45, 0.9}, xticklabels={}, yticklabels={}, xmin=-1, xmax=1, ymin=0, ymax=0.9]
\addplot []
graphics [xmin=-1,xmax=1,ymin=0,ymax=0.9] { 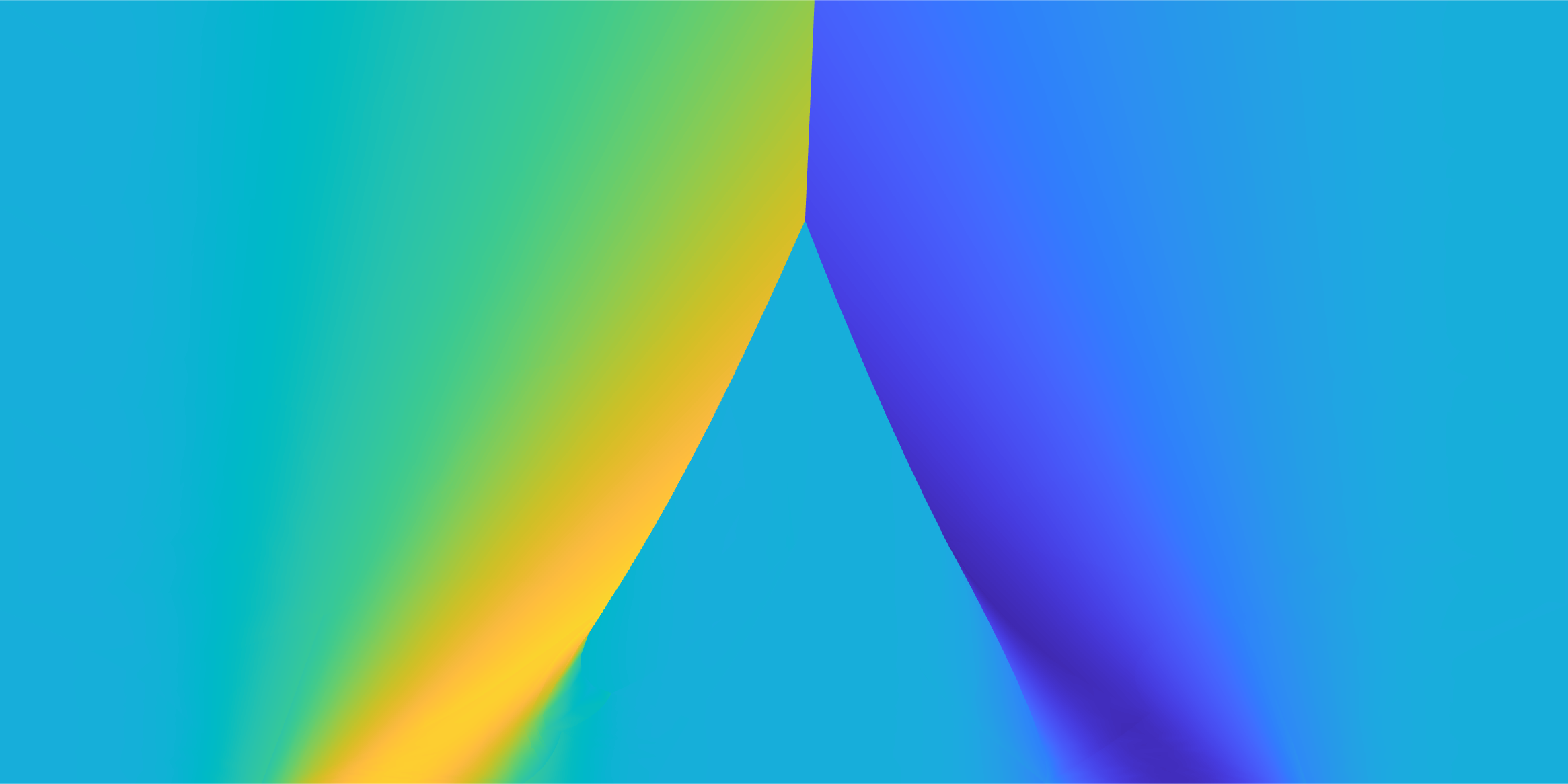};

\nextgroupplot[axis equal image, width=0.5\textwidth, xtick={-1, 0, 1}, ytick={0, 0.45, 0.9}, xticklabels={}, yticklabels={0, 0.45, 0.9}, ylabel={$t$}, xmin=-1, xmax=1, ymin=0, ymax=0.9]
\addplot []
graphics [xmin=-1,xmax=1,ymin=0,ymax=0.9] { 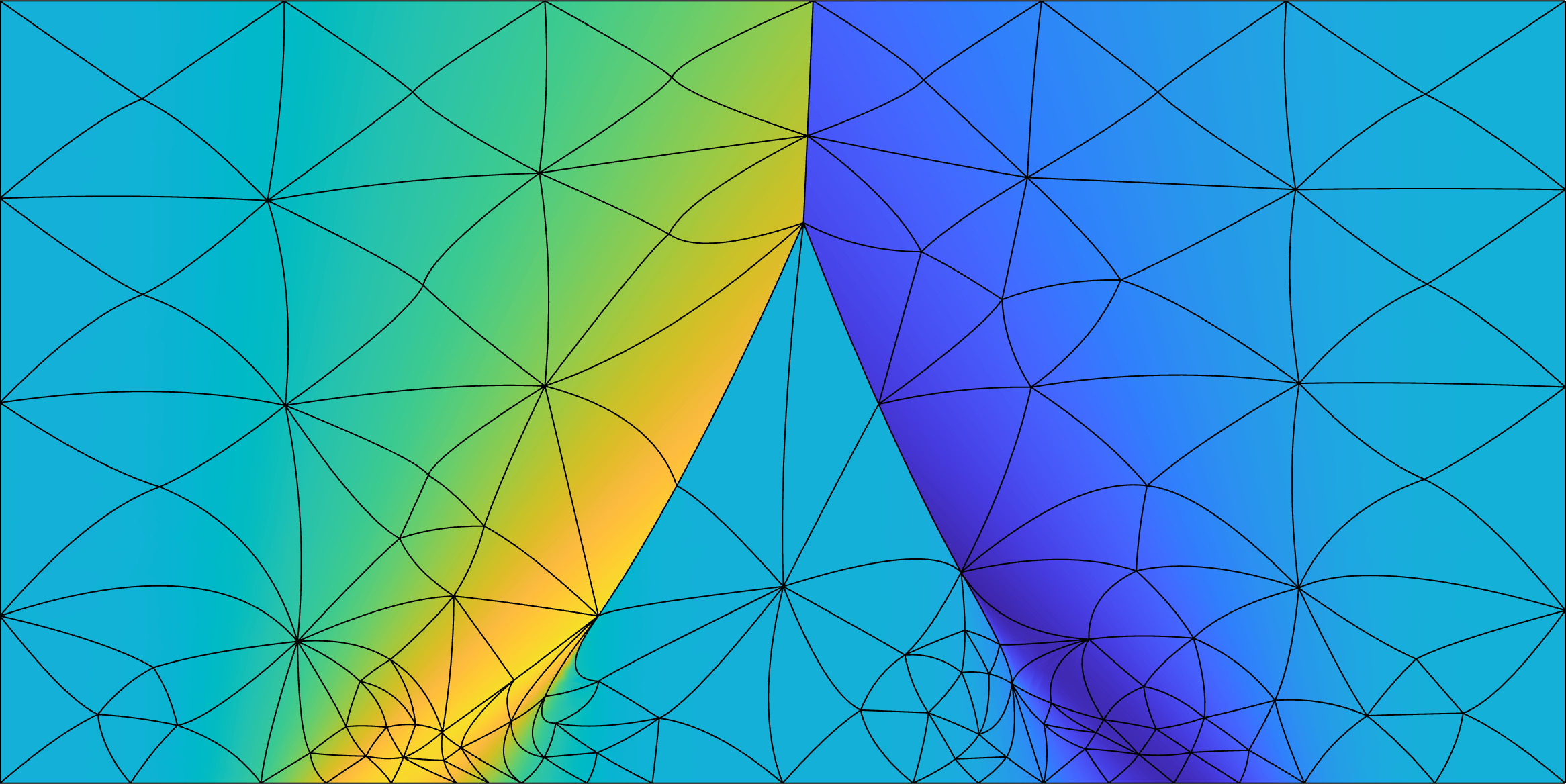};

\nextgroupplot[axis equal image, width=0.5\textwidth, xtick={-1, 0, 1}, ytick={0, 0.45, 0.9}, xticklabels={}, yticklabels={}, xmin=-1, xmax=1, ymin=0, ymax=0.9]
\addplot []
graphics [xmin=-1,xmax=1,ymin=0,ymax=0.9] { 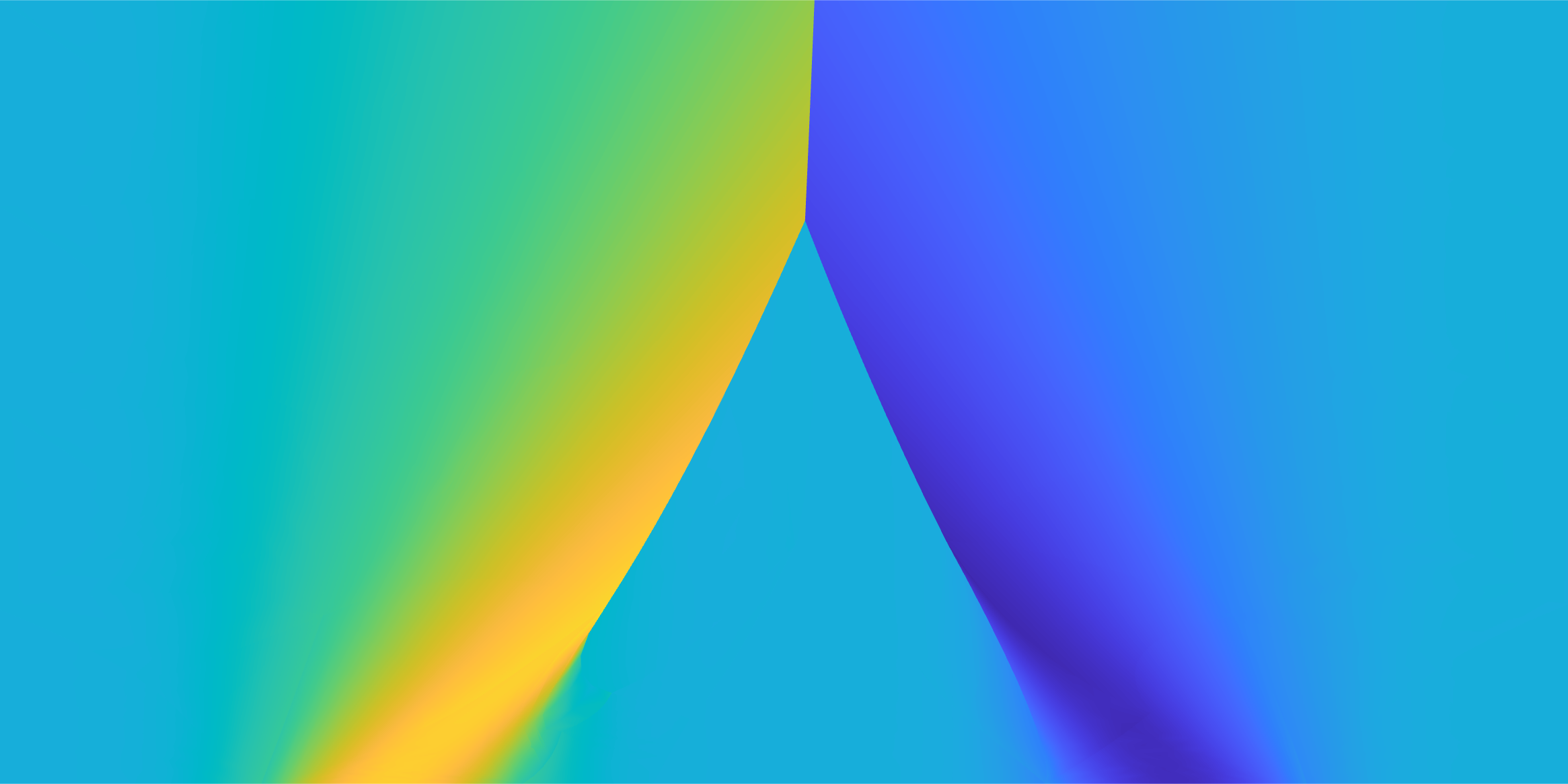};

\nextgroupplot[axis equal image, width=0.5\textwidth, xtick={-1, -0.5, 0, 0.5, 1}, ytick={0, 0.45, 0.9}, xticklabels={-1, -0.5, 0, 0.5, 1}, yticklabels={0, 0.45, 0.9}, xlabel={$x$}, ylabel={$t$}, xmin=-1, xmax=1, ymin=0, ymax=0.9]
\addplot []
graphics [xmin=-1,xmax=1,ymin=0,ymax=0.9] { 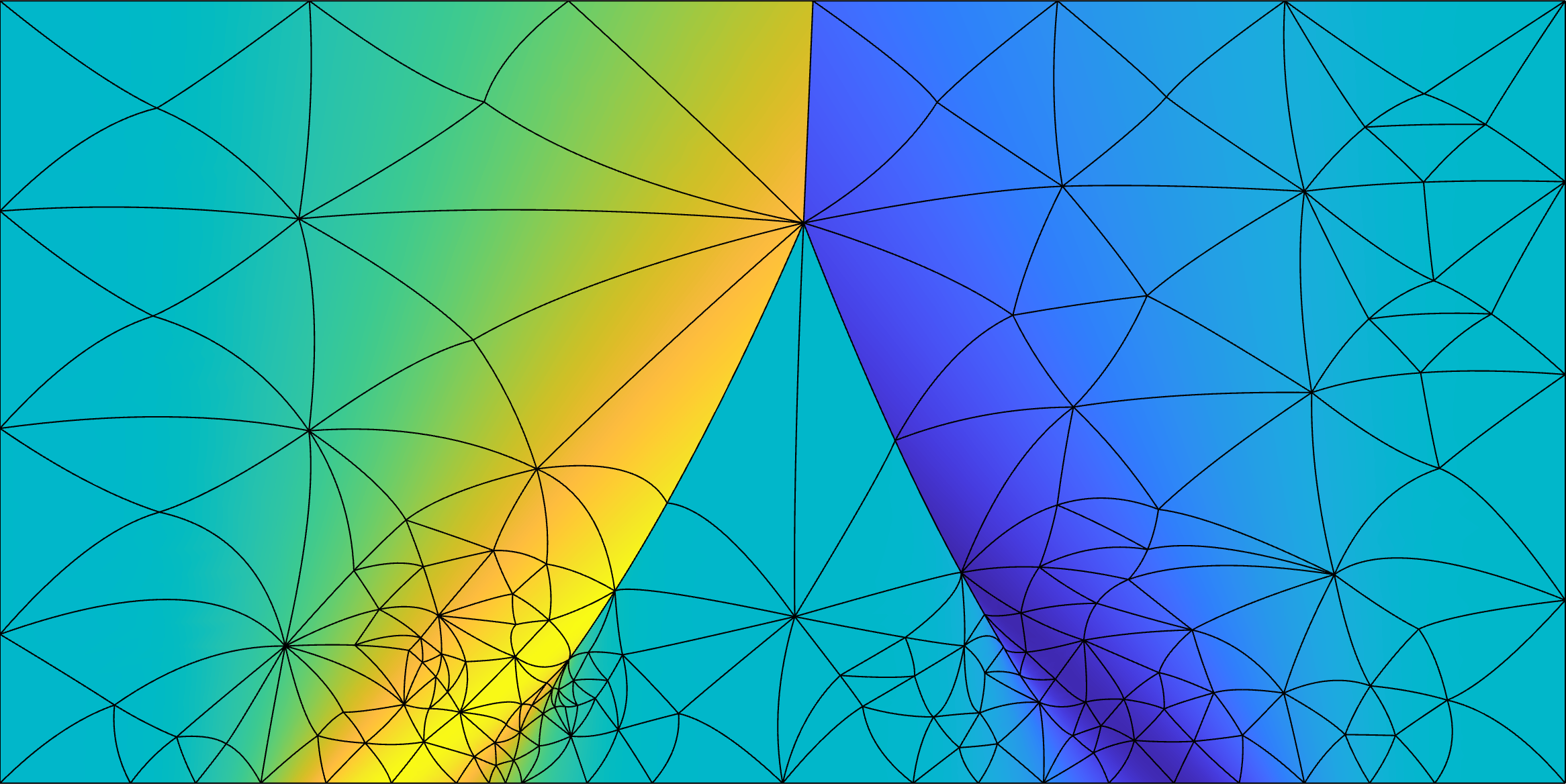};

\nextgroupplot[axis equal image, width=0.5\textwidth, xtick={-1, -0.5, 0, 0.5, 1}, ytick={0, 0.45, 0.9}, xticklabels={-1, -0.5, 0, 0.5, 1}, yticklabels={}, xlabel={$x$}, xmin=-1, xmax=1, ymin=0, ymax=0.9]
\addplot []
graphics [xmin=-1,xmax=1,ymin=0,ymax=0.9] { 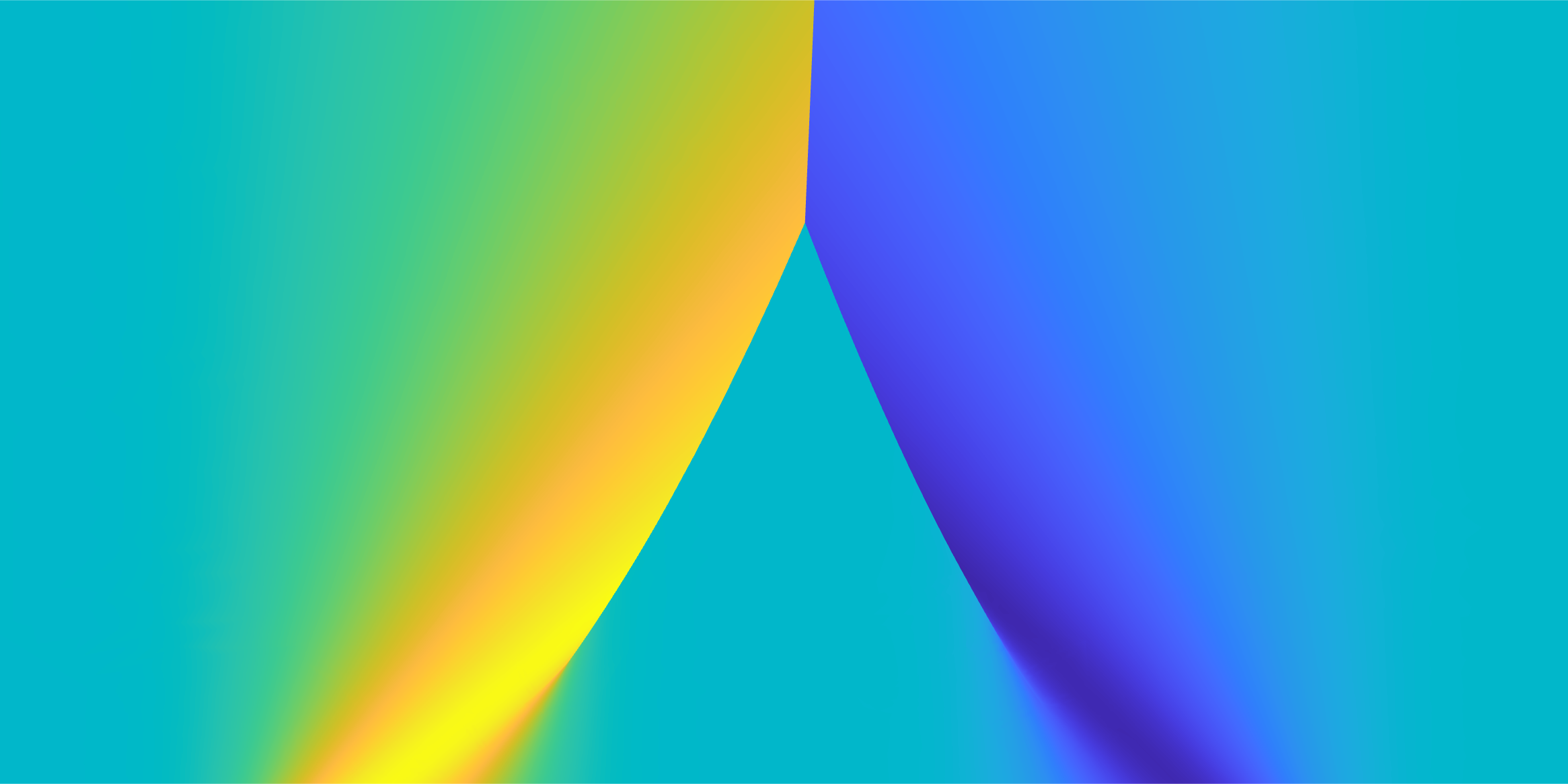};

\end{groupplot}\end{tikzpicture}}
	\colorbarMatlabParula{-1}{-0.5}{0}{0.5}{1}
 	\caption{HOIST solutions for the inviscid Burgers' problem with shock formation including the initial mesh and flow solution (\textit{top row}), as well as the mesh and flow solution after one (\textit{second row}), two (\textit{third row}), four (\textit{fourth row}), and eight (\textit{fifth row}) rounds of adaptive mesh refinement.}
 	\label{fig:iburg_shkfrm_slabs_1_thru_2}
\end{figure}

\begin{figure}
	\centering
 	\raisebox{-0.5\height}{\input{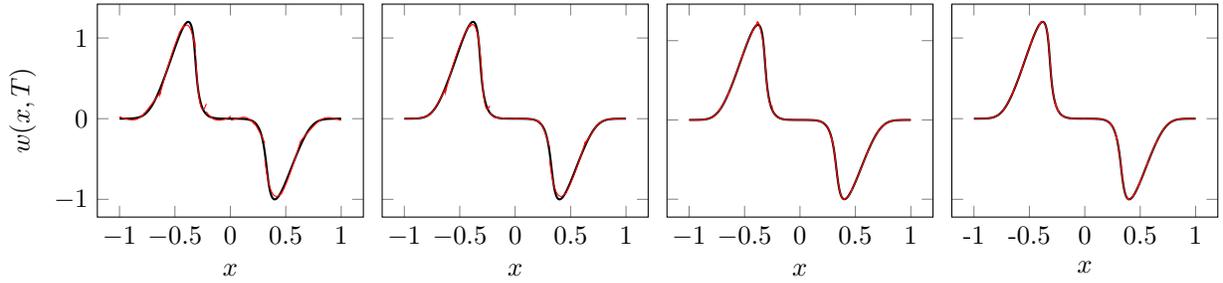}}
 	\caption{HOIST solution (\ref{line:u11_solid})
  for the inviscid Burgers' problem with shock formation at temporal
  slices $t = 0.1$ (just before shock formation) after one (\textit{left}),
  two (\textit{middle left}), four (\textit{middle right}), and eight
  (\textit{right}) rounds of adaptive mesh refinement. A reference solution
  computed with a highly refined second-order finite volume method
  is shown in (\ref{line:ucfd1}).}
 	\label{fig:iburg_shkfrm_slices}
\end{figure}

To study adaptive mesh refinement, we consider a single slab $S = 1$ with a mesh of $96$
quadratic ($p = q = 2$) space-time triangles. Space-time singularities form as the continuous features steepen into shocks. The adaptive refinement procedure increases resolution in these regions of the domain with minimal refinement in regions where the flow is well-resolved after eight rounds of mesh refinement (Figure~\ref{fig:iburg_shkfrm_slabs_1_thru_2}). The converged, fully adapted mesh contains $207$ triangles, and compares well to a reference solution computed with a highly refined second-order finite volume scheme at a temporal slice immediately before shock formation (Figure~\ref{fig:iburg_shkfrm_slices}).

	
\subsubsection{Propagating, curved shock in two dimensions}
\label{sec:numexp:burg:2d}

Finally, we consider Burgers' equation (\ref{eqn:burg}) ($\beta = (1,0)$)
with a quadratic initial condition \cite{2021_shi_methlines}
\begin{equation}
 w_0 : (x_1, x_2) \mapsto
 \begin{cases}
  (\frac{4}{3}(x_1 + 0.75))(0.5 - 2(x_2^2 - 0.25)) & x \in (-1,0)\times (-0.5,0.5) \\
  0 & \text{else}
 \end{cases}
\end{equation}
over the two-dimensional spatial domain $\Omega_x \coloneqq (-1, 1)^2$
and temporal domain $\Tcal \coloneqq (0, 1.5)$. The initial condition is a
straight-sided shock with a quadratically varying conserved variable, which causes the shock to curve as it propagates.

\begin{figure}
	\centering
 	\raisebox{-0.5\height}{\begin{tikzpicture}
\begin{groupplot} [
group style={group size = 3 by 2, horizontal sep = 0.05cm, vertical sep = 0.5cm},
title style={at={(current bounding box.north west)}, anchor=west}]
\nextgroupplot[axis equal image, width=0.4\textwidth, xticklabels={}, yticklabels={}, ylabel={}, xlabel={}, xmin=-1, xmax=1, ymin=0, ymax=1.5]
\addplot []
graphics [xmin=-1,xmax=1,ymin=0,ymax=1.5] { 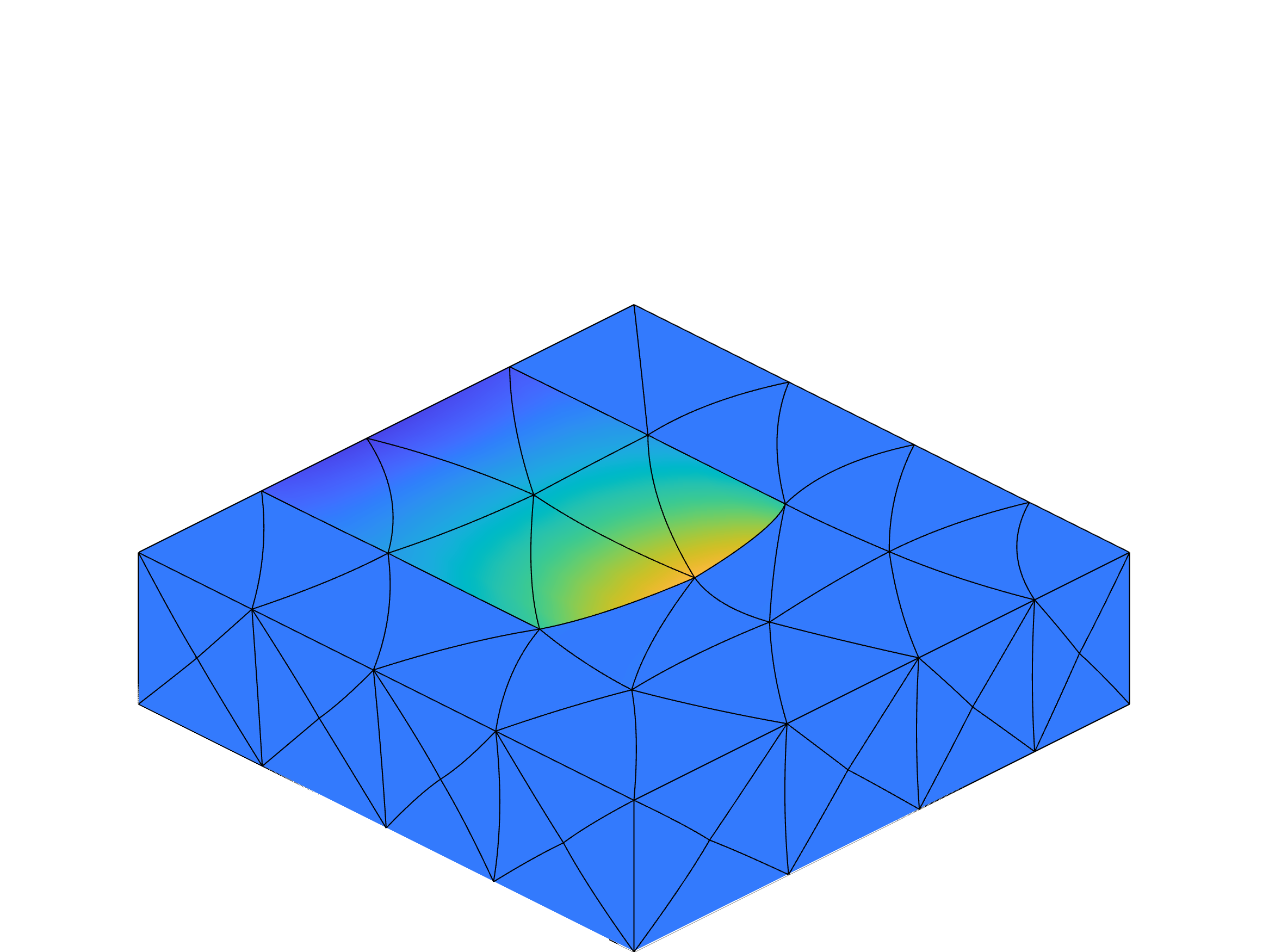};

\nextgroupplot[axis equal image, width=0.4\textwidth, xticklabels={}, yticklabels={}, xmin=-1, xmax=1, ymin=0, ymax=1.5]
\addplot []
graphics [xmin=-1,xmax=1,ymin=0,ymax=1.5] { 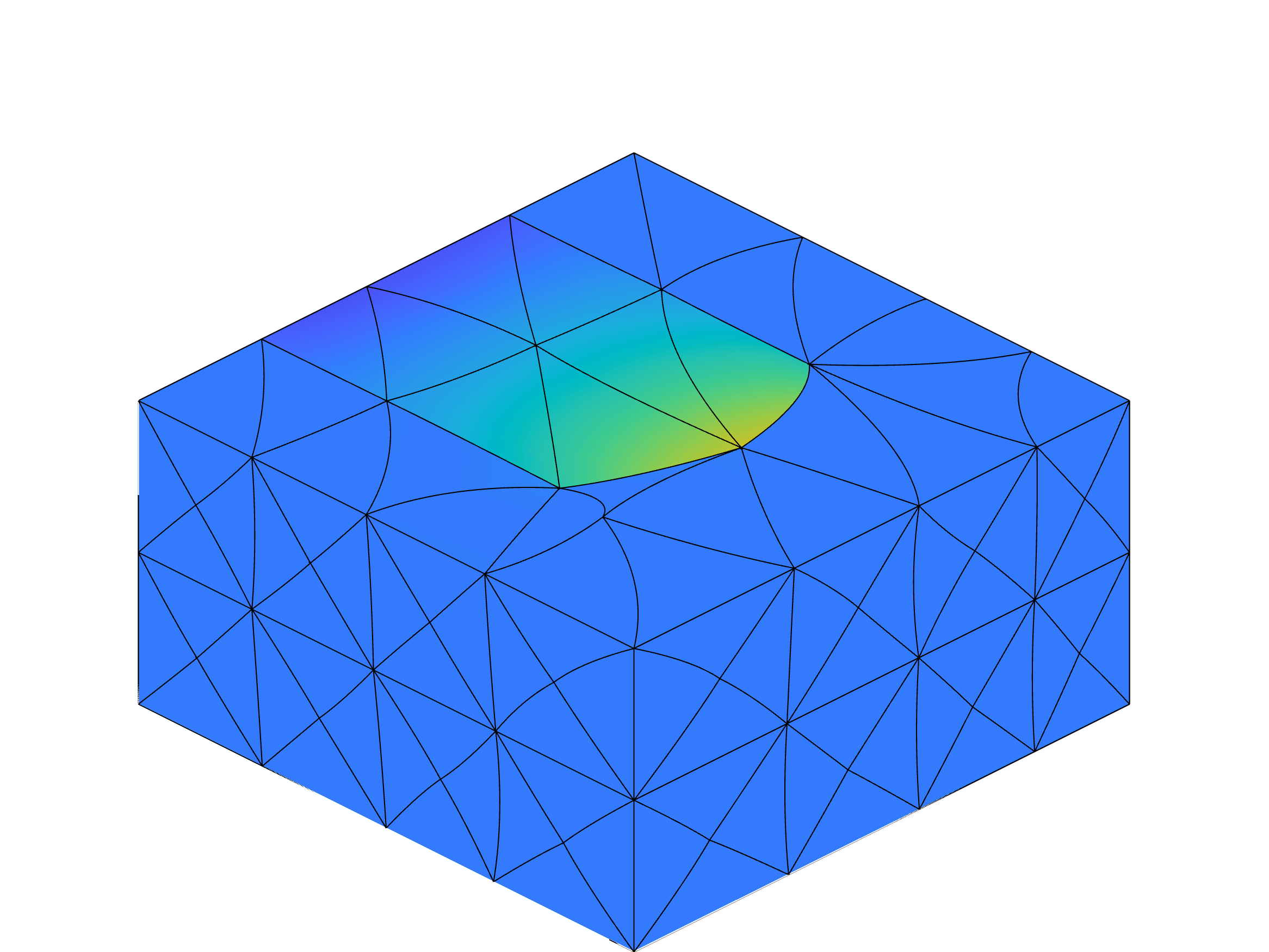};

\nextgroupplot[axis equal image, width=0.4\textwidth, xticklabels={}, yticklabels={}, ylabel={}, xlabel={}, xmin=-1, xmax=1, ymin=0, ymax=1.5]
\addplot []
graphics [xmin=-1,xmax=1,ymin=0,ymax=1.5] { 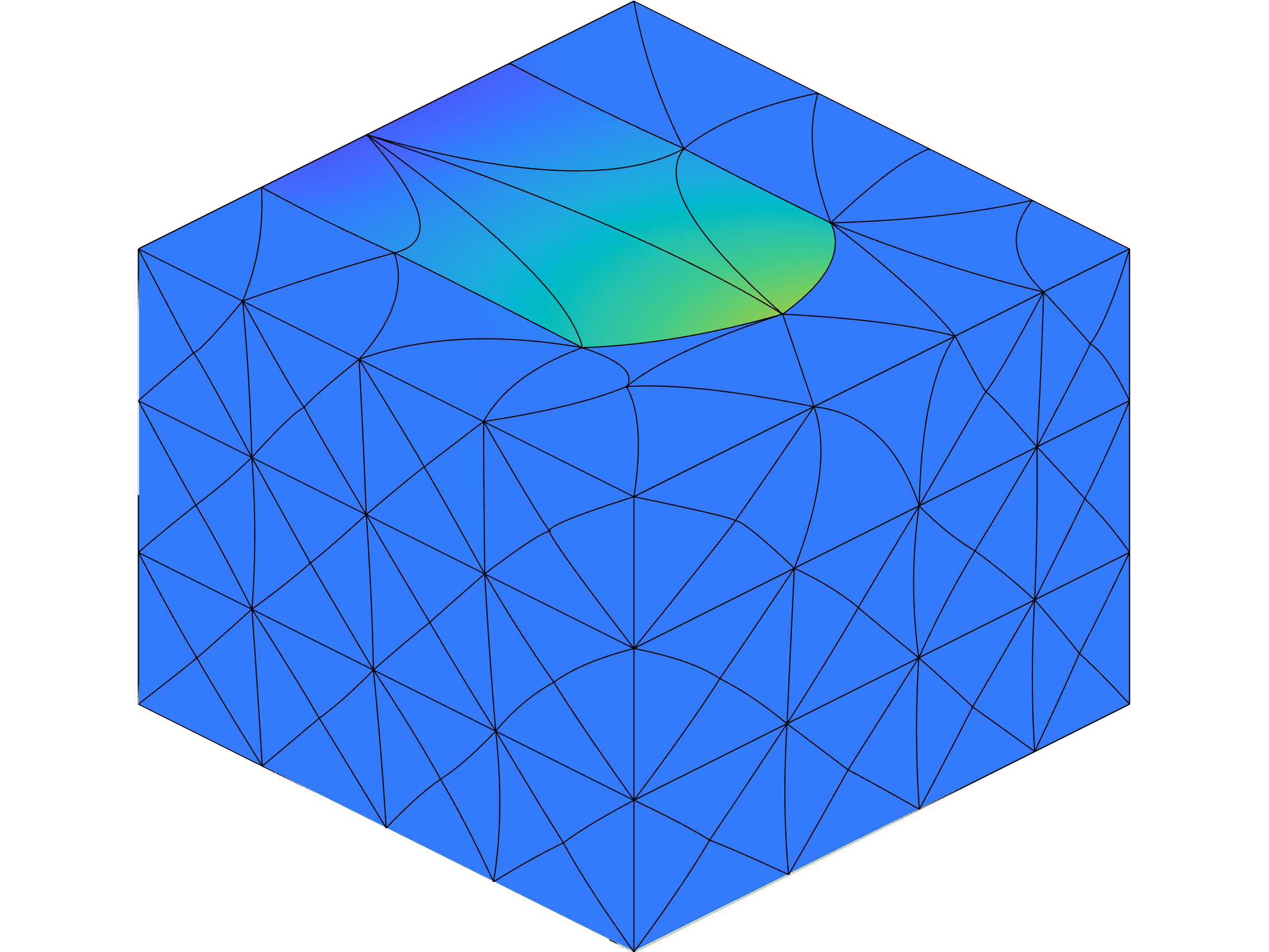};

\nextgroupplot[axis equal image, width=0.4\textwidth, xticklabels={}, yticklabels={}, ylabel={}, xlabel={}, xmin=-1, xmax=1, ymin=0, ymax=1.5]
\addplot []
graphics [xmin=-1,xmax=1,ymin=0,ymax=1.5] { 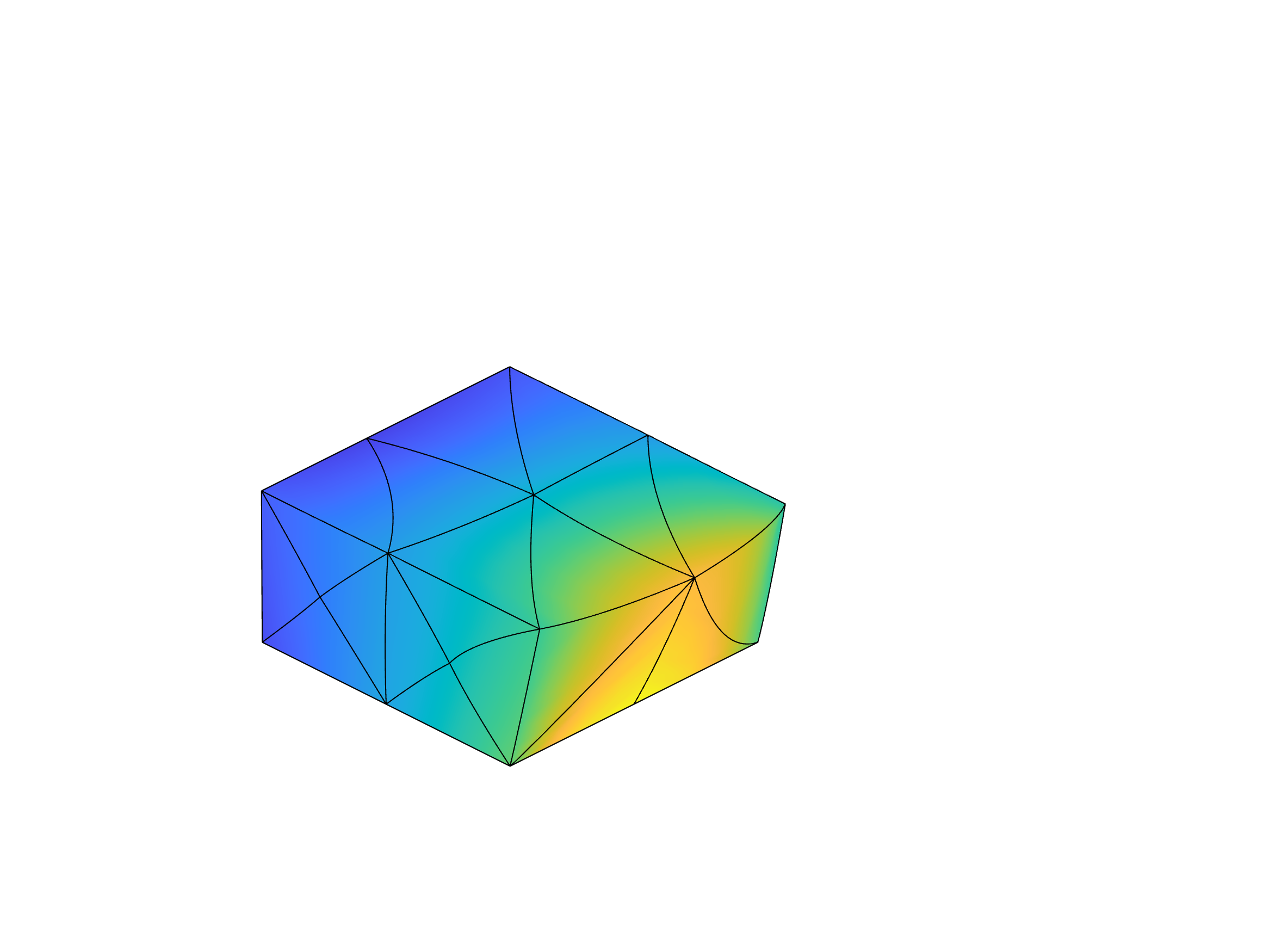};

\nextgroupplot[axis equal image, width=0.4\textwidth, xticklabels={}, yticklabels={}, ylabel={}, xlabel={}, xmin=-1, xmax=1, ymin=0, ymax=1.5]
\addplot []
graphics [xmin=-1,xmax=1,ymin=0,ymax=1.5] { 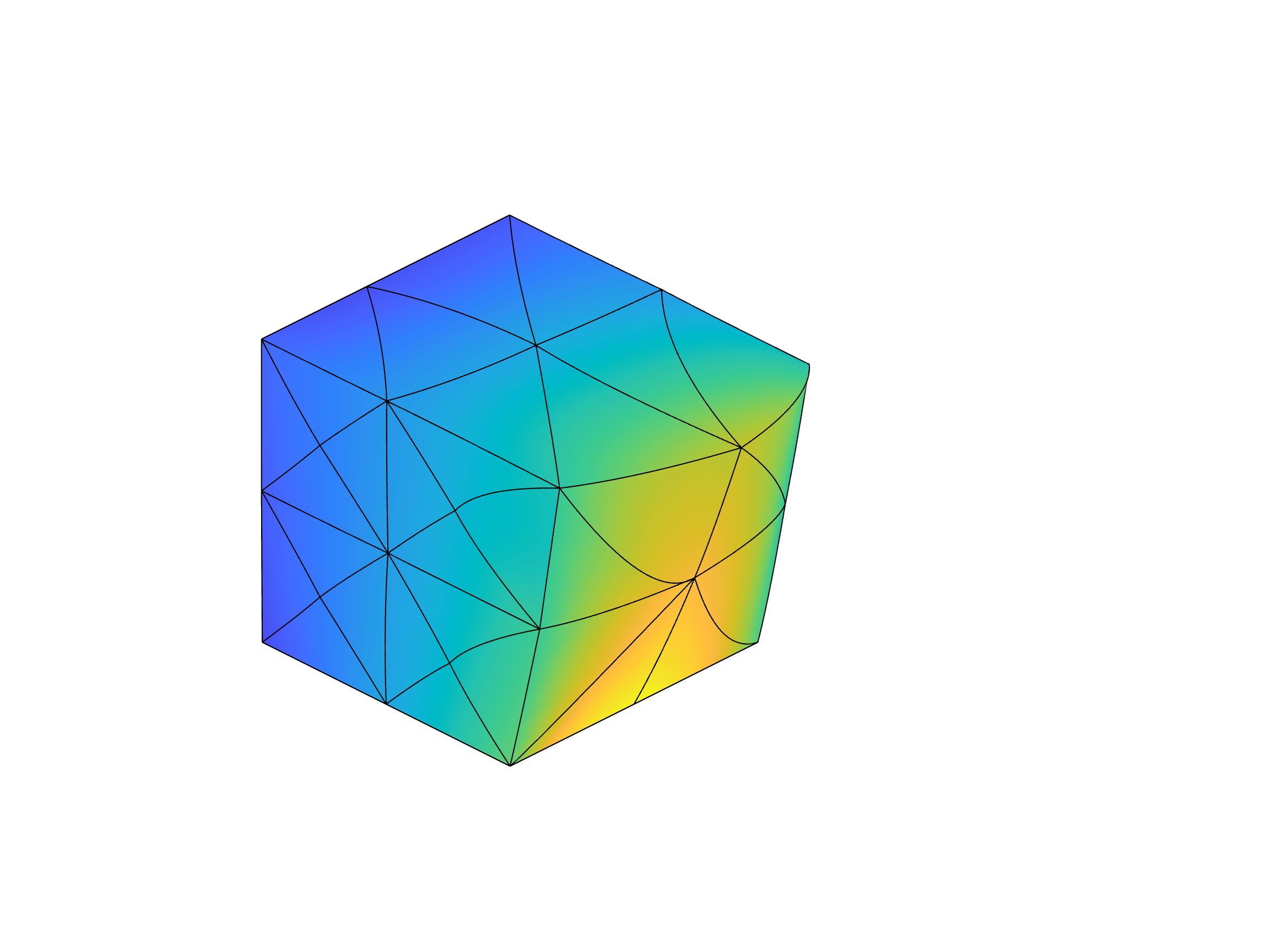};

\nextgroupplot[axis equal image, width=0.4\textwidth, xticklabels={}, yticklabels={}, ylabel={}, xlabel={}, xmin=-1, xmax=1, ymin=0, ymax=1.5]
\addplot []
graphics [xmin=-1,xmax=1,ymin=0,ymax=1.5] { 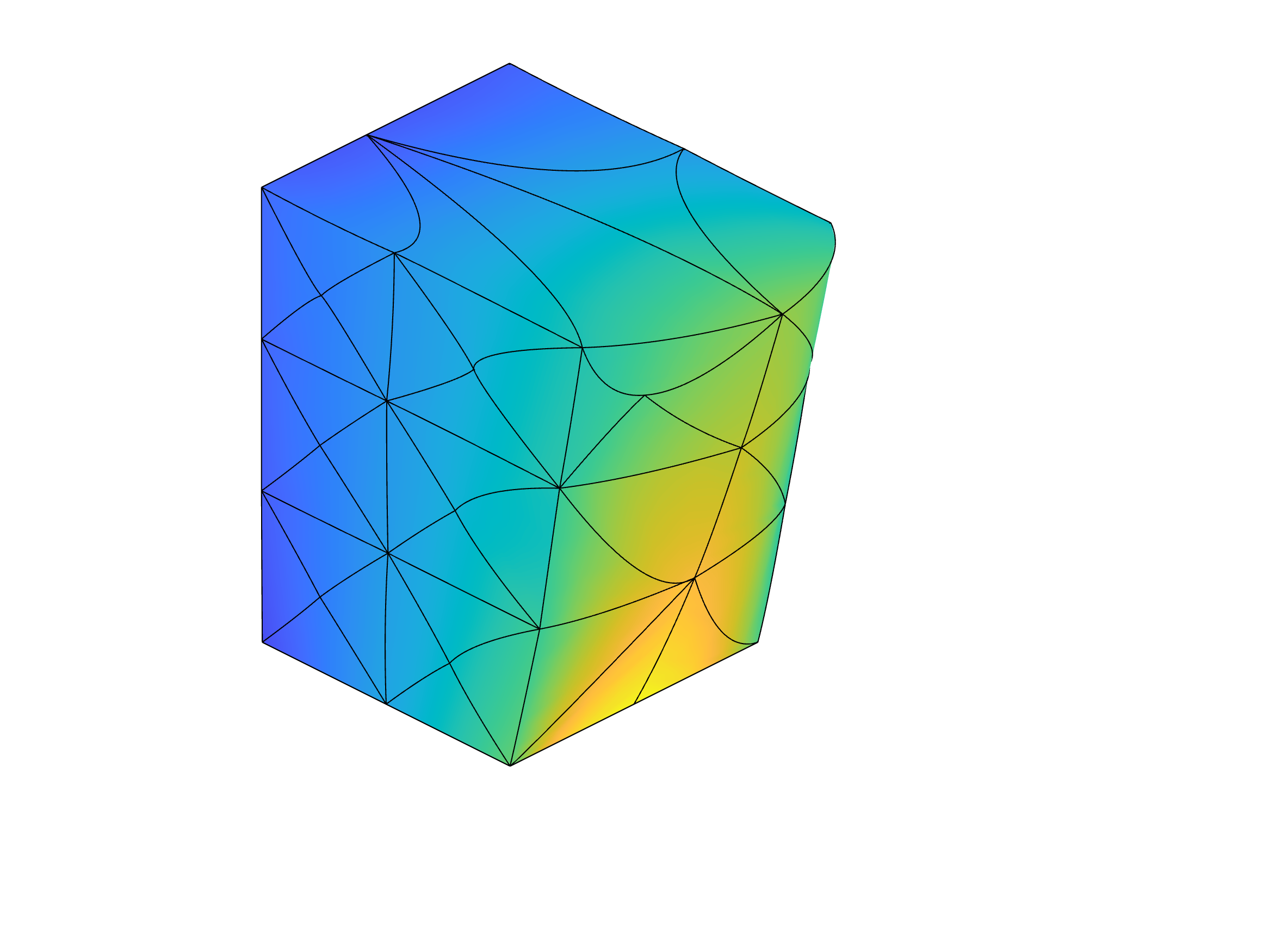};

\end{groupplot}\end{tikzpicture}}
	\colorbarMatlabParula{0}{0.25}{0.5}{0.75}{1}
 	\caption{Slab-base HOIST solution to the two-dimensional inviscid Burgers' problem.}
 	\label{fig:iburg2D}
\end{figure}

We use the slab-based HOIST method with $S = 3$ slabs. The spatial domain is discretized
with $32$ triangular elements, which are extruded and split to form an initial slab with
$448$ quadratic ($p = q = 2$) space-time tetrahedral elements.
As time evolves, the initially straight-sided shock begins to obtain a curved profile,
which is captured by the $q=2$ elements (Figure~\ref{fig:iburg2D}) and the $p=2$ solution
approximation accurately resolves the conserved variable.

\subsection{Shallow water equations}
Next we consider flow of an incompressible, inviscid, non-heat conducting fluid subject
to gravitational forces, but neglecting vertical acceleration through a domain
$\Omega_x\subset\Rbb^{d'}$ over a time interval $\Tcal\subset\Rbb_{\geq 0}$,
which is governed by the shallow water equations (SWE)
\begin{equation} \label{eqn:swe}
\begin{aligned}
\pder{}{t}\rho(x,t) + \pder{}{x_j} \rho(x,t) v_j(x,t) &= 0 \\
\pder{}{t}(\rho(x, t)v_i(x,t)) + \pder{}{x_j}(\rho(x,t)v_i(x,t)v_j(x,t)+\delta_{ij}\rho(x,t)^2/2) &= g \rho(x, t) \pder{h}{x_i}(x)
\end{aligned}
\end{equation}
for $i = 1, ..., d'$. The geopotential $\rho(\cdot,t) : \Omega_x \rightarrow \Rbb_{>0}$
and fluid velocity $v(\cdot,t) : \Omega_x \rightarrow \Rbb^d$ are implicitly defined as
the solution to (\ref{eqn:swe}). The geopotential is related to the total depth of the
fluid $H(\cdot,t) : \Omega_x \rightarrow \Rbb_{>0}$ as
\begin{equation}
 \rho(x,t) = g H(x,t), \qquad H(x,t) = \eta(x,t) + h(x)
\end{equation}
where $g \in \Rbb_{>0}$ is the gravitational acceleration and the total depth of the
fluid is written as the summation of the channel bed depth
$h : \Omega_x \rightarrow \Rbb_{>0}$ and the height
of the free surface $\eta(\cdot,t) : \Omega_x \rightarrow \Rbb_{>0}$ .
The shallow water equations are cast as a conservation law
of the form (\ref{eqn:gen_cons_law}) and the projected inviscid flux Jacobian
and its eigenvalue decomposition are given in \ref{app:projjac}. From
this information, the transformed space-time version of the conservation law
(including Riemann solver-based numerical flux) follows systematically from
Section~\ref{sec:govern:sptm}-\ref{sec:govern:transf}.

We consider the canonical one-dimensional ($d' = 1$) dam break probem with
$g = 1$ and no bathymetry ($h(x) = 0$) over $\Omega_x \coloneqq (-5, 5)$ and time
interval $\Tcal \coloneqq (0, 6)$, which is a Riemann probem for the SWE (\ref{eqn:swe})
with initial condition
\begin{equation}
  \eta (x, 0) =
  \begin{cases}
    3 & x \in [-5, 0) \\
    1 & x \in (0, 5 ]
  \end{cases}, \qquad
  v (x, 0) = 0.
\end{equation}
At both boundaries, we assume a non-penetration condition, i.e.,
$v(x,t) \cdot \eta_x(x) = 0$, where $\eta_x \in \Sbb_{d'}$ is the outward
unit normal to the boundary. As such, this problem features two waves---a
shock and rarefaction---that propagate from the origin and reflect off
the boundaries. This numerical experiment will demonstrate the slab-based
space-time HOIST method for the SWE and its ability to resolve wave reflections.

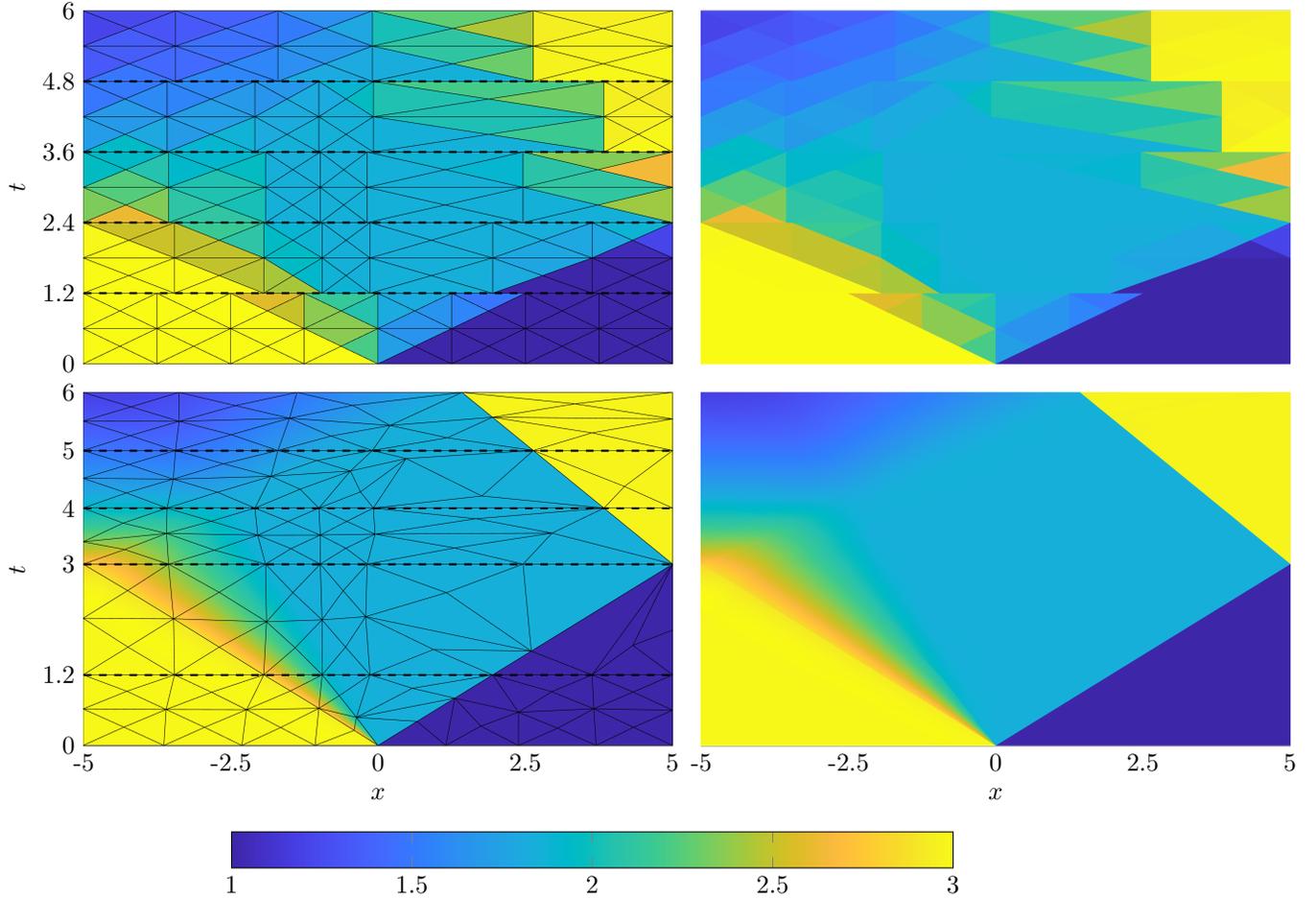
\begin{figure}
\centering
\begin{tikzpicture}
\begin{groupplot} [
group style={group size = 2 by 2, horizontal sep = 0.4cm, vertical sep = 0.4cm},
title style={at={(current bounding box.north west)}, anchor=west}]
\nextgroupplot[axis equal image, width=0.59\textwidth, xtick={-5, -2.5, 0, 2.5, 5}, ytick={0, 1.2, 2.4, 3.6, 4.8, 6}, xticklabels={}, yticklabels={0, 1.2, 2.4, 3.6, 4.8, 6}, ylabel={$t$}, xmin=-5, xmax=5, ymin=0, ymax=6]
\addplot []
graphics [xmin=-5,xmax=5,ymin=0,ymax=6] { 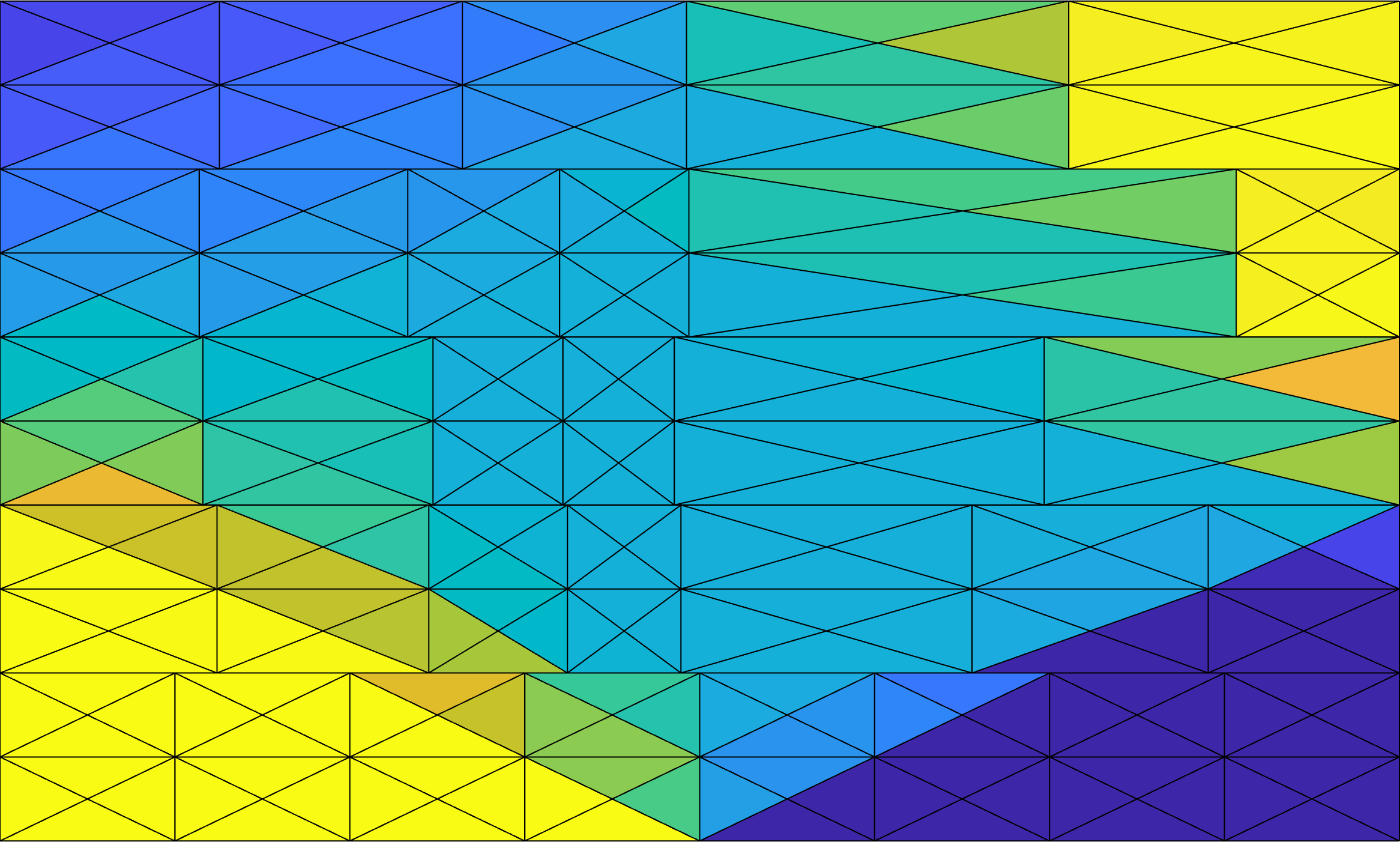};

\addplot [black, thick, dashed]
coordinates {
(-5.00000000e+00,  1.20000000e+00)
( 5.00000000e+00,  1.20000000e+00)};\label{line:swe_slab1/2}

\addplot [black, thick, dashed]
coordinates {
(-5.00000000e+00,  2.40000000e+00)
( 5.00000000e+00,  2.40000000e+00)};\label{line:slab2/3}

\addplot [black, thick, dashed]
coordinates {
(-5.00000000e+00,  3.60000000e+00)
( 5.00000000e+00,  3.60000000e+00)};\label{line:slab3/4}

\addplot [black, thick, dashed]
coordinates {
(-5.00000000e+00,  4.80000000e+00)
( 5.00000000e+00,  4.80000000e+00)};\label{line:slab4/5}

\nextgroupplot[axis equal image, width=0.59\textwidth, xtick={-5, -2.5, 0, 2.5, 5}, ytick={0, 1, 2}, xticklabels={}, yticklabels={}, xmin=-5, xmax=5, ymin=0, ymax=6]
\addplot []
graphics [xmin=-5,xmax=5,ymin=0,ymax=6] { 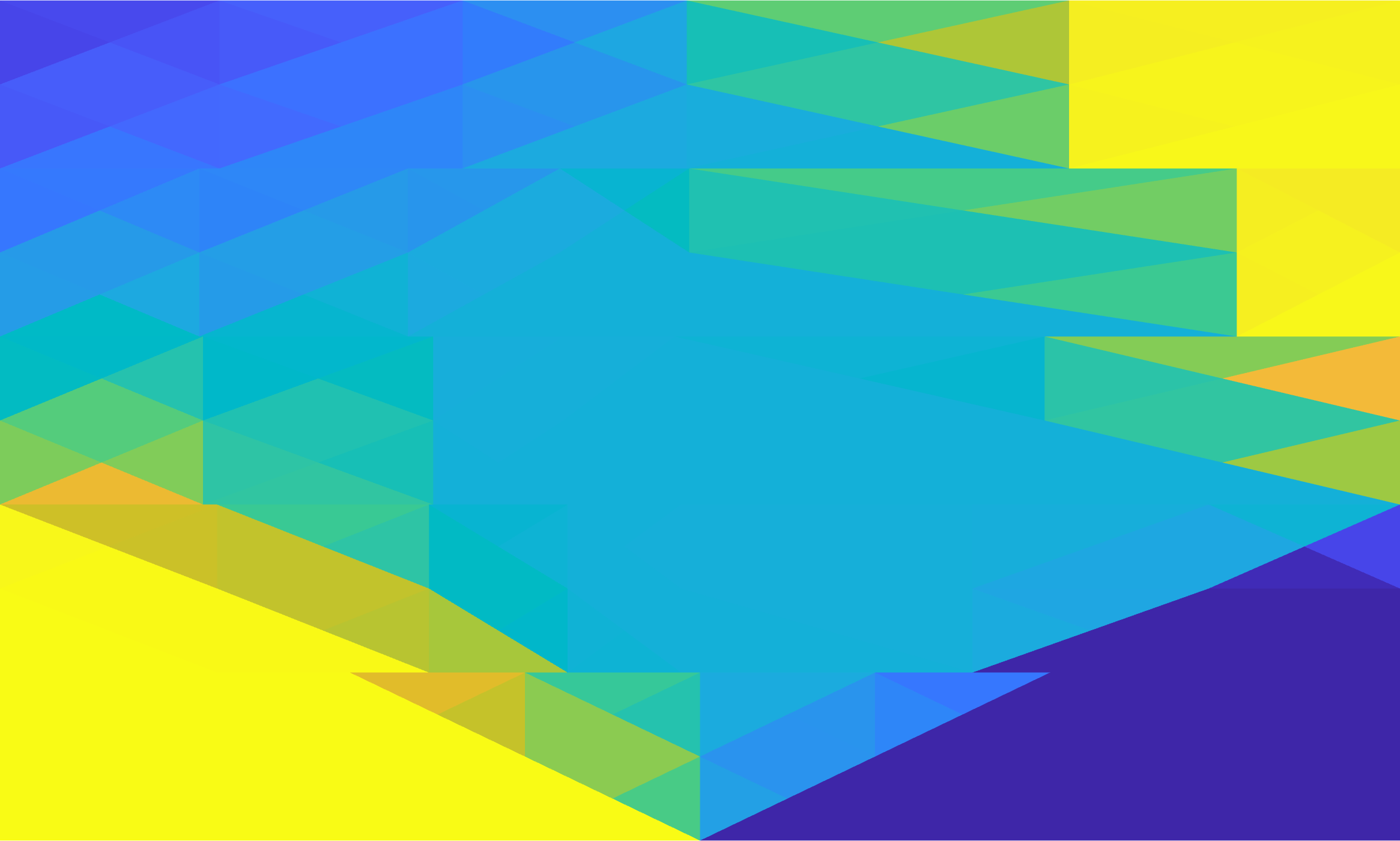};

\nextgroupplot[axis equal image, width=0.59\textwidth, xtick={-5, -2.5, 0, 2.5, 5}, ytick={0, 1.2, 3.082, 4.034, 5.0123, 6}, xticklabels={-5, -2.5, 0, 2.5, 5}, yticklabels={0, 1.2, 3, 4, 5, 6}, xlabel={$x$}, ylabel={$t$}, xmin=-5, xmax=5, ymin=0, ymax=6]
\addplot []
graphics [xmin=-5,xmax=5,ymin=0,ymax=6] { 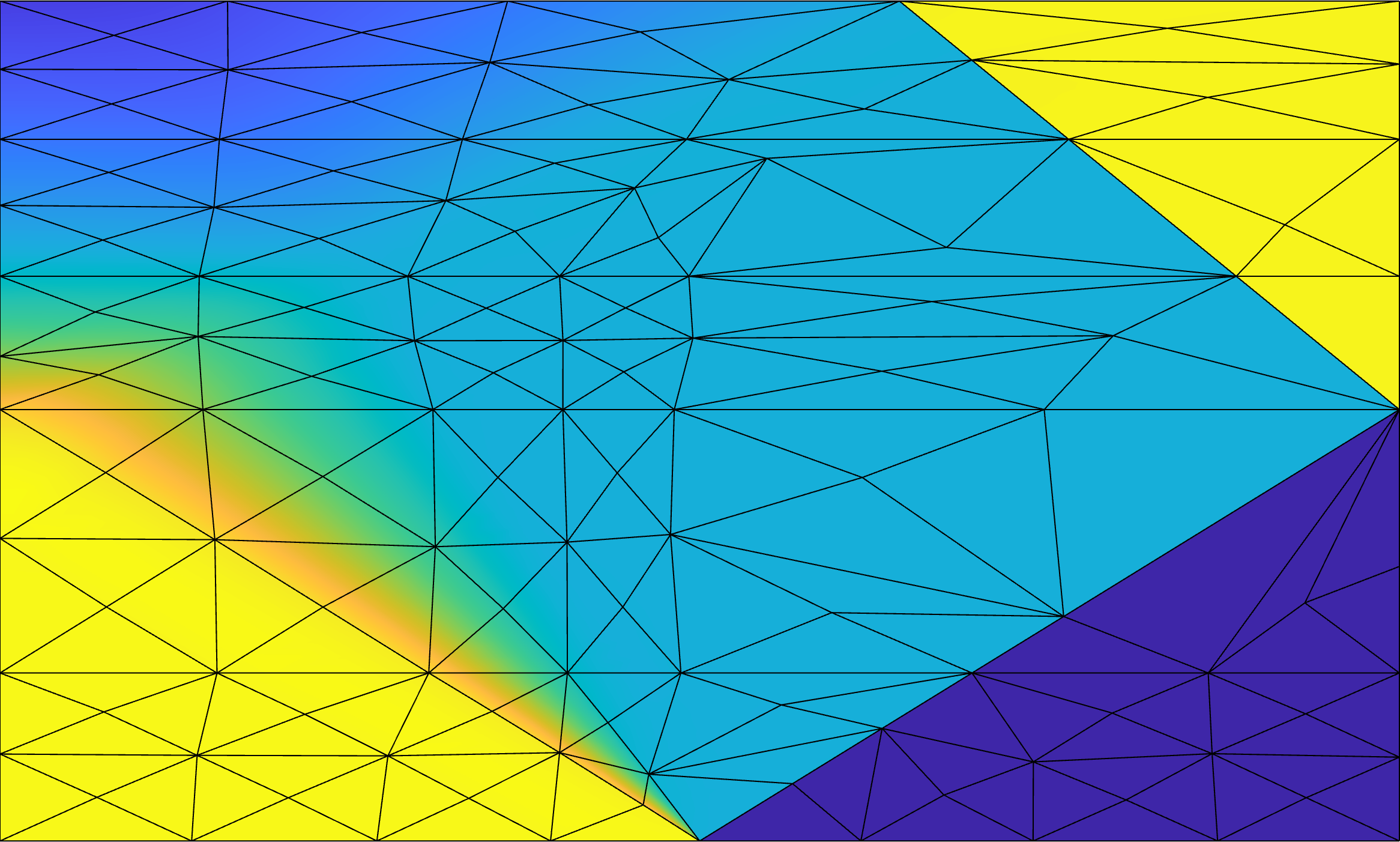};

\addplot [black, thick, dashed]
coordinates {
(-5.00000000e+00,  1.20000000e+00)
( 5.00000000e+00,  1.20000000e+00)};\label{line:slab1/2}

\addplot [black, thick, dashed]
coordinates {
(-5.00000000e+00,  3.08200000e+00)
( 5.00000000e+00,  3.08200000e+00)};\label{line:slab2/3}

\addplot [black, thick, dashed]
coordinates {
(-5.00000000e+00,  4.03400000e+00)
( 5.00000000e+00,  4.03400000e+00)};\label{line:slab3/4}

\addplot [black, thick, dashed]
coordinates {
(-5.00000000e+00,  5.01230000e+00)
( 5.00000000e+00,  5.01230000e+00)};\label{line:slab4/5}

\nextgroupplot[axis equal image, width=0.59\textwidth, xtick={-5, -2.5, 0, 2.5, 5}, ytick={0, 1, 2}, xticklabels={-5, -2.5, 0, 2.5, 5}, yticklabels={}, xlabel={$x$}, xmin=-5, xmax=5, ymin=0, ymax=6]
\addplot []
graphics [xmin=-5,xmax=5,ymin=0,ymax=6] { 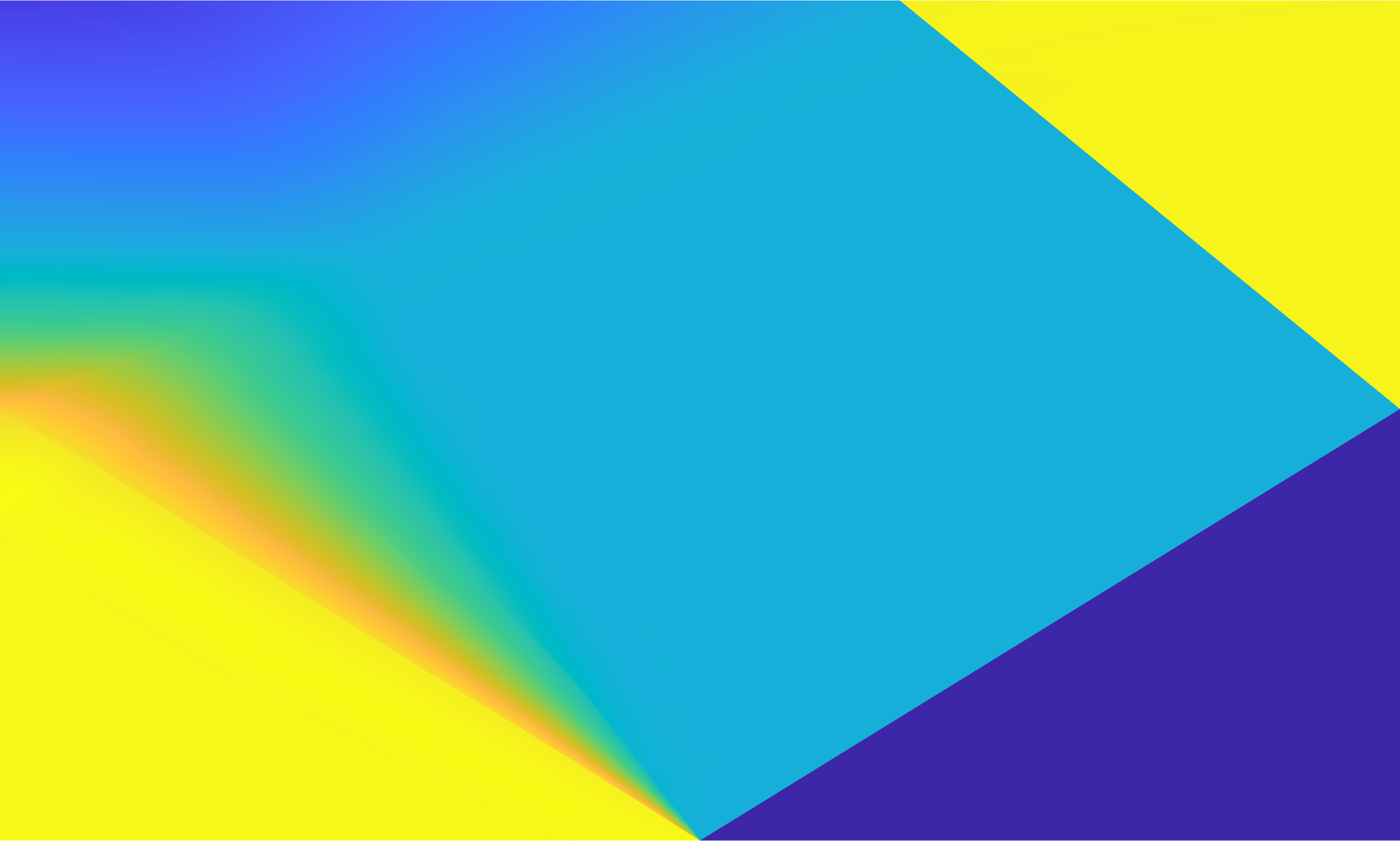};

\end{groupplot}\end{tikzpicture}
\colorbarMatlabParula{1}{1.5}{2}{2.5}{3}
\caption{Slab-based HOIST solution (geopotential) of the SWE including the initial (\textit{top row}) and converged (\textit{bottom row}) mesh and flow solution. The interface between each slab is denoted by (\ref{line:swe_slab1/2}) in both the initial solution and tracked solution.}
\label{fig:swe_refl_slabs_1_thru_6}
\end{figure}

This problem is discretized using $5$ slabs and solved using the space-time HOIST method.
The first slab is generated from a spatial mesh with $8$ elements, which are extruded and
split to produce $32$ space-time linear ($q = 1$) triangles with a quadratic ($p = 2$)
solution approximation. As shown in Figure \ref{fig:swe_refl_slabs_1_thru_6}, in the first slab, the shock and head/tail of the rarefaction are
successfully tracked as they propagate away from the origin. In the second slab,
the space-time points at which the shock and rarefaction reflect off the wall are
tracked by means of the variable temporal boundary (Section~\ref{sec:ist:dommap}),
i.e., the second slab is stretched in time such that these intersections are represented
perfectly. Because of this, the second slab covers more of the temporal domain
than the first slab while still accurately approximating the flow field. In
subsequent slabs, the shock and rarefaction are reflected and propagate in the
reflected directions until the end of the temporal domain is reached by the
fifth slab.

\begin{remark}
Figure \ref{fig:swe_refl_slabs_1_thru_6}  (\textit{top row}) shows the initial solution and mesh for each slab in the space-time solution. The initial mesh is guaranteed to be conforming to the converged mesh and solution from the previous slab (Figure \ref{fig:swe_refl_slabs_1_thru_6}  (\textit{bottom row})), but the initial meshes for each slab interface are not guaranteed
to be conforming (as shown). The initial meshes are never used simultaneously so conformity is irrelevant. This remark also applies to Figures \ref{fig:sod_slabs_1_thru_2} and \ref{fig:shu_osher_initial_level1}.
\end{remark}

%

\subsection{Euler equations of gasdynamics}
\label{sec:numexp:euler}
Finally, we consider flow of an inviscid, compressible fluid through a domain
$\Omega_x \subset \Rbb^{d'}$, which are governed by the Euler equations
\begin{equation} \label{eqn:euler}
\begin{aligned}
\pder{}{t}\rho(x,t) + \pder{}{x_j}\left(\rho(x,t) v_j(x,t)\right) &= 0 \\
\pder{}{t}(\rho(x,t)v_i(x,t)) + \pder{}{x_j}(\rho(x,t) v_i(x,t) v_j(x,t) + P(x,t) \delta_{ij} ) &= 0 \\
\pder{}{t}(\rho(x,t)E(x,t)) + \pder{}{x_j}([\rho(x,t) E(x,t)+P(x,t)]v_j(x,t)) &= 0 
\end{aligned}
\end{equation}
for all $x \in \Omega_x$ and $t \in \mathcal{T} \subset \Rbb$, where $i = 1, ..., d'$
and summation is implied over the repeated index $j = 1, ..., d'$.  The density of the
fluid $\rho : \Omega_x \times \mathcal{T} \rightarrow \Rbb_{>0}$, the velocity of the
fluid $v_i : \Omega_x \times \mathcal{T} \rightarrow \Rbb$ in the $x_i$ direction for
$i = 1,..., d'$, and the total energy of the fluid
$E : \Omega_x \times \mathcal{T} \rightarrow \Rbb_{>0}$ are implicitly defined
as the solution of (\ref{eqn:euler}). For a calorically ideal fluid, the pressure of
the fluid, $P : \Omega_x \times \mathcal{T} \rightarrow \Rbb_{>0}$ is related to the
energy via the ideal gas law
\begin{equation}
P = (\gamma - 1) ( \rho E - \rho v_i v_i/2)
\end{equation}
where $\gamma \in \Rbb_{>0}$ is the ratio of specific heats.
The Euler equations are cast as a conservation law
of the form (\ref{eqn:gen_cons_law}) and the projected inviscid flux Jacobian
and its eigenvalue decomposition are given in \ref{app:projjac}. From
this information, the transformed space-time version of the conservation law
(including Riemann solver-based numerical flux) follows systematically from
Section~\ref{sec:govern:sptm}-\ref{sec:govern:transf}.

\subsubsection{Sod's shock tube}
\label{sec:numexp:euler:sod}
Sod's shock tube is a Riemann problem for the Euler equations that models
an idealized shock tube where the membrane separating a high pressure region
from a low pressure one is instantaneously removed. This is a commonly used
validation problem since it has an analytical solution that features a shock
wave, rarefaction wave, and contact discontinuity. The flow domain is
$\Omega_x\coloneqq(0, 1)$, the time domain is $\Tcal\coloneqq(0,1)$, and
the initial condition is given in terms of the density, velocity, and
pressure as
\begin{equation}
 \rho(x,0) = \begin{cases} 1 & x<0.5 \\ 0.125 & x \geq 0.5 \end{cases}, \quad
  v(x,0) = 0, \quad
  P(x,0) = \begin{cases} 1 & x<0.5 \\ 0.1 & x \geq 0.5, \end{cases}.
\end{equation}
Wall boundary conditions are prescribed at both $x=0$ and $x=1$, i.e.,
$v(x,t) \cdot \eta_x(x) = 0$, where $\eta_x \in \Sbb_{d'}$ is the outward
unit normal to the boundary, which means any waves that reach the boundary
will reflect back into the domain. The solution of this problem contains
three waves (shock, contact, rarefaction) that emanate from $x=0.5$ at $t=0$
and propagate until the shock and rarefaction reflect off the walls. The shock
and contact eventually interact creating a pass through, reflection, and intermediate
that continue to interact with one another and the right wall. Eventually the
reflected rarefaction interacts with the shock causing it to curve.
The time domain considered in this work is five times longer than
the one usually associated with this problem \cite{toro2013riemann} to induce
these interactions. The purpose of this numerical experiment is to demonstrate
the slab-based space-time HOIST method is able to track various types of non-smooth
features and locate triple points without seeding with \textit{a priori} information.

\begin{figure}
\centering
 \raisebox{-0.5\height}{\begin{tikzpicture}
\begin{groupplot} [
group style={group size = 2 by 2, horizontal sep = 0.4cm, vertical sep = 0.4cm},
title style={at={(current bounding box.north west)}, anchor=west}]
\nextgroupplot[axis equal image, width=0.55\textwidth, xtick={-5, -2.5, 0, 2.5, 5}, ytick={0, 1, 2, 3, 4, 5.5, 10}, xticklabels={}, yticklabels={0, 0.1, 0.2, 0.3, 0.4, 0.55, 1}, ylabel={$t$}, xmin=-5, xmax=5, ymin=0, ymax=10]
\addplot []
graphics [xmin=-5,xmax=5,ymin=0,ymax=10] { 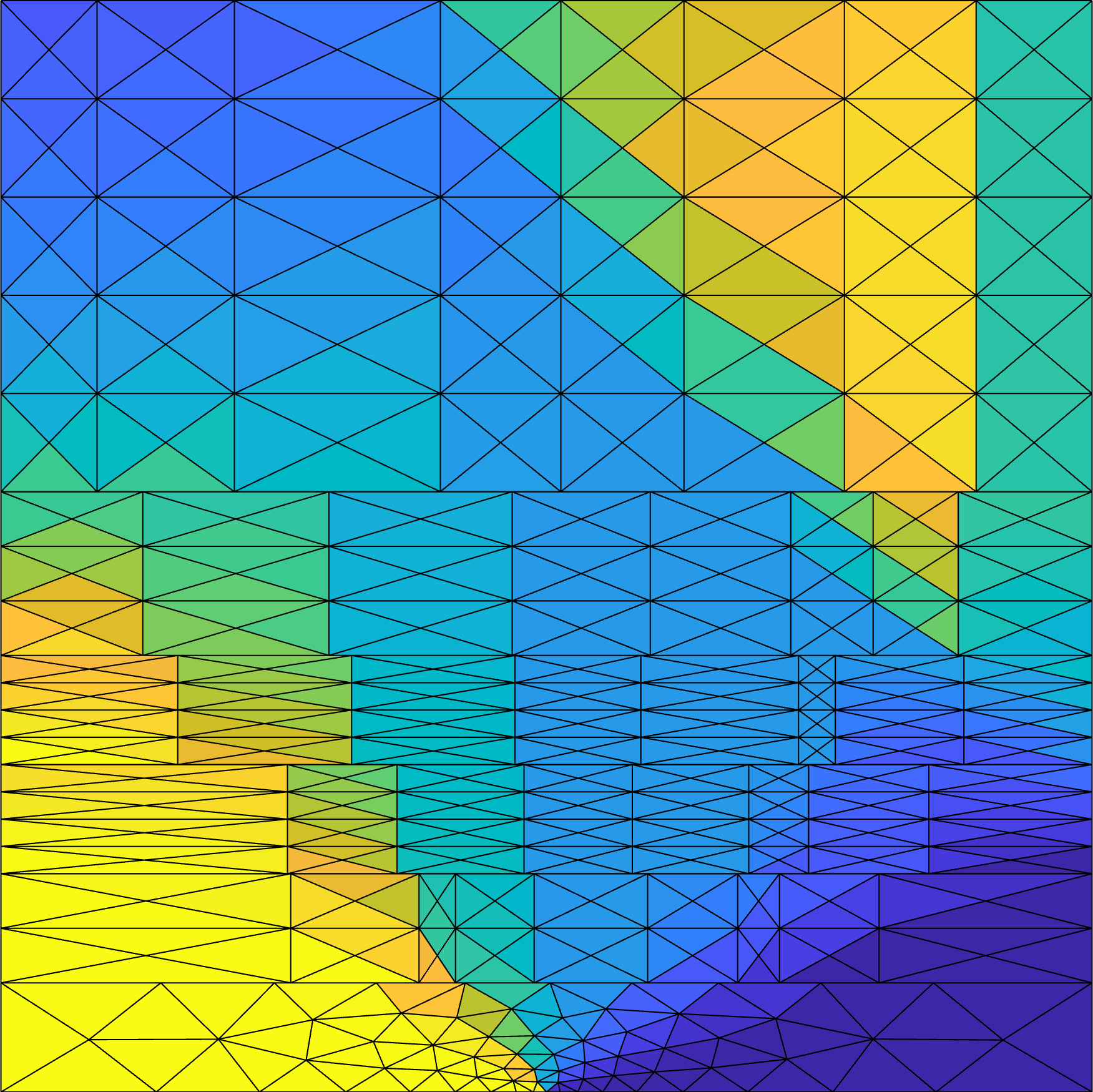};

\addplot [black, thick, dashed]
coordinates {
(-5.00000000e+00,  1.00000000e+00)
( 5.00000000e+00,  1.00000000e+00)};\label{line:sod_slab1/2}

\addplot [black, thick, dashed]
coordinates {
(-5.00000000e+00,  2.00000000e+00)
( 5.00000000e+00,  2.00000000e+00)};\label{line:slab2/3}

\addplot [black, thick, dashed]
coordinates {
(-5.00000000e+00,  3.00000000e+00)
( 5.00000000e+00,  3.00000000e+00)};\label{line:slab3/4}

\addplot [black, thick, dashed]
coordinates {
(-5.00000000e+00,  4.00000000e+00)
( 5.00000000e+00,  4.00000000e+00)};\label{line:slab4/5}

\addplot [black, thick, dashed]
coordinates {
(-5.00000000e+00,  5.50000000e+00)
( 5.00000000e+00,  5.50000000e+00)};\label{line:slab5/6}

\nextgroupplot[axis equal image, width=0.55\textwidth, xtick={-5, 0, 5}, ytick={0, 10}, xticklabels={}, yticklabels={}, xmin=-5, xmax=5, ymin=0, ymax=10]
\addplot []
graphics [xmin=-5,xmax=5,ymin=0,ymax=10] { 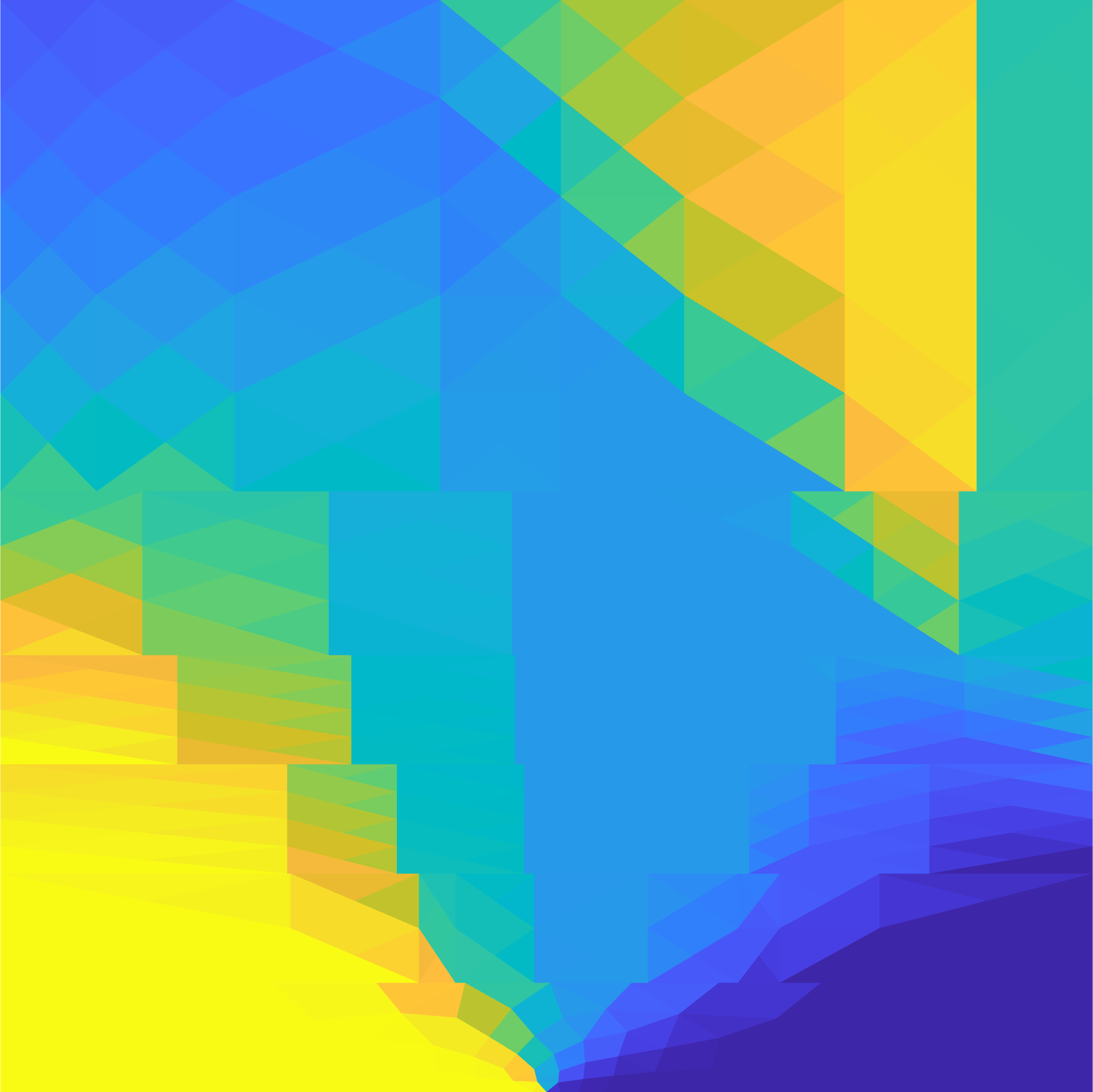};

\nextgroupplot[axis equal image, width=0.55\textwidth, xtick={-5, -2.5, 0, 2.5, 5}, ytick={0, 1, 2, 2.85, 4.1, 5.6, 10}, xticklabels={0, 0.25, 0.5, 0.75, 1}, yticklabels={0, 0.1, 0.2, 0.285, 0.41, 0.56, 1}, xlabel={$x$}, ylabel={$t$}, xmin=-5, xmax=5, ymin=0, ymax=10]
\addplot []
graphics [xmin=-5,xmax=5,ymin=0,ymax=10] { 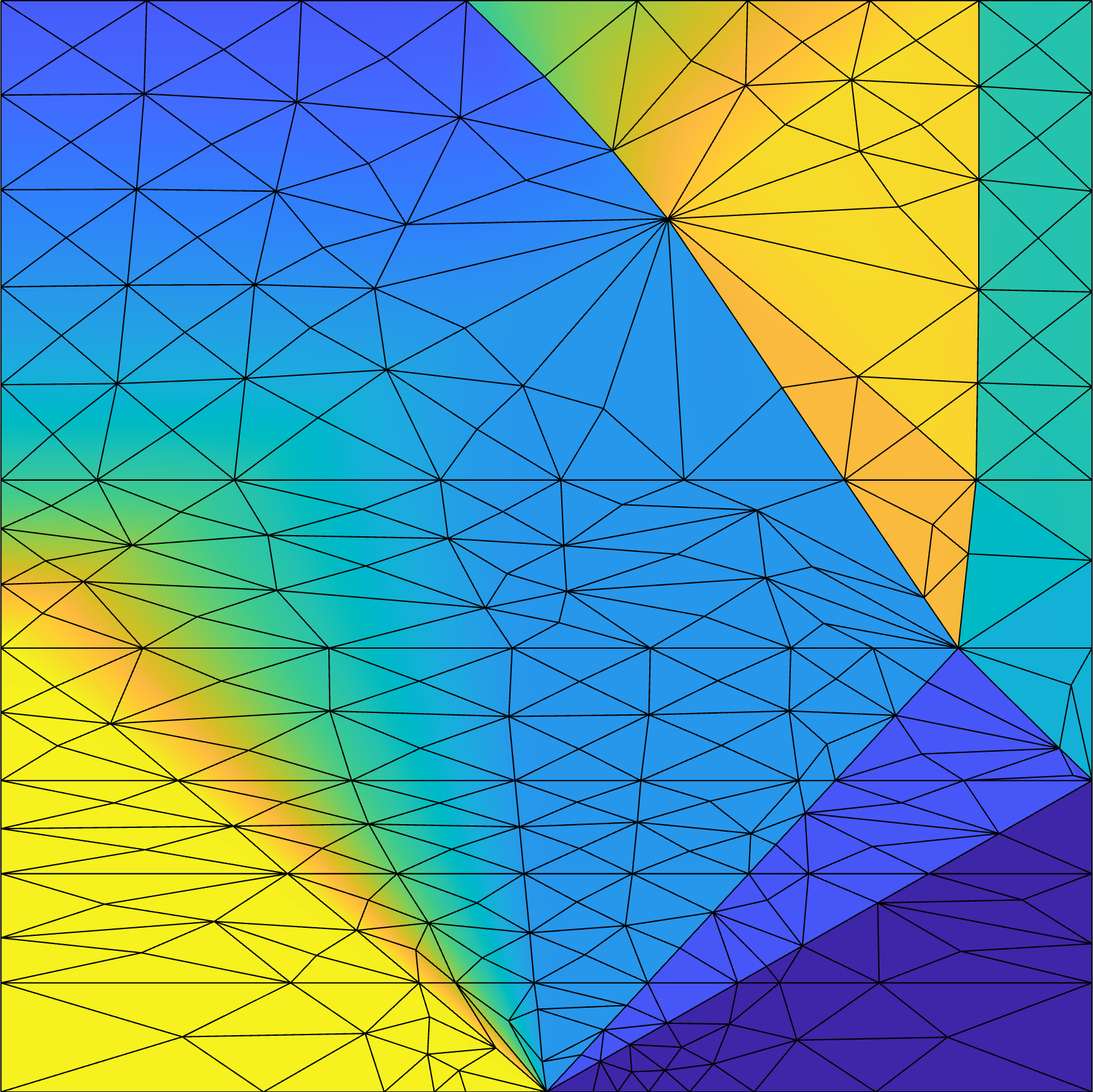};

\addplot [black, thick, dashed]
coordinates {
(-5.00000000e+00,  1.00000000e+00)
( 5.00000000e+00,  1.00000000e+00)};\label{line:slab1/2}

\addplot [black, thick, dashed]
coordinates {
(-5.00000000e+00,  2.00000000e+00)
( 5.00000000e+00,  2.00000000e+00)};\label{line:slab2/3}

\addplot [black, thick, dashed]
coordinates {
(-5.00000000e+00,  2.85400000e+00)
( 5.00000000e+00,  2.85400000e+00)};\label{line:slab3/4}

\addplot [black, thick, dashed]
coordinates {
(-5.00000000e+00,  4.06800000e+00)
( 5.00000000e+00,  4.06800000e+00)};\label{line:slab4/5}

\addplot [black, thick, dashed]
coordinates {
(-5.00000000e+00,  5.60700000e+00)
( 5.00000000e+00,  5.60700000e+00)};\label{line:slab5/6}

\nextgroupplot[axis equal image, width=0.55\textwidth, xtick={-5, -2.5, 0, 2.5, 5}, ytick={0, 10}, xticklabels={0, 0.25, 0.5, 0.75, 1}, yticklabels={}, xlabel={$x$}, xmin=-5, xmax=5, ymin=0, ymax=10]
\addplot []
graphics [xmin=-5,xmax=5,ymin=0,ymax=10] { 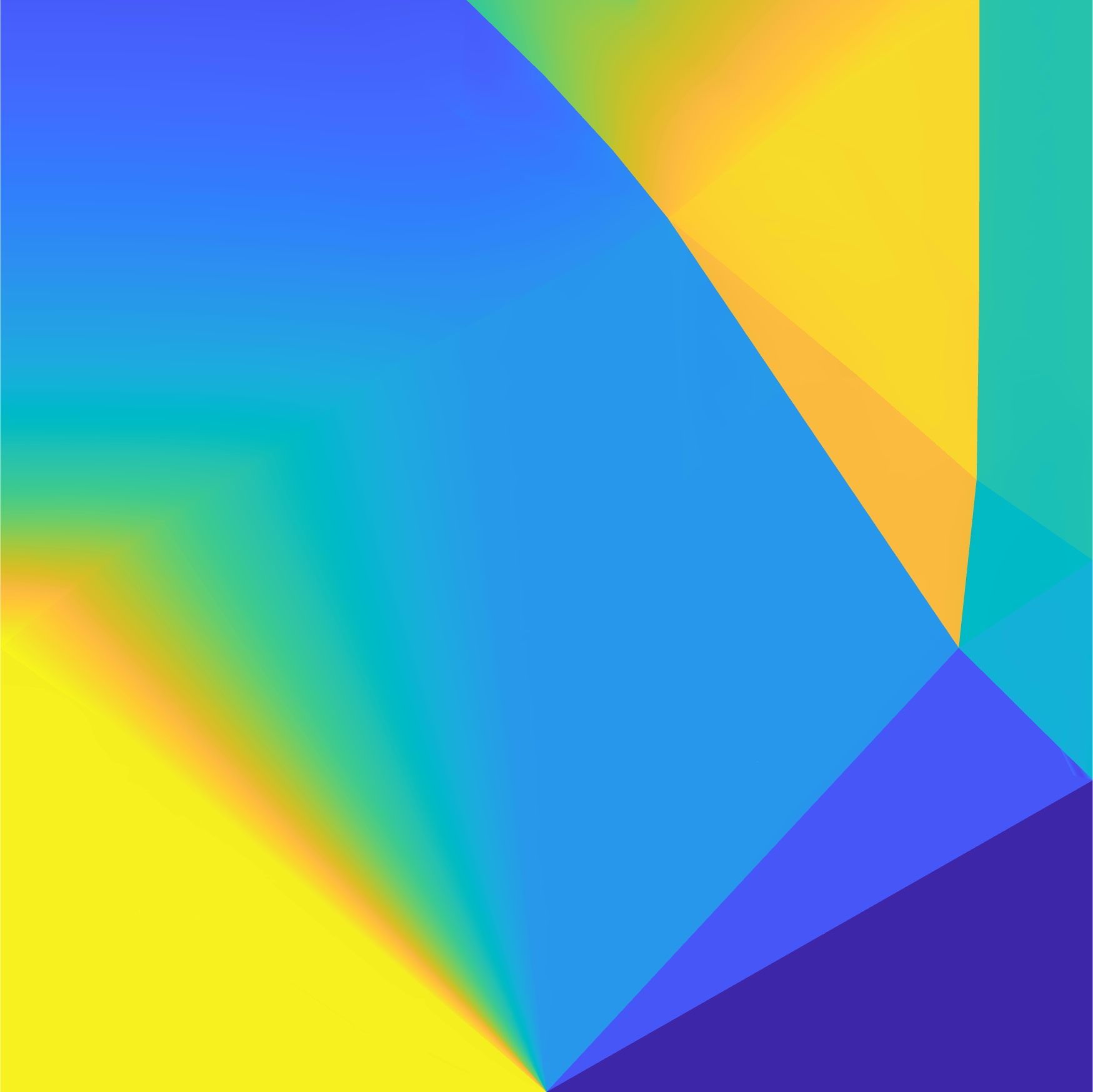};

\addplot [red, thin, solid]
coordinates {
(-5.00000000e+00,  2.00000000e+00)
( 5.00000000e+00,  2.00000000e+00)};\label{line:slice4}

\addplot [red, thin, solid]
coordinates {
(-5.00000000e+00,  3.40000000e+00)
( 5.00000000e+00,  3.40000000e+00)};\label{line:slice1}

\addplot [red, thin, solid]
coordinates {
(-5.00000000e+00,  4.30000000e+00)
( 5.00000000e+00,  4.30000000e+00)};\label{line:slice2}

\addplot [red, thin, solid]
coordinates {
(-5.00000000e+00,  5.40000000e+00)
( 5.00000000e+00,  5.40000000e+00)};\label{line:slice3}

\end{groupplot}\end{tikzpicture}}
 \colorbarMatlabParula{0.125}{0.25}{0.5}{0.75}{1}
 \caption{Slab-based HOIST solution (density) of Sod's shock tube including the initial (\textit{top row}) and converged (\textit{bottom row}) mesh and flow solution. The interface between each slab is denoted by (\ref{line:sod_slab1/2}) in both the initial solution and tracked solution. Slices through the domain indicated by (\ref{line:slice1}) at times near shock interactions can be seen in Figure \ref{fig:sod_slices}.}
 \label{fig:sod_slabs_1_thru_2}
\end{figure}
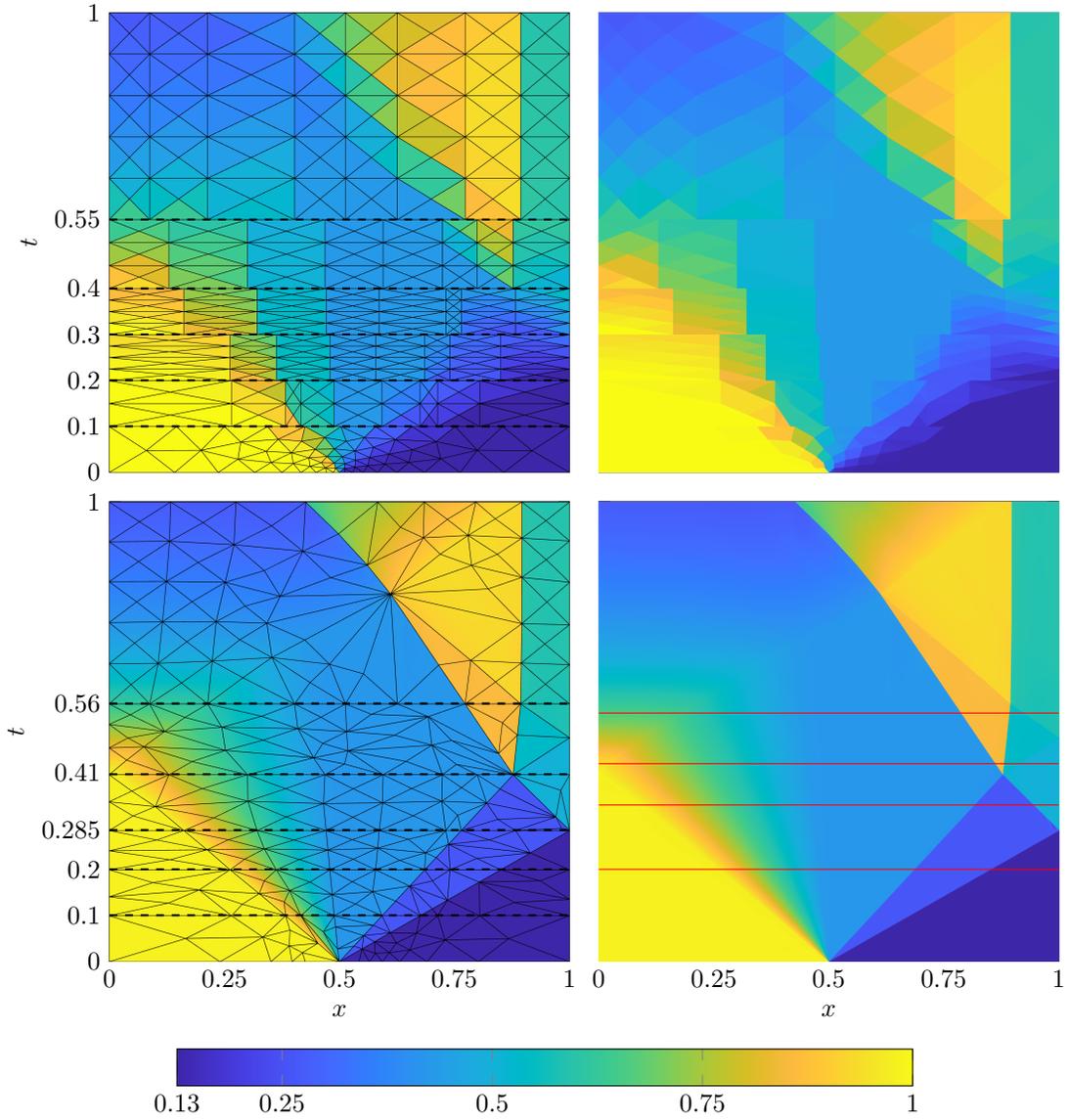

\begin{figure}
\centering
 \raisebox{-0.5\height}{\input{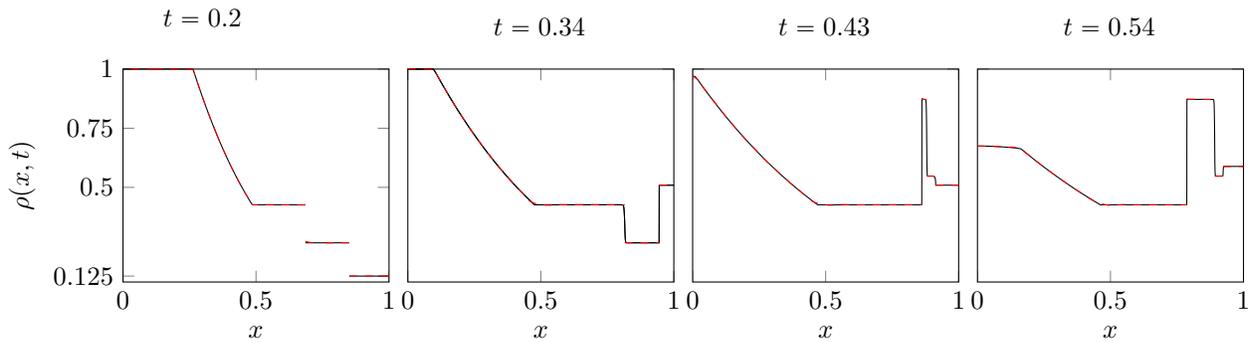}}
 \caption{Slab-based HOIST solution slices (density) (\ref{line:ist_sol_sod}) of Sod's shock tube compared against the reference solution (\ref{line:analytical_sol_sod}) at times indicated by (\ref{line:slice1}) in Figure \ref{fig:sod_slabs_1_thru_2}.}
 \label{fig:sod_slices}
\end{figure}

This problem is discretized with $S = 6$ slabs and solved using the space-time HOIST method.
An unstructured space-time mesh of $65$ linear ($q = 1$) simplices with quadratic solution
approximation ($p = 2$) is used for the first slab. Subsequent space-time slabs are
generated using the approach described in Section~\ref{sec:slab} whereby the spatial mesh
at the end of a time slab is extruded and split into simplices. All of the solution
features described above are tracked by the slab-based HOIST method
(Figure~\ref{fig:sod_slabs_1_thru_2}) in various slabs.
The shock, contact, and head/tail of the rarefaction are all tracked in the first slab,
which shows the method can locate multiple non-smooth features intersecting at a generalized
triple point. This is demonstrated again at later times when various shock-shock
and shock-contact interactions occurs. Furthermore, the variable temporal boundary
allows the method to track shock (top of slab 3) and rarefaction (top of slab 4)
reflections off the walls, as well as the shock-contact interaction (top of slab 4). Figure \ref{fig:sod_slices} shows several spatial solutions at given points in time for the HOIST method and a reference solution (second-order finite
volume method with superbee limiter on a grid with over $10^5$ elements and
fourth-order Runge-Kutta time integration with adaptive time steps). The slices in time display relevant shock dynamics and the ability of HOIST to represent the shock fronts as perfect discontinuities both before and after wave interactions.

\begin{figure}
\centering
 \raisebox{-0.5\height}{\begin{tikzpicture}
\begin{groupplot} [
group style={group size = 2 by 1, horizontal sep = 0.4cm, vertical sep = 0.4cm},
title style={at={(current bounding box.north west)}, anchor=west}]
\nextgroupplot[axis equal image, width=0.5\textwidth, xtick={0, 0.25, 0.5, 0.75, 1}, ytick={0, 0.4397}, xticklabels={0, 0.25, 0.5, 0.75, 1}, yticklabels={0, 0.44}, xlabel={$x$}, ylabel={$t$}, xmin=0, xmax=1, ymin=0, ymax=0.439713]
\addplot []
graphics [xmin=0,xmax=1,ymin=0,ymax=0.406756] { 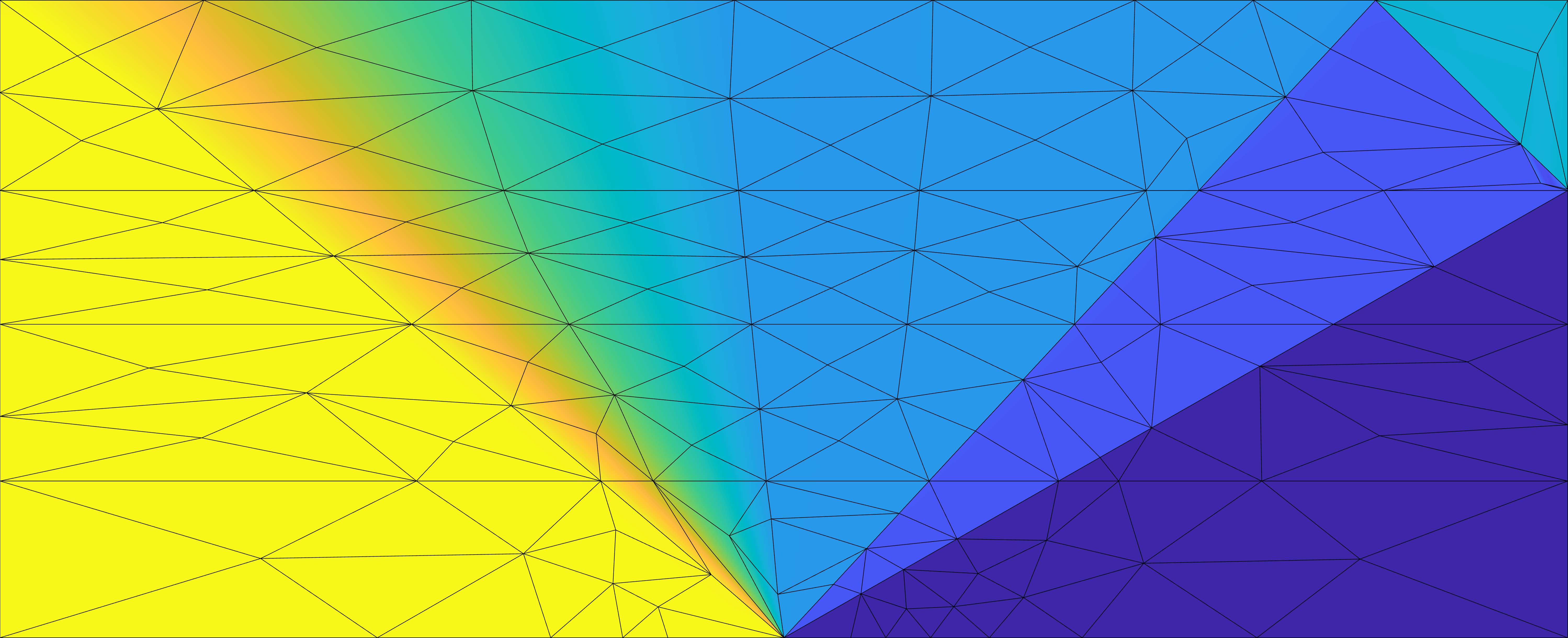};

\addplot [red, thin, dashed]
coordinates {
( 0.00000000e+00,  4.28000000e-01)
( 1.02040816e-02,  4.19263776e-01)
( 2.04081633e-02,  4.10527551e-01)
( 3.06122449e-02,  4.01791327e-01)
( 4.08163265e-02,  3.93055102e-01)
( 5.10204082e-02,  3.84318878e-01)
( 6.12244898e-02,  3.75582653e-01)
( 7.14285714e-02,  3.66846429e-01)
( 8.16326531e-02,  3.58110204e-01)
( 9.18367347e-02,  3.49373980e-01)
( 1.02040816e-01,  3.40637755e-01)
( 1.12244898e-01,  3.31901531e-01)
( 1.22448980e-01,  3.23165306e-01)
( 1.32653061e-01,  3.14429082e-01)
( 1.42857143e-01,  3.05692857e-01)
( 1.53061224e-01,  2.96956633e-01)
( 1.63265306e-01,  2.88220408e-01)
( 1.73469388e-01,  2.79484184e-01)
( 1.83673469e-01,  2.70747959e-01)
( 1.93877551e-01,  2.62011735e-01)
( 2.04081633e-01,  2.53275510e-01)
( 2.14285714e-01,  2.44539286e-01)
( 2.24489796e-01,  2.35803061e-01)
( 2.34693878e-01,  2.27066837e-01)
( 2.44897959e-01,  2.18330612e-01)
( 2.55102041e-01,  2.09594388e-01)
( 2.65306122e-01,  2.00858163e-01)
( 2.75510204e-01,  1.92121939e-01)
( 2.85714286e-01,  1.83385714e-01)
( 2.95918367e-01,  1.74649490e-01)
( 3.06122449e-01,  1.65913265e-01)
( 3.16326531e-01,  1.57177041e-01)
( 3.26530612e-01,  1.48440816e-01)
( 3.36734694e-01,  1.39704592e-01)
( 3.46938776e-01,  1.30968367e-01)
( 3.57142857e-01,  1.22232143e-01)
( 3.67346939e-01,  1.13495918e-01)
( 3.77551020e-01,  1.04759694e-01)
( 3.87755102e-01,  9.60234694e-02)
( 3.97959184e-01,  8.72872449e-02)
( 4.08163265e-01,  7.85510204e-02)
( 4.18367347e-01,  6.98147959e-02)
( 4.28571429e-01,  6.10785714e-02)
( 4.38775510e-01,  5.23423469e-02)
( 4.48979592e-01,  4.36061224e-02)
( 4.59183673e-01,  3.48698980e-02)
( 4.69387755e-01,  2.61336735e-02)
( 4.79591837e-01,  1.73974490e-02)
( 4.89795918e-01,  8.66122449e-03)
( 5.00000000e-01, -7.50000000e-05)};\label{line:rare_line}

\nextgroupplot[axis equal image, width=0.5\textwidth, xtick={0, 0.25, 0.5, 0.75, 1}, ytick={0, 1}, xticklabels={0, 0.25, 0.5, 0.75, 1}, yticklabels={}, xlabel={$x$}, xmin=0, xmax=1, ymin=0, ymax=0.439713]
\addplot []
graphics [xmin=0,xmax=1,ymin=0,ymax=0.439713] { 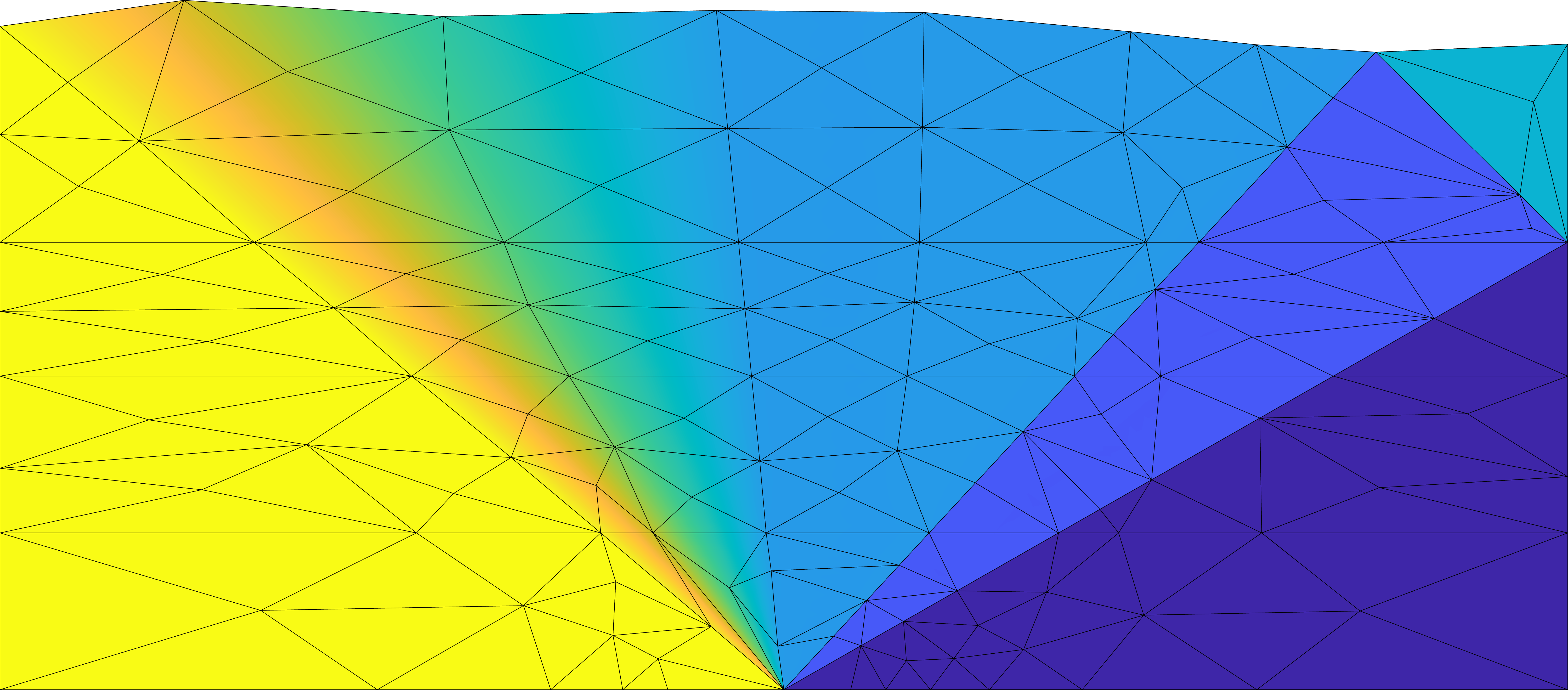};

\addplot [red, thin, dashed]
coordinates {
( 0.00000000e+00,  4.28000000e-01)
( 1.02040816e-02,  4.19263776e-01)
( 2.04081633e-02,  4.10527551e-01)
( 3.06122449e-02,  4.01791327e-01)
( 4.08163265e-02,  3.93055102e-01)
( 5.10204082e-02,  3.84318878e-01)
( 6.12244898e-02,  3.75582653e-01)
( 7.14285714e-02,  3.66846429e-01)
( 8.16326531e-02,  3.58110204e-01)
( 9.18367347e-02,  3.49373980e-01)
( 1.02040816e-01,  3.40637755e-01)
( 1.12244898e-01,  3.31901531e-01)
( 1.22448980e-01,  3.23165306e-01)
( 1.32653061e-01,  3.14429082e-01)
( 1.42857143e-01,  3.05692857e-01)
( 1.53061224e-01,  2.96956633e-01)
( 1.63265306e-01,  2.88220408e-01)
( 1.73469388e-01,  2.79484184e-01)
( 1.83673469e-01,  2.70747959e-01)
( 1.93877551e-01,  2.62011735e-01)
( 2.04081633e-01,  2.53275510e-01)
( 2.14285714e-01,  2.44539286e-01)
( 2.24489796e-01,  2.35803061e-01)
( 2.34693878e-01,  2.27066837e-01)
( 2.44897959e-01,  2.18330612e-01)
( 2.55102041e-01,  2.09594388e-01)
( 2.65306122e-01,  2.00858163e-01)
( 2.75510204e-01,  1.92121939e-01)
( 2.85714286e-01,  1.83385714e-01)
( 2.95918367e-01,  1.74649490e-01)
( 3.06122449e-01,  1.65913265e-01)
( 3.16326531e-01,  1.57177041e-01)
( 3.26530612e-01,  1.48440816e-01)
( 3.36734694e-01,  1.39704592e-01)
( 3.46938776e-01,  1.30968367e-01)
( 3.57142857e-01,  1.22232143e-01)
( 3.67346939e-01,  1.13495918e-01)
( 3.77551020e-01,  1.04759694e-01)
( 3.87755102e-01,  9.60234694e-02)
( 3.97959184e-01,  8.72872449e-02)
( 4.08163265e-01,  7.85510204e-02)
( 4.18367347e-01,  6.98147959e-02)
( 4.28571429e-01,  6.10785714e-02)
( 4.38775510e-01,  5.23423469e-02)
( 4.48979592e-01,  4.36061224e-02)
( 4.59183673e-01,  3.48698980e-02)
( 4.69387755e-01,  2.61336735e-02)
( 4.79591837e-01,  1.73974490e-02)
( 4.89795918e-01,  8.66122449e-03)
( 5.00000000e-01, -7.50000000e-05)};\label{line:rare_line}

\end{groupplot}\end{tikzpicture}}
 \caption{Slab-based HOIST solution (first four slabs) to Sod's shock tube using the translating boundary (\textit{left}) and a more general parametrization of the upper temporal boundary (\textit{right}). The true rarefaction trajectory is identified with (\ref{line:rare_line}). Colorbar in Figure~\ref{fig:sod_slabs_1_thru_2}.}
 \label{fig:sod_boundary_variable}
\end{figure}

The reflected shock and contact intersect at nearly (but not exactly) the same time as the rarefaction reflects off the wall. For the coarse mesh used, the translating temporal slab is insufficient to track both features, and the HOIST method chose to track the shock-contact interaction because it is the stronger feature. Thus, the interaction of the rarefaction with the left wall is an approximation (Figure~\ref{fig:sod_boundary_variable}). Two options to resolve this issue are to use a: 1) finer mesh (possibly with $h$-adaptation) or 2) more flexible parametrization of the upper temporal boundary. To demonstrate the flexibility of HOIST, we consider the latter by allowing free temporal motion of the nodes on the top of the fourth space-time slab, similar to parametrization of shock boundaries (Section~\ref{sec:ist:shkbc}). This is dangerous because causality may be violated if the relative motion of two adjacent nodes is too several, in which case the solution in the current slab would depend on the solution in the next slab. However, for this problem, the approach works robustly without causality violated and allows the rarefaction reflection and shock-contact interaction to be tracked simultaneously on the coarse grid (Figure~\ref{fig:sod_boundary_variable}). Future work will push this approach further with additional constraints to enforce causality.

\subsubsection{Shu-Osher problem}
\label{sec:numexp:euler:shuosher}
The Shu-Osher problem is a one-dimensional idealization of shock-turbulence interaction
in which a shock propagates into a density field with artificial fluctuations. This
phenomena is modeled by the ideal gas one-dimensional Euler equations, i.e.,
(\ref{eqn:euler}) with $d'=1$, with $\gamma = 1.4$. The flow domain is
$\Omega_x \coloneqq (-5, 5)$ and temporal domain is $\Tcal \coloneqq (0, T)$
with $T = 2$. The initial condition for the flow is
\begin{equation*}
 (\rho(x,0), v(x,0), P(x,0)) =
  \begin{cases}
   (3.857143, 2.629369, 10.3333) & x < -4 \\
   (1+0.2\sin(5x), 0, 1) & x \geq -4
  \end{cases}
\end{equation*}
with a supersonic inlet at $x = -5$ that prescribes the density, velocity, and pressure
\begin{equation*}
  (\rho(-5,t), v(-5,t), P(-5,t)) = (3.857143, 2.629369, 10.3333)
\end{equation*}
and a solid wall at $x = 5$. This problem corresponds to a Mach $M = 3$ shock moving
into a field with a small density (or entropy) disturbance. As time progresses, the
flucuations that have passed through the discontinuity become more oscillatory and
eventually steepen into shocks. The final time ($T = 2$) is chosen larger than often
seen in the literature \cite{2021_shi_methlines,2019_corrigan_MGDICEunsteady,shu1989efficient, johnsen2010assessment} to expose these challenging features.
The purpose of this numerical experiment is to demonstrate the slab-based space-time
HOIST method can simultaneously resolve strong shocks and weak shocks, shock formation,
and fluctuations that can easily be washed out by numerical or artificial dissipation.
This problem will also demonstrate the benefit of using adaptive mesh refinement
(Section~\ref{sec:ist:amr}) and shock boundary conditions (Section~\ref{sec:ist:shkbc}).

\begin{figure}
 \centering
 \raisebox{-0.5\height}{\begin{tikzpicture}
\begin{groupplot} [
group style={group size = 1 by 4, horizontal sep = 0.05cm, vertical sep = 0.05cm},
title style={at={(current bounding box.north west)}, anchor=west}]
\nextgroupplot[axis equal image, width=1\textwidth, xtick={-5, -4, -2.5, 0, 2.5, 5}, ytick={0.5, 1, 1.5, 2}, xticklabels={}, ylabel={$t$}, xmin=-5, xmax=5, ymin=0, ymax=2]
\addplot []
graphics [xmin=-5,xmax=5,ymin=0,ymax=2] { 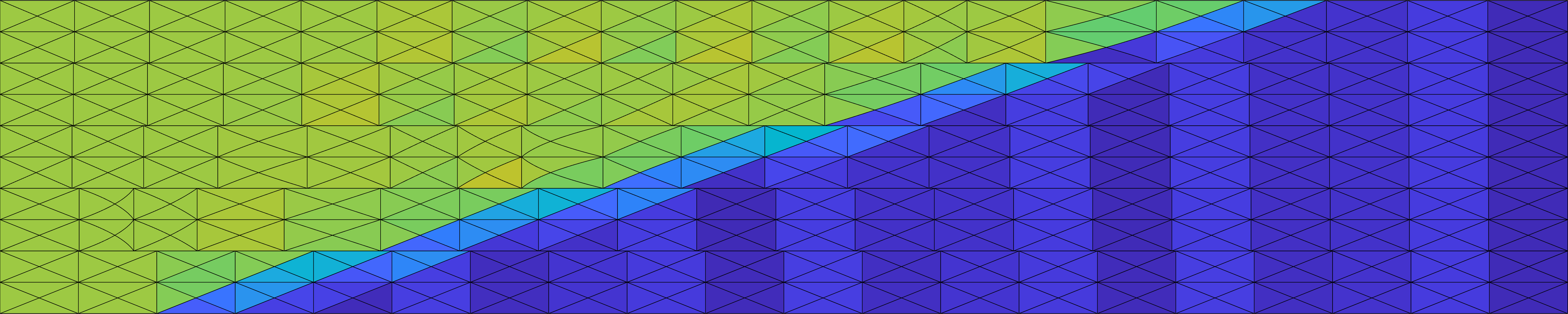};

\addplot [black, thick, dashed]
coordinates {
(-5.00000000e+00,  4.00000000e-01)
( 5.00000000e+00,  4.00000000e-01)};\label{line:shu_slab1/2}

\addplot [black, thick, dashed]
coordinates {
(-5.00000000e+00,  8.00000000e-01)
( 5.00000000e+00,  8.00000000e-01)};\label{line:slab2/3}

\addplot [black, thick, dashed]
coordinates {
(-5.00000000e+00,  1.20000000e+00)
( 5.00000000e+00,  1.20000000e+00)};\label{line:slab3/4}

\addplot [black, thick, dashed]
coordinates {
(-5.00000000e+00,  1.60000000e+00)
( 5.00000000e+00,  1.60000000e+00)};\label{line:slab4/5}

\nextgroupplot[axis equal image, width=1\textwidth, xtick={-5, -4, 0, 5}, ytick={0.5, 1, 1.5, 2}, xticklabels={}, ylabel={$t$}, xmin=-5, xmax=5, ymin=0, ymax=2]
\addplot []
graphics [xmin=-5,xmax=5,ymin=0,ymax=2] { 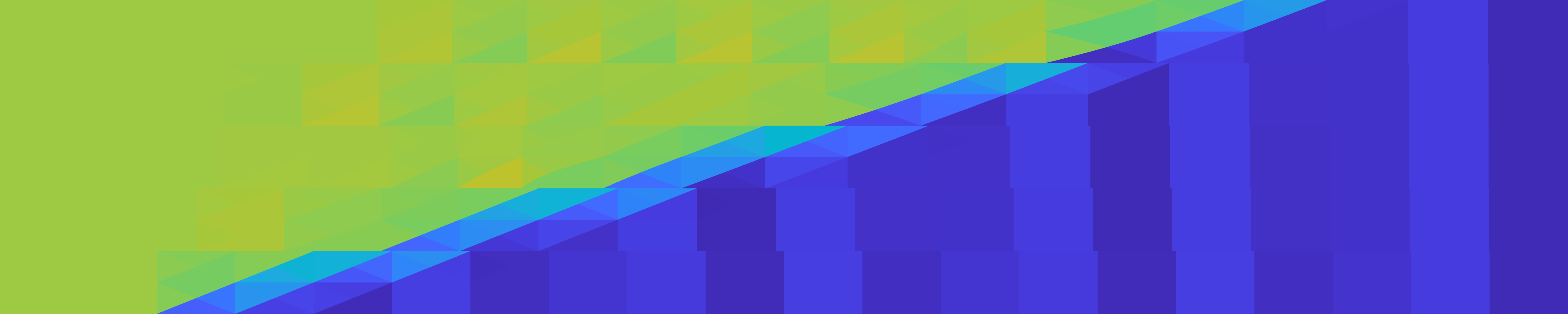};

\nextgroupplot[axis equal image, width=1\textwidth, xtick={-5, -4, 0, 5}, ytick={0.5, 1, 1.5, 2}, xticklabels={}, ylabel={$t$}, xmin=-5, xmax=5, ymin=0, ymax=2]
\addplot []
graphics [xmin=-5,xmax=5,ymin=0,ymax=2] { 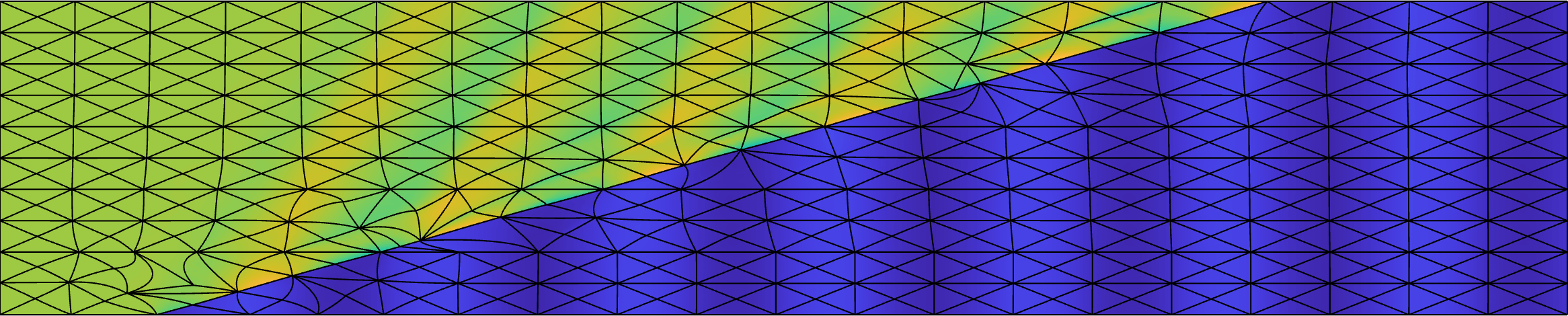};

\addplot [black, thick, dashed]
coordinates {
(-5.00000000e+00,  4.00000000e-01)
( 5.00000000e+00,  4.00000000e-01)};\label{line:slab1/2}

\addplot [black, thick, dashed]
coordinates {
(-5.00000000e+00,  8.00000000e-01)
( 5.00000000e+00,  8.00000000e-01)};\label{line:slab2/3}

\addplot [black, thick, dashed]
coordinates {
(-5.00000000e+00,  1.20000000e+00)
( 5.00000000e+00,  1.20000000e+00)};\label{line:slab3/4}

\addplot [black, thick, dashed]
coordinates {
(-5.00000000e+00,  1.60000000e+00)
( 5.00000000e+00,  1.60000000e+00)};\label{line:slab4/5}

\nextgroupplot[axis equal image, width=1\textwidth, xtick={-5, -4, -2.5, 0, 2.5, 5}, ytick={0, 0.5, 1, 1.5, 2}, xlabel={$x$}, ylabel={$t$}, xmin=-5, xmax=5, ymin=0, ymax=2]
\addplot []
graphics [xmin=-5,xmax=5,ymin=0,ymax=2] { 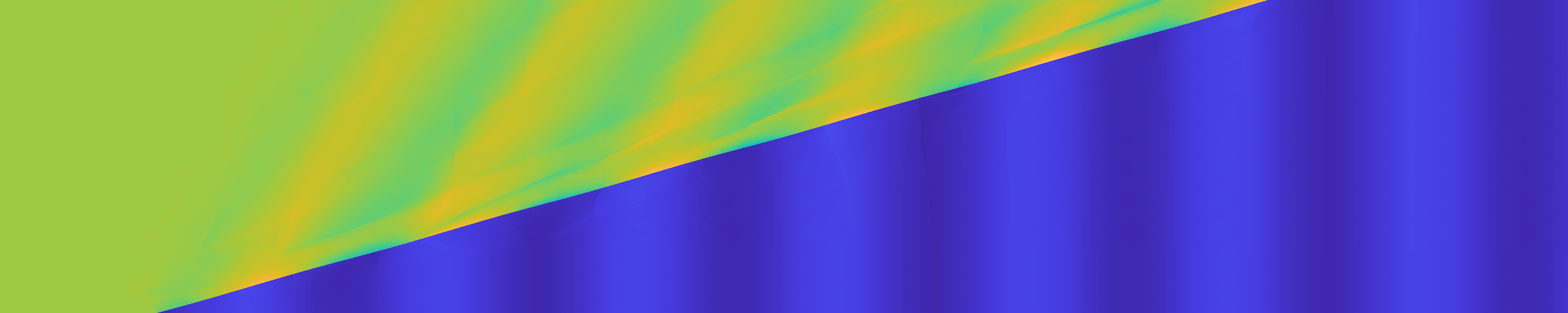};

\end{groupplot}\end{tikzpicture}}
 \colorbarMatlabParula{0.75}{2}{3}{4}{5.2}
 \caption{Slab-based HOIST solution (density) for the Shu-Osher problem, including the initial (\textit{top two rows}) and converged (\textit{bottom two rows}) mesh and flow solution. The interface between each slab is denoted by (\ref{line:shu_slab1/2}) in both the initial solution and tracked solution.}
 \label{fig:shu_osher_initial_level1}
\end{figure}

\begin{figure}
 \centering
 \raisebox{-0.5\height}{\begin{tikzpicture}
\begin{groupplot} [
group style={group size = 1 by 4, horizontal sep = 0.05cm, vertical sep = 0.05cm},
title style={at={(current bounding box.north west)}, anchor=west}]
\nextgroupplot[axis equal image, width=1\textwidth, xtick={-5, -4, 0, 5}, ytick={0.5, 1, 1.5, 2}, xticklabels={}, ylabel={$t$}, xmin=-5, xmax=5, ymin=0, ymax=2]
\addplot []
graphics [xmin=-5,xmax=5,ymin=0,ymax=2] { 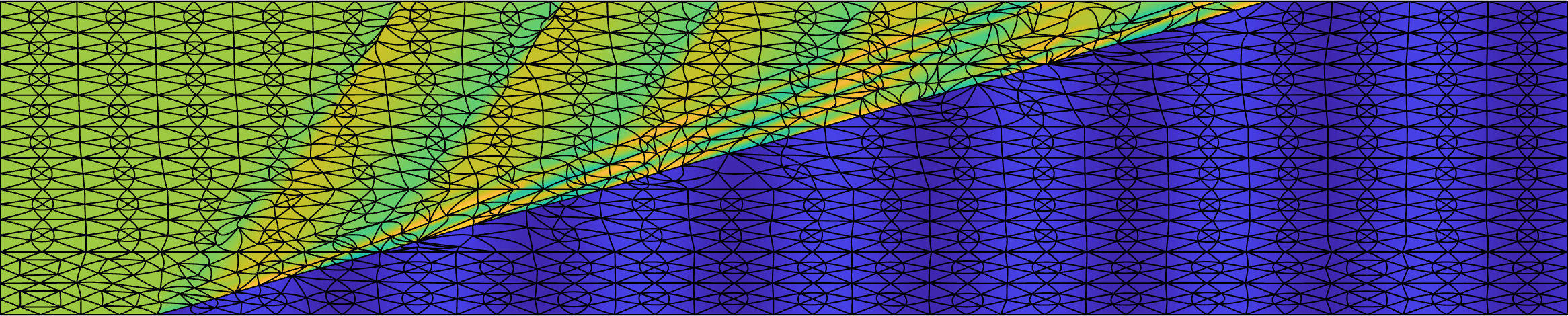};

\nextgroupplot[axis equal image, width=1\textwidth, xtick={-5, -4, 0, 5}, ytick={0.5, 1, 1.5, 2}, xticklabels={}, ylabel={$t$}, xmin=-5, xmax=5, ymin=0, ymax=2]
\addplot []
graphics [xmin=-5,xmax=3.061,ymin=0,ymax=2] { 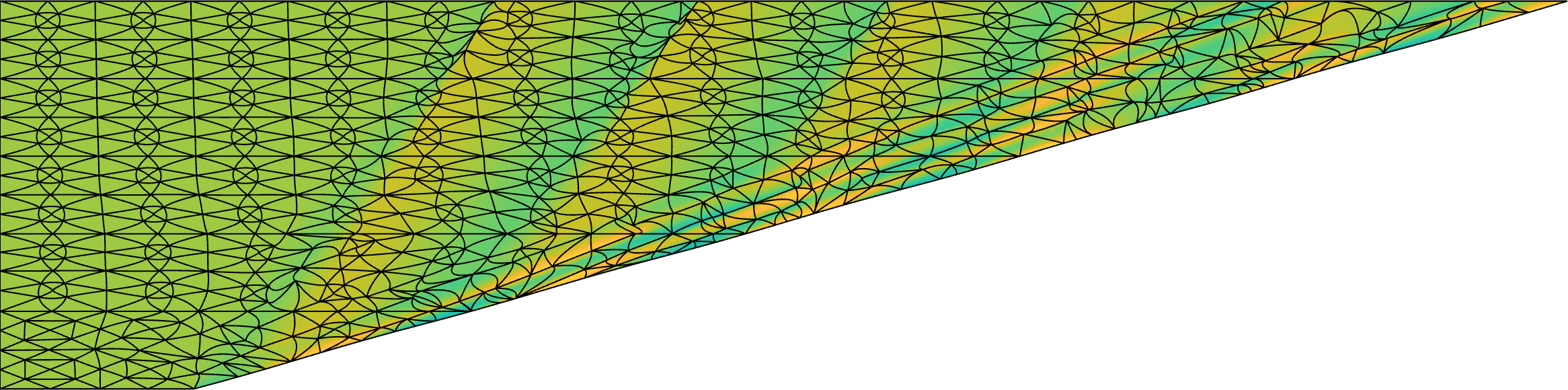};

\nextgroupplot[axis equal image, width=1\textwidth, xtick={-5, -4, 0, 5}, ytick={0.5, 1, 1.5, 2}, xticklabels={}, ylabel={$t$}, xmin=-5, xmax=5, ymin=0, ymax=2]
\addplot []
graphics [xmin=-5,xmax=3.061,ymin=0,ymax=2] { 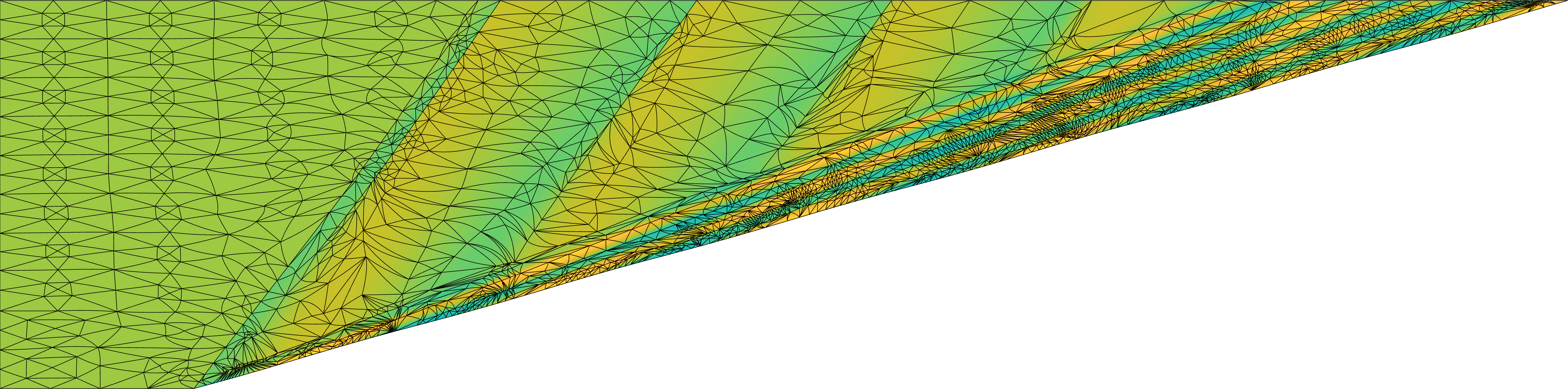};

\nextgroupplot[axis equal image, width=1\textwidth, xtick={-5, -4, -2.5, 0, 2.5, 5}, ytick={0, 0.5, 1, 1.5, 2}, xlabel={$x$}, ylabel={$t$}, xmin=-5, xmax=5, ymin=0, ymax=2]
\addplot []
graphics [xmin=-5,xmax=3.061,ymin=0,ymax=2] { 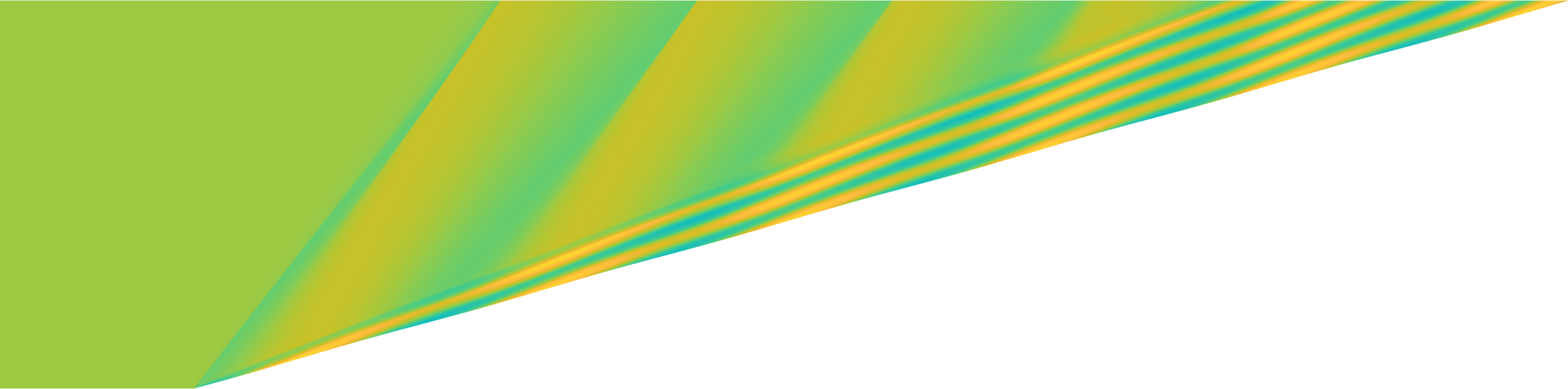};

\end{groupplot}\end{tikzpicture}}
 \caption{HOIST solution for the Shu-Osher problem after uniform mesh refinement (\textit{top row}), removing the region of the mesh with known solution (\textit{second row}), and using adaptive mesh refinement (\textit{bottom two rows}). Colorbar in Figure~\ref{fig:shu_osher_initial_level1}.}
 \label{fig:shu_osher_full_slabs}
\end{figure}

\begin{figure}
 \centering
 \raisebox{-0.5\height}{\input{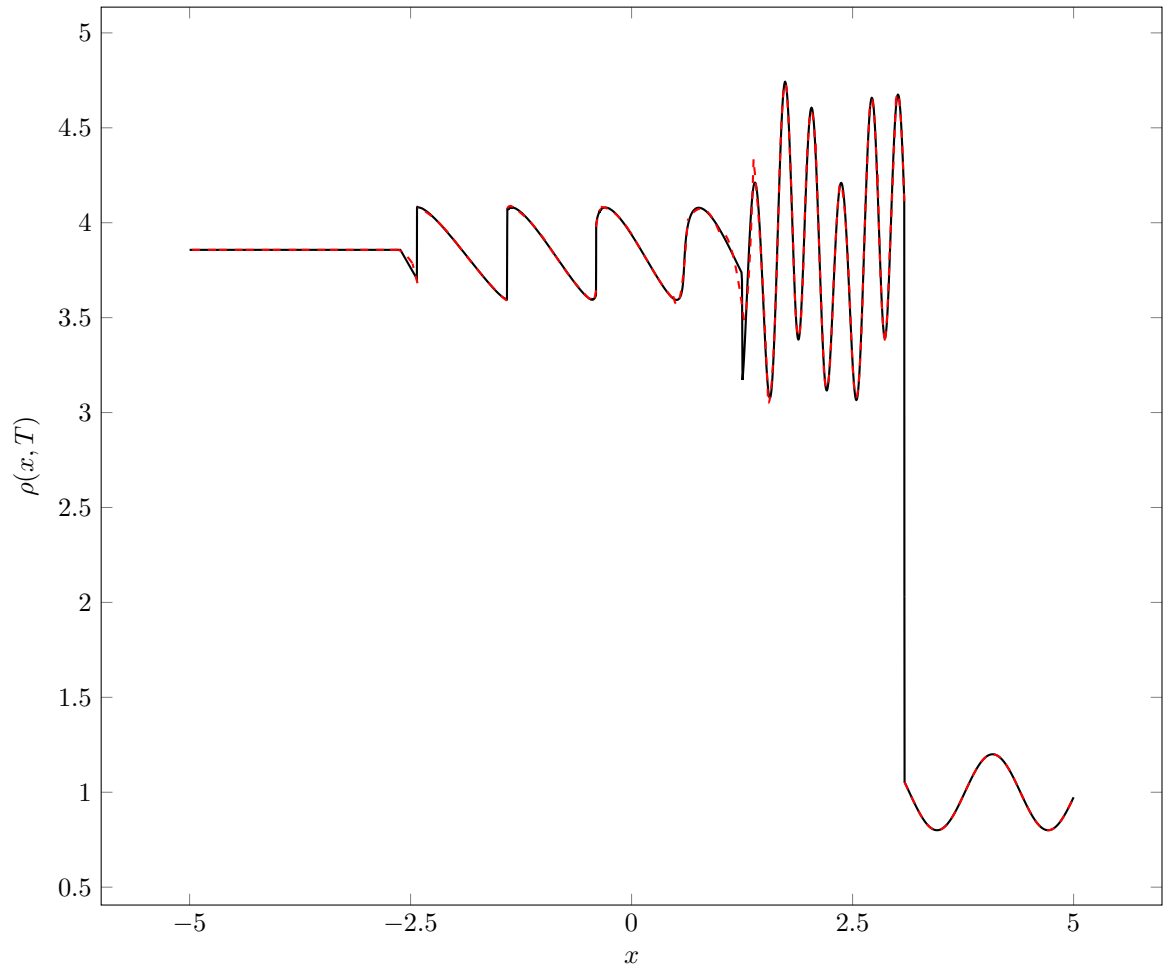}}
 \caption{Time slice of HOIST solution (density) on adapted mesh (\ref{line:ist_sol}) and a reference solution (\ref{line:ref_sol}) at $T = 2$.}
 \label{fig:shu_osher_slice}
\end{figure}

This problem is initially solved using the slab-based space-time HOIST method
with $S = 5$ slabs. The initial slab is generated by extruding a spatial mesh
of $20$ elements and splitting the quadrilaterals to produce $160$ quadratic
($p = q = 2$) space-time simplices. Subsequent slabs are generated using the
approach in Section~\ref{sec:slab}. The purpose of this initial simulation
is to approximately track the main discontinuity so shock boundary conditions
can be used. As such, the HOIST method is terminated once a weak optimality
tolerance ($\epsilon = 10^{-8}$) is satisfied or after $200$ SQP iterations. This
leads to accurate representation of the main shock; however, the oscillations
and weaker shocks are not resolved (Figure~\ref{fig:shu_osher_initial_level1}).
Next, the mesh is uniformly refined and the region of the domain upstream of the shock,
where the solution is
known analytically, is removed, which reduces the elements in the mesh from
$3184$ elements to $1500$ elements. Finally, the space-time HOIST
method is applied to this reduced problem with adaptive mesh refinement
and shock boundary conditions, which allows the finer features (shock
formation and post-shock oscillation) to be represented accurately
(Figure~\ref{fig:shu_osher_full_slabs}). The HOIST solution at the final time $T = 2$ shows
excellent agreement with a reference solution (second-order finite
volume method with superbee limiter on a grid with over $5\times 10^4$ and
fourth-order Runge-Kutta time integration with adaptive time steps)
everywhere except a few spatial locations (Figure~\ref{fig:shu_osher_slice}). The relative $L^2$ error norm of the HOIST solution compared to the reference solution at final time $T=2$ is $3.97 \times 10^{-4}$.

\section{Conclusion}
\label{sec:conclude}
High Order Implicit Shock Tracking (HOIST) is a method that approximates solutions of
conservation laws by a DG discretization of a conservation to represent non-smooth
features with inter-element jumps, while smooth regions of the flow are approximated
using high-order polynomials, which eliminates the need for stabilization techniques
(e.g., limiting and artificial viscosity) commonly used for shock capturing.
Shock-aligned meshes are produced as the solution of an optimization problem
over the DG state vector and nodal coordinates of the grid.

In this work, we extend the HOIST method to time-dependent problems using 
a slab-based space-time approach. This is achieved by reformulating a
time-dependent conservation law in $d'$ dimensions as steady conservation
law in $d'+1$ dimensions and applying the HOIST method of
\cite{2020_zahr_HOIST,huang2022robust}.
By systematically building a transformed space-time conservation law
from a time-dependent conservation law, the flux function, source term,
and numerical flux function can be explicitly expressed in terms of the
well-known spatial counterparts (Section~\ref{sec:govern}). To avoid
computations over the entire time domain and unstructured mesh generation
in $d'+1$ dimensions, we introduce a general procedure to generate conforming,
simplex-only meshes of space-time slabs in such a way that preserves features
(e.g., curved elements, refinement regions) from previous time slabs
(Section~\ref{sec:slab}). In addition to avoiding global temporal coupling,
space-time slabs also reduce the difficulty of the shock tracking problem,
especially for problems with complex interactions that evolve over time.

While the space-time HOIST method is formulated exactly as the original HOIST
method \cite{2020_zahr_HOIST} over the space-time slab,
several practical adaptations proved beneficial in the space-time setting
including 1) a translating temporal boundary (Section~\ref{sec:ist:dommap}),
2) nondimensionalization of a space-time slab (Section~\ref{sec:ist:scaling}),
3) adaptive mesh refinement (Section~\ref{sec:ist:amr}), and
4) shock boundary conditions (Section~\ref{sec:ist:shkbc}). Several
numerical experiments demonstrate the proposed space-time HOIST
method provides highly accurate solutions of conservation laws
on coarse space-time grids and is sufficiently robust to accurately
track and resolve striaght-sided non-smooth features (rarefactions,
contact discontinuities, shock waves) curved shocks, shock formation,
shock-shock and shock-boundary interactions, and triple points.

Future research directions include extending this methods to more complex
and higher dimensional problems. This will require parallelization of the
HOIST solver, likely via domain decomposition, and suitable preconditioners
for the linearize HOIST optimality system. Additionally, the extrude-then-split
approach proposed in this work applies in any dimension $d$; however, become
inefficient in $d=3$ and impractical in $d = 4$ because of the larger number
of simplices required for a conforming split of a prism ($N_3 = 14$, $N_4 = 58$).
Several approaches should be investigated to circumvent this: 1) a global
splitting algorithm like \cite{2015_wang_sptm} to reduce the $N_d$ factor
and produce a conforming simplex-only mesh, 2) more efficient local splittings
that do not guarantee conforming simplex-only meshes with non-conforming
extension of the HOIST method, and 3) generalization of HOIST for
non-simplex grids so (extruded) block or prism meshes can be used
directly for shock tracking.


\appendix

\section{Projected inviscid flux Jacobians}
\label{app:projjac}
In this section, we formulate the three conservation laws considered in this work
in the generic form of (\ref{eqn:gen_cons_law}) and present the inviscid flux
Jacobians and their eigenvalue decompositions. The Jacobians and their eigenvalue
decompositions were derived by hand for Burgers' equation and using
Maple\textsuperscript{TM} for the other equations.

\subsection{Inviscid Burgers' equation}
For the inviscid Burgers' equation (\ref{eqn:burg}), the solution vector, inviscid flux
function, and source term take the form
\begin{equation}
 U_x(x,t) = w(x,t), \qquad \Fcal_x(W_x) = (W_x^2/2)\beta^T, \qquad \Scal_x(W_x) = 0,
\end{equation}
where the projected inviscid flux Jacobian and its eigenvalue decomposition are given
by
\begin{equation}
 B_x(W_x, \eta_x) = (\beta\cdot\eta_x)W_x, \qquad
 V_x(W_x, \eta_x) = V_x(W_x, \eta_x)^{-1} = 1, \qquad
 \Lambda_x(W_x, \eta_x) = (\beta\cdot\eta_x)W_x.
\end{equation}

\subsection{Shallow water equations}
For the shallow water equations (\ref{eqn:swe}), the solution vector, inviscid flux
function, and source term take the form
\begin{equation}
 U_x(x,t) = \begin{bmatrix} \rho(x,t) \\ \rho(x,t) v(x,t) \end{bmatrix},
 \qquad
 \Fcal_x(W_x) = \begin{bmatrix} \check\rho\check{v}^T \\ \check\rho\check{v}\check{v}^T + (\check\rho^2/2) I_{d'}\end{bmatrix},
 \qquad
 \Scal_x(W_x) = \begin{bmatrix} 0 \\ g \check\rho \nabla h(x) \end{bmatrix}
\end{equation}
where $W_x = (\check\rho, \check\rho\check{v})$. The projected inviscid flux Jacobian
is given by
\begin{equation}
 B_x(W_x, \eta_x) = \begin{bmatrix} 0 & \eta_x^T \\ \check\rho\eta_x - \check{v}_{\eta_x}\check{v} & \check{v}\eta_x^T + \check{v}_{\eta_x}I_{d'} \end{bmatrix},
\end{equation}
where $\check{v}_{\eta_x} = \check{v}\cdot \eta_x$. In the $d'=1$ case, the
eigenvalue decomposition of the projected inviscid flux Jacobian is
\begin{equation}
\begin{aligned}
 \Lambda_x(W_x, \eta_x) &= \begin{bmatrix} \check{v}_{\eta_x}-\sqrt{\check\rho} & \\ & \check{v}_{\eta_x}+\sqrt{\check\rho} \end{bmatrix} \\
 V_x(W_x, \eta_x) &= \begin{bmatrix} 1 & 1 \\ \check{v}-\sqrt{\check\rho}\eta_x & \check{v}+\sqrt{\check\rho}\eta_x \end{bmatrix} \\
 V_x(W_x, \eta_x)^{-1} &= \frac{1}{2\check\rho} \begin{bmatrix} \check{v}_{\eta_x} + \sqrt{\check\rho} & -\eta_x^T \\ -\check{v}_{\eta_x} + \sqrt{\check\rho} & \eta_x^T \end{bmatrix}
\end{aligned}
\end{equation}
whereas in the $d'=2$ case, it takes the form
\begin{equation}
\begin{aligned}
 \Lambda_x(W_x, \eta_x) &= \begin{bmatrix} \check{v}_{\eta_x}-\sqrt{\check\rho} & & \\ & \check{v}_{\eta_x} & \\ & & \check{v}_{\eta_x}+\sqrt{\check\rho} \end{bmatrix} \\
 V_x(W_x, \eta_x) &= \begin{bmatrix} 1 & 0 & 1 \\ \check{v}-\sqrt{\check\rho}\eta_x & g(\eta_x) & \check{v}+\sqrt{\check\rho}\eta_x \end{bmatrix} \\
 V_x(W_x, \eta_x)^{-1} &= \frac{1}{2\check\rho} \begin{bmatrix} \check{v}_{\eta_x} + \sqrt{\check\rho} & -\eta_x^T \\ -2\sqrt{\check\rho}(\check{v}\cdot g(\eta_x)) & 2\sqrt{\check\rho}g(\eta_x)^T \\ -\check{v}_{\eta_x} + \sqrt{\check\rho} & \eta_x^T \end{bmatrix}
\end{aligned}
\end{equation}
where $g(\eta_x) = (-(\eta_x)_2, (\eta_x)_1)$.

\subsection{Euler equations of gasdynamics}
For the Euler equations with ideal gas law (\ref{eqn:euler}), the solution vector,
inviscid flux function, and source term take the form
\begin{equation}
 U_x(x,t) = \begin{bmatrix} \rho(x,t) \\ \rho(x,t) v(x,t) \\ \rho(x,t)E(x,t) \end{bmatrix},
 \qquad
 \Fcal_x(W_x) = \begin{bmatrix} \check\rho\check{v}^T \\ \check\rho\check{v}\check{v}^T + \check{P} I_{d'} \\ (\check\rho\check{E} + \check{P})\check{v}^T \end{bmatrix},
 \qquad
 \Scal_x(W_x) = 0
\end{equation}
where $W_x = (\check\rho, \check\rho\check{v}, \check\rho\check{E})$ and
$\check{P} = (\gamma-1)(\check\rho\check{E}-\check\rho\check{v}\cdot\check{v}/2)$.
The projected inviscid flux Jacobian is
\begin{equation}
B_x(W_x,\eta_x) = 
 \begin{bmatrix}
 0 & \eta_x^T & 0 \\
 \left(\frac{\check\gamma-1}{2}\right)\norm{\check{v}}^2\eta_x-\check{v}_{\eta_x} \check{v} & \check{v}_{\eta_x} I_{d'} + \check{v}\eta_x^T - (\check\gamma-1)\eta_x\check{v}^T& (\check\gamma-1)\eta_x \\
 \left[\left(\frac{\check\gamma-1}{2}\right)\norm{\check{v}}^2-\check{H}\right]\check{v}_{\eta_x} & \check{H}\eta_x^T-(\check\gamma-1)\check{v}_{\eta_x}\check{v}^T & \check\gamma\check{v}_{\eta_x}
\end{bmatrix}
\end{equation}
where $\check{v}_{\eta_x} = \check{v}\cdot \eta_x$ and
$\check{H} = (\check\rho\check{E}+\check{P})/\check\rho$.
The eigenvalues and right eigenvectors of the projected inviscid Jacobian are given by
\begin{equation} \label{eqn:euler-evd1}
\begin{aligned}
 \Lambda(W_x,\eta_x) =
 \begin{bmatrix}
  \check{v}_{\eta_x}-\check{c} & & \\
  & \check{v}_{\eta_x} I_{d'} & \\
  & & \check{v}_{\eta_x}+\check{c}
 \end{bmatrix}, \qquad
 V_x(W_x,\eta_x) = 
 \begin{bmatrix}
  1 & \eta_x^T & 1 \\
  \check{v}-\check{c} \eta_x & (\check{v}-\eta_x)\eta_x^T+I_{d'} & \check{v}+\check{c}\eta_x \\
  \check{H}-\check{v}_{\eta_x}\check{c} & \check{v}^T+\left(\norm{\check{v}}^2/2-\check{v}_{\eta_x}\right)\eta_x^T & \check{H}+\check{v}_{\eta_x}\check{c}
 \end{bmatrix},
\end{aligned}
\end{equation}
and the left eigenvectors are
\begin{equation} \label{eqn:euler-evd3}
 V_x(W_x,\eta_x)^{-1} = 
 \frac{\check\gamma-1}{2\check{c}^2}
 \begin{bmatrix}
 \norm{\check{v}}^2/2+\frac{\check{v}_{\eta_x}\check{c}}{\check\gamma-1} & -\check{v}^T-\frac{\check{c}}{\check\gamma-1}\eta_x^T & 1 \\
 \frac{2\check{c}^2}{\check\gamma-1}(\check{v}_{\eta_x}\eta_x+\eta_x-\check{v})-\norm{\check{v}}^2\eta_x & 2\eta_x\check{v}^T + \frac{2\check{c}^2}{\check\gamma-1}(I_{d'}-\eta_x\eta_x^T)& -2 \eta_x \\
 \norm{\check{v}}^2/2-\frac{\check{v}_{\eta_x}\check{c}}{\check\gamma-1} & -\check{v}^T+\frac{\check{c}}{\check\gamma-1}\eta_x^T & 1
 \end{bmatrix},
\end{equation}
where $\check{c} = \sqrt{\gamma\check{P}/\check\rho}$.

\section*{Acknowledgments}
This work is supported by AFOSR award numbers FA9550-20-1-0236,
FA9550-22-1-0002, FA9550-22-1-0004, and ONR award number
N00014-22-1-2299. The content of this publication does not
necessarily reflect the position or policy of any of these
supporters, and no official endorsement should be inferred.

\bibliographystyle{plain}
\bibliography{ist2023_jcp}

\end{document}